\documentclass{memo-l}

\def\myqed{}
\usepackage[numbers]{natbib}
\usepackage{amsmath, amssymb,stmaryrd}
\usepackage{mathtools}
\usepackage{type1cm}
\usepackage{tensor}
\usepackage{enumerate}
\usepackage{graphicx}
\usepackage[czech,english]{babel}
\usepackage{mdwtab}
\usepackage[OT1]{fontenc}
\usepackage[all]{xy}
\newcommand{\sig}{\lambda}
\newcommand{\leb}{\lambda}
\newcommand{\spect}{\hspace{1pt}\mathsf{sp}}
\newcommand{\entails}{\vdash}

\newcommand{\bbbn}{\mathbb{N}}
\newcommand{\bbbz}{\mathbb{Z}}
\newcommand{\bbbr}{\mathbb{R}}

\newcommand{\Sym}[1]{\mathfrak S_{#1}}
\newcommand{\clos}{\rm Clos}
\newcommand{\mybx}[2]{{\begin{array}{c}{#1}\\{#2}\end{array}}}
\newcommand{\myb}[2]{\mybx{#1}{\text{with }#2}}

\newcommand{\shtm}{\,\widetilde{\triangledown}\,}

\usepackage{color}
\definecolor{shadecolor}{gray}{.85}%
\definecolor{tintedcolor}{gray}{.80}%
\usepackage{framed}%
\definecolor{mytintedcolor}{gray}{.95}%
\makeatletter
\newdimen\svparindent
\newenvironment{mytinted}{%
 \begin{samepage} \MakeFramed {\FrameRestore}}%
{\endMakeFramed\end{samepage}}
{\endlist\end{mytinted}\egroup}
\makeatother
\newcommand{\rss}{relational sample space}
\newcommand{\Rss}{Relational sample space}
\newcommand{\rps}{modeling}
\newcommand{\Rps}{Modeling}
\newcommand{\Type}[2]{\tensor*[^{#1}_{\scriptscriptstyle\mathbf #2}]{\sf\textstyle Tp}{}}

\newcommand{\ERCagreement}{{\begin{minipage}{.56\textwidth}This paper is part of a project that has received funding from the European Research Council (ERC) under the European Union's Horizon 2020 research and innovation programme (grant agreement No 810115 -- {\sc Dynasnet}). \end{minipage}\hfill\begin{minipage}{.33\textwidth}\includegraphics[width=\textwidth]{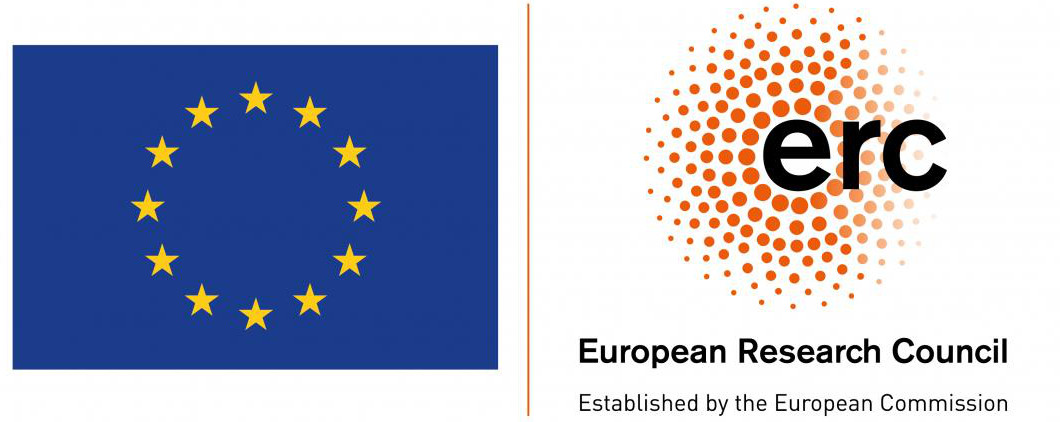}\end{minipage}\hfill}}

\newtheorem{theorem}{Theorem}[chapter]
\newtheorem{proposition}[theorem]{Proposition}
\newtheorem{corollary}[theorem]{Corollary}
\newtheorem{lemma}[theorem]{Lemma}

\theoremstyle{definition}
\newtheorem{definition}[theorem]{Definition}
\newtheorem{example}[theorem]{Example}

\theoremstyle{remark}
\newtheorem{remark}[theorem]{Remark}
\newtheorem{fact}[theorem]{Fact}
\newtheorem{problem}{Problem}[chapter]
\newtheorem{conjecture}{Conjecture}[chapter]

\numberwithin{section}{chapter}
\numberwithin{equation}{chapter}

\makeindex

\begin{document}

\frontmatter

\title{A Unified Approach to Structural Limits and Limits of Graphs with Bounded Tree-depth}


\author{Jaroslav Ne{\v s}et{\v r}il}
\address{Jaroslav Ne{\v s}et{\v r}il\\
Computer Science Institute of Charles University (IUUK and ITI)\\
   Malostransk\' e n\' am.25, 11800 Praha 1, Czech Republic}
\email{nesetril@iuuk.mff.cuni.cz}
\thanks{\ERCagreement}

\author{Patrice Ossona~de~Mendez}
\address{Patrice~Ossona~de~Mendez\\
Centre d'Analyse et de Math\'ematiques Sociales (CNRS, UMR 8557)\\
  190-198 avenue de France, 75013 Paris, France}
\email{pom@ehess.fr}

\date{October 20th, 2013}




\begin{abstract}
In this paper we introduce a general framework for the study of 
limits of relational structures and graphs in particular, which is 
based on a combination of model theory and (functional) analysis. We show how the various
approaches to graph limits fit to this framework and that they naturally appear 
as ``tractable cases'' of a general theory.
As an outcome of this, we provide extensions
of known results. We believe that this puts these into a broader context.
The second part of the paper is devoted to the study of sparse structures. First, we consider limits
of  structures with bounded diameter connected components and  
we prove
that in this case the convergence can be ``almost'' studied component-wise.
We also propose the structure of limit objects for convergent sequences of sparse structures.
Eventually, we consider the specific case of limits of colored rooted trees with bounded height
and of graphs with bounded tree-depth, motivated by
their role as ``elementary bricks'' these graphs play in decompositions of sparse
graphs, and give an explicit construction of a limit object in this case. This
limit object is a graph built on a standard probability space with the
property that every first-order definable set of tuples is measurable. 
This is an example of the general concept of {\em modeling} we introduce
here. Our example is also the first ``intermediate class'' with explicitly defined
limit structures where the inverse problem has been solved.
\end{abstract}

\maketitle

\tableofcontents


\mainmatter
\chapter{Introduction}
\label{sec:intro}
To facilitate the study of the asymptotic properties of finite graphs (and more generally of finite structures) in a sequence $G_1, G_2,\ldots, G_n,\ldots$, it
is natural to introduce notions of {\em structural convergence}. By structural convergence, we mean that
we are interested in the characteristics of a typical vertex (or group of vertices) in the graph $G_n$, as $n$ grows to infinity.
This convergence can be concisely expressed by various means. We note two main directions:
\begin{itemize}
	\item the convergence of the sampling distributions;
\item the convergence with respect to a metric in the space of structures (such as the cut metric).
\end{itemize}

Also, sampling from a limit structure  may also be used to define a sequence convergent to the
limit structure.

All these directions lead to a rich theory which originated in a probabilistic context by Aldous~\cite{Aldous1981} and Hoover~\cite{Hoover1979}
(see also the monograph of Kallenberg~\cite{Kallenberg2005} and the survey of Austin~\cite{Austin2008}) 
and, independently, in the study of random graph
processes, and in analysis of properties of random (and quasirandom) graphs (in turn motivated among others by statistical physics 
\cite{Borgs20081801,Borgs2012,Lov'asz2006}).  This development is nicely documented in  the recent monograph
 of Lov\' asz \cite{LovaszBook}.

The asymptotic properties of large graphs are studied also in the context of decision problems as exemplified e.g.\ by structural graphs theory, \cite{D,Sparsity}.
However it seems that the existential approach typical for decision problems, structural graph theory and model theory on the one side and the counting approach typical for statistics and probabilistic approach on the other side
have little in common and lead to different directions: on the one side to study, say, definability of various classes and the properties of the homomorphism order and on the other side, say, properties of partition functions.
It has been repeatedly stated that these two extremes are somehow incompatible and lead to different areas of study 
(see e.g.~\cite{borgchaylovas06:_count_graph_homom,HN}).
In this paper we take a radically different approach which unifies these both extremes.

  We propose here a model which  is a mixture of the analytic, model theoretic and algebraic approach. It is also a mixture of existential and probabilistic approach.
Precisely, our  approach 
is based on the {\em Stone pairing} $\langle \phi,G\rangle$ of a first-order formula $\phi$ 
(with set of free variables ${\rm Fv}(\phi)$) and 
a graph $G$, which 
is defined by  the  following expression
$$
\langle \phi,G\rangle=
\frac{|\{(v_1,\dots,v_{|{\rm Fv}(\phi)|})\in G^{|{\rm Fv}(\phi)|}: G\models\phi(v_1,\dots,v_{|{\rm Fv}(\phi)|})\}|}
{|G|^{|{\rm Fv}(\phi)|}}.
$$

Stone pairing induces a notion of convergence: a sequence of graphs $(G_n)_{n\in\bbbn}$ is {\em ${\rm FO}$-convergent} if, for every first order
formula $\phi$ (in the language of graphs), the values $\langle \phi, G_n\rangle$ converge
as $n\rightarrow\infty$.
In other words, $(G_n)_{n\in\bbbn}$ is ${\rm FO}$-convergent if
the probability that a formula $\phi$ is satisfied by the graph $G_n$ with a random assignment of vertices 
of $G_n$ to the free variables of $\phi$ converges as $n$ grows to infinity.
We also consider analogously defined $X$-convergence, where $X$ is a fragment of ${\rm FO}$.

Our main result is that this model of FO-convergence is a suitable model for the analysis of limits of sparse graphs (and particularly of graphs with bounded tree depth).
This fits to a broad context of recent research. 

For graphs, and more generally for finite structures, there is a class dichotomy: 
{\em nowhere dense} and {\em somewhere dense} \cite{ND_characterization,ND_logic}. Each
class of graphs falls in one of these two categories. Somewhere dense class
 ${\mathcal C}$ may be characterised by saying that there exists a (primitive positive) 
 FO interpretation of all graphs into them. Such class   ${\mathcal C}$ is
 inherently a class of dense graphs. 
In the theory of nowhere dense structures~\cite{Sparsity} there are two extreme conditions related to sparsity: bounded degree and bounded diameter.
Limits of bounded degree graphs have been studied thoroughly \cite{Benjamini2001}, and this setting has been partially extended to sparse graphs with far away large degree vertices \cite{Lyons2005}. The class of graphs with bounded diameter is considered in Section~\ref{sec:decompose} (and leads to a difficult analysis of componentwise convergence). This analysis provides a first-step for the study of limits of graphs with bounded {\em tree-depth}.
Classes of graphs with bounded tree-depth can be defined by logical terms as well as combinatorially in various ways;
the most concise definition is perhaps that a class of graphs has  bounded tree depth if and only if 
the maximal length of a path in every $G$ in the class is bounded by a constant. 
Graphs with bounded tree-depth play also the role of building blocks of graphs in a nowhere
dense class (by means of  {\em low tree-depth decompositions}~\cite{Taxi_tdepth,POMNI,Sparsity}).
So the solution of limits for graphs with bounded tree depth presents a step (and perhaps provides a road map) in solving the limit problem for sparse graphs.

We propose here a new type of measurable structure, called {\em modeling}, which extends the notion of graphing, and which we believe is a good candidate for limit objects of sequence of graphs in a nowhere dense class.  The convergence of
graphs with bounded tree depth is analysed in detail and this leads to a construction of a modeling limits for those sequences of graphs  where all members of the sequence have uniformly bounded tree depth (see Theorem~\ref{thm:limtd}). Moreover, we characterize modelings which are limits of graphs with bounded tree-depth.

There is more to this than meets the eye:  We prove that if  ${\mathcal C}$ is a monotone class of graphs such that every ${\rm FO}$-convergent sequence
has a modeling limit then the class ${\mathcal C}$ is nowhere dense (see Theorem \ref{thm:modnd0}). This shows the natural limitations to modeling ${\rm FO}$-limits.
To create a proper model for bounded height trees we have to introduce the
model in a greater generality and it appeared that our approach relates and in
most cases generalizes, by properly choosing a fragment $X$ of ${\rm FO}$,  all existing models of graph limits.
For instance, for the fragment $X$ of all existential first-order formulas, $X$-convergence means
that the probability that a structure has a particular extension property converges.
Our approach is encouraged by the deep connections to the {\em four notions} of
convergence which have been proposed to study graph limits in different contexts.

The ultimate goal of the study of structural limits is to provide (as effectively as possible) limit objects themselves: we would like to find an object which will induce the limit distribution and encode the convergence. %

For dense graphs Lov\' asz and Szegedy managed to unveil the essential notion of 
 a {\em graphon}, which exactly fits their notion of convergence: In this representation the limit \cite{Lov'asz2006,Borgs20081801}
is a symmetric Lebesgue measurable function $W : [0, 1]^2 \rightarrow [0, 1]$ called a {\em graphon} and every graphon is the limit of a sequence of graphs.
Such a representation is of course not unique, in the sense that different graphons may define
the same graph limit, but equivalence of graphons is well understood \cite{Borgs2012a,pre05504139}. 
A connection between graph limits and
de Finetti's theorem for exchangeable arrays (and the early works of Aldous \cite{Aldous1981}, Hoover \cite{Hoover1979} 
and Kallenberg \cite{Kallenberg2005}) has been established, see e.g. Diaconis and Janson \cite{pre05504139}. Note that representation of  graph limits by graphons extend (in a non-trivial way) to regular hypergraphs \cite{ElekSze,Zhao2014} and, more generally, to relational structures \cite{Aroskar2012,aroskar2014limits}.
  
A representation of the limit for our second example of bounded degree graphs is a {\em measurable graphing} (notion introduced by Adams \cite{Adams1990} in the
context of Ergodic theory), that is a standard Borel space with a measure $\mu$ and $d$ measure preserving Borel
involutions. The existence of such a representation has been made explicit by Elek \cite{Elek2007b}, 
and relies on the works of Benjamini \cite{Benjamini2001} and Gaboriau \cite{Gaboriau2005}. Graphing representation is not unique, but the equivalence of graphings (called {\em local equivalence}) can be characterized by means of  bi-local isomorphism  \cite{LovaszBook}. However, it is a difficult open problem, known as Aldous--Lyons conjecture, whether every graphing is the limit of some sequence of finite graphs (see Conjecture~\ref{conj:AL}).

Both of these models of convergence are particular cases of our general approach. 
One of the main issue of our general approach is to determine a representation of ${\rm FO}$-limits as
measurable graphs. A natural limit object is a standard probability space
$(V,\Sigma,\mu)$ together with a graph with vertex set $V$ and edge set $E$, with
the property that every first-order definable subset of a power of $V$ is
measurable. This leads to the notion of {\rss} and to the notion of {\em modeling}.
This notion seems to be particularly suitable for sparse graphs (and in the full generality only for sparse graphs, see Theorem~\ref{thm:modnd0}).
We shall see that modelings inherit most of the nice properties of graphings and that open problems on graphings can be generalized to open problems on modelings (in particular the Aldous--Lyons conjecture mentioned above). It is open which type of limit object could be considered for the general (sparse and dense) case, which would generalize graphons and graphings.

In this paper, we shed a new light on all these constructions by an approach inspired by
functional analysis. The preliminary material and our framework are introduced in 
Sections~\ref{sec:def} and~\ref{sec:stone}.
The general approach presented in the first sections of this paper leads
to several new results. Let us mention a sample of such results.

 Central to the theory of graph limits  stand random graphs
  (in the Erd{\H o}s-R{\'e}nyi model, where each edge is present  with a given probability $p$, independently of the other edges \cite{ErdHos1960}): a sequence of random graphs with increasingly many vertices 
 and fixed edge probability $0<p<1$  
 is almost surely convergent to the constant graphon $p$
  \cite{Lov'asz2006}. On the other hand, it follows from the work of
    Erd{\H o}s and R{\'e}nyi \cite{asym} and the work of
Glebskii,  Kogan,  Liagonkii and  Talanov \cite{Glebskii1969a},  Fagin \cite{Fagin1976}
   that such a sequence is almost surely elementarily convergent to an ultra-homogeneous graph, 
  called the {\em Rado graph}. 
  We prove that these
  two facts --- elementary convergence to the  Rado graph and convergence to a constant graphon --- together with the quantifier elimination property of ultra-homogeneous graphs, imply that
   a sequence of random graphs with increasing   order and fixed edge probability
   $0<p<1$ is almost surely ${\rm FO}$-convergent, see Section~\ref{sec:homog}.
   (However, we know that this limit cannot be either a random-free graphon or a modeling, see Theorem~\ref{thm:modnd0})

 We shall  prove that a sequence of bounded degree graphs $(G_n)_{n\in\bbbn}$ with $|G_n|\rightarrow\infty$ is ${\rm FO}$-convergent
  if and only if it is both convergent in the sense of Benjamini-Schramm and in the sense of elementary convergence. The limit can still be represented by a graphing, see
  Sections~\ref{sec:BS} and~\ref{sec:BS2}.

 For the general case we prove that the limit of an ${\rm FO}$-convergent sequence of graphs is a probability measure on
the Stone space of the Boolean algebra of first-order formulas, which is invariant under the action
of the symmetric group $\Sym{\omega}$ on this space, see Section~\ref{sec:stone}. This representation theorem holds generally and it is the basis of our approach. 
Fine interplay of
these notions is depicted on Table~\ref{tab:trsf}.

\begin{table}[ht]
\begin{tabular}[C]{|c|c|}
\multicolumn{1}{c}{Boolean algebra $\mathcal B(X)$}&
\multicolumn{1}{c}{Stone Space $S(\mathcal B(X))$}\\
\hlx{v[0,1,2]hvvv}
Formula $\phi$&Continuous function $f_\phi$\\
\hlx{vvvhvvv}
Vertex $v$&``Type'' of vertices $T$\\
\hlx{vvvhvvv}
Graph $G$&statistics of types \\
\hlx{v}
&=probability measure $\mu_G$\\
\hlx{vvvhvvv}
$\langle\phi,G\rangle$&$\displaystyle\int f_\phi(T)\ {\rm d}\mu_G(T)$\\
\hlx{vvvhvvv}
$X$-convergent $(G_n)$&weakly convergent $\mu_{G_n}$\\
\hlx{vvvhvvv}
$\Gamma={\rm Aut}(\mathcal B(X))$&$\Gamma$-invariant measure\\ 
\hlx{vvvh}
\end{tabular}
\caption{Some correspondences}
\label{tab:trsf}
\end{table}

Graph limits (in the sense of Lov\'asz et al.) --- and more generally hypergraph limits --- 
  have been studied by Elek and Szegedy \cite{ElekSze} through the introduction 
   of a measure on the ultraproduct of the graphs in the sequence
   (via Loeb measure construction, see \cite{Loeb1975}). 
  The fundamental theorem of ultraproducts proved 
  by {\L}o\'s \cite{Los1955a} implies that the ultralimit of a sequence of graphs is 
  (as a measurable graph) an ${\rm FO}$-limit. Thus in this non-standard setting we get
  ${\rm FO}$-limits (almost) for free see \cite{CMUC}. 
  However this very general construction has several major drawbacks in an analytical context: it involves 
  countably many measures (which are not simply product measures) and non-separable sigma algebras, while major tools from analysis rely on  Borel product measures on standard Borel spaces (like for graphings).

We believe that the approach taken in this paper is  natural and that it enriches
the existing notions of limits by several natural notions of $X$-convergence (such as elementary, quantifier-free, and local convergences), and gives the whole area a new perspective, which we explain in the next section. 
In a sense we proceed dually to homomorphism counting (see e.g. \cite{borgchaylovas06:_count_graph_homom, LovaszBook}).
We do not view $\langle\phi,G\rangle$  as a ``$\phi$ test''  for $G$  but rather as  
a ``$G$ test'' for $\phi$: A graph defines an operator 
on the Boolean algebra of all ${\rm FO}$-formulas  (or on the sub-algebra induced by a fragment $X \subset {\rm FO}$).
It also presents %
a promising approach
to more general intermediate classes (see the final comments).
\vskip 5mm

\section{Main Definitions and Results}
\label{sec:def}
If we consider relational structures with signature $\sig$, the symbols of the relations and constants
in $\sig$ define the  non-logical symbols of the vocabulary of the first-order language ${\rm FO}(\sig)$
associated to $\sig$-structures. Notice that if $\sig$ is  countable then ${\rm FO}(\sig)$ is countable.
The symbols of variables will be assumed to be taken from 
a countable 
set $\{x_1,\dots,x_n,\dots\}$ indexed by $\bbbn$.
Let $u_1,\dots,u_k$ be terms. 
The set of used free variables of a formula $\phi$ will be denoted by ${\rm Fv}(\phi)$
(by saying that a variable $x_i$ is ``used'' in $\phi$ we mean that $\phi$ is not logically equivalent
to a formula in which $x_i$ does not appear).
The formula $\phi_{x_{i_1},\dots,x_{i_k}}(u_1,\dots,u_k)$ denotes the formula obtained by substituting simultaneously
the term $u_j$ to the free occurrences of $x_{i_j}$ for $j=1,\dots,k$.
 In the sake of simplicity, we will denote
by $\phi(u_1,\dots,u_k)$ the substitution $\phi_{x_1,\dots,x_k}(u_1,\dots,u_k)$. 

A {\em relational structure} $\mathbf A$ with signature $\sig$ is defined by its {\em domain} (or {\em universe}) $A$
and relations with names and arities as defined in $\sig$. In the following we will denote relational structures
by bold face letters $\mathbf A,\mathbf B,\dots$ and their domains by the corresponding light face letters
$A,B,\dots$ 

 The key to our approach are the following two definitions.
\begin{mytinted} 
 \begin{definition}[Stone pairing]
 \label{def:bracket}
 Let $\sig$ be a signature, let $\phi\in{\rm FO}(\sig)$ be a first-order formula with 
 free variables $x_1,\dots,x_p$ and
 let $\mathbf A$ be a finite $\sig$-structure. 
 
 Put
$$
\Omega_\phi(\mathbf A)=\{(v_1,\dots,v_p)\in A^p:\ \mathbf
A\models\phi(v_1,\dots,v_p)\}.
$$
 
 We define the {\em Stone pairing} of $\phi$ and $\mathbf A$ by
 \begin{equation}
 \label{eq:def}
 \langle \phi,\mathbf A\rangle=\frac{|\Omega_\phi(\mathbf A)|}{|A|^p}.
 \end{equation}
\end{definition}
 \end{mytinted}
 
 In other words, $\langle\phi,\mathbf A\rangle$ is the probability that 
 $\phi$ is satisfied in $\mathbf A$ when we interpret the $p$ free variables
 of $\phi$ by $p$ vertices of $G$ chosen randomly, uniformly and independently.
	Also, $\Omega_\phi(\mathbf A)$ is interpreted as the {\em solution set} of $\phi$ in $\mathbf A$.

Note that  
in the case of a sentence $\phi$ (that is a formula with no free variables, thus $p=0$), the
definition of the Stone pairing reduces to
$$\langle \phi,\mathbf A\rangle=\begin{cases}
1,&\text{if }\mathbf A\models\phi;\\
0,&\text{otherwise.}
\end{cases}$$

\begin{mytinted}
\begin{definition}[FO-convergence]
A sequence $(\mathbf A_n)_{n\in\bbbn}$ of finite $\sig$-structures is {\em {\rm FO}-convergent} if,
for every formula $\phi\in{\rm FO}(\sig)$ the sequence
$(\langle\phi,\mathbf A_n\rangle)_{n\in\bbbn}$ is (Cauchy) convergent.
\end{definition}
\end{mytinted}

In other words, a sequence $(\mathbf A_n)_{n\in\bbbn}$ is FO-convergent if the sequence of
mappings $\langle\,\cdot\,,\mathbf A_n\rangle:{\rm FO}(\sig)\rightarrow [0,1]$ is pointwise-convergent. 

The interpretation of the Stone pairing as a probability 
suggests to extend this view to more general $\sig$-structures which will be our candidates for limit objects.

\begin{mytinted}
\begin{definition}[\Rss]
\label{def:rss}
A {\em \rss} is a relational structure $\mathbf A$ (with signature $\sig$) 
with extra structure:
The domain $A$ of $\mathbf A$ of a sample model is a standard Borel space
(with Borel $\sigma$-algebra $\Sigma_{\mathbf A}$)
 with the property that every subset of $A^p$ that is first-order definable 
 in ${\rm FO}(\sig)$
  is measurable (in $A^p$ with respect to the product $\sigma$-algebra). For brevity we shall use the same letter $\mathbf A$
  for  structure and relational sample space.
\end{definition}
\end{mytinted}

In other words, if $\mathbf A$ is a {\rss} then 
for every integer $p$ and every $\phi\in{\rm FO}(\sig)$
with $p$ free variables 
it holds that
$\Omega_\phi(\mathbf A)\in\Sigma_{\mathbf A}^p$.

\begin{mytinted}
\begin{definition}[\Rps]
\label{def:rps}
A {\em \rps} $\mathbf A$ is a 
{\rss} $\mathbf A$ equipped with a probability measure
(denoted $\nu_{\mathbf A}$).
By the abuse of symbols the modeling will be denoted by $\mathbf A$ (with $\sigma$-algebra $\Sigma_{\mathbf A}$
and corresponding measure $\nu_{\mathbf A}$).
A $\rps$ with signature $\sig$ is a {\em $\sig$-\rps}.

\end{definition}
\end{mytinted}

\begin{remark}
We take time for some comments on the above definitions:
\begin{itemize}
  \item According to Kuratowski's isomorphism theorem, the domains of {\rss}s are Borel-isomorphic to either
$\bbbr$, $\bbbz$, or a finite space.
\item {\em Borel graphs} (in the sense of Kechris et al. \cite{Kechris1999}) are
generally not {\rps}s (in our sense) as Borel graphs are only required to have a measurable adjacency relation.  
\item By equipping its domain with the discrete $\sigma$-algebra,  
 every finite $\sig$-struc\-ture defines 
a {\rss}. Considering the uniform probability measure on this space
then canonically defines a uniform {\rps}.
\item It follows immediately from Definition~\ref{def:rss} that
any {\em $k$-rooting} of a {\rss} is a {\rss}.
\end{itemize}
\end{remark}

We can extend the definition of Stone pairing from finite structures to {\rps}s as follows.

\begin{mytinted} 
 \begin{definition}[Stone pairing for modeling]
 \label{def:bracket2}
 Let $\sig$ be a signature, let $\phi\in{\rm FO}(\sig)$ be a first-order formula with 
 free variables $x_1,\dots,x_p$ and
 let $\mathbf A$ be a $\sig$-{\rps}. 
 
 We can define the {\em Stone pairing} of $\phi$ and $\mathbf A$ by
 \begin{equation}
 \label{eq:def2}
 \langle \phi,\mathbf A\rangle=\int_{x\in A^p}
 1_{\Omega_\phi(\mathbf A)}(x)\,{\rm d}\nu_{\mathbf A}^{p}(x).
 \end{equation}
\end{definition}
 \end{mytinted}

Note that the definition of a {\rps} is simply tailored to make the 
expression~\eqref{eq:def2} meaningful.
Based on this definition, {\rps}s can sometimes be used as a representation of the limit
of an FO-convergent sequence of finite $\sig$-structures. 

\begin{mytinted}
\begin{definition}
\label{def:modFOlim}
A {\rps} $\mathbf L$ is a {\em \rps\ {\rm FO}-limit} of an FO-convergent sequence $(\mathbf A_n)_{n\in\bbbn}$
of finite $\sig$-structures if $\langle\phi\,,\mathbf A_n\rangle$ converges pointwise
to  $\langle\,\phi\,,\mathbf L\rangle$ for every first order formula $\phi$.
\end{definition}
\end{mytinted}

As we shall see in Lemma~\ref{lem:wu}, a {\rps} FO-limit of an FO-convergent sequence $(\mathbf A_n)_{n\in\bbbn}$
of finite $\sig$-structures is necessarily weakly uniform
(meaning that all the singletons of the limit have the same measure). It follows that if a {\rps}
$\mathbf L$ is a {\rps} FO-limit then $L$ is either finite or uncountable. 

We shall see that not every FO-convergent sequence of finite relational structures admits a 
{\rps} FO-limit. In particular we prove (see Theorem~\ref{thm:modnd}):
\begin{mytinted}
\begin{theorem}
\label{thm:modnd0}
Let $\mathcal C$ be a monotone class of finite graphs, such that every ${\rm FO}$-convergent sequence
of graphs in $\mathcal C$ has a modeling ${\rm FO}$-limit. Then the class $\mathcal C$ is nowhere dense.
\end{theorem}
\end{mytinted}

Recall that a class of graphs is {\em monotone} if it is closed under the operation of taking a subgraph, and that
a monotone class of graphs $\mathcal C$ is {\em nowhere dense} if, for every integer $p$, there exists
an integer $N(p)$ such that the $p$-th subdivision of the complete graph $K_{N(p)}$ on $N(p)$ vertices 
does not belong to $\mathcal C$ (see~\cite{ND_logic,ND_characterization,Sparsity}).

However, we conjecture that the theorem above expresses exactly when modeling FO-limits exist:
\begin{mytinted}
\begin{conjecture}
\label{conj:nd}
If $(G_n)_{n\in\bbbn}$ is an ${\rm FO}$-convergent sequence
of graphs and if $\{G_n: n\in\bbbn\}$ is a nowhere dense class, then 
the sequence $(G_n)_{n\in\bbbn}$ has a modeling ${\rm FO}$-limit.
\end{conjecture}
\end{mytinted}

As a first step, we prove that {\rps} FO-limits
exist in two particular cases, which form in a certain sense the building blocks of nowhere dense classes.

\begin{mytinted}
\begin{theorem}
\label{thm:modlim}
Let $C$ be a integer.
\begin{enumerate}
  \item Every FO-convergent sequence of graphs with maximum degree at most $C$ has a {\rps} FO-limit;
  \item Every FO-convergent sequence of rooted trees with height at most $C$ has a {\rps} FO-limit.
\end{enumerate}
\end{theorem}
\end{mytinted} 

The first item will be derived from the graphing representation of limits of Benjamini-Schramm convergent
sequences of graphs with bounded maximum degree with no major difficulties. 
Recall that a {\em graphing} \cite{Adams1990} is a Borel graph $G$ such that the following 
{\em Intrinsic Mass Transport Principle} (IMTP) holds: 
\begin{mytinted}
$$
\forall A,B\qquad\int_A{\rm deg}_B(x)\,{\rm d}x=\int_B{\rm deg}_A(y){\rm d}y,
$$
\end{mytinted}
\noindent where the quantification is on all measurable subsets of vertices, and where
${\rm deg}_B(x)$ (resp.\  ${\rm deg}_A(y)$) denotes the degree in $B$ (resp.\  in $A$) of
the vertex $x$ (resp.\ of the vertex $y$). In other words, the Mass Transport Principle
states that if we count the edges between sets $A$ and $B$ by summing up the degrees in $B$ of vertices
in $A$ or by summing up the degrees in $A$ of vertices in $B$, we should get the same result.

\begin{theorem}[Elek \cite{Elek2007b}]
The Benjamini-Schramm limit of a bounded degree graph sequence can be represented by a graphing.
\end{theorem}

A full characterization
of the limit objects in this case is not known, and is related to the following conjecture.
\begin{mytinted} 
\begin{conjecture}[Aldous, Lyons \cite{Aldous2006}]
\label{conj:AL}
Every graphing is the Benjamini-Schramm limit of a bounded degree graph sequence.

Equivalently, every unimodular distribution on rooted countable graphs with bounded degree is
the Benjamini-Schramm limit of a bounded degree graph sequence.
\end{conjecture}
\end{mytinted} 

We conjecture that a similar condition could characterize modeling ${\rm FO}$-limits
of sequences of graphs with bounded degree. In this more general setting, we have to add a
new condition, namely to have the {\em finite model property}. Recall that an infinite structure $\mathbf L$ has 
the finite model property if every sentence satisfied by $\mathbf L$ has a finite model.

\begin{mytinted} 
\begin{conjecture}
\label{conj:ALFO}
A modeling is the Benjamini-Schramm limit of a bounded degree graph sequence if 
and only if it is a graph
with bounded degree, is weakly uniform, it satisfies both the Intrinsic Mass Transport Principle,
 and it has the finite model property.
\end{conjecture}
\end{mytinted} 

When handling infinite degrees, we do not expect to be able to keep the Intrinsic Mass Transport Principle as is.
If a sequence of finite graphs is ${\rm FO}$-convergent 
to some modeling $\mathbf L$ then we require the following condition to hold, 
which we call {\em Finitary Mass Transport Principle} (FMTP):
\begin{mytinted}
For every measurable subsets of vertices $A$ and $B$, if it holds that 
${\rm deg}_B(x)\geq a$ for every $x\in A$ and ${\rm deg}_A(y)\leq b$ for every $y\in B$ then
$a\,\nu_{\mathbf L}(A)\leq b\,\nu_{\mathbf L}(B)$.
\end{mytinted}
Note that in the case of modelings with bounded degrees, the Finitary Mass Transport Principle is
equivalent to the Intrinsic Mass Transport Principle. Also note that the above equation holds necessarily
when $A$ and $B$ are first-order definable, according to the convergence of the Stone pairings and the fact
that the Finitary Mass Transport Principle obviously holds for finite graphs.

The second item of Theorem~\ref{thm:modlim} will be
quite difficult to establish and is the main result of this paper. 
We formulate it together with the inverse theorem as follows:

\begin{mytinted} 
\begin{theorem}
Every sequence of finite rooted colored trees with height at most $C$ has a modeling ${\rm FO}$-limit that is
a rooted colored tree with height at most $C$, is weakly uniform, and satisfies the Finitary Mass Transport Principle.

Conversely, every %
rooted colored tree modeling with height at most $C$
that satisfies the Finitary Mass Transport Principle is the FO-limit of a sequence of finite rooted colored trees. 
\end{theorem}
\end{mytinted}

By Theorem~\ref{thm:modnd0},  {\rps} FO-limits do not exist in general.
However, we have a general representation of the limit of an FO-convergent 
sequence of $\sig$-structures by means of a probability distribution on a
compact Polish space $S_\sig$ defined from ${\rm FO}(\sig)$ using Stone duality:

\begin{mytinted}
\begin{theorem}
\label{thm:rep1}
Let $\sig$ be a fixed (possibly finite) countable signature.
Then there exist two mappings
 $\mathbf A\mapsto\mu_{\mathbf A}$  and $\phi\mapsto K(\phi)$ such that 
\begin{itemize}
  \item $\mathbf A\mapsto\mu_{\mathbf A}$ is an injective mapping from the class of finite $\sig$-structures to the 
  space of regular probability measures on  $S_\sig$,
  \item $\phi\mapsto K(\phi)$ is a mapping from ${\rm FO}(\sig)$ to the set of the clopen subsets of $S_\sig$,
\end{itemize}
 such that for every finite $\sig$-structure $\mathbf A$ and every first-order formula
$\phi\in{\rm FO}(\sig)$ the following equation holds:
$$
\langle\phi,\mathbf A\rangle=\int_{S_\sig}1_{K(\phi)}\,{\rm d}\mu_{\mathbf A}. 
$$
\end{theorem}
\end{mytinted}

(To prevent risks of notational ambiguity, we shall use $\mu$ as root symbol for measures on Stone spaces
and keep $\nu$ for measures on {\rps}s.)

Consider an FO-convergent sequence $(\mathbf A_n)_{n\in\bbbn}$.
Then the pointwise convergence of $\langle\,\cdot\,,\mathbf A_n\rangle$ translates as a
weak $\ast$-convergence of the measures $\mu_{\mathbf A_n}$ and we get:
 
\begin{mytinted}
\begin{theorem}
\label{thm:rep2}
A sequence $(\mathbf A_n)_{n\in\bbbn}$ of finite $\sig$-structures is FO-convergent if and only if
the sequence $(\mu_{\mathbf A_n})_{n\in\bbbn}$ is weakly $\ast$-convergent.\\
Moreover, if $\mu_{\mathbf A_n}\Rightarrow\mu$ then for every first-order formula $\phi\in{\rm FO}(\sig)$ the following equation holds:
$$\int_{S_\sig}1_{K(\phi)}\,{\rm d}\mu=\lim_{n\rightarrow\infty}\langle\phi,\mathbf A_n\rangle.
$$
\end{theorem}
\end{mytinted}
These last two Theorems are established in the next section as a warm up for our general theory.

\chapter{General Theory}
\section{Limits as Measures on Stone Spaces}
\label{sec:stone}
In order to prove the representation theorems Theorem~\ref{thm:rep1} and Theorem~\ref{thm:rep2}, 
we first need to prove a general representation for additive functions on Boolean algebras.
 
\subsection{Representation of Additive Functions} 
 Recall that a {\em Boolean algebra} $B=(B,\wedge,\vee,\neg,0,1)$ is an algebra  
with two binary operations $\vee$ and $\wedge$, a unary operation $\neg$
 and two elements $0$ and $1$, such that
$(B,\vee,\wedge)$ is a complemented distributive lattice with minimum $0$ and maximum $1$. 
The two-elements Boolean algebra is denoted $\mathbf 2$. 
 
To a Boolean algebra $B$ is associated a topological space, denoted $S(B)$, whose points 
are the ultrafilters on $B$  (or equivalently the homomorphisms $B\rightarrow\mathbf 2$).
The topology on $S(B)$ is generated by a sub-basis consisting of all sets 
    $$K_B(b)=\{ x \in S(B): b \in x\},$$
where $b\in B$. When the considered Boolean algebra will be clear from context
we shall omit the subscript and write $K(b)$ instead of $K_B(b)$.   

A topological space is a {\em Stone space} if it is Hausdorff, compact,
and has a basis of clopen subsets. Boolean algebras and Stone spaces are equivalent 
as formalized by Stone representation theorem \cite{Stone1936}, which states
(in the language of category theory)
that there is a duality between the category of Boolean algebras
(with homomorphisms) and the category of Stone spaces (with continuous functions).
This justifies calling $S(B)$ the {\em Stone space} of the Boolean algebra $B$.
The two contravariant functors defining this duality are denoted by $S$ and $\Omega$
and defined as follows:

For every homomorphism $h:A\rightarrow B$ between two Boolean algebra,
we define the map $S(h):S(B)\rightarrow S(A)$ by $S(h)(g)=g\circ h$
(where points of $S(B)$ are identified with homomorphisms $g:B\rightarrow\mathbf 2$).
Then for every homomorphism $h:A\rightarrow B$, the map
$S(h):S(B)\rightarrow S(A)$ is a continuous function.

Conversely, for every continuous function $f:X\rightarrow Y$
between two Stone spaces, define the map $\Omega(f):\Omega(Y)\rightarrow\Omega(X)$
by $\Omega(f)(U)=f^{-1}(U)$ (where elements of $\Omega(X)$ are identified with
 clopen sets of $X$). Then for every continuous function 
 $f:X\rightarrow Y$, the map $\Omega(f):\Omega(Y)\rightarrow\Omega(X)$ is a 
 homomorphism of Boolean algebras.

We denote by $K=\Omega\circ S$ one of the two natural isomorphisms
defined by the duality. Hence, for a Boolean algebra $B$, 
$K(B)$ is the set algebra $\{K_B(b): b\in B\}$, and this 
algebra is isomorphic to $B$. 

An ultrafilter of a Boolean algebra $B$ can be considered as a finitely additive measure, 
for which every subset has either measure $0$ or $1$.
Because of the equivalence of the notions of Boolean algebra and of set algebra, we define
the space ${\rm ba}(B)$ as the space of all bounded 
additive functions $f:B\rightarrow\bbbr$.
Recall that a function $f:B\rightarrow\bbbr$ is {\em additive} 
if for all $x,y\in B$ the following implication holds
$$
x\wedge y=0\quad\Longrightarrow\quad f(x\vee y)=f(x)+f(y).
$$
The space ${\rm ba}(B)$ is a Banach space for the norm
$$
\|f\|_{{\rm ba}(B)}=\sup_{x\in B} f(x) - \inf_{x\in B} f(x).
$$
(Recall that the ba space of an algebra of sets $\Sigma$ is the Banach space consisting of all bounded
 and finitely additive measures on $\Sigma$ with the total variation norm.)

Let $V(B)$ be the normed vector space (of so-called {\em simple functions}) 
generated by the indicator functions 
of the clopen sets (equipped with supremum norm).
The indicator function of the clopen set $K(b)$ (for some $b\in B$) is denoted by
 $\mathbf 1_{K(b)}$.

\begin{lemma}
\label{lem:baV}
The space ${\rm ba}(B)$ is the topological dual of $V(B)$.
\end{lemma}
\begin{proof}
One can identify ${\rm ba}(B)$ with the space 
${\rm ba}(K(B))$ of finitely additive measures defined on the set algebra
$K(B)$. As a vector space, ${\rm ba}(B)\approx{\rm ba}(K(B))$ is then 
clearly the (algebraic) dual of the normed vector space $V(B)$.
  
The pairing of a function $f\in {\rm ba}(B)$ and
a vector $X=\sum_{i=1}^n a_i \mathbf 1_{K(b_i)}$ is defined by
$$
[f,X]=\sum_{i=1}^n a_i f(b_i).
$$
That $[f,X]$ does not depend on a particular choice of a decomposition of $X$
follows from the additivity of $f$. We include a short proof for completeness:
Assume $\sum_i \alpha_i\mathbf 1_{K(b_i)}=\sum_i \beta_i\mathbf 1_{K(b_i)}$. As
for every $b,b'\in B$ it holds that $f(b)=f(b\wedge b')+f(b\wedge\neg b')$
and $\mathbf 1_{K(b)}=\mathbf 1_{K(b\wedge b')}+\mathbf 1_{K(b\wedge \neg b')}$
we can express the two sums as
$\sum_j \alpha_j'\mathbf 1_{K(b_j')}=\sum_j \beta_j'\mathbf 1_{K(b_j')}$ (where
$b_i'\wedge b_j'=0$ for every $i\neq j$), with
$\sum_i \alpha_i f(b_i)=\sum_j \alpha_j' f(b_j')$ and
$\sum_i \beta_i f(b_i)=\sum_j \beta_j' f(b_j')$.
As $b_i'\wedge b_j'=0$ for every $i\neq j$, for $x\in K(b_j')$ 
it holds that  $\alpha_j'=X(x)=\beta_j'$. Hence $\alpha_j'=\beta_j'$ for every $j$.
Thus  $\sum_i \alpha_i f(b_i)=\sum_i \beta_i f(b_i)$.

Note that $X\mapsto [f,X]$ is indeed continuous. Thus ${\rm ba}(B)$ is
the topological dual of $V(B)$. 
\myqed\end{proof}
\begin{lemma}
\label{lem:densv}
The vector space $V(B)$ is dense in $C(S(B))$ (with the uniform norm).
\end{lemma}
\begin{proof}
Let $f\in C(S(B))$ and let $\epsilon>0$. 
For $z\in f(S(B))$ let $U_z$ be the preimage by $f$ of the open ball $B_{\epsilon/2}(z)$ of $\bbbr$ centered in $z$.
As $f$ is continuous, $U_z$ is a open set of $S(B)$. As $\{K(b): b\in B\}$ is a basis of the topology of $S(B)$, 
$U_z$ can be expressed as a union $\bigcup_{b\in\mathcal F(U_z)}K(b)$. It follows that
$\bigcup_{z\in f(S(B))}\bigcup_{b\in\mathcal F(U_z)}K(b)$ is a covering of $S(B)$ by open sets. As $S(B)$ is compact,
there exists a finite subset $\mathcal F$ of $\bigcup_{z\in f(S(B))}\mathcal F(U_z)$ that covers $S(B)$.
Moreover, as for every $b,b'\in B$ it holds that $K(b)\cap K(b')=K(b\wedge b')$ and $K(b)\setminus K(b')=K(b\wedge\neg b')$
it follows that we can assume that there exists a finite family $\mathcal F'$ such that
$S(B)$ is covered by open sets $K(b)$ (for $b\in \mathcal F'$) and such that for every $b\in \mathcal F'$ there exists
$b'\in\mathcal F$ such that $K(b)\subseteq K(b')$. In particular, it follows that for every $b\in\mathcal F'$, $f(K(b))$
is included in an open ball of radius $\epsilon/2$ of $\bbbr$. For each $b\in\mathcal F'$ choose a point $x_b\in S(B)$ such 
that $b\in x_b$. Now define
$$
\hat f=\sum_{b\in\mathcal F'}f(x_b)\mathbf 1_{K(b)}
$$   
Let $x\in S(B)$. Then there exists $b\in \mathcal F'$ such that $x\in K(b)$. Thus
$$|f(x)-\hat f(x)|=|f(x)-f(x_b)|<\epsilon.$$
Hence $\|f-\hat{f}\|_\infty<\epsilon$. 
\myqed\end{proof}

\begin{lemma}
\label{lem:barca}
Let $B$ be a Boolean algebra, let ${\rm ba}(B)$ be the Banach space of bounded additive real-valued functions 
equipped with the norm $\|f\|=\sup_{b\in B} f(b)-\inf_{b\in B} f(b)$, let $S(B)$ be the Stone space associated 
to $B$ by the Stone representation theorem, 
and let ${\rm rca}(S(B))$ 
be the Banach space of the regular countably additive measure on $S(B)$ equipped with
the total variation norm.

Then the mapping $C_K: {\rm rca}(S(B))\rightarrow {\rm ba}(B)$ defined by $C_K(\mu)=\mu\circ K$ is an isometric isomorphism.
In other words, $C_K$ is defined by
$$C_K(\mu)(b)=\mu(\{x\in S(B):\ b\in x\})$$ (considering that the points of $S(B)$ are the ultrafilters on $B$).
\end{lemma}
\begin{proof}

According to Lemma~\ref{lem:baV}, the Banach space ${\rm ba}(B)$ is the topological dual of $V(B)$ 
and as $V(B)$ is dense in $C(S(B))$ (according to Lemma~\ref{lem:densv})
we deduce that ${\rm ba}(B)$ can 
be identified with the continuous dual of $C(S(B))$.
By Riesz representation theorem, the topological dual of $C(S(B))$ is the space ${\rm rca}(S(B))$ of
regular countably additive measures on $S(B)$.
From these observations follows the equivalence of ${\rm ba}(B)$ and
${\rm rca}(S(B))$.

This equivalence is easily made explicit, leading to the conclusion that the mapping
$C_K: {\rm rca}(S(B))\rightarrow {\rm ba}(B)$ defined by 
$C_K(\mu)=\mu\circ K$ is an isometric isomorphism.
\myqed\end{proof}

Note also that, similarly, the 
restriction of $C_K$ to the space ${\rm Pr}(S(B))$ of all (regular) probability measures on $S(B)$ 
is an isometric isomorphism of ${\rm Pr}(S(B))$ and the 
subset ${\rm ba}_1(B)$ of ${\rm ba}(B)$ of all non-negative additive functions
$f$ on $B$ such that $f(1)=1$.

Recall that given a measurable function $f:X\rightarrow Y$ (where $X$ and $Y$ are measurable spaces), the
{\em pushforward}\label{def:pushforward} $f_*(\mu)$ of a measure $\mu$ on $X$ is the measure on $Y$ defined by
$f_*(\mu)(A)=\mu(f^{-1}(A))$ (for every measurable set $A$ of $Y$).  
Note that if $f$ is a continuous function and if $\mu$ is a regular measure on $X$,
then the pushforward measure $f_*(\mu)$ is a regular measure on $Y$.
By similarity with the definition of $\Omega(f):\Omega(Y)\rightarrow\Omega(X)$
(see above definition) we denote by $\Omega_*(f)$ the mapping
from ${\rm rca}(X)$ to ${\rm rca}(Y)$ defined by
$(\Omega_*(f))(\mu)=f_*(\mu)$.

All the functors defined above are consistent in the sense that
if $h:A\rightarrow B$ is a homomorphism and $f\in{\rm ba}(B)$
then
$$\Omega_*(S(h))(\mu_f)\circ K_A=f\circ h.$$ 

A standard notion of convergence in ${\rm rca}(S(B))$ 
(as the continuous dual of $C(S(B))$) is the weak $*$-convergence: 
a sequence $(\mu_n)_{n\in\bbbn}$ of measures
is convergent if, for every $f\in C(S(B))$ the sequence $\int f(x)\,{\rm d}\mu_n(x)$ is convergent.
Thanks to the density of $V(B)$
 this convergence translates as pointwise convergence in ${\rm ba}(B)$ as follows:
a sequence $(g_n)_{n\in\bbbn}$ of functions in ${\rm ba}(B)$ is convergent if, for every $b\in B$ the sequence $(g_n(b))_{n\in\bbbn}$ is convergent.
As ${\rm rca}(S(B))$ is complete, so is ${\rm rca}(B)$. Moreover, it is easily checked that 
${\rm ba}_1(B)$ is closed in ${\rm ba}(B)$.

In a more concise way, we can write, for a sequence $(f_n)_{n\in\bbbn}$ of functions in ${\rm ba}(B)$ and
for the corresponding sequence $(\mu_{f_n})_{n\in\bbbn}$ of regular measures on $S(B)$:
$$
f_n\rightarrow f\text{ pointwise}\qquad\iff\qquad \mu_{f_n}\Rightarrow\mu_f.
$$
 
We now apply this classical machinery to structures and models.

\subsection{Basics of Model Theory and Lindenbaum--Tarski Algebras}
We denote by ${\mathcal B}({\rm FO}(\sig))$ the equivalence classes of ${\rm FO}(\sig)$ defined by logical equivalence.
The (class of) unsatisfiable formulas
(resp.\ of tautologies) will be designated by $0$ (resp.\ $1$). Then, ${\mathcal B}({\rm FO}(\sig))$ gets a natural structure
of Boolean algebra (with minimum $0$, maximum $1$, infimum $\wedge$, supremum $\vee$, and complement $\neg$).
This algebra is called the {\em Lindenbaum--Tarski algebra} of ${\rm FO}(\sig)$.
Notice that all the Boolean algebras ${\rm FO}(\sig)$ for countable $\sig$ are isomorphic, as there exists 
only one countable atomless Boolean algebra up to isomorphism (see \cite{Hodges1993}).

For an integer $p\geq 1$, the fragment ${\rm FO}_p(\sig)$ of ${\rm FO}(\sig)$ contains
first-order formulas $\phi$ such that ${\rm Fv}(\phi)\subseteq\{x_1,\dots,x_p\}$. 
The fragment ${\rm FO}_0(\sig)$ of ${\rm FO}(\sig)$ contains 
first-order 
formulas without free variables (that is {\em sentences}).
 
We check that the permutation group $\Sym{p}$ on $[p]$ acts on ${\rm FO}_p(\sig)$ by 
$\sigma\cdot\phi=\phi(x_{\sigma(1)},\dots,x_{\sigma(p)})$ and that each permutation indeed define
an automorphism of $\mathcal B({\rm FO}_p(\sig))$.
Similarly, the group $\Sym{\omega}$ of permutations on $\bbbn$ 
acts on ${\rm FO}(\sig)$
and $\mathcal B({\rm FO}(\sig))$.
Note that ${\rm FO}_0(\sig)\subseteq\dots\subseteq{\rm FO}_p(\sig)\subseteq{\rm FO}_{p+1}(\sig)\subseteq\dots\subseteq{\rm FO}(\sig)$.
Conversely, let ${\rm rank}(\phi)=\max\{i: x_i\in{\rm Fv}(\phi)\}$. Then we have a natural projection
$\pi_p:{\rm FO}(\sig)\rightarrow{\rm FO}_p(\sig)$ defined by
$$
\pi_p(\phi)=\begin{cases}
\phi&\text{if }{\rm rank}(\phi)\leq p\\
\exists x_{p+1}\,\exists x_{p+2}\,\dots\,\exists x_{{\rm rank}(\phi)}\,\phi&\text{otherwise}
\end{cases}
$$

An {\em elementary class} (or {\em axiomatizable class}) 
$\mathcal C$ of $\sig$-structures is 
a class consisting of all $\sig$-structures satisfying a fixed consistent first-order theory $T_{\mathcal C}$.
Denoting by $\mathcal I_{T_{\mathcal C}}$ the ideal of all first-order formulas in $\mathcal L$ that are 
provably false from axioms in $T_{\mathcal C}$, 
The Lindenbaum--Tarski algebra $\mathcal B({\rm FO}(\sig),T_{\mathcal C})$ 
associated to the theory $T_{\mathcal C}$ of $\mathcal C$ is
the quotient Boolean algebra
$\mathcal B({\rm FO}(\sig),T_{\mathcal C})=\mathcal B({\rm FO}(\sig))/\mathcal I_{T_{\mathcal C}}$.
As a set, $\mathcal B({\rm FO}(\sig),T_{\mathcal C})$ is simply the quotient of 
${\rm FO}(\sig)$ by logical equivalence modulo $T_{\mathcal C}$.

As we consider countable languages, $T_{\mathcal C}$ is at most countable and
it is easily checked that $S(\mathcal B({\rm FO}(\sig),T_{\mathcal C}))$ 
is homeomorphic to the compact subspace of 
$S(\mathcal B({\rm FO}(\sig)))$ defined as $\{T\in S(\mathcal B({\rm FO}(\sig))): T\supseteq T_{\mathcal C}\}$. 
Note that, for instance, $S(\mathcal B({\rm FO}_0(\sig),T_{\mathcal C}))$
is a clopen set of $S(\mathcal B({\rm FO}_0(\sig)))$ if and only if 
$\mathcal C$ is {\em finitely axiomatizable} (or a {\em basic} elementary class), that is
 if $T_{\mathcal C}$ can be chosen to be a single sentence. These explicit correspondences 
 are particularly useful to our setting.

\subsection{Stone Pairing Again}
 We add a few comments to Definition~\ref{def:bracket2}. Note first that this definition
 is consistent in the sense that for every  {\rps} $\mathbf A$ and for
 every formula $\phi\in{\rm FO}(\sig)$ with $p$ free
 variables can be considered as a formula with $q\geq p$ free variables with
 $q-p$ unused variables, we have
 
 $$\int_{A^q}1_{\Omega_\phi(\mathbf A)}(x)\,{\rm d}\nu_{\mathbf A}^q(x)
 =\int_{A^p}1_{\Omega_\phi(\mathbf A)}(x)\,{\rm d}\nu_{\mathbf A}^p(x).
 $$
 
It is immediate that for every formula $\phi$ it holds that 
$\langle\neg\phi,\mathbf A\rangle=1-\langle\phi,\mathbf A\rangle$. Moreover, 
if $\phi_1,\dots,\phi_n$ are formulas, then by de Moivre's formula, the following equation holds:
 
 $$\langle\bigvee_{i=1}^n\phi_i,\mathbf A\rangle=\sum_{k=1}^n
 (-1)^{k+1}\biggl(\sum_{1\leq i_1<\dots<i_k\leq n}\langle
 \bigwedge_{j=1}^k\phi_{i_j},\mathbf A\rangle\biggr).$$
 
In particular, if $\phi_1,\dots,\phi_k$ are {\em mutually exclusive}
 (meaning that $\phi_i\wedge\phi_j=0$)
 then the following equation holds:
$$\langle\bigvee_{i=1}^k\phi_i,\mathbf A\rangle=\sum_{i=1}^k\langle \phi_i,\mathbf A\rangle.$$
 
It follows that for every fixed {\rps} $\mathbf A$, the mapping $\phi\mapsto\langle\phi,\mathbf A\rangle$ is additive
(i.e.\ $\langle\,\cdot\,,\mathbf A\rangle\in{\rm ba}(\mathcal B({\rm FO}(\sig)))$):
$$
\phi_1\wedge\phi_2=0\quad\Longrightarrow\quad \langle\phi_1\vee\phi_2,\mathbf A\rangle=\langle\phi_1,\mathbf A\rangle
+\langle\phi_2,\mathbf A\rangle.
$$ 
 
The Stone pairing is antimonotone: 
\begin{mytinted}
Let $\phi,\psi\in{\rm FO}(\sig)$.
 For every {\rps} $\mathbf A$ the following implication holds:
$$\phi\entails\psi\quad\Longrightarrow\quad\langle\phi,G\rangle\geq\langle\psi,G\rangle.$$
\end{mytinted}
However, even if $\phi$ and $\psi$ are sentences and $\langle\phi,\,\cdot\,\rangle\geq \langle\psi,\,\cdot\,\rangle$
on finite $\sig$-structures, this does not imply in general that $\phi\entails\psi$: let $\theta$ be a 
sentence with only infinite models and let $\phi$ be a sentence with only finite models.
On finite $\sig$-structures it holds that $\langle\phi\vee\theta,\,\cdot\,\rangle=\langle\phi,\,\cdot\,\rangle$
 although $\phi\vee\theta\nvdash\phi$ (as witnessed by an infinite model of $\theta$).  

Nevertheless, inequalities between Stone pairing that are valid for finite $\sig$-structures will 
of course still hold at the limit.
For instance, for $\phi_1,\phi_2\in{\rm FO}_1(\sig)$, for $\zeta\in{\rm FO}_2(\sig)$, and for $a,b\in\bbbn$ define
the first-order sentence $B(a,b,\phi_1,\phi_2,\zeta)$ expressing that for every vertex $x$ such that $\phi_1(x)$ holds
there exist at least $a$ vertices $y$ such that $\phi_2(y)\wedge\zeta(x,y)$ holds and that for every vertex $y$ such
that $\phi_2(x)$ holds there exist at most $b$ vertices $x$ such that $\phi_1(x)\wedge\zeta(x,y)$ holds.
Then it is easily checked that for every finite $\sig$-structure $\mathbf A$ the following implication holds:
$$
\mathbf A\models B(a,b,\phi_1,\phi_2,\zeta)\quad\Longrightarrow\quad a\langle\phi_1,\mathbf A\rangle\leq b\langle\phi_2,\mathbf A\rangle.
$$  
For example, if a finite directed graph is such that every arc connects a vertex with out-degree $2$ to a vertex
with in-degree $1$, it is clear that the probability that a random vertex has out-degree $2$ is half the probability that
a random vertex has in-degree $1$. 
 
 Now we come to important twist and the basic of our approach.
The Stone pairing $\langle\,\cdot\,,\,\cdot\,\rangle$ can be considered from both sides: On the right side
the functions of type $\langle\phi,\,\cdot\,\rangle$ are a generalization of the homomorphism density functions \cite{borgchaylovas06:_count_graph_homom}:
$$
t(F,G)=\frac{|{\rm hom}(F,G)|}{|G|^{|F|}}
$$
(these functions correspond to $\langle\phi,G\rangle$ for Boolean conjunctive queries $\phi$ and a graph $G$).
Also the density function used in \cite{Benjamini2001} to measure the probability that the ball of radius $r$
rooted at a random vertex as a given isomorphism type may be expressed as a function $\langle\phi,\,\cdot\,\rangle$.
Note again that we follow here, in a sense, a dual approach:  we consider for fixed $\mathbf A$ the 
function $\langle\,\cdot\,,\mathbf A\rangle$, which
is an additive function on $\mathcal B({\rm FO}(\sig))$ with the following properties:
\begin{itemize}
  \item $\langle\,\cdot\,,\mathbf A\rangle\geq 0$ and $\langle 1,\mathbf A\rangle=1$;
  \item $\langle\sigma\cdot\phi,\mathbf A\rangle=\langle\phi,\mathbf A\rangle$ for every $\sigma\in \Sym{\omega}$;
  \item if ${\rm Fv}(\phi)\cap{\rm Fv}(\psi)=\emptyset$, then 
  $\langle\phi\wedge\psi,\mathbf A\rangle=\langle\phi,\mathbf A\rangle\ \langle\psi,\mathbf A\rangle$.
\end{itemize}

Thus $\langle\,\cdot\,,\mathbf A\rangle$ is, for a given $\mathbf A$, an operator on the 
class of first-order formulas.

We now can apply Lemma~\ref{lem:barca} to derive a representation by means of a 
regular measure on a Stone space.
The fine structure and interplay of additive functions, Boolean functions, and dual spaces can be 
used effectively if we consider finite $\sig$-structures as probability spaces as we did
when we considered finite $\sig$-structures as a particular case of Borel models.

The following two theorems generalize Theorems~\ref{thm:rep1} and~\ref{thm:rep2} mentioned in Section~\ref{sec:def}.
\begin{mytinted}
\begin{theorem}
\label{thm:genmod}
Let $\sig$ be a signature, let $\mathcal B({\rm FO}(\sig))$ be the Lindenbaum--Tarski algebra
of ${\rm FO}(\sig)$, let $S(\mathcal B({\rm FO}(\sig)))$ be the associated Stone space,
and let ${\rm rca}(S(\mathcal B({\rm FO}(\sig))))$
be the Banach space of the regular countably additive measures on $S(\mathcal B({\rm FO}(\sig)))$.
Then:
\begin{enumerate}
  \item There is a mapping from the class of $\sig$-{\rps}
  to ${\rm rca}(S(\mathcal B({\rm FO}(\sig))))$, which maps a {\rps} $\mathbf A$
to the unique regular measure $\mu_{\mathbf A}$ such that for every $\phi\in{\rm FO}(\sig)$
the following equation holds:
$$
\langle\phi,\mathbf A\rangle=\int_{S(\mathcal B({\rm FO}(\sig)))}\mathbf 1_{K(\phi)}\, {\rm d}\mu_{\mathbf A},
$$
where $\mathbf 1_{K(\phi)}$ is the indicator function of $K(\phi)$ in $S(\mathcal B({\rm FO}(\sig)))$.
Moreover, this mapping is injective of finite $\sig$-structures.
\item A sequence $(\mathbf A_n)_{n\in\bbbn}$ of finite $\sig$-structures is {\rm FO}-convergent if and only if
the sequence $(\mu_{\mathbf A_n})_{n\in\bbbn}$ is weakly converging in ${\rm rca}(S(\mathcal B({\rm FO}(\sig))))$;
\item If $(\mathbf A_n)_{n\in\bbbn}$ is an {\rm FO}-convergent sequence of finite $\sig$-structures
then the weak limit $\mu$
of $(\mu_{\mathbf A_n})_{n\in\bbbn}$ is such that for every $\phi\in {\rm FO}(\sig)$ the following equation holds:
$$
\lim_{n\rightarrow\infty}\langle\phi,\mathbf A_n\rangle=\int_{S(\mathcal B({\rm FO}(\sig)))}\mathbf 1_{K(\phi)}\, {\rm d}\mu.
$$
\end{enumerate}
\end{theorem}
\end{mytinted}
\begin{proof}
The proof follows from Lemma~\ref{lem:barca}, considering the additive functions
$\langle\,\cdot\,,\mathbf A\rangle$. 

Let $\mathbf A$ be a finite $\sig$-structure.
As $\mu_{\mathbf A}$ allows one to recover the complete theory of $\mathbf A$ and as $\mathbf A$ is finite,
the mapping $\mathbf A\mapsto \mu_{\mathbf A}$ is injective.
\myqed\end{proof}

It is important to  consider fragments of ${\rm FO}(\sig)$ to define a weaker notion
of convergence. This allows us to capture limits of dense graphs too.

\begin{mytinted}
\begin{definition}[$X$-convergence]
Let $X$ be a fragment of ${\rm FO}(\sig)$.
A sequence $(\mathbf A_n)_{n\in\bbbn}$ of finite $\sig$-structures is {\em $X$-convergent} if
$\langle\phi,\mathbf A_n\rangle$ is convergent for every  $\phi\in X$.
\end{definition}
\end{mytinted}

In the particular case that $X$ is a Boolean sub-algebra of $\mathcal B({\rm FO}(\sig))$ we can apply all above methods and
in this context we can extend Theorem~\ref{thm:genmod}.

\begin{mytinted}
\begin{theorem}
\label{thm:genmodx}
Let $\sig$ be a signature, and let $X$ be a fragment of ${\rm FO}(\sig)$ defining
a Boolean algebra $\mathcal B(X)\subseteq\mathcal B({\rm FO}(\sig))$.
Let $S(\mathcal B(X))$ be the associated Stone space, 
and let ${\rm rca}(S(\mathcal B(X)))$ 
be the Banach space of the regular countably additive measure on $S(\mathcal B(X))$.
Then:
\begin{enumerate}
  \item The canonical injection $\iota^X:\mathcal B(X)\rightarrow\mathcal B({\rm FO}(\sig))$
  defines by duality a continuous projection $p^X:S(\mathcal B({\rm FO}(\sig)))\rightarrow S(\mathcal B(X))$;
  The pushforward $p^X_*\mu_{\mathbf A}$ 
  of the measure $\mu_{\mathbf A}$ associated to a {\rps} $\mathbf A$ (see Theorem~\ref{thm:genmod}) 
  is the unique regular measure on $S(\mathcal B(X))$ such that:
$$
\langle\phi,\mathbf A\rangle=\int_{S(\mathcal B(X))}\mathbf 1_{K(\phi)}\, {\rm d}p^X_*\mu_{\mathbf A},
$$
where $\mathbf 1_{K(\phi)}$ is the indicator function of $K(\phi)$ in $S(\mathcal B(X))$.
\item A sequence $(\mathbf A_n)_{n\in\bbbn}$ of finite $\sig$-structures is $X$-convergent if and only if
the sequence $(p^X_*\mu_{\mathbf A_n})_{n\in\bbbn}$ is weakly converging in ${\rm rca}(S(\mathcal B(X)))$;
\item If $(\mathbf A_n)_{n\in\bbbn}$ is an $X$-convergent sequence of finite $\sig$-structures 
then the weak limit $\mu$
of $(p^X_*\mu_{\mathbf A_n})_{n\in\bbbn}$ is such that for every $\phi\in X$ the following equation holds:
$$
\lim_{n\rightarrow\infty}\langle\phi,\mathbf A_n\rangle=\int_{S(\mathcal B(X))}\mathbf 1_{K(\phi)}\, {\rm d}\mu.
$$
\end{enumerate}
\end{theorem}
\end{mytinted}
\begin{proof}
If $X$ is closed under conjunction, disjunction and negation, thus
defining a Boolean algebra $\mathcal B(X)$, then the inclusion of $X$ in ${\rm FO}(\sig)$
translates as a canonical injection $\iota$ from $\mathcal B(X)$ to $\mathcal B({\rm FO}(\sig))$.
By Stone duality, the injection $\iota$ corresponds to a continuous projection
$p:S(\mathcal B({\rm FO}(\sig)))\rightarrow S(\mathcal B(X))$.
As every measurable function, this continuous projection also transports measures by
pushforward: the projection $p$ transfers the 
measure $\mu$ on $S(\mathcal B({\rm FO}(\sig)))$ to  $S(\mathcal B(X))$ as the pushforward measure 
$p_*\mu$ defined by the identity $p_*\mu(Y)=\mu(p^{-1}(Y))$, which holds for every measurable subset 
$Y$ of $S(\mathcal B(X))$. 

The proof follows from Lemma~\ref{lem:barca}, considering the additive functions
$\langle\,\cdot\,,\mathbf A\rangle$. 
\myqed\end{proof}

We can also consider a notion of convergence restricted to $\sig$-structures satisfying a fixed axiom.
\begin{mytinted}
\begin{theorem}
\label{thm:genmodxc}
Let $\sig$ be a signature, and let $X$ be a fragment of ${\rm FO}(\sig)$ defining
a Boolean algebra $\mathcal B(X)\subseteq\mathcal B({\rm FO}(\sig))$.
Let $S(\mathcal B(X))$ be the associated Stone space, 
and let ${\rm rca}(S(\mathcal B(X)))$ 
be the Banach space of the regular countably additive measure on $S(\mathcal B(X))$.

Let $\mathcal C$ be a basic elementary class defined by a single 
axiom $\Psi\in X\cap {\rm FO}_0$, and 
let $\mathcal I_{\Psi}$ be the principal ideal of $\mathcal B(X)$ 
generated by $\neg\Psi$.

Then:
\begin{enumerate}
  \item The Boolean algebra obtained by taking the quotient
  of $X$ equivalence modulo $\Psi$ is the quotient Boolean algebra
  $\mathcal B(X,\Psi)=\mathcal B(X)/\mathcal I_{\Psi}$.
  Then $S(\mathcal B(X,\Psi))$ is homeomorphic to the clopen subspace 
  $K(\Psi)$ of $S(\mathcal B(X))$.
  
  If $\mathbf A\in\mathcal C$ is a finite $\sig$-structure
  then the support of the measure $p^X_*\mu_{\mathbf A}$ associated to $\mathbf A$ 
  (see Theorem~\ref{thm:genmodx}) is included
  in $K(\Psi)$ and for every $\phi\in X$ the following equation holds:
$$
\langle\phi,\mathbf A\rangle=\int_{K(\Psi)}\mathbf 1_{K(\phi)}\, {\rm d}p^X_*\mu_{\mathbf A}.
$$
\item A sequence $(\mathbf A_n)_{n\in\bbbn}$ of finite $\sig$-structures of $\mathcal C$ is $X$-convergent
 if and only if
the sequence $(p^X_*\mu_{\mathbf A_n})_{n\in\bbbn}$ is weakly converging in ${\rm rca}(S(\mathcal B(X,\Psi)))$;
\item If $(\mathbf A_n)_{n\in\bbbn}$ is an $X$-convergent sequence of finite $\sig$-structures in $\mathcal C$
then the weak limit $\mu$
of $(p^X_*\mu_{\mathbf A_n})_{n\in\bbbn}$ is such that for every $\phi\in X$ the following equation holds:
$$
\lim_{n\rightarrow\infty}\langle\phi,\mathbf A_n\rangle=\int_{K(\Psi)}\mathbf 1_{K(\phi)}\, {\rm d}\mu.
$$
\end{enumerate}
\end{theorem}
\end{mytinted}
\begin{proof}
 The quotient algebra $\mathcal B(X,\Psi)=\mathcal B(X)/\mathcal I_{\Psi}$ 
 is isomorphic to the sub-Boolean algebra $\mathcal B'$ of
  $\mathcal B$ of all (equivalence classes of) formulas $\phi\wedge\Psi$ for $\phi\in X$.
  To this isomorphism corresponds by duality
  the identification of $S(\mathcal B(X,\Psi))$ with the clopen subspace 
  $K(\Psi)$ of $S(\mathcal B(X))$.
\myqed\end{proof}

The situation expressed by these theorems is summarized in the following diagram.
$$
\xy
\xymatrix@C=10mm{
\mathcal B({\rm FO}(\sig))\ar@{<->}[d]&&\mathcal B(X)\ar@{<->}[d]\ar[ll]_(.45){\text{canonical injection}}&
\mathcal B'\ar[l]_(.4){\text{inclusion}}\ar@{<->}[d]\ar@{<->}[rr]^(.45){isomorphism}&&\mathcal B(X,\Psi)\ar@{<->}[d]\\
S(\mathcal B({\rm FO}(\sig)))\ar[rr]^(.55){\text{projection } p^X}&&S(\mathcal B(X))&
K(\Psi)\ar[l]_(.45){\text{inclusion}}&&S(\mathcal B(X,\Psi))\ar@{<->}[ll]_(.55){\text{homeomorphism}}\\
\mu\ar[rr]^{\text{pushforward}}&&p^X_*\mu\ar[r]^{\text{restriction}}&p^X_*\mu
}
\endxy
$$
The essence of our approach is that we follow a dual path: we view a graph $G$ as an operator
on first-order formulas through Stone pairing $\langle\,\cdot\,,G\rangle$.
%
\subsection{Limit of Measures Associated to Finite Structures}
\label{sec:limfin}
We consider a signature $\sig$ and fragment ${\rm FO}_p$ of ${\rm FO}(\sig)$.
Let $(\mathbf A_n)_{n\in\bbbn}$ be an $X$-convergent sequence of $\sig$-structures, 
let $\mu_{\mathbf A_n}$ be the measure on $S(\mathcal B(X))$ associated to $\mathbf A_n$, and 
let $\mu$ be the weak limit of $\mu_{\mathbf A_n}$.

\begin{fact}
As we consider countable languages only,  $S(\mathcal B({\rm FO}_p))$ is a Radon space
and thus for every (Borel) probability measure $\mu$ on $S(\mathcal B({\rm FO}_p))$, any measurable set outside
the support of $\mu$ has zero $\mu$-measure.
\end{fact}

\begin{definition}
\label{def:pure}
Let
$\pi$ be the natural projection
$S(\mathcal B({\rm FO}_p))\rightarrow S(\mathcal B({\rm FO}_0))$.

A measure $\mu$ on $S(\mathcal B({\rm FO}_p))$ is {\em pure} if
$|\pi({\rm Supp}(\mu))|=1$. The unique element $T$
of $\pi({\rm Supp}(\mu))$ is then called the {\em complete theory} of $\mu$.
\end{definition}

\begin{mytinted}
\begin{remark}
\label{rem:pure}
Consider ${\rm FO}_p$ or ${\rm FO}$ convergence.
Every measure $\mu$ that is the weak limit of some sequence of measures associated to finite
structures is pure and its complete theory has the finite model property.
\end{remark}
\end{mytinted}
Indeed, if a sequence $(\mathbf A_n)_{n\in\bbbn}$ of finite structures is ${\rm FO}_p$ or ${\rm FO}$-convergent it is in particular ${\rm FO}_0$-convergent. It follows that  if $\mu_{\mathbf A_n}$ weakly converges to $\mu$ then 
$\pi(\mu)$ is concentrated on the complete theory $T$ of the elementary limit of $(\mathbf A_n)_{n\in\bbbn}$ (thus $\mu$ is pure) and as $T$ is the complete theory of the elementary limit of finite structures it  has the finite model property.

\begin{definition}
For $T\in S(\mathcal B({\rm FO}_p),\psi,\phi\in {\rm FO}_p$, and $\beta\in{\rm FO}_{2p}$ define
\begin{align*}
{\rm deg}_\psi^{\beta+}(T)&=\begin{cases}
k&\text{if }T\ni(\exists^{=k}(y_1,\dots,y_p)\,\beta(x_1,\dots,x_p,y_1,\dots,y_p)\wedge\psi(y_1,\dots,y_p))\\
\infty&\text{ otherwise.}\\
\end{cases}\\
{\rm deg}_\phi^{\beta-}(T)&=\begin{cases}
k&\text{if }T\ni(\exists^{=k}(x_1,\dots,x_p)\,\phi(x_1,\dots,x_p)\wedge\beta(x_1,\dots,x_p,y_1,\dots,y_p))\\
\infty&\text{ otherwise.}
\end{cases}
\end{align*}
\end{definition}

Denote by $\xi_k$ the formula $\exists^{=k}(\mathbf y)\,\beta(\mathbf x,\mathbf y)\wedge\psi(\mathbf y)$  (where $\mathbf x=(x_1,\dots,x_p)$ and $\mathbf y=(y_1,\dots,y_p)$)
 then 
 for every finite structure $\mathbf A$ 
 it holds that $\mathbf A\models (\forall\mathbf x)\ \neg\xi_k(\mathbf x,\mathbf y)$ if $k>|A|^p$. Thus ${\rm deg}_\psi^{\beta^+}=\sum_{k=1}^{|A|^p}	 \mathbf 1_{K(\xi_k)}$ and 
\begin{align*}
\int_{K(\phi)}{\rm deg}_\psi^{\beta^+}(T)\,{\rm d}\mu_{\mathbf A}(T)
	&=	\adjustlimits\int_{K(\phi)} \sum_{k=1}^{|A|^p}\mathbf 1_{K(\xi_k)}\,{\rm d}\mu_{\mathbf A}(T)\\
		&=\sum_{k=1}^{|A|^p}\langle\xi_k\wedge\phi,\mathbf A\rangle\\
		&=\frac{1}{|A|^p}\sum_{\mathbf v\in \phi(\mathbf A)}\Bigl|\{\mathbf w\in\psi(\mathbf A): \mathbf A\models\beta(\mathbf v,\mathbf w)\}\Bigr|\\
		&=\frac{1}{|A|^p}\Bigl|\{((\mathbf v,\mathbf w)\in \phi(\mathbf A)\times \psi(\mathbf A):\ \mathbf A\models\beta(\mathbf v,\mathbf w)\}\Bigr|
		\intertext{and, similarly we get}
\int_{K(\psi)}{\rm deg}_\phi^{\beta^-}(T)\,{\rm d}\mu_{\mathbf A}(T)
&=\frac{1}{|A|^p}\sum_{\mathbf w\in \psi(\mathbf A)}\Bigl|\{\mathbf v\in\phi(\mathbf A): \mathbf A\models\beta(\mathbf v,\mathbf w)\}\Bigr|\\
		&=\frac{1}{|A|^p}\Bigl|\{((\mathbf v,\mathbf w)\in \phi(\mathbf A)\times \psi(\mathbf A):\ \mathbf A\models\beta(\mathbf v,\mathbf w)\}\Bigr|.
\end{align*}

Thus if $\mu$ is a measure associated to a finite structure then for every $\phi,\psi\in {\rm FO}_p$ the following equation holds:
$$\int_{K(\phi)}{\rm deg}_\psi^{\beta^+}(T)\,{\rm d}\mu(T)=\int_{K(\psi)}{\rm deg}_\phi^{\beta-}(T)\,{\rm d}\mu(T).$$

Hence for every measure $\mu$ that is the weak limit of some sequence of measures associated to finite
structures the following property holds:

\begin{mytinted}
\begin{center}
	{\bf General Finitary Mass Transport Principle} (GFMTP)
\end{center}
For every $\phi,\psi\in{\rm FO}_p$, every $\beta\in{\rm FO}_{2p}$, and all integers $a,b$ that are such that
\begin{align*}
\forall T\in K(\phi)\quad&\ {\rm deg}_\psi^{\beta^+}(T)\geq a\\
\forall T\in K(\psi)\quad&\ {\rm deg}_\phi^{\beta-}(T)\leq b
\end{align*}
the following inequality holds:
$$a\,\mu(K(\phi))\leq b\,\mu(K(\psi)).$$
\end{mytinted}
Of course, similar statement holds as well for the projection 
of $\mu$ on $S(\mathcal B({\rm FO}_q))$ for $q<p$.
In the case of digraphs, when $p=1$ and $\beta(x_1,x_2)$ is existence
of an arc from $x_1$ to $x_2$, we shall write 
${\rm deg}_\psi^+$ and ${\rm deg}_\phi^-$ instead of
${\rm deg}_\psi^{\beta+}$ and ${\rm deg}_\phi^{\beta-}$. 
(In the case of graphs, we have  ${\rm deg}_\psi^+={\rm deg}_\psi^-={\rm deg}_\psi$.)
Thus the following property holds.

\begin{mytinted}
\begin{center}
	{\bf Finitary Mass Transport Principle} (FMTP)
\end{center}
For every $\phi,\psi\in{\rm FO}_1$, and all integers $a,b$ that are such that
\begin{align*}
	\forall T\in K(\phi)\quad&\ {\rm deg}_\psi^+(T)\geq a\\
\forall T\in K(\psi)\quad&\ {\rm deg}_\phi^-(T)\leq b
\end{align*}
the following inequality holds:
$$a\,\mu(K(\phi))\leq b\,\mu(K(\psi)).$$
\end{mytinted}

GFMTP and FMTP will play a key role in the analysis of modeling limits.
\section{Convergence, Old and New}
\label{sec:oldnew}
As we have seen above, there are many possible notions of convergence for a sequence $(\mathbf A_n)_{n\in\bbbn}$ of  finite $\sig$-structures.
As we considered $\sig$-structures defined with a countable signature $\sig$, the Boolean algebra
$\mathcal B({\rm FO}(\sig))$ is countable. It follows that 
the Stone space $S(\mathcal B({\rm FO}(\sig)))$ is a  Polish space, and thus
(with the Borel $\sigma$-algebra) it is a standard Borel space. Hence every probability distribution  turns
$S(\mathcal B({\rm FO}(\sig)))$ into a standard probability space.
However, the fine structure of $S(\mathcal B({\rm FO}(\sig)))$ is complex 
and we have no simple description of this space.

${\rm FO}$-convergence is of course the most restrictive notion of convergence and it seems (at least at the first glance)
that this is perhaps too much to ask, as we may encounter many particular difficulties and specific cases.
But we shall exhibit later classes for which ${\rm FO}$-convergence is captured --- for special basic elementary
classes of structures --- by $X$-convergence for a small fragment $X$ of ${\rm FO}$.

At this time it is natural to ask whether one can consider 
fragments
whose corresponding Boolean algebras are not sub-Boolean algebras of $\mathcal B({\rm FO}(\lambda))$
and still have a description of the limit 
of a converging sequence as a probability measure on a nice
measurable space.
There is obviously a case where this is possible: when the 
convergence of $\langle\phi,\mathbf A_n\rangle$ for every $\phi$ in 
a fragment $X$ implies the convergence of $\langle\psi,\mathbf A_n\rangle$ 
for every $\psi$ in the minimum Boolean algebra containing $X$.
We prove now that this is for instance the case when $X$ is a fragment
closed under conjunction.

For a Boolean algebra $B$ and a subset $X$ of $B$ we denote by 
$B[X]$ the Boolean sub-algebra of $B$ {\em generated by } $X$, that
is the sub-algebra of $B$ whose elements  can be expressed as a finite combination of elements of $X$, using the Boolean operations (in $B$).
We shall need the following preliminary lemma:
\begin{lemma}
\label{lem:meetbasis}
Let $B$ be a Boolean algebra and let $X\subseteq B$ be closed under $\wedge$ and such
that $X$ {\em generates} $B$ (i.e.\ such that $B[X]=B$).

Then $\{\mathbf 1_b:\ b\in X\}\cup\{\mathbf 1\}$ 
(where $\mathbf 1$ is the constant function with value $1$) 
includes a basis of the vector space $V(B)$ generated by the whole
set $\{\mathbf 1_b:\ b\in B\}$.
\end{lemma}
\begin{proof}
Let $b\in B$. As $X$ generates $B$ there exist
$b_1,\dots,b_k\in X$ and a Boolean function $F$
such that $b=F(b_1,\dots,b_k)$.
As $\mathbf 1_{x\wedge y}=\mathbf 1_x\,\mathbf 1_y$ and
$\mathbf 1_{\neg x}=\mathbf 1-\mathbf 1_x$ there exists a polynomial
$P_F$ such that $\mathbf 1_b=P_F(\mathbf 1_{b_1},\dots,\mathbf 1_{b_k})$.
For $I\subseteq [k]$, the monomial $\prod_{i\in I}\mathbf 1_{b_i}$
rewrites as $\mathbf 1_{b_I}$ where $b_I=\bigwedge_{i\in I}b_i$.
It follows that $\mathbf 1_b$ is a linear combination of the functions
$\mathbf 1_{b_I}$  ($I\subseteq [k]$) which belong to $X$ if $I\neq\emptyset$ (as $X$ is 
closed under $\wedge$ operation) and equals $\mathbf 1$, otherwise.
\myqed\end{proof}
 
\begin{proposition}
\label{lem:extBA}
Let $X$ be a fragment of ${\rm FO}(\sig)$ closed under (finite) 
conjunction --- thus defining a {\em meet semilattice} of $\mathcal B({\rm FO}(\sig))$ --- 
and let ${\mathcal B}(X)$ be the sub-Boolean algebra of $\mathcal B({\rm FO}(\sig))$ generated by $X$.
Let $\overline{X}$ be the fragment of ${\rm FO}(\sig)$ consisting of all formulas with equivalence
class in ${\mathcal B}(X)$.

Then $X$-convergence is equivalent to $\overline{X}$-convergence. 
\end{proposition}
\begin{proof}
Let $\Psi\in \overline{X}$. 
According to Lemma~\ref{lem:meetbasis}, there exist
$\phi_1,\dots,\phi_k\in X$ and $\alpha_0,\alpha_1,\dots,\alpha_k\in\bbbr$ such that

$$
\mathbf 1_\Psi=\alpha_0\mathbf 1+\sum_{i=1}^k\alpha_i\mathbf 1_{\phi_i}.
$$
Let $\mathbf A$ be a $\sig$-structure, let $\Omega=S({\mathcal B}(X))$ and let $\mu_{\mathbf A}\in{\rm rca}(\Omega)$ 
be the associated measure. Then
$$
\langle\Psi,\mathbf A\rangle=\int_{\Omega}\mathbf 1_\Psi\,{\rm d}\mu_{\mathbf A}
=\int_{\Omega}\bigl(\alpha_0\mathbf 1+\sum_{i=1}^k\alpha_i\mathbf 1_{\phi_i}\bigr)\,{\rm d}\mu_G
=\alpha_0+\sum_{i=1}^k\alpha_i\langle\phi_i,\mathbf A\rangle.
$$

It follows that if $(\mathbf A_n)_{n\in\bbbn}$ is an $X$-convergent sequence,
the sequence \linebreak $(\langle\psi,\mathbf A_n\rangle)_{n\in\bbbn}$ converges for every
$\psi\in\overline{X}$, that is $(\mathbf A_n)_{n\in\bbbn}$ is $\overline{X}$-convergent.
\myqed\end{proof}

Now we  demonstrate the expressive power of $X$-convergence by relating it to the main types of
convergence of graphs studied previously: 
\begin{enumerate}
  \item the notion of {\em dense graph limit} \cite{Borgs2005,Lov'asz2006};
  \item the notion of {\em bounded degree graph limit} \cite{Benjamini2001,Aldous2006};
  \item the notion of {\em elementary limit} derived from two important
  results in first-order logic, namely G\"odel's completeness theorem and the
  compactness theorem.
\end{enumerate}

These standard notions of graph limits, which have inspired this work, correspond
to special fragments of ${\rm FO}(\gamma)$, where $\gamma$ is the signature of graphs.
In the remainder of this section, we shall
only consider undirected graphs, thus we shall omit making mention of their signature in
the notations as well as the axioms defining the basic elementary class of undirected graphs.

\subsection{{\rm L}-convergence and ${\rm QF}$-convergence}
\label{sec:Lconv}
Recall that a sequence \linebreak $(G_n)_{n\in\bbbn}$ of graphs
 is {\em {\rm L}-convergent} if 
 $$t(F,G_n)=\frac{{\rm hom}(F,G_n)}{|G_n|^{|F|}}$$
 converges for every fixed (connected) graph $F$, where ${\rm hom}(F,G)$ denotes the number
 of homomorphisms of $F$ to $G$ \cite{Lov'asz2006,Borgs20081801,Borgs2012}.

It is a classical observation that homomorphisms between finite structures can be expressed by
Boolean conjunctive queries \cite{Chandra1977}.
We denote by ${\rm HOM}$ the fragment of ${\rm FO}$ consisting of formulas formed
by conjunction of atoms. For instance, the formula
$$(x_1\sim x_2)\wedge (x_2\sim x_3)\wedge (x_3\sim x_4)\wedge (x_4\sim x_5)\wedge (x_5\sim x_1)$$
belongs to ${\rm HOM}$ and it expresses that $(x_1,x_2,x_3,x_4,x_5)$ form a homomorphic image of $C_5$.
Generally, to a finite graph $F$ we associate the canonical formula $\phi_F\in{\rm HOM}$ defined by
$$\phi_F:=\bigwedge_{ij\in E(F)}(x_i\sim x_j).$$
Then, for every graph $G$ the following equation holds:

$$
\langle\phi_F,G\rangle=\frac{{\rm hom}(F,G)}{|G|^{|F|}}=t(F,G).
$$

Thus {\rm L}-convergence is equivalent to ${\rm HOM}$-convergence.
According to Proposition~\ref{lem:extBA}, ${\rm HOM}$-convergence is equivalent
to $\overline{\rm HOM}$-convergence. It is easy to see that $\overline{\rm HOM}$ is 
the fragment ${\rm QF}^-$ of quantifier free formulas that do not use equality.
We prove now that ${\rm HOM}$-convergence is actually equivalent to 
${\rm QF}$-convergence, where ${\rm QF}$ is the fragment of all quantifier free formulas.
Note that ${\rm QF}$ is a proper fragment of the fragment ${\rm FO}^{\rm local}$ of local formulas (that is of formulas whose satisfaction only depends on a fixed neighborhood of the free variables, see Section~\ref{sec:BS} for a formal definition).  
\begin{mytinted}
\begin{theorem}
\label{thm:LQF}
Let $(G_n)$ be a sequence of finite graphs such that
$\lim_{n\rightarrow\infty}|G_n|=\infty$.

Then the following conditions are equivalent:
\begin{enumerate}
  \item the sequence $(G_n)$ is {\rm L}-convergent; 
  \item the sequence $(G_n)$ is ${\rm QF}^-$-convergent;
  \item the sequence $(G_n)$ is ${\rm QF}$-convergent;   
\end{enumerate} 
\end{theorem}
\end{mytinted}
\begin{proof}
As {\rm L}-convergence is equivalent to ${\rm HOM}$-convergence and
as ${\rm HOM}\subset {\rm QF}^-\subset {\rm QF}$, it is sufficient
to prove that {\rm L}-convergence implies ${\rm QF}$-convergence.
 
Assume $(G_n)$ is {\rm L}-convergent. 
The inclusion/exclusion principle implies that for
every finite graph $F$ the density of induced subgraphs isomorphic
 to $F$ converges too. Define

$${\rm dens}(F,G_n)=\frac{(\#F\subseteq_i G_n)}{|G_n|^{|F|}}.$$

Then ${\rm dens}(F,G_n)$ is a converging sequence for each $F$.

Let $\theta$ be a quantifier free formula with ${\rm Fv}(\theta)\subseteq [p]$.
We first consider all possible cases of equalities between the free variables.
For a partition $\mathcal P=(I_1,\dots,I_k)$ of $[p]$, we define
$|\mathcal P|=k$ and $s_{\mathcal P}(i)=\min I_i$ (for $1\leq i\leq |\mathcal P|$).
Consider the formula 
$$\zeta_{\mathcal P}:=\bigwedge_{i=1}^{|\mathcal
P|}\biggl(\bigwedge_{j\in I_i}(x_j=x_{s_{\mathcal P}(i)})\ \wedge\ 
\bigwedge_{j=i+1}^{|\mathcal P|}(x_{s_{\mathcal P}(j)}\neq x_{s_{\mathcal
P}(i)})\biggr).$$
Then $\theta$ is logically equivalent to
\begin{equation*}
(\bigwedge_{i\neq j}(x_i\neq
x_j)\wedge\theta)\ \vee\bigvee_{\mathcal P:|\mathcal P|<p}\zeta_{\mathcal P}\wedge \theta_{\mathcal
P}(x_{s_{\mathcal P}(1)},\dots,x_{s_{\mathcal P}(|\mathcal P|)}).
\end{equation*}

Note that all the formulas in the disjunction are mutually exclusive.
Also
$\bigwedge_{i\neq j}(x_i\neq
x_j)\wedge\theta$ may be expressed as a disjunction of mutually exclusive
terms:
 
$$\bigwedge_{i\neq j}(x_i\neq
x_j)\wedge\theta=\bigvee_{F\in\mathcal
F}\theta_{F}',$$
where $\mathcal F$ is a finite family of finite graphs $F$ and where 
$G\models\theta_{F}'(v_1,\dots,v_p)$ if and only if the mapping $i\mapsto
 v_i$ is an isomorphism from $F$ to $G[v_1,\dots,v_p]$.

It follows that for every graph $G$ it holds that
\begin{align*}
\langle\theta,G\rangle&=\sum_{F\in\mathcal
F}\langle\theta_{F}',G\rangle+
\sum_{\mathcal P:|\mathcal P|<p}\langle\zeta_{\mathcal
P}\wedge \theta_{\mathcal P}(x_{s_{\mathcal
P}(1)},\dots,x_{s_{\mathcal
P}(|\mathcal P|)}),G\rangle\\
&=\sum_{F\in\mathcal
F}\langle\theta_{F}',G\rangle+\sum_{\mathcal
P:|\mathcal P|<p}|G|^{|\mathcal P|-p}\langle\theta_{\mathcal P},G\rangle\\
&=\sum_{F\in\mathcal F}\frac{1}{p!}\sum_{\sigma\in\Sym{p}}\frac{|\{(v_1,\dots,v_p):\
G\models\theta_F'(v_{\sigma(1)},\dots,v_{\sigma(p)})\}|}
{|G|^{p}}+ O(|G|^{-1})\\
&=\sum_{F\in\mathcal F}\frac{{\rm
Aut}(F)}{p!}\,{\rm dens}(F,G)+ O(|G|^{-1}).\\
\end{align*}
Thus $\langle\theta,G_n\rangle$ converge for every quantifier free formula
$\theta$. Hence $(G_n)$ is {\rm QF}-convergent.
\myqed\end{proof}
Notice that the condition that $\lim_{n\rightarrow\infty}|G_n|$ is necessary
as witnessed by the sequence $(G_n)$ where $G_n$ is $K_1$ if $n$ is odd and
$2K_1$ if $n$ is even. The sequence is obviously {\rm L}-convergent, but not QF
convergent as witnessed by the formula $\phi(x,y): x\neq y$, which has density 
$0$ in $K_1$ and $1/2$ in $2K_1$.

\begin{remark}
The Stone space of the fragment ${\rm QF}^-$ has a simple description.
 Indeed, a homomorphism $h:\mathcal B({\rm QF}^-)\rightarrow\mathbf 2$ is 
determined by its values on the formulas $x_i\sim x_j$ and any mapping from this subset of formulas to $\mathbf 2$ extends
(in a unique way) to a homomorphism of $\mathcal B({\rm QF}^-)$ to $\mathbf 2$.
Thus the points of $S(\mathcal B({\rm QF}^-))$ can be identified with the mappings from $\binom{\bbbn}{2}$ to $\{0,1\}$ that 
is to the graphs on $\bbbn$.
Hence the considered measures $\mu$ are probability measures of graphs on $\bbbn$ that have the property that
they are invariant under the natural action of $\Sym{\omega}$ on $\bbbn$. Such random graphs on $\bbbn$ are called
{\em infinite exchangeable random graphs}. For more on infinite exchangeable random graphs and graph limits, 
see e.g.~\cite{Austin2008,pre05504139}.
\end{remark}  

\subsection{BS-convergence and ${\rm FO}^{\rm local}$-convergence}
\label{sec:BS}
The class of graphs with maximum degree at most $D$ (for some integer $D$)
has received much attention. Specifically, the notion of {\em local weak convergence} 
of bounded degree graphs was introduced in
\cite{Benjamini2001}, which is called here {\em BS-convergence}:

A {\em rooted graph} is a pair $(G, o)$, where $o\in V(G)$.
 An {\em isomorphism} of rooted graph $\phi: (G, o)\rightarrow (G', o')$
  is an isomorphism of the underlying
graphs which satisfies $\phi(o) = o'$. 
Let $D\in\bbbn$. 
Let ${\mathcal G}_D$ denote the collection of
all isomorphism classes of connected 
rooted graphs with maximal degree at most $D$.
For the sake of simplicity, we denote elements of $\mathcal G_D$ simply as graphs.
 For $(G, o)\in{\mathcal G}_D$ and
$r\geq 0$ let $B_G(o, r)$ 
denote the subgraph of $G$ spanned by the vertices at 
distance at most $r$ from
$o$. If $(G, o), (G', o') \in{\mathcal G}_D$ and $r$
 is the largest integer such that
$(B_G(o, r), o)$ is
rooted-graph isomorphic to
$(B_{G'}(o', r), o')$, then set 
$\rho((G, o), (G', o'))= 1/r$, say. Also take 
$\rho((G, o), (G, o))=0$.
 Then $\rho$ is metric on ${\mathcal G}_D$. 
 Let $\mathfrak M_D$ denote the space of all probability measures
 on ${\mathcal G}_D$  that
are measurable with respect to the Borel $\sigma$-field of 
$\rho$. Then $\mathfrak M_D$ is endowed
with the topology of weak convergence, and is compact in this topology.

A sequence $(G_n)_{n\in\bbbn}$ of finite connected graphs with maximum degree at most $D$ is 
{\em BS-convergent} if, for every integer $r$ and every rooted connected graph 
$(F,o)$ with maximum degree at most $D$ the following limit exists:
$$
\lim_{n\rightarrow\infty}\frac{|\{v: B_{G_n}(v,r)\cong (F,o)\}|}{|G_n|}.
$$

This notion of limits leads to the definition of a limit object
as a probability measure on ${\mathcal G}_D$ \cite{Benjamini2001}.

To relate BS-convergence to $X$-convergence, we shall consider the fragment 
of local formulas:

Let $r\in\bbbn$. A formula $\phi\in{\rm FO}_p$ is {\em $r$-local} if, for
every graph $G$ and every $v_1,\dots,v_p\in G^p$ the following equivalence holds:
$$
G\models\phi(v_1,\dots,v_p)\quad\iff\quad
G[N_r(v_1,\dots,v_p)]\models\phi(v_1,\dots,v_p),$$
where $G[N_r(v_1,\dots,v_p)]$ denotes the subgraph of $G$ induced by all the vertices at (graph) distance
at most $r$ from one of $v_1,\dots,v_p$ in $G$.

A formula $\phi$ is {\em local} if it is $r$-local for some $r\in\bbbn$;
the fragment ${\rm FO}^{\rm local}$ is the set of all
local formulas in ${\rm FO}$.
Notice that if $\phi_1$ and $\phi_2$ are local formulas, so are $\phi_1\wedge\phi_2$,
$\phi_1\vee\phi_2$ and $\neg\phi_1$. It follows that the quotient of ${\rm FO}^{\rm local}$ by
the relation of logical equivalence defines a sub-Boolean algebra
$\mathcal B({\rm FO}^{\rm local})$ of $\mathcal B({\rm FO})$.
For $p\in\bbbn$ we further define ${\rm FO}_p^{\rm local}={\rm FO}^{\rm local}\cap {\rm FO}_p$.  

\begin{mytinted}
\begin{theorem}
\label{thm:BS}
Let $(G_n)$ be a sequence of finite graphs with maximum degree $d$, with
$\lim_{n\rightarrow\infty}|G_n|=\infty$. 

Then the following properties are equivalent:
\begin{enumerate}
  \item the sequence $(G_n)_{n\in\bbbn}$ is BS-convergent;
  \item the sequence $(G_n)_{n\in\bbbn}$ is ${\rm FO}_1^{\rm
local}$-convergent;
  \item the sequence $(G_n)_{n\in\bbbn}$ is ${\rm FO}^{\rm
local}$-convergent.
\end{enumerate}
\end{theorem}
\end{mytinted}
\begin{proof}
If $(G_n)_{n\in\bbbn}$ is ${\rm FO}^{\rm local}$-convergent, it is ${\rm FO}^{\rm
local}_1$-convergent;

If $(G_n)_{n\in\bbbn}$ is ${\rm FO}_1^{\rm local}$-convergent then it is 
BS-convergent as for any finite rooted graph $(F,o)$, testing whether the
the ball of radius $r$ centered at a vertex $x$ is isomorphic to
$(F,o)$ can be formulated by a local first order formula.

Assume $(G_n)_{n\in\bbbn}$ is BS-convergent.
As we consider graphs with maximum degree $d$, there are only
finitely many isomorphism types for the balls of radius $r$ centered at a
vertex. It follows that any local formula $\xi(x)$ with a single variable can be
expressed as the conjunction of  a finite number of (mutually exclusive)
formulas $\xi_{(F,o)}(x)$, which in turn correspond to subgraph testing.
It follows that BS-convergence implies ${\rm FO}_1^{\rm local}$-convergence.

Assume $(G_n)_{n\in\bbbn}$ is ${\rm FO}_1^{\rm local}$-convergent and
let $\phi\in{\rm FO}_p^{\rm local}$ be an $r$-local formula.
Let $\mathcal F_\phi$ be the set of all $p$-tuples 
$((F_1,f_1),\dots,(F_p,f_p))$ of rooted connected graphs
with maximum degree at most $d$ and radius (from the root) at most $r$ such that
$\bigcup_i F_i\models \phi(f_1,\dots,f_p)$.

Then, for every graph $G$ the sets

$$\Omega_\phi(G)=\{(v_1,\dots,v_p):\ G\models\phi(v_1,\dots,v_p)\}$$
and
$$\biguplus_{((F_1,f_1),\dots,(F_p,f_p))\in\mathcal F_\phi}\prod_{i=1}^p
\{v:\ G\models\theta_{(F_i,f_i)}(v)\}$$
differ by at most $O(|G|^{p-1})$ elements.
Indeed, according to the definition of an $r$-local formula, the $p$-tuples
$(x_1,\dots,x_p)$ belonging to exactly one of these sets are such that there
exists $1\leq i<j\leq p$ such that ${\rm dist}(x_i,x_j)\leq 2r$.

It follows that
$$
\langle\phi,G\rangle=\bigl(\sum_{((F_i,f_i))_{1\leq i\leq p}\in\mathcal
F_\phi}\, \prod_{i=1}^p\,\langle\theta_{(F_i,f_i)},G\rangle\bigr)+O(|G|^{-1}).
$$

It follows that ${\rm FO}^{\rm local}_1$-convergence (hence BS-convergence)
implies full ${\rm FO}^{\rm local}$-convergence.
\myqed\end{proof}

\begin{remark}
\label{rem:BS}
According to this proposition and Theorem~\ref{thm:genmodxc}, the BS-limit of a sequence of graphs
with maximum degree at most $D$ corresponds to a probability measure
on $S(\mathcal B({\rm FO}_1^{\rm local}))$ whose support is include in
the clopen set $K(\zeta_D)$,  where $\zeta_D$ is the sentence
expressing that the maximum degree is at most $D$.
The Boolean algebra $\mathcal B({\rm FO}_1^{\rm local})$ is isomorphic
to the Boolean algebra defined by the fragment 
$X\subset {\rm FO}_0(\sig_1)$ 
of sentences for rooted graphs
that are local with respect to the root (here, $\sig_1$ denotes the signature of graphs
augmented by one symbol of constant).  
According to this locality, any two countable rooted graphs $(G_1,r_1)$ 
and $(G_2,r_2)$,
the trace of the complete theories of $(G_1,r_1)$ and $(G_2,r_2)$ on $X$
 are the same if and only if the (rooted) connected component $(G_1',r_1)$ of $(G_1,r_1)$ 
 containing the root $r_1$ is elementary equivalent to the (rooted) connected component $(G_2',r_2)$ of $(G_2,r_2)$ 
 containing the root $r_2$. As isomorphism and elementary equivalence
are equivalent for  countable connected graphs with bounded degrees (see Lemma~\ref{lem:locfin})
it is easily checked that $K_X(\zeta_D)$ is homeomorphic to $\mathcal G_D$.
Hence our setting (based on a very different and dual approach) 
leads essentially to the same limit object as \cite{Benjamini2001}
for BS-convergent sequences. 
\end{remark}
\subsection{Elementary-convergence and ${\rm FO}_0$-convergence}
\label{sec:elem}
We already mentioned that ${\rm FO}_0$-convergence is nothing but elementary
convergence.  Elementary convergence is implicitly part of the classical model theory.
Although we only consider graphs here, the definition and results indeed generalize
to general $\sig$-structures We now reword the notion of elementary convergence:

A sequence $(G_n)_{n\in\bbbn}$ is {\em elementarily convergent} if, for every
sentence $\phi\in{\rm FO}_0$, there exists a integer $N$ such that 
either all the graphs $G_n$ $(n\geq N)$ satisfy $\phi$ or none of them do.

Of course, the limit object (as a graph) is not unique in general and formally, the
limit of an elementarily convergent sequence of graphs is an elementary
class defined by a complete theory.

Elementary convergence is also the backbone of all the $X$-convergences we consider
in this paper. The ${\rm FO}_0$-convergence is induced by an easy ultrametric defined
on equivalence classes of elementarily equivalent graphs. Precisely,
two (finite or infinite) graphs $G_1,G_2$ are {\em elementarily
equivalent} (denoted $G_1\equiv G_2$) if, for every sentence $\phi$ the following equivalence holds:

$$G_1\models\phi\ \iff\ G_2\models\phi.$$

In other words, two graphs are elementarily equivalent if they satisfy the same sentences.

A weaker (parametrized) notion of equivalence will be crucial:
two graphs $G_1,G_2$ are {\em $k$-elementarily
equivalent} (denoted $G_1\equiv^k G_2$) if, for every sentence $\phi$ 
with quantifier rank at most $k$ it holds that
$G_1\models\phi\ \iff\ G_2\models\phi$.

It is easily checked that for every two graphs $G_1,G_2$ the following equivalence holds:
$$
G_1\equiv G_2\quad\iff\quad (\forall k\in\bbbn)\ G_1\equiv^k G_2.
$$
For every fixed $k\in\bbbn$, checking whether two graphs $G_1$ and $G_2$
are $k$-elementarily equivalent can be done using the so-called
Ehrenfeucht-Fra{\"\i}ss{\'e} game.

From the notion of $k$-elementary equivalence naturally derives
a pseudometric ${\rm dist}_0(G_1,G_2)$: 
$$
{\rm dist}_0(G_1,G_2)=\begin{cases}
0&\text{if }G_1\equiv G_2\\
\min\{2^{-{\rm qrank}(\phi)}: (G_1\models\phi)\wedge(G_2\models\neg\phi)\}&\text{otherwise}
\end{cases}
$$

\begin{proposition}
The metric space of countable graphs (up to elementary equivalence) with
ultrametric ${\rm dist}_0$ is compact.
\end{proposition}
\begin{proof}
This is a direct consequence of the compactness theorem for first-order logic
(a theory has a model if and only if every finite subset of it has a model) and 
of the downward L\"owenheim-Skolem theorem (if a theory has a model and the
language is countable then the theory has a countable model).
\myqed\end{proof}

Note that not every countable graph is (up to elementary equivalence) the
limit of a sequence of finite graphs. A graph $G$ that is a limit of a
sequence finite graphs is said to have the {\em finite model property}, as
such a graph is characterized by the property that every finite set of sentences
satisfied by $G$ has a finite model (which does not imply that $G$ is
elementarily equivalent to a finite graph). As proved by Trakhtenbrot \cite{Trakhtenbrot1950}
the set of finitely satisfiable sentences is not decidable and 
deciding wether a given theory has a finite model is usually an extremely  difficult problem (see for instance Example~\ref{ex:triangle_free}).

\begin{example}
A {\em ray} is not an elementary limit of finite graphs as it contains exactly one 
vertex of degree $1$ and all the other vertices have degree $2$, what can be expressed 
in first-order logic but is satisfied by no finite graph.
However, the union of two rays is an elementary limit from the sequence $(P_n)_{n\in\bbbn}$ of paths
of order $n$.
\end{example}

Although two finite graphs are elementary equivalent if and only if they
are isomorphic, this property does not holds in general for countable graphs.
For instance, the union of a ray and a line is elementarily equivalent to a ray.
However we shall make use of the equivalence of isomorphisms and elementary equivalences 
for rooted connected countable locally finite graphs, which we prove now
for completeness.

\begin{lemma}
\label{lem:locfin}
Let $(G,r)$ and $(G',r')$ be two rooted connected countable graphs. 

If $G$ is locally finite then $(G,r)\equiv (G',r')$ if and only if
  $(G,r)$ and $(G',r')$ are isomorphic. 
\end{lemma}
\begin{proof}
If two rooted graphs are isomorphic they are obviously elementarily equivalent.
Assume that $(G,r)$ and $(G',r')$ are elementarily equivalent.
Enumerate the vertices of $G$ in a way that distance to the root is not
decreasing. Using $n$-back-and-forth equivalence (for all $n\in\bbbn$), one
builds a tree of partial isomorphisms of the subgraphs induced by the $n$
first vertices, where ancestor relation is restriction.
This tree is infinite and has only finite degrees. Hence, by K{\H o}nig's lemma,
it contains an infinite path.
It is easily checked that it defines an isomorphism from $(G,r)$ to $(G',r')$ as
these graphs are connected.
\myqed\end{proof}

Fragments of ${\rm FO}_0$ allow to define convergence notions, which are
weaker than elementary convergence. The hierarchy of the convergence schemes defined 
by sub-algebras of $\mathcal B({\rm FO}_0)$ is as strict as one could expect. Precisely,
if $X\subset Y$ are two sub-algebras of $\mathcal B({\rm FO}_0)$ then $Y$-convergence
is strictly stronger than $X$-convergence --- meaning that there exists graph sequences that are 
$X$-convergent but not $Y$-convergent --- if and only if there exists a sentence $\phi\in Y$ such that
for every sentence $\psi\in X$, there exists a (finite) graph $G$ disproving $\phi\leftrightarrow\psi$. 
 
We shall see that the special case of elementary convergent
sequences is of particular
importance. Indeed, every limit measure is a Dirac measure concentrated on a single point
of $S(\mathcal B({\rm FO}_0))$. 
This point is the complete theory of the elementary limit of the considered sequence.
This limit can be represented by a finite or countable graph. 
As ${\rm FO}$-convergence (and any ${\rm FO}_p$-convergence) implies ${\rm FO}_0$-convergence, the 
support of a limit measure $\mu$ corresponding to an ${\rm FO}_p$-convergent sequence (or to 
an ${\rm FO}$-convergent sequence) is such that ${\rm Supp}(\mu)$ projects to a single point
of $S(\mathcal B({\rm FO}_0))$.

Finally, let us remark that all the results of this section can be readily formulated and
proved for $\sig$-structures.

\section{Combining Fragments}
\label{sec:combine}
\subsection{The ${\rm FO}_p$ Hierarchy}
\label{sec:fop}
When we consider ${\rm FO}_p$-convergence of finite $\sig$-structures 
for finite a signature $\sig$, 
 the space $S(\mathcal B({\rm FO}_p(\sig)))$
can be given the following ultrametric ${\rm dist}_p$ (compatible with the topology of 
$S(\mathcal B({\rm FO}_p(\sig)))$):
Let $T_1,T_2\in S(\mathcal B({\rm FO}_p(\sig)))$ (where the points of $S(\mathcal B({\rm FO}_p(\sig)))$ are identified
with ultrafilters on $\mathcal B({\rm FO}_p(\sig))$). Then
$$
{\rm dist}_p(T_1,T_2)=\begin{cases}
0&\text{if }T_1=T_2\\
2^{-\min\{{\rm qrank}(\phi):\ \phi\in T_1\setminus T_2\}}&\text{otherwise}
\end{cases}
$$ 

This ultrametric has several other nice properties:
\begin{itemize}
  \item actions of $\Sym{p}$ on $S(\mathcal B({\rm FO}_p(\sig)))$ are isometries:
  $$\forall \sigma\in \Sym{p}\ \forall T_1,T_2\in S(\mathcal B({\rm FO}_p(\sig)))\quad 
  {\rm dist}_p(\sigma\cdot T_1,\sigma\cdot T_2)={\rm dist}_p(T_1,T_2);
   $$
  \item projections $\pi_p$ are contractions:
  $$\forall q\geq p\ \forall T_1,T_2\in S(\mathcal B({\rm FO}_q(\sig)))\quad
  {\rm dist}_p(\pi_p(T_1),\pi_p(T_2))\leq {\rm dist}_q(T_1,T_2).
  $$
\end{itemize}

We prove that there is a natural isometric
embedding $\eta_p:S(\mathcal B({\rm FO}_p(\sig)))\rightarrow S(\mathcal B({\rm FO}(\sig)))$.
This may be seen as follows:
for an ultrafilter $X\in S(\mathcal B({\rm FO}_p(\sig)))$, consider 
the filter $X^+$ on $\mathcal B({\rm FO}(\sig))$ generated by $X$ and all the formulas
$x_i=x_{i+1}$ (for $i\geq p$). This filter is an ultrafilter: for every
sentence $\phi\in{\rm FO}(\sig)$, let $\widetilde{\phi}$ be the
sentence obtained from $\phi$ by replacing each free occurrence of a variable $x_q$ with $q>p$ by $x_p$. It is
clear that $\phi$ and $\widetilde{\phi}$ are equivalent modulo the theory $T_p=\{(x_i=x_{i+1}):\ i\geq p\}$. 
As either $\widetilde{\phi}$ or $\neg\widetilde{\phi}$ belongs to $X$, either $\phi$ or
$\neg\phi$ belongs to $\eta_p(X)$. Moreover, we deduce easily from the fact
that $\widetilde{\phi}$ and $\phi$ have the same quantifier rank that 
if $q\geq p$ then $\pi_q\circ\eta_p$ is an
isometry. Finally, let us note that $\pi_p\circ\eta_p$ is the
identity of $S(\mathcal B({\rm FO}_p(\sig)))$.

Let $\sig_p$ be the signature $\sig$ augmented by $p$ symbols of constants $c_1,\dots,c_p$.
There is a natural isomorphism of Boolean algebras $\nu_p:{\rm FO}_p(\sig)\rightarrow {\rm FO}_0(\sig_p)$,
which replaces the free occurrences of the variables $x_1,\dots,x_p$  
in a formula $\phi\in{\rm FO}_p$ by the corresponding symbols of constants
$c_1,\dots,c_p$, so that the following equation holds, for every {\rps} $\mathbf A$, for every $\phi\in{\rm FO}_p$ and
every $v_1,\dots,v_p\in A$:

$$
\mathbf A\models\phi(v_1,\dots,v_p)\quad\iff\quad (\mathbf
A,v_1,\dots,v_p)\models\nu_p(\phi).$$

This mapping induces an isometric isomorphism of the metric spaces \linebreak $(S(\mathcal B({\rm FO}_p(\sig))),{\rm dist}_p)$
and $(S(\mathcal B({\rm FO}_0(\sig_p))),{\rm dist}_0)$.
Note that the Stone space \linebreak $S(\mathcal B({\rm FO}_0(\sig_p)))$ associated to the Boolean algebra $\mathcal B({\rm FO}_0(\sig_p))$
 is the space of all complete theories of $\sig_p$-structures.
 In particular, points of $S(\mathcal B({\rm FO}_p(\sig))$
can be represented (up to elementary equivalence) by countable $\sig$-structures with $p$ special points.
All these transformations may seem routine but they need to be carefully formulated and checked.

We can test whether the distance
${\rm dist}_p$ of two theories $T$ and $T'$ is smaller than $2^{-n}$ by means of an
Ehrenfeucht-Fra{\"\i}ss{\'e} game: Let $\nu_p(T)=\{\nu_p(\phi): \phi\in T\}$ and, similarly, let
$\nu_p(T')=\{\nu_p(\phi): \phi\in T'\}$. Let $(\mathbf A,v_1,\dots,v_p)$ be a model of $T$ and let
$(\mathbf A',v_1',\dots,v_p')$ be a model of $T'$. Then the following equivalence holds:
$$
{\rm dist}_p(T,T')<2^{-n}\quad\iff\quad (\mathbf A,v_1,\dots,v_p)\equiv^n (\mathbf A',v_1',\dots,v_p').
$$
Recall that the $n$-rounds {\em Ehrenfeucht-Fra{\"\i}ss{\'e} game} on two $\sig$-structures $\mathbf A$ and $\mathbf A'$,
denoted ${\rm EF}(\mathbf A,\mathbf A',n)$ is the perfect information game with two players --- the Spoiler and the Duplicator --- 
defined as follows:
The game has $n$ rounds and each round has two parts. At each round, the Spoiler first chooses one of $\mathbf A$ 
and $\mathbf A'$ and
accordingly selects either a vertex $x\in A$ or a vertex $y\in A'$. Then, the Duplicator selects a vertex in the other 
$\sig$-structure. At the end of the $n$ rounds, $n$ vertices have been selected from each structure: $x_1,\dots,x_n$ in $A$
and $y_1,\dots,y_n$ in $A'$ ($x_i$ and $y_i$ corresponding to vertices $x$ and $y$ selected during the $i$th round).
The Duplicator wins if the substructure induced by the selected vertices are order-isomorphic (i.e.\ $x_i\mapsto y_i$ is 
an isomorphism of $\mathbf A[\{x_1,\dots,x_n\}]$ and $\mathbf A'[\{y_1,\dots,y_n\}]$). As there are no hidden moves and no draws, one of the two players
has a winning strategy, and we say that that player wins ${\rm EF}(\mathbf A,\mathbf A',n)$.
The main property of this game is the following equivalence, due to Fra{\"\i}ss{\'e} \cite{Fraiss'e1950,Fraiss'e1953}
and Ehrenfeucht \cite{Ehrenfeucht1961}: The duplicator wins ${\rm EF}(\mathbf A,\mathbf A',n)$ if and only if 
$\mathbf A\equiv^n \mathbf A'$.
In our context this translates to the following equivalence:
$$
{\rm dist}_p(T,T')<2^{-n}\ \iff\ \text{Duplicator wins }{\rm EF}((\mathbf A,v_1,\dots,v_p),(\mathbf A',v_1',\dots,v_p'),n).
$$
As ${\rm FO}_0\subset {\rm FO}_1\subset\dots\subset {\rm FO}_p\subset {\rm FO}_{p+1}\subset\dots\subset{\rm FO}=\bigcup_i {\rm FO}_i$,
the fragments ${\rm FO}$ form a hierarchy of more and more restrictive notions of convergence. In particular,
${\rm FO}_{p+1}$-convergence implies ${\rm FO}_{p}$-convergence and ${\rm FO}$-convergence is equivalent to
${\rm FO}_{p}$ for all $p$. If a sequence $(\mathbf A_n)_{n\in\bbbn}$ is ${\rm FO}_p$-convergent then for every $q\leq p$
the ${\rm FO}_q$-limit of $(\mathbf A_n)_{n\in\bbbn}$ is a measure $\mu_q\in{\rm rca}(S(\mathcal B({\rm FO}_q)))$, which is
the pushforward of $\mu_p$ by the projection $\pi_q$ (more precisely, by the restriction of $\pi_q$ to
$S(\mathcal B({\rm FO}_p))$):
$$
\mu_q=(\pi_q)_*(\mu_p).
$$

\subsection{${\rm FO}^{\rm local}$ and  Locality}

${\rm FO}$-convergence can be reduced to the conjunction
of elementary convergence and ${\rm FO}^\text{local}$-convergence, which we call {\em local convergence}.
This is a consequence of a result, which we recall now:

\begin{theorem}[Gaifman locality theorem \cite{Gaifman1982}]
\label{thm:gaifman}
For every first-order formula $\phi(x_1,\dots,x_n)$ there exist integers
$t$ and $r$ such that $\phi$ is equivalent to a Boolean combination of $t$-local
formulas $\xi_s(x_{i_1},\dots,x_{i_s})$ and sentences of the form

\begin{equation}
\label{eq:gl}
\exists
y_1\dots\exists y_m\biggl(\bigwedge_{1\leq i<j\leq m}{\rm
dist}(y_i,y_j)>2r\wedge\bigwedge_{1\leq i\leq
m}\psi(y_i)\biggr)
\end{equation}

\noindent where $\psi$ is $r$-local. Furthermore, we can choose
 $$r\leq 7^{{\rm qrank}(\phi)-1},\ t\leq
(7^{{\rm qrank}(\phi)-1}-1)/2,\ m\leq n+{\rm qrank}(\phi),$$
 and, if
$\phi$ is a sentence, only sentences~\eqref{eq:gl} occur in the Boolean
combination. Moreover, these sentences can be chosen with quantifier rank
at most $q({\rm qrank}(\phi))$, for some fixed function $q$.  
\end{theorem}

From this theorem and the following folklore technical result 
will follow the claimed decomposition of
${\rm FO}$-convergence into elementary and local convergence.
\begin{lemma}
\label{lem:supba}
Let $B$ be a Boolean algebra, let $A_1$ and $A_2$ be sub-Boolean algebras of $B$, and let
$b\in B[A_1\cup A_2]$ be a Boolean combination of elements from $A_1$ and $A_2$.
Then $b$ can be written as 
$$b=\bigvee_{i\in I} x_i\wedge y_i,$$
where $I$ is finite, $x_i\in A_1$, $y_i\in A_2$, and
for every $i\neq j$ in $I$ it holds that 
$(x_i\wedge y_i)\wedge (x_j\wedge y_j)=0$.
\end{lemma}
\begin{proof}
Let $b=F(u_1,\dots,u_a,v_1,\dots,v_b)$ with $u_i\in A_1$ ($1\leq i\leq a$) and $v_j\in A_2$
($1\leq j\leq b$) where $F$ is a Boolean combination. By using iteratively Shannon's expansion, we can write
$F$ as 
$$F(u_1,\dots,u_a,v_1,\dots,v_b)=\bigvee_{(X_1,X_2,Y_1,Y_2)\in\mathcal F}\ (\bigwedge_{i\in X_1}u_i\wedge
\bigwedge_{i\in X_2}\neg u_i\wedge \bigwedge_{j\in Y_1}v_j\wedge\bigwedge_{j\in Y_2}\neg v_j),
$$
where $\mathcal F$ is a subset of the quadruples $(X_1,X_2,Y_1,Y_2)$ such that $(X_1,X_2)$ is a partition of 
$[a]$ and $(Y_1,Y_2)$ is a partition of $[b]$. For a quadruple $Q=(X_1,X_2,Y_1,Y_2)$, define
$x_Q=\bigwedge_{i\in X_1}u_i\wedge\bigwedge_{i\in X_2}\neg u_i$ and 
$y_Q=\bigwedge_{j\in Y_1}v_j\wedge\bigwedge_{j\in Y_2}\neg v_j$. Then for every $Q\in\mathcal F$ it holds that 
$x_Q\in A_1, y_Q\in A_2$, for every $Q\neq Q'\in\mathcal F$ it holds that  $x_Q\wedge y_Q\wedge x_{Q'}\wedge y_{Q'}=0$, and
we have $b=\bigvee_{Q\in\mathcal F}\ x_Q\wedge y_Q$. 
\myqed\end{proof}

\begin{mytinted}
\begin{theorem}
\label{thm:fole}
Let $(\mathbf A_n)$ be a sequence of finite $\sig$-structures. Then $(\mathbf A_n)$ is ${\rm FO}$-convergent if
and only if it is both ${\rm FO}^{\rm local}$-convergent and ${\rm FO}_0$-convergent.
Precisely, $(\mathbf A_n)$ is ${\rm FO}_p$-convergent if
and only if it is both ${\rm FO}_p^{\rm local}$-convergent and ${\rm FO}_0$-convergent.
\end{theorem}
\end{mytinted}
\begin{proof}
Assume $(\mathbf A_n)_{n\in\bbbn}$ is both ${\rm FO}^{\rm local}_p$-convergent and ${\rm FO}_0$-convergent and
let $\phi\in{\rm FO}_p$.
According to Theorem~\ref{thm:gaifman}, there exist integers
$t$ and $r$ such that $\phi$ is equivalent to a Boolean combination of $t$-local
formula $\xi(x_{i_1},\dots,x_{i_s})$ and of sentences.
As both ${\rm FO}^{\rm local}$ and ${\rm FO}_0$ define a sub-Boolean algebra of
$\mathcal B({\rm FO})$, according to Lemma~\ref{lem:supba}, $\phi$ can be written
as $\bigvee_{i\in I}\psi_i\wedge\theta_i$, where $I$ is finite, $\psi_i\in{\rm FO}^{\rm local}$,
$\theta_i\in{\rm FO_0}$, and $\psi_i\wedge\theta_i\wedge\psi_j\wedge\theta_j=0$ if $i\neq j$.
Thus for every finite $\sig$-structure $\mathbf A$ the following equation holds:
$$
\langle\phi,\mathbf A\rangle=\sum_{i\in I}\langle\psi_i\wedge\theta_i,\mathbf A\rangle.
$$
As $\langle\,\cdot\,,\mathbf A\rangle$ is additive and
$\langle\theta_i,\mathbf A\rangle\in\{0,1\}$ we have
 $\langle\psi_i\wedge\theta_i,\mathbf A\rangle=\langle\psi_i,\mathbf A\rangle\,\langle\theta_i,\mathbf A\rangle$.
 Hence 
$$
\langle\phi,\mathbf A\rangle=\sum_{i\in I}\langle\psi_i,\mathbf A\rangle\,\langle\theta_i,\mathbf A\rangle.
$$
Thus if $(\mathbf A_n)_{n\in\bbbn}$ is both ${\rm FO}^{\rm local}_p$-convergent and ${\rm FO}_0$-convergent
then $(\mathbf A_n)_{n\in\bbbn}$ is ${\rm FO}_p$-convergent.
\myqed\end{proof}

Similarly points of $S(\mathcal B({\rm FO}_p(\sig))$
can be represented (up to elementary equivalence) by countable $\sig$-structures with $p$ special points, and 
points of \linebreak $S(\mathcal B({\rm FO}_p^{\rm local}(\sig))$
can be represented by countable $\sig$-structures with $p$ special points
such that every connected component contains at least one special point.
In particular, points of $S(\mathcal B({\rm FO}_1^{\rm local}(\sig))$
can be represented by rooted connected countable $\sig$-structures.

\begin{figure}[ht]
\begin{center}
  \includegraphics[width=.8\textwidth]{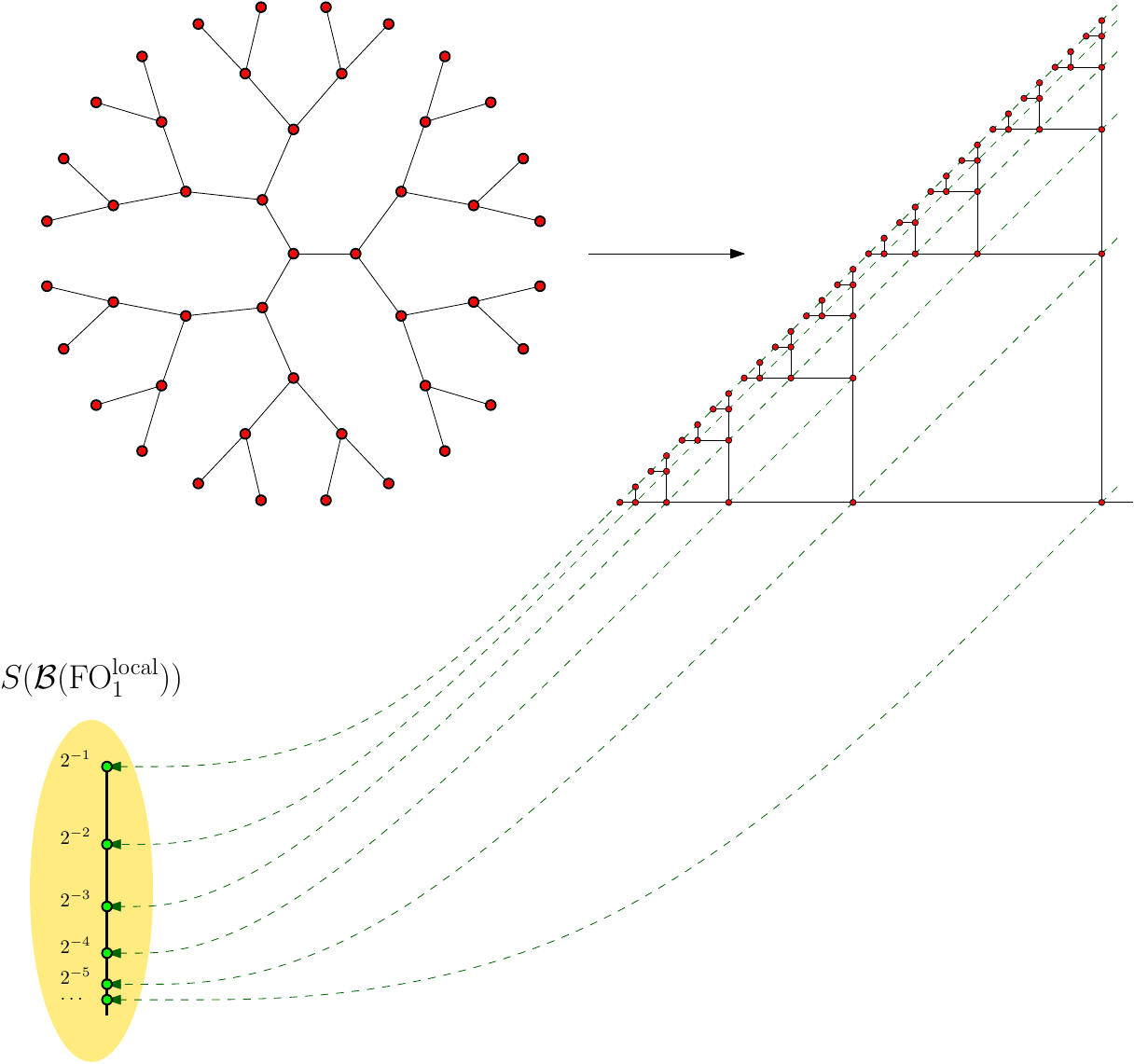}
  \caption{An outline of the local limit of a sequence of trees}
  \label{fig:3tree}
\end{center}
\end{figure}

Also, the structure of an ${\rm FO}_2^{\rm local}$-limit of graphs can be outlined
by considering that points of $S(\mathcal B({\rm FO}_2^{\rm local}))$ as countable graphs
with two special vertices $c_1$ and $c_2$, such that every connected component contains
at least one of $c_1$ and $c_2$. Let $\mu_2$ be the limit probability 
measure on $S(\mathcal B({\rm FO}_2^{\rm local}))$
for an ${\rm FO}_2^{\rm local}$-convergent sequence $(G_n)_{n\in\bbbn}$, let 
$\pi_1$
be the standard projection 
of $S(\mathcal B({\rm FO}_2^{\rm local}))$ into $S(\mathcal B({\rm FO}_1^{\rm local}))$, and let $\mu_1$
be the pushforward of $\mu_2$ by $\pi_1$.
We construct a measurable
graph $\hat G$ as follows: the vertex set of $\hat G$ is the support ${\rm Supp}(\mu_1)$ of $\mu_1$.
Two vertices $x$ and $y$ of $\hat G$ are adjacent if there exists $x'\in \pi_1^{-1}(x)$ and 
$y'\in \pi_1^{-1}(y)$
such that (considered as ultrafilters of $\mathcal B({\rm FO}_2^{\rm local})$) it holds that:
\begin{itemize}
  \item $x_1\sim x_2$ belongs to both $x'$ and $y'$,
  \item the transposition $\tau_{1,2}$ exchanges $x'$ and $y'$ (i.e.\ $y'=\tau_{1,2}\cdot x'$).
\end{itemize}
The vertex set of $\hat G$ is of course endowed with a structure of a probability space
(as a measurable subspace of $S(\mathcal B({\rm FO}_1^{\rm local}))$ equipped with the probability
measure $\mu_1$).
In the case of bounded degree graphs, the obtained graph $\hat G$ is the {\em graph of graphs} 
introduced in~\cite{Lovasz2009}. Notice that this graph may have loops.
An example of such a graph is shown Fig.~\ref{fig:3tree}.
\subsection{Component-Local Formulas}
\label{sec:piloc}
It is sometimes possible to reduce \linebreak ${\rm FO}^{\rm local}$ to a smaller fragment.
This is in particular the case when connected components of the considered structures can be identified 
by some first-order formula. Precisely:

\begin{definition}
\label{def:pirel}
Let $\sig$ be a signature and let $T$ be a theory of $\sig$-structures.
A binary relation $\varpi\in\sig$ is a {\em component relation} in $T$ 
if $T$ entails that $\varpi$ is an equivalence relation such that 
for every $k$-ary relation $R\in\sig$ with $k\geq 2$ the following property holds:
$$T\models (\forall x_1,\dots,x_k)\ \Bigl(R(x_1,\dots,x_k)\rightarrow\bigwedge_{1\leq i<j\leq k} \varpi(x_i,x_j)\Bigr).$$

A local  formula $\phi$ with $p$ free variables is {\em $\varpi$-local}  
if $\phi$ is equivalent (modulo $T$) to
$\phi\wedge\bigwedge_{x_i,x_j\in{\rm Fv}(\phi)}\varpi(x_i,x_j)$.
\end{definition}

In presence of a component relation, it is possible to reduce from ${\rm FO}^{\rm local}$ to the fragment
of $\varpi$-local formulas, thanks to the following result.

\begin{mytinted}
\begin{lemma}
\label{lem:sl}
Let $\varpi$ be a component relation in a theory $T$.
For every local formula $\phi$ with quantifier rank $r$ there exist $\varpi$-local formulas $\xi_{i,j}\in{\rm FO}_{q_{i,j}}^{\rm local}$ 
($1\leq i\leq n$, $j\in I_i$) with quantifier rank at most $r$ and permutations $\sigma_i$ of $[p]$ ($1\leq i\leq n$) such that
for each $1\leq i\leq n$, $\sum_{j\in I_i} q_{i,j}=p$
and, for every model $\mathbf A$ of $T$ the following equation holds:
$$\Omega_\phi(\mathbf A)=\biguplus_{i=1}^n F_{\sigma_i}\Bigl(
\prod_{j\in I_i}\Omega_{\xi_{i,j}}(\mathbf A)\Bigr),$$
where 
$F_{\sigma_i}(X)$ performs a permutation of the coordinates according to $\sigma_i$.
\end{lemma}
\end{mytinted}
\begin{proof}
First note that if two $\varpi$-local formulas $\phi_1$ and $\phi_2$ share a free variable then 
$\phi_1\wedge\phi_2$ is $\varpi$-local.
For this obvious fact, we deduce that if
 $\psi_1,\dots,\psi_n$ are $\varpi$-local formulas in ${\rm FO}_p$, 
then there is a partition $\tau$ and a permutation $\sigma$ of $[p]$
such that for every $\sig$-structure $\mathbf A$
 the following equation holds:

$$\Omega_{\bigwedge_{i=1}^n\psi_i}(\mathbf A)=F_\sigma\Bigl(\prod_{P\in \tau}\Omega_{\bigwedge_{i\in
P}\psi_i}(\mathbf A)\Bigr),$$
where each $\bigwedge_{i\in P}\psi_i$ is $\varpi$-local, and $F_\sigma:A^p\rightarrow A^p$ is defined by
$$F_\sigma(X)=\{(v_{\sigma(1)},\dots,v_{\sigma(p)}):\ (v_1,\dots,v_p)\in X\}.$$

For a partition $\tau$ of $[p]$ we denote by $\zeta_\tau$ the 
conjunction of $\varpi(x_{i},x_{j})$ for every $i,j$ belonging to a same part
and of $\neg\varpi(x_{i},x_{j})$ for every $i,j$ belonging to different parts.
Then, for any two distinct partitions
$\tau$ and $\tau'$, the formula $\zeta_\tau\wedge\zeta_\tau'$ is never
satisfied; moreover $\bigvee_\tau \zeta_\tau$ is always satisfied. Thus
for every local formula $\phi$ the following equation holds:
 $$
\phi=\bigvee_{\tau} (\zeta_\tau\wedge\phi)=\bigoplus_{\tau}(\zeta_\tau\wedge\phi)
$$
(where only the partitions $\tau$ for which $\zeta_\tau\wedge\phi\neq 0$ have to be considered).

We denote by $\Lambda_\tau$ the formula $\bigwedge_{P\in\tau}\bigwedge_{i,j\in P}\varpi(x_i,x_j)$.
Obviously the following equation holds:
$$\Lambda_\tau=\bigoplus_{\tau'\geq\tau}\zeta_\tau,$$
where $\oplus$ stands for the exclusive disjunction ($a\oplus b=(a\wedge\neg n)\vee(\neg a\wedge b)$) and $\tau'\geq\tau$ means
that $\tau'$ is a partition of $[p]$, which is coarser than $\tau$. Then there exists (by M\"obius inversion or immediate induction) a function $M$ from the set of the partitions of $[p]$ to the powerset of the set of partitions of $[p]$ such that for every partition $\tau$ of $[p]$ the following equation holds:
$$\zeta_\tau=\bigoplus_{\tau'\in M(\tau)}\Lambda_{\tau'}.$$
Hence
$$\phi=\bigoplus_{\tau}\bigoplus_{\tau'\in M(\tau)}\Lambda_{\tau'}\wedge\phi.$$

It follows that $\phi$ is a Boolean combination of formulas $\Lambda_{\tau}\wedge\phi$, for partitions
$\tau$ such that $\zeta_\tau\wedge\phi\neq 0$ (as $\zeta_\tau\wedge\phi\neq 0$ and $\tau'\geq\tau$ imply
$\zeta_{\tau'}\wedge\phi\neq 0$). Each formula $\Lambda_{\tau}\wedge\phi$ is itself a Boolean combination of
$\varpi$-local formulas. Putting this in standard form (exclusive disjunction of conjunctions) and gathering in the conjunctions
the $\varpi$-local formulas whose set of free variables intersect, we get that there exist
families $\mathcal F_\tau$ of $\varpi$-local formulas $\varphi_{P}$ ($P\in\tau$) 
with free variables ${\rm Fv}(\varphi_{P})=\{x_j: j\in P\}$ such that 
$$\phi=\bigvee_\tau\bigvee_{\varphi\in \mathcal F_\tau} \bigwedge_{P\in\tau} \varphi_{P},$$
where the disjunction is exclusive.

Hence, considering adequate permutations $\sigma_\tau$ of $[p]$
the following equation holds:

$$\Omega_\phi(\mathbf A)=\biguplus_\tau\biguplus_{\varphi\in \mathcal F_\tau}F_{\sigma_\tau}\Bigl(
\prod_{P\in\tau}\Omega_{\tilde{\varphi}_P}(\mathbf A)\Bigr),$$
which is the requested form.

Note that the fact that ${\rm qrank}(\xi_{i,j})\leq {\rm qrank}(\phi)$ is obvious as we did not introduce any quantifier
in our transformations.
\end{proof}
As a consequence, we get the the desired:
\begin{corollary}
Let $\varpi$ be a component relation in a theory $T$ and let $(\mathbf A_n)_{n\in\bbbn}$ be a sequence of models of $T$.
Then the sequence $(\mathbf A_n)_{n\in\bbbn}$ is ${\rm FO}^{\rm local}$-convergent if and only if it is
${\rm FO}^{\mbox{\scriptsize $\varpi$-local}}$-convergent.
\end{corollary}
\subsection{Sequences with Homogeneous Elementary Limit}
\label{sec:homog}
Elementary convergence is an important aspect of ${\rm FO}$-convergence and
we shall see that in several contexts, ${\rm FO}$-convergence can be reduced 
to the conjunction of elementary convergence and $X$-convergence (for some suitable fragment $X$).

\begin{figure}[htp]
\begin{center}
  \includegraphics[width=.7\textwidth]{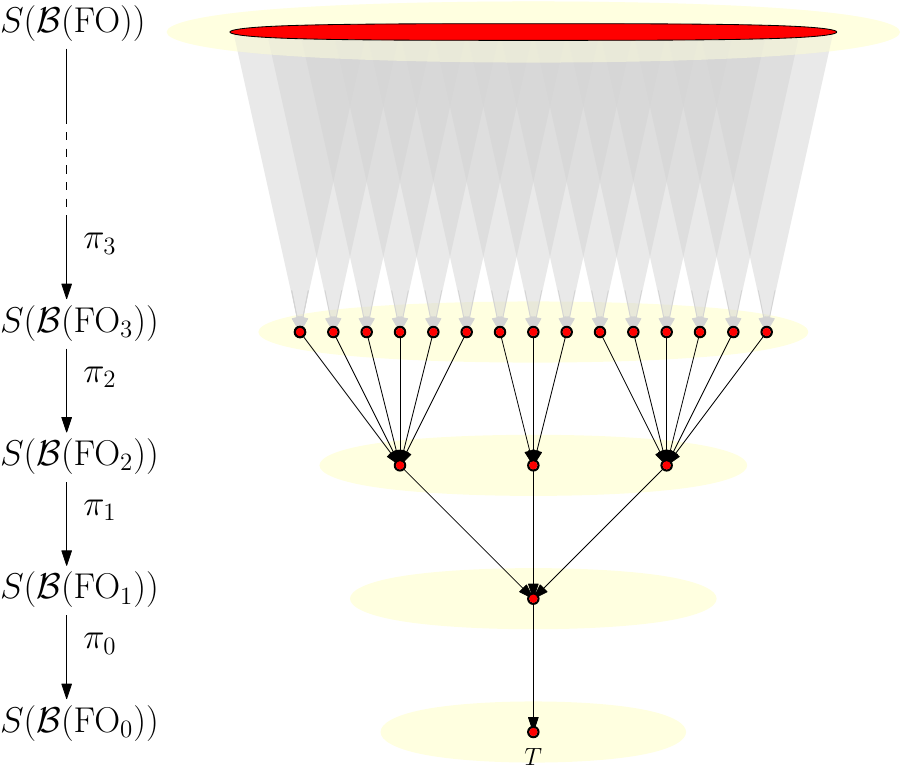}
  \caption{Ultrafilters projecting to an $\omega$-categorical theory}
  \label{fig:ultra}
\end{center}
\end{figure}

In some special cases, the limit (as a countable structure) will be unique. This means
that some particular complete theories have exactly one countable model (up to isomorphism).
Such complete theories are called {\em $\omega$-categorical}.
Several properties are known to be equivalent to $\omega$-categoricity. For instance, for a complete 
theory $T$ the following statements are equivalent:
\begin{itemize}
\item $T$ is $\omega$-categorical;
  \item for every  $p\in\bbbn$,
  the Stone space $S(\mathcal B({\rm FO}_p(\sig),T))$ is finite (see Fig.~\ref{fig:ultra}); 
  \item every countable model $\mathbf A$ of $T$ has an {\em oligomorphic} automorphism group, what means
 that for every $n\in\bbbn$, $A^n$ has finitely many orbits under the action of ${\rm Aut}(\mathbf A)$.  
\end{itemize}

A theory $T$ is said to have {\em quantifier elimination} if, for every $p$ and every formula
$\phi\in{\rm FO}_p(\sig)$ there exists $\widetilde\phi\in{\rm QF}_p(\sig)$ such that
$T\models \phi\leftrightarrow\widetilde\phi$. 
If a theory (in the language of  relational structures with given finite signature, like the language of graphs)
 has quantifier elimination then it is
$\omega$-categorical. Indeed, for every $p$, there exists only finitely many quantifier free
formulas with $p$ free variables hence (up to equivalence modulo $T$) only finitely many formulas
with $p$ free variables. The unique countable model of a complete theory $T$ (in the language of  relational structures with given finite signature) with quantifier elimination 
is {\em ultra-homogeneous}, what means that every partial isomorphism of finite induced substructures extends 
as a full automorphism. In the context of relational structures with given finite signature, 
the property of having a countable ultra-homogeneous model is equivalent to the
property of having quantifier elimination. 
We provide a proof of this 
folklore result in the context of graphs in order to illustrate these notions.
\begin{lemma}
\label{lem:UHQE}
Let $T$ be a complete theory of graphs with
 no finite model.

Then $T$ has quantifier elimination if and only if some (equivalently, every) countable
model of $T$ is ultra-homogeneous.
\end{lemma}
\begin{proof}
Assume that $T$ has an ultra-homogeneous countable model $G$.
Let $(a_1,\dots,a_p)$, $(b_1,\dots,b_p)$ be $p$-tuples of vertices of $G$.
Assume that $a_i\mapsto b_i$ is an isomorphism between
$G[a_1,\dots,a_p]$ and $G[b_1,\dots,b_p]$. Then, as $G$ is ultra-homogeneous, there exists
an automorphism $f$ of $G$ such that $f(a_i)=b_i$ for every $1\leq i\leq p$. As 
the satisfaction of a first-order formula is invariant under the action of the automorphism group,
for every formula $\phi\in{\rm FO}_p$ the following equivalence holds:
$$G\models\phi(a_1,\dots,a_p)\quad\iff\quad G\models\phi(b_1,\dots,b_p).$$
Consider a maximal set $\mathcal F$ of $p$-tuples $(v_1,\dots,v_p)$ of $G$
such that $G\models\phi(v_1,\dots,v_p)$ and no two $p$-tuples induce isomorphic (ordered) induced subgraphs.
Obviously $|\mathcal F|=2^{O(p^2)}$ is finite. Moreover, each $p$-tuple
$\vec v=(v_1,\dots,v_p)$ defines a quantifier free formula $\eta_{\vec v}$
with $p$ free variables such that 
$G\models \eta_{\vec v}(x_1,\dots,x_p)$ if and only if $x_i\mapsto v_i$ is an isomorphism between
$G[x_1,\dots,x_p]$ and $G[v_1,\dots,v_p]$. Hence the following property holds:
$$
G\models \phi\leftrightarrow \bigvee_{\vec v\in\mathcal F}\eta_{\vec v}.
$$
In other words, $\phi$ is equivalent (modulo $T$) to the quantifier free formula
$\widetilde\phi=\bigvee_{\vec v\in\mathcal F}\eta_{\vec v}$, that is: $T$ has quantifier elimination.

Conversely, assume that $T$ has quantifier elimination. As notice above, $T$ is $\omega$-categorical thus
has a unique countable model. Assume $(a_1,\dots,a_p)$ and $(b_1,\dots,b_p)$ are $p$-tuples
of vertices such that $f:a_i\mapsto b_i$ is a partial isomorphism.
Assume that $f$ does not extend into an automorphism of $G$. 
Let $(a_1,\dots,a_q)$ be a tuple of vertices of $G$ of maximal length
such that there exists $b_{p+1},\dots,b_q$ such that $a_i\mapsto b_i$ is 
a partial isomorphism. Let $a_{q+1}$ be a vertex distinct from $a_1,\dots,a_q$.
Let $\phi(x_1,\dots,x_q)$ be the formula
\begin{equation*}
\begin{split}
\bigwedge_{a_i\sim a_j}(x_i\sim x_j)\wedge\bigwedge_{\neg(a_i\sim a_j)}\neg(x_i\sim x_j)\wedge\bigwedge_{1\leq i\leq q} \neg(x_i=x_j)\\
\quad\wedge\ (\exists y)\ \biggl(\bigwedge_{a_i\sim a_{q+1}}(x_i\sim y)\wedge\bigwedge_{\neg(a_i\sim a_{q+1})}\neg(x_i\sim y)\wedge\bigwedge_{1\leq i\leq q} \neg(x_i=y)\biggr)
\end{split}
\end{equation*}
As $T$ has quantifier elimination, there exists a quantifier free formula
$\widetilde\phi$ such that
$T\models\phi\leftrightarrow\widetilde\phi$. 
As $G\models\phi(a_1,\dots,a_q)$ (witnessed by $a_{q+1}$) it holds that 
$G\models\widetilde\phi(a_1,\dots,a_q)$ hence $G\models\widetilde\phi(b_1,\dots,b_q)$
(as $a_i\mapsto b_i, 1\leq i\leq q$ is a partial isomorphism) thus
 $G\models\phi(b_1,\dots,b_q)$. It follows that there exists $b_{q+1}$
 such that $a_i\mapsto b_i, 1\leq i\leq q+1$ is a partial isomorphism, contradicting
 the maximality of $(a_1,\dots,a_q)$.
\myqed
\end{proof}

When a sequence of graphs is elementarily convergent to an ultra-homogeneous graph
(i.e.\ to a complete theory with quantifier elimination), we shall prove that 
${\rm FO}$-convergence reduces to ${\rm QF}$-convergence. 
This later mode of convergence is of particular interest as it is equivalent to {\rm L}-convergence. 

\begin{mytinted}
\begin{lemma}
\label{lem:FOQF}
Let $(G_n)_{n\in\bbbn}$ be sequence of graphs that converges elementarily to some
ultra-homogeneous graph $\hat G$. Then the following properties are equivalent:
\begin{itemize}
  \item the sequence $(G_n)_{n\in\bbbn}$ is ${\rm FO}$-convergent;
  \item the sequence $(G_n)_{n\in\bbbn}$ is ${\rm QF}$-convergent;
  \item the sequence $(G_n)_{n\in\bbbn}$ is {\rm L}-convergent.
\end{itemize}
\end{lemma}
\end{mytinted}
\begin{proof}
As ${\rm FO}$-convergence implies ${\rm QF}$-convergence we only have to prove the opposite direction.
Assume that the sequence $(G_n)_{n\in\bbbn}$ is ${\rm QF}$-convergent.
According to Lemma~\ref{lem:UHQE}, for every formula $\phi\in{\rm FO}_p$ there exists a quantifier free formula
$\widetilde\phi\in{\rm QF}_p$ such that $\hat G\models \phi\leftrightarrow\widetilde\phi$ (i.e.
${\rm Th}(\hat G)$ has quantifier elimination). As $\hat G$ is an elementary limit of the
sequence $(G_n)_{n\in\bbbn}$ there exists some $N$ such that for every $n\geq N$ it holds that 
 $G_n\models \phi\leftrightarrow\widetilde\phi$. It follows that for every $n\geq N$ it holds that 
 $\langle \phi,G_n\rangle=\langle\widetilde\phi,G_n\rangle$ hence
 $\lim_{n\rightarrow\infty}\langle \phi,G_n\rangle$ exists. Thus
 the sequence $(G_n)_{n\in\bbbn}$ is ${\rm FO}$-convergent.
 \myqed
\end{proof}

There are not so many countable ultra-homogeneous graphs.
\begin{theorem}[Lachlan and Woodrow \cite{Lachlan1980}]
 Every infinite countable ultrahomogeneous undirected
graph is isomorphic to one of the following:
\begin{itemize}
  \item the disjoint union of $m$ complete graphs of size $n$, where $m, n \leq\omega$ and at least one
of $m$ or $n$ is $\omega$, (or the complement of it);
\item the generic graph for the class of all countable graphs not containing $K_n$ for a given
$n\geq 3$ (or the complement of it);
\item the Rado graph R (the generic graph for the class of all countable graphs).
\end{itemize}
\end{theorem}

Among them, the Rado graph $R$
(also called ``the random graph'')
 is characterized by the {\em extension property}:
for every finite disjoint subsets of vertices $A$ and $B$ of $R$ there exists a vertex
$z$ of $R-A-B$ such that $z$ is adjacent to every vertex in $A$ and to no vertex in $B$.
We deduce for instance the following application of Lemma~\ref{lem:FOQF}.

\begin{example}
It is known \cite{Blass1981,Bollobas1981} that for every fixed $k$, Paley graphs 
of sufficiently large order satisfy the $k$-extension property hence the sequence of Paley graphs
converge elementarily to the Rado graph. Moreover, Paley graphs is a standard example of quasi-random 
graphs~\cite{Chung1989}, and the sequence of Paley graphs is {\rm L}-convergent to the $1/2$-graphon.
Thus, according to Lemma~\ref{lem:FOQF}, the sequence of Paley graphs is ${\rm FO}$-convergent.
\end{example}

%

We now relate more precisely the extension property with quantifier
elimination.

\begin{definition}
Let $k\in\bbbn$. A graph $G$ has the {\em $k$-extension property}
if, for every disjoint subsets of vertices $A,B$ of $G$ with  
size $k$ there exists a vertex $z$ not in $A\cup B$ that is
adjacent to every vertex in $A$ and to no vertex in $B$. In other words,
$G$ has the $k$-extension property if $G$ satisfies the sentence
$\Upsilon_k$ below:
\begin{equation*}
\begin{split}
(\forall x_1,\dots,x_{2k})\quad 
&\biggl(\bigwedge_{1\leq i<j\leq 2k}\neg(x_i=x_j)\\
&\qquad\rightarrow\quad
(\exists z)\ \bigwedge_{i=1}^{2k}\neg(x_i=z)\wedge
\bigwedge_{i=1}^{k}(x_i\sim z)\wedge\bigwedge_{i=k+1}^{2k}\neg(x_i\sim z)
\biggr).
\end{split}
\end{equation*}
\end{definition}

\begin{lemma}
\label{lem:extp}
Let $G$ be a graph and let $p,r$ be integers.

If $G$ has the $(p+r)$-extension property then every
formula $\phi$ with $p$ free variables and quantifier rank $r$
is equivalent, in $G$, with a quantifier free formula. 
\end{lemma}
\begin{proof}
Let $\phi$ be a formula with $p$ free variables and quantifier rank $r$.
Let $(a_1,\dots,a_p)$ and $(b_1,\dots,b_p)$ be two $p$-tuples of
vertices of $G$ such that $a_i\mapsto b_i$ is a partial isomorphism.
The $(p+r)$-extension properties allows to easily play a $r$-turns 
back-and-forth game between $(G,a_1,\dots,a_p)$ and $(G,b_1,\dots,b_p)$,
 thus proving that  $(G,a_1,\dots,a_p)$ and $(G,b_1,\dots,b_p)$ are
 $r$-equivalent. It follows that $G\models\phi(a_1,\dots,a_p)$ if
 and only if $G\models\phi(b_1,\dots,b_p)$. Following the lines 
 of Lemma~\ref{lem:UHQE}, we deduce that there exists a quantifier
 free formula $\widetilde\phi$ such that $G\models\phi\leftrightarrow\widetilde\phi$. 
\myqed\end{proof}

We now prove that random graphs converge elementarily to the countable random graphs.
\begin{lemma}
\label{lem:ERado}
Let $1/2>\delta>0$. Assume that for every positive integer $n\geq 2$
and every $1\leq i<j\leq n$, $p_{n,i,j}\in[\delta,1-\delta]$. Assume
that for each $n\in\bbbn$, $G_n$ is a random graph on $[f(n)]$ where $f(n)\geq n$, and
where $i$ and $j$ are adjacent with probability $p_{n,i,j}$ (all 
these events being independent). Then the sequence $(G_n)_{n\in\bbbn}$ almost
surely converges elementarily to the Rado graph. 
\end{lemma}
\begin{proof}
Let $p\in\bbbn$ and let $\alpha=\delta(1-\delta)$.
 The probability that $G_n\models\Upsilon_p$ is at least
$1-(1-\alpha^p)^{f(n)}$. It follows that for $N\in\bbbn$ 
the probability that all the graphs $G_n$ ($n\geq N$) satisfy
$\Upsilon_p$ is at least $1-\alpha^{-p}(1-\alpha^p)^{f(N)}$.
According to Borel-Cantelli lemma,  the probability that 
$G_n$ does not satisfy $\Upsilon_p$ infinitely many is zero.
As this holds for every integer $p$, it follows that, 
with high probability,
 every elementarily converging subsequence
of $(G_n)_{n\in\bbbn}$ converges to the Rado graph hence,
with high probability, $(G_n)_{n\in\bbbn}$ converges  
elementarily to the Rado graph.
\myqed\end{proof}
Thus we get:
\begin{mytinted}
\begin{theorem}
\label{thm:rand}
Let $0<p<1$ and let $G_n\in G(n,p)$ be independent random graphs with edge
probability $p$. Then $(G_n)_{n\in\bbbn}$ is almost surely ${\rm FO}$-convergent. 
\end{theorem}
\end{mytinted}
\begin{proof}
This is an immediate consequence of Lemma~\ref{lem:FOQF}, Lemma~\ref{lem:ERado} and
the easy fact that $(G_n)_{n\in\bbbn}$ is almost surely ${\rm QF}$-convergent. 
\myqed\end{proof}

\begin{mytinted}
\begin{theorem}
For every $\phi\in{\rm FO}_p$ there exists a polynomial
$P_\phi\in\bbbz[X_1,\dots,X_{\binom{p}{2}}]$ such that for every
sequence $(G_n)_{n\in\bbbn}$ of finite graphs that converges elementarily
to the Rado graph the following holds:

If $(G_n)_{n\in\bbbn}$ is {\rm L}-convergent to some graphon $W$ then
$$
\lim_{n\rightarrow\infty}\langle\phi,G_n\rangle=\idotsint P_\phi((W_{i,j}(x_i,x_j))_{1\leq i<j\leq p})\,
{\rm d}x_1\dots{\rm d}x_p.
$$
\end{theorem}
\end{mytinted}
 \begin{proof}
 Assume the sequence $(G_n)_{n\in\bbbn}$ is elementarily convergent to the Rado graph and that it is
 {\rm L}-convergent to some graphon $W$.
 
According to Lemma~\ref{lem:UHQE}, there exists a quantifier
 free formula $\widetilde\phi$ such that  
 
 $$G\models (\forall x_1\dots x_p)\
 \phi(x_1,\dots,x_p)\leftrightarrow\widetilde\phi(x_1,\dots,x_p)$$
 (hence $\Omega_\phi(G)=\Omega_{\widetilde\phi}(G)$) holds when $G$ is
 the Rado graph.  
 As $(G_n)_{n\in\bbbn}$ is elementarily convergent to the Rado graph, this sentence holds for all but finitely many graphs $G_n$.
  Thus for all but finitely many $G_n$ it holds that 
  $\langle\phi,G_n\rangle=\langle\widetilde\phi,G_n\rangle$.
  Moreover, according to Lemma~\ref{lem:FOQF}, the sequence $(G_n)_{n\in\bbbn}$ is ${\rm FO}$-convergent and
  thus the following equation holds:
  $$
\lim_{n\rightarrow\infty}\langle\phi,G_n\rangle=\lim_{n\rightarrow\infty}\langle\widetilde\phi,G_n\rangle.
$$  
 
 By using an inclusion/exclusion argument and the general form of the density of homomorphisms of fixed target
 graphs to a graphon we deduce that there exists a  polynomial
$P_\phi\in\bbbz[X_1,\dots,X_{\binom{p}{2}}]$ (which depends only on $\phi$) such that
$$
\lim_{n\rightarrow\infty}\langle\widetilde\phi,G_n\rangle=\idotsint P_\phi((W_{i,j}(x_i,x_j))_{1\leq i<j\leq p})\,
{\rm d}x_1\dots{\rm d}x_p.
$$ 
The theorem follows.
 \end{proof}

Although elementary convergence to Rado graph seems quite
a natural assumption for graphs which are neither too sparse nor
too dense, elementary convergence to other ultra-homogeneous graphs
may be problematic. 
\begin{example}\label{ex:triangle_free}
Cherlin \cite{Cherlin2011} posed 
the problem whether there is a finite $k$-saturated 
triangle-free graph, for each $k\in\bbbn$, where 
a triangle free graph is called {\em $k$-saturated}
 if for every set $S$ of at most $k$ vertices, and for 
 every independent subset $T$ of $S$, there exists a
 vertex adjacent to each vertex of $T$ and to no 
 vertex of $S-T$.
 In other words, Cherlin asks whether the generic 
 countable triangle-free graph has the finite model 
 property, that is if it is an elementary limit
 of a sequence of finite graphs. 
 See \cite{hubivcka2015universal} for more ultra-homogeneous structures defined by forbidden substructures.
\end{example}

It is possible to extend Lemma~\ref{lem:FOQF} to 
sequences of graph having a non ultra-homogeneous elementary
limit if we restrict ${\rm FO}$ to a smaller fragment. 
An example is the following:

\begin{example}
A graph $G$ is {\em IH-Homogeneous} \cite{Cameron2006} if every partial finite isomorphism
extends into an endomorphism.
Let ${\rm PP}$ be the fragment of ${\rm FO}$ that consists 
into {\em primitive positive} formulas, that is formulas formed using adjacency, equality, conjunctions and
existential quantification only, and let ${\rm BA}({\rm PP})$ be the minimum sub-Boolean algebra of 
${\rm FO}$ containing ${\rm PP}$. 

Following the lines of Lemma~\ref{lem:FOQF} and using
 Theorem~\ref{thm:LQF} and Proposition~\ref{lem:extBA}, one proves that
 if a sequence of graphs $(G_n)_{n\in\bbbn}$ converges elementarily to 
 some IH-homogeneous infinite countable graph then
  $(G_n)_{n\in\bbbn}$ is ${\rm BA}({\rm PP})$-convergent if and only if
it is ${\rm QF}$-convergent.
\end{example}
\subsection{${\rm FO}$-convergence of Graphs with Bounded Maximum Degree}

We now consider how full ${\rm FO}$-convergence differs to BS-convergence
for sequence of graphs with maximum degree at most $D$.
As a corollary of Theorems~\ref{thm:fole} and~\ref{thm:BS} we have:
\begin{mytinted}
\begin{corollary}
\label{cor:BSE}
A sequence $(G_n)$ of finite graphs with maximum degree at most $d$ such that
$\lim_{n\rightarrow\infty}|G_n|=\infty$ is
${\rm FO}$-convergent if and only if it is both 
BS-convergent and elementarily convergent.
\end{corollary}
\end{mytinted}
\section{Interpretation Schemes}
In the process of this research we discovered the increasing role played by interpretations. They are described in this section.
\subsection{Continuous Functions and Interpretations}
Let $X$ and $Y$ be fragments of ${\rm FO}(\kappa)$ and ${\rm FO}(\sig)$, respectively. 
Let $f:S(\mathcal B(X))\rightarrow S(\mathcal B(Y))$.
A function 
$f:S(\mathcal B(X))\rightarrow S(\mathcal B(Y))$ is continuous if and only if the inverse image of an open subset
of $S(\mathcal B(Y))$ is an open subset of $S(\mathcal B(X))$. In the case of Stone spaces (where clopen subsets generates the topology), we can further restrict our attention to clopen subsets: $f$ will be continuous if the inverse image of a clopen subset is a clopen subset. In other words, $f$ is continuous if there exists $f^*:\mathcal B(Y)\rightarrow\mathcal B(X)$, such that for every $\phi\in Y$, the following equation holds:
$$
f^{-1}(K(\phi))=K(f^*(\phi)).
$$
It follows that if $f$ is continuous then
for every $X$-convergent sequence $(\mathbf A_n)_{n\in\bbbn}$, the sequence $(f(\mathbf A_n))_{n\in\bbbn}$ is 
$Y$-convergent.
Note that $f^*$ will be a homomorphism from $\mathcal B(Y)$ to $\mathcal B(X)$, and that the duality between $f$ and $f^*$ is
nothing more the duality between Stone spaces and Boolean algebras.

The above property can be sometimes restated in terms of definable sets in structures.
For a fragment $X$ of ${\rm FO}$ and a relational structure $\mathbf A$, a subset
$F\subseteq A^p$ is {\em $X$-definable} if there exists a formula $\phi\in X$ with free variables
$x_1,\dots,x_p$ such that
$$F=\Omega_\phi(\mathbf A)=\{(v_1,\dots,v_p)\in A^p:\ \mathbf A\models\phi(v_1,\dots,v_p)\}.
$$

Let $\mathbf A$ be a $\kappa$-structure, let $\mathbf B$ be a $\sig$-structure, and let
$g:A^k\rightarrow B$ be surjective. 
Assume that there exists a function
$g^*:Y\rightarrow X$ such that for every $\phi\in Y$ with free variables
$x_1,\dots,x_p$ ($p\geq 0$), and every $v_{i,j}\in A$  ($1\leq i\leq p$, $1\leq j\leq k$) the following holds:

\begin{gather*}
\mathbf B\models \phi(g(v_{1,1},\dots,v_{1,k}),\dots,g(v_{p,1},\dots,v_{p,k}))\\
\iff\\
\mathbf A\models g^*(\phi)(v_{1,1},\dots,v_{1,k},\dots,v_{p,1},\dots,v_{p,k})
\end{gather*}

\noindent then $g^*$ is a homomorphism, and thus it defines a continuous function from $S(\mathcal B(X))$ to
$S(\mathcal B(Y))$. Note that the above formula can be restated as
$$\Omega_{g^*(\phi)}(\mathbf A)=\hat g^{-1}(\Omega_{\phi}(\mathbf B)),$$
 where
$$\hat g((v_{1,1},\dots,v_{p,k}))=(g(v_{1,1},\dots,v_{1,k}),\dots,g(v_{p,1},\dots,v_{p,k})).$$
In other words, the inverse image of a $Y$-definable set of $\mathbf B$
is an $X$-definable set of $\mathbf A$.

When $X={\rm FO}(\kappa)$ and $Y={\rm FO}(\sig)$, the property that the inverse image of a first-order definable
set of $\mathbf B$ is a first-order definable set of $\mathbf A$ leads to the model theoretical notion of interpretation
(without parameters) of $\mathbf B$ in $\mathbf A$. We recall now the formal definition of an interpretation.
\begin{definition}[Interpretation]
An {\em interpretation} of $\mathbf B$ in $\mathbf A$ {\em with parameters} (or {\em without parameters}, respectively) with {\em exponent} $k$ is a surjective map from a subset of $A^k$ onto $B$ such that the inverse image
 of every set $X$ definable in $\mathbf B$ by a first-order formula without parameters 
is definable in $\mathbf A$ by a first-order formula with parameters (or without parameters, respectively).
\end{definition}

\subsection{Interpretation Schemes}
The main drawback of interpretations is that they only concerns two specific structures $\mathbf A$ and $\mathbf B$.
However, it is frequent that interpretations naturally generalize to a family of interpretations of $\sig$-structures 
in $\kappa$-structures with the same associated homomorphism of Boolean algebras. Moreover, this homomorphism is
uniquely defined by the way it transforms each relation in $\sig$ (including equality) into a formula in $\kappa$ and by the formula which defines the domain of the $\kappa$-structures.
This can be formalized as follows.

\begin{definition}[Interpretation Scheme]
\label{def:interp}
Let $\kappa,\sig$ be signatures, where $\sig$ has 
$q$ relational symbols $R_1,\dots,R_q$ with respective arities $r_1,\dots,r_q$.

An {\em interpretation scheme} ${\mathsf I}$ of $\sig$-structures
in $\kappa$-structures is defined by an integer $k$ --- the {\em exponent} of the interpretation scheme --- a formula
$E\in{\rm FO}_{2k}(\kappa)$, a formula $\theta_0\in{\rm FO}_k(\kappa)$, and
 a formula 
$\theta_i\in{\rm FO}_{r_ik}(\kappa)$ for each symbol $R_i\in\sig$, such that:
\begin{itemize}
  \item the formula $E$ defines an equivalence relation
  of $k$-tuples;
   \item each formula $\theta_i$ is compatible with $E$, in the sense that
   for
   every $0\leq i\leq q$ the following property holds:
    $$
  \bigwedge_{1\leq j\leq r_i}\,E(\mathbf x_j,\mathbf y_j)\quad\entails\quad
  \theta_i(\mathbf x_1,\dots,\mathbf x_{r_i})\leftrightarrow\theta_i(\mathbf
  y_1,\dots,\mathbf y_{r_i}),
   $$
   where $r_0=1$, boldface $\mathbf x_j$ and $\mathbf y_j$ represent
   $k$-tuples of free variables, and 
   where $\theta_i(\mathbf x_1,\dots,\mathbf x_{r_i})$ stands for
 $\theta_i(x_{1,1},\dots,x_{1,k},\dots,x_{r_i,1},\dots,x_{r_i,k})$.
    
\end{itemize}

For a $\kappa$-structure $\mathbf A$, we denote by $\mathsf{I}(\mathbf A)$ the
$\sig$-structure $\mathbf B$ defined as follows:
 \begin{itemize}
   \item the domain $B$ of $\mathbf B$ is the subset
   of the $E$-equivalence classes $[\mathbf x]\subseteq A^k$  of the tuples $\mathbf x=(x_1,\dots,x_k)$ 
  such that $\mathbf A\models \theta_0(\mathbf x)$;
   \item for each $1\leq i\leq q$ and every 
   $\mathbf v_1,\dots,\mathbf v_{s_i}\in A^{kr_i}$ such that
   $\mathbf A\models\theta_0(\mathbf v_j)$ (for every $1\leq j\leq r_i$) the following equivalence
   holds:
   $$
   \mathbf B\models R_i([\mathbf v_1],\dots,[\mathbf v_{r_i}])\quad\iff\quad
   \mathbf A\models \theta_i(\mathbf v_1,\dots,\mathbf v_{r_i}).
   $$  
 \end{itemize}
\end{definition}

From the standard properties of model theoretical interpretations
(see, for instance \cite{Lascar2009} p.~180), we state the following: if
$\mathsf I$ is an interpretation of $\sig$-structures in $\kappa$-structures,
then there exists a mapping
$\tilde{\mathsf I}:{\rm FO}(\sig)\rightarrow {\rm FO}(\kappa)$ (defined by
means of the formulas $E,\theta_0,\dots,\theta_q$ above) such that
for every $\phi\in {\rm FO}_p(\sig)$, and every $\kappa$-structure $\mathbf A$,
 the following property holds (while letting $\mathbf B=\mathsf I(\mathbf A)$
 and identifying elements of $B$ with the corresponding equivalence classes of $A^k$):

For every $[\mathbf v_1],\dots,[\mathbf v_p]\in B^{p}$ (where $\mathbf v_i=(v_{i,1},\dots,v_{i,k})\in A^k$)
the following equivalence holds:
$$
  \mathbf B\models \phi([\mathbf v_1],\dots,[\mathbf v_p])\quad\iff\quad \mathbf
  A\models \tilde{\mathsf I}(\phi)(\mathbf v_1,\dots,\mathbf v_p).$$
It directly follows from the existence of the mapping $\tilde{\mathsf I}$ that
an interpretation scheme ${\mathsf I}$ of $\sig$-structures in
$\kappa$-structures defines a continuous mapping from $S(\mathcal B({\rm
FO}(\kappa)))$ to $S(\mathcal B({\rm FO}(\sig)))$.
Thus, interpretation schemes have the following general property:
\begin{mytinted}
\begin{proposition}
\label{prop:intlim}
Let $\mathsf I$ be an interpretation scheme  of $\sig$-structures
in $\kappa$-structures.

Then, if a sequence $(\mathbf A_n)_{n\in\bbbn}$ of finite $\kappa$-structures
is ${\rm FO}$-convergent then the sequence 
$(\mathsf{I}(\mathbf A_n))_{n\in\bbbn}$ of (finite) $\sig$-structures
is ${\rm FO}$-convergent. 
\end{proposition}
\end{mytinted}

We shall be mostly interested in very specific and simple types of interpretation schemes.

\begin{definition}
Let $\kappa,\sig$ be signatures.
A {\em basic interpretation scheme} ${\mathsf I}$ of $\sig$-structures
in $\kappa$-structures  with {\em exponent} $k$ is defined by a formula 
$\theta_i\in{\rm FO}_{kr_i}(\kappa)$ for each symbol $R_i\in\sig$ with arity
$r_i$. 

For a $\kappa$-structure $\mathbf A$, we denote by $\mathsf{I}(\mathbf A)$ the
structure with domain $A^k$ such that, for every $R_i\in\sig$ with arity $r_i$ and every
$\mathbf v_1,\dots,\mathbf v_{r_i}\in A^k$ the following equivalence holds:
$$
\mathsf I(\mathbf A)\models R_i(\mathbf v_1,\dots,\mathbf v_{r_i})\quad\iff\quad
\mathbf A\models\theta_i(\mathbf v_1,\dots,\mathbf v_{r_i}).
$$
\end{definition}

It is immediate that every basic interpretation scheme $\mathsf I$ defines
a mapping $\tilde{\mathsf I}:{\rm FO}(\sig)\rightarrow{\rm FO}(\kappa)$ such
that for every $\kappa$-structure $\mathbf A$, every $\phi\in{\rm FO}_p(\sig)$,
and every $\mathbf v_1,\dots,\mathbf v_p\in A^k$ the following equivalence holds:
$$
\mathsf I(\mathbf A)\models \phi(\mathbf v_1,\dots,\mathbf v_{p})\quad\iff\quad
\mathbf A\models\tilde{\mathsf I}(\phi)(\mathbf v_1,\dots,\mathbf v_{p})
$$
and
$$
{\rm qrank}(\tilde{\mathsf I}(\phi))\leq k({\rm qrank}(\phi)+\max_i{\rm qrank}(\theta_i)).
$$

It follows that for every $\kappa$-structure $\mathbf A$, every $\phi\in{\rm FO}_p(\sig)$,
the following equation holds:
$$\Omega_{\phi}(\mathsf I(\mathbf A))=\Omega_{\tilde{\mathsf I}(\phi)}(\mathbf A).$$
In particular, if $\mathbf A$ is a finite structure, the following equation holds:
$$\langle\phi,\mathsf I(\mathbf A)\rangle=\langle\tilde{\mathsf I}(\phi),\mathbf A\rangle.$$

\chapter{Modelings for Sparse Structures}
\section{Relational Samples Spaces}
\label{sec:rss}
The notion of {\em \rss} is a strenghtening of the one of relational structure, where it is required that the domain
shall be endowed with a suitable structure of  a (nice) measurable space. 

\subsection{Definition and Basic Properties}
\begin{definition}
Let $\sig$ be a signature.
A {\em $\sig$-\rss} is a $\sig$-structure $\mathbf A$, 
whose domain $A$ is a standard Borel space
 with the property that every first-order definable subset of $A^p$ is measurable.  
Precisely, for every integer $p$, and every $\phi\in{\rm FO}_p(\sig)$,
denoting
$$
\Omega_\phi(\mathbf A)=\{(v_1,\dots,v_p)\in A^p:\ \mathbf
A\models\phi(v_1,\dots,v_p)\},
$$
it holds that $\Omega_\phi(\mathbf A)\in\Sigma_{\mathbf A}^p$,
where $\Sigma_{\mathbf A}$ is the
Borel $\sigma$-algebra of $A$.
\end{definition}
Note, that in the case of graphs, every {\rss} 
is a Borel graph (that is a graph whose vertex set is a standard Borel space and whose edge set is Borel), but the
converse is  not true.

\begin{mytinted}
\begin{lemma}
\label{lem:rss}
Let $\sig$ be a signature, let $\mathbf A$ be a 
$\sig$-structure, whose domain $A$ is a standard Borel 
space with $\sigma$-algebra $\Sigma_{\mathbf A}$.

Then the following conditions are equivalent:
\begin{enumerate}[(a)]
  \item\label{enum:rss1} $\mathbf A$ is a $\sig$-{\rss};
  \item\label{enum:rss2} for every integer $p\geq 0$ and every $\phi\in{\rm
  FO}_p(\sig)$, it holds that $\Omega_\phi(\mathbf A)\in\Sigma_{\mathbf A}^p$; 
\item\label{enum:rss3} for every integer $p\geq 1$ and every $\phi\in{\rm
FO}_p^{\rm local}(\sig)$, it holds that $\Omega_\phi(\mathbf A)\in\Sigma_{\mathbf A}^p$;
\item\label{enum:rss4} for all integers $p,q\geq 0$, every $\phi\in{\rm
FO}_{p+q}(\sig)$, and every $a_1,\dots,a_q\in A^q$ the set
   $$
   \{(v_1,\dots,v_p)\in A^p:\ \mathbf
A\models\phi(a_1,\dots,a_q,v_1,\dots,v_p)\}
$$
belongs to $\Sigma_{\mathbf A}^p$.
\end{enumerate}
\end{lemma}
\end{mytinted}
\begin{proof}
Items (\ref{enum:rss1}) and (\ref{enum:rss2}) are equivalent by definition.
Also we obviously have the implications 
(\ref{enum:rss4}) $\Rightarrow$ (\ref{enum:rss2}) $\Rightarrow$
(\ref{enum:rss3}). That (\ref{enum:rss3}) $\Rightarrow$ (\ref{enum:rss2})
is a direct consequence of Gaifman locality theorem, and
the implication (\ref{enum:rss2}) $\Rightarrow$ (\ref{enum:rss4})
is a direct consequence of Fubini-Tonelli theorem.
 \myqed\end{proof}

\begin{lemma}
\label{lem:cc}
Let $\mathbf A$ be a {\rss}, let $a\in A$, and let 
$\mathbf A_a$ be the connected component of $\mathbf A$ containing $a$.
 
Then $\mathbf A_a$ has a measurable domain and, equipped with the
$\sigma$-algebra of the Borel sets of $A$ included in $A_a$, it is
a {\rss}.
\end{lemma}
\begin{proof}
Let $\phi\in{\rm FO}_p^{\rm local}$ and let 
$$X=\{(v_1,\dots,v_p)\in A_a^p:\ \mathbf A_a\models\phi(v_1,\dots,v_p)\}.$$
As $\phi$ is local, there is an integer $D$ such that the satisfaction of $\phi$
only depends on the $D$-neighborhoods of the free variables. 

For every integer $n\in\bbbn$, denote by $B(\mathbf A,a,n)$ the substructure
of $\mathbf A$ induced by all vertices at distance at most $n$ from $a$. 
By the locality of $\phi$, for every $v_1,\dots,v_p$ at distance at most $n$
from $a$ the following equivalence holds:
$$
\mathbf A_a\models\phi(v_1,\dots,v_p)\quad\iff\quad
B(\mathbf A,a,n+D)\models\phi(v_1,\dots,v_p).
$$
However, it is easily checked that there is a local first-order formula
$\varphi_n\in{\rm FO}_{p+1}^{\rm local}$ such that for every 
$v_1,\dots,v_p$ the following equivalence holds:
$$
B(\mathbf A,a,n+D)\models\phi(v_1,\dots,v_p)\wedge\bigwedge_{i=1}^p{\rm
dist}(a,v_i)\leq n
\quad\iff\quad \mathbf A\models\varphi_n(a,v_1,\dots,v_p).$$

By Lemma~\ref{lem:rss}, it follows that  the set 
$X_n=\{(v_1,\dots,v_n)\in A:\ \mathbf A\models\varphi_n(a,v_1,\dots,v_p)\}$ is
measurable. As $X=\bigcup_{n\in\bbbn}X_n$, we deduce that 
$X$ is measurable (with respect to $\Sigma_{\mathbf A}^p$).
In particular, $A_a$ is a Borel subset of $A$ hence
$A_a$, equipped with the $\sigma$-algebra $\Sigma_{\mathbf A_a}$ of the Borel
sets of $A$ included in $A_a$, is a standard Borel set. Moreover, it is immediate that a subset
of $A_a^p$ belongs to  $\Sigma_{\mathbf A_a}^p$ if and only if it belongs to
 $\Sigma_{\mathbf A}^p$. Hence, every subset of $A_a^p$ defined  
 by a local formula is measurable with respect to $\Sigma_{\mathbf A_a}^p$.
By Lemma~\ref{lem:rss}, it follows that $\mathbf A_a$ is a
{\rss}.
\myqed\end{proof}

\subsection{Interpretations of Relational Sample Spaces}
An elementary interpretation with parameter amounts to 
distinguishing a single element, the parameter, by adding a new unary symbol to the signature (e.g.\ representing a root).

\begin{lemma}
\label{lem:mark}
Let $\mathbf A$ be a $\sig$-{\rss}, let $\sig^+$ be the signature obtained
from $\sig$ by adding a new unary symbol $M$ and let $\mathbf A^+$ be
obtained from $\mathbf A$ by marking a single $a\in A$ (i.e. $a$ is the only 
element $x$ of $A^+=A$ such that $\mathbf A^+\models M(x)$).

Then $\mathbf A^+$ is a \rss.
\end{lemma}
\begin{proof}
Let $\phi\in{\rm FO}_p(\sig^+)$. There exists $\phi'\in{\rm FO}_{p+1}(\sig)$ such
that for every $x_1,\dots,x_p\in A$ the following equivalence holds:
$$
\mathbf A^+\models\phi(x_1,\dots,x_p)\quad\iff\quad
\mathbf A\models\phi(a,x_1,\dots,x_p).
$$
According to Lemma~\ref{lem:rss}, the set of all
$(x_1,\dots,x_p)$ such that $\mathbf A\models\phi(a,x_1,\dots,x_p)$
is measurable. It follows that $\mathbf A^+$ is a \rss.
\myqed
\end{proof}

\begin{mytinted}
\begin{lemma}
\label{lem:interp}
Every injective first-order interpretation (with or without parameters) of a {\rss} is a {\rss}.

Precisely, if $f$ is an injective first-order interpretation of a $\sig$-structure $\mathbf B$ in a 
$\kappa$-{\rss} $\mathbf A$ and if we define
$$\Sigma_{\mathbf B}=\{X\subseteq B: f^{-1}(X)\in\Sigma_{\mathbf A}^k\},$$
then $(\mathbf B,\Sigma_{\mathbf B})$ is a \rss.
\end{lemma}
\end{mytinted}
\begin{proof}
According to Lemma~\ref{lem:mark}, we can first mark all the parameters and reduce to the case where
the interpretation has no parameters.

Let $D$ be the domain of $f$. 
As $B$ is first-order definable in $\mathbf B$, $D$ is first-order definable in $\mathbf A$ hence
$D\in\Sigma_{\mathbf A}^k$. 
Then $D$ is a Borel sub-space of $A^k$. As $f$ is a bijection from $D$ to $B$, we deduce
that $(B,\Sigma_{\mathbf B})$ is a standard Borel space.

Moreover, as the inverse image of every first-order definable set of $\mathbf B$ is first-order definable in $\mathbf A$, we
deduce that $(\mathbf B,\Sigma_{\mathbf B})$ is a $\sig$-\rss. 
\end{proof}

\subsection{Disjoint union}
Let $\mathbf H_i$ be $\sig$-{\rss}s for $i\in I\subseteq\bbbn$.
We define the {\em disjoint union} 
$${\mathbf H}=\coprod_{i\in I}\mathbf H_i$$
of the $\mathbf H_i$'s as the relational structure, which is the disjoint union of the $\mathbf H_i$'s
endowed with the $\sigma$-algebra  $\Sigma_{\mathbf H}=\{\bigcup_i X_i:\ X_i\in\Sigma_{\mathbf H_i}\}$.

\begin{mytinted}
\begin{lemma}
\label{lem:unionrss}
Let $\mathbf H_i$ be $\sig$-{\rss}s for $i\in I\subseteq\bbbn$. Then
$\mathbf H=\coprod_{i\in I}\mathbf H_i$ is a $\sig$-\rss, in which every
$H_i$ is measurable.
\end{lemma}
\end{mytinted}
\begin{proof}
%
We consider the signature $\sig^+$ obtained from $\sig$ by adding a new binary
relation $\varpi$, and the basic interpretation scheme $\mathsf I_1$ of
$\sig^+$-structures in $\sig$-structures corresponding to the addition of the
new relation $\varpi$ by the formula $\theta_\varpi=1$. This means that for
every $\sig$-structure $\mathbf A$ it holds that $\mathsf I_1(\mathbf
A)\models(\forall x,y)\ \varpi(x,y)$.
Let $\mathbf H_i^+=\mathsf I_1(\mathbf H_i)$. 

Let $\mathbf H^+=\coprod_{i\in I}\mathbf H_i^+$.
Clearly, $\Sigma_{\mathbf H^+}=\Sigma_{\mathbf H}$ and $(H,\Sigma_{\mathbf H})$ is a standard Borel space.
Moreover, by construction, each $H_i$ is measurable.

Let $\phi\in{\rm FO}_{p}(\sig)$. 
First notice that for every $(v_1,\dots,v_p)\in H^{p+q}$ (which is also $(H^+)^{p+q}$) it holds that 
$\Omega_\phi(\mathbf H)=\Omega_\phi(\mathbf H^+)$, that is:
$$\mathbf H\models\phi(v_1,\dots,v_p)\quad\iff\quad\mathbf H^+\models\phi(v_1,\dots,v_p).$$
It follows from Lemma~\ref{lem:sl} that the set $\Omega_\phi(\mathbf H^+)$ may
be obtained by Boolean operations, products, and coordinate permutations from sets defined by $\varpi$-local formulas
(which we introduced in Section~\ref{sec:piloc}).
As all these operations preserve measurability, we can assume that $\phi$ is $\varpi$-local.
 Then $\Omega_\phi(\mathbf H^+)$ is the union
 of the sets $\Omega_\phi(\mathbf H_i)$.
All these sets are measurable (as $\mathbf H_i$ is a {\rps}) thus their union is measurable 
 (by construction of $\Sigma_{\mathbf H}$). It follows that $\mathbf H^+$ is a {\rss}, and so
 is $\mathbf H$ (every first-order definable set of $\mathbf H$ is first-order definable in $\mathbf H^+$).
 \myqed\end{proof}
\section{Modelings}
\label{sec:rps}
We introduced a notion of limit objects --- called modelings --- for sequences of sparse graphs and structures, which is a natural generalization of graphings.
These limit objects are defined by considering a probability measure on a \rss.
In this section, we show that the most we can expect
is that {\rps}s are limit objects for sequence of sparse structures, and we conjecture that an unavoidable qualitative jump
occurs for notions of limit structures, which coincides with the nowhere dense/somewhere dense frontier (see Conjecture~\ref{conj:nd}). 

\subsection{Definition and Basic Properties}
Recall Definitions~\ref{def:rps}
and~\ref{def:bracket2}: a {\em $\sig$-\rps} $\mathbf A$ is a 
$\sig$-{\rss} equipped with a probability measure (denoted $\nu_{\mathbf A}$), and
the {\em Stone pairing} of $\phi\in{\rm FO}(\sig)$ and a $\sig$-modeling $\mathbf A$ is
$\langle \phi,\mathbf A\rangle=\nu_{\mathbf A}^p(\Omega_\phi(\mathbf A))$.
Notice that it follows (by Fubini's theorem) that it holds that
 
 \begin{align*}
 \langle \phi,\mathbf A\rangle&=
\int_{\mathbf x\in
 A^p}1_{\Omega_\phi(\mathbf A)}(\mathbf x)\,{\rm d}\nu_{\mathbf A}^p(\mathbf
 x)\\
&=\idotsint \mathbf 1_{\Omega_\phi(\mathbf A)}(x_1,\dots,x_p)\ 
{\rm d}\nu_{\mathbf A}(x_1)\,\dots\,{\rm d}\nu_{\mathbf A}(x_p).
 \end{align*}

Then, generalizing Definition~\ref{def:modFOlim}, we extend the notion of $X$-convergence to modelings:

\begin{mytinted}
\begin{definition}[{\rps} $X$-limit]
Let $X$ be a fragment of ${\rm FO}(\sig)$.

If an $X$-convergent sequence $(\mathbf A_n)_{n\in\bbbn}$ of $\sig$-{\rps}s
satisfies 
$$(\forall\phi\in X)\quad\langle\phi,\mathbf
L\rangle=\lim_{n\rightarrow\infty}\langle\phi,\mathbf A_n\rangle $$
for some $\sig$-{\rps} $\mathbf L$, then we say
that $\mathbf L$ is  a  {\em {\rps} $X$-limit}
of $(\mathbf A_n)_{n\in\bbbn}$.
\end{definition}
\end{mytinted}

A $\sig$-{\rps} $\mathbf A$ is {\em weakly uniform} if all the
singletons of $A$ have the same measure. Clearly, every finite $\sig$-structure
$\mathbf A$ can be identified with the weakly uniform modeling obtained by
considering the discrete topology on $A$. This identification is clearly
consistent with our definition of the Stone pairing of a formula and a {\rps}.

In the case where a {\rps} $\mathbf A$ has an infinite domain, the condition
for $\mathbf A$ to be weakly uniform is equivalent to the condition for
$\nu_{\mathbf A}$ to be atomless. This property is usually fulfilled by
{\rps} $X$-limits of sequences of finite structures.

\begin{lemma}
\label{lem:wu}
Let $X$ be a fragment of ${\rm FO}$ that includes ${\rm FO}_0$ and
the formula $(x_1=x_2)$.
Then every {\rps} $X$-limit of weakly uniform
modelings is weakly uniform.
\end{lemma}
\begin{proof}
Let $\phi$ be the formula $(x_1=x_2)$. Notice that for every finite 
$\sig$-structure $\mathbf A$
it holds that $\langle\phi,\mathbf A\rangle=1/|A|$ and
that for every infinite weakly uniform $\sig$-structure it holds that 
$\langle\phi,\mathbf A\rangle=0$. 

Let $\mathbf L$ be a {\rps} $X$-limit of a sequence
$(\mathbf A_n)_{n\in\bbbn}$.
Assume $\lim_{n\rightarrow\infty}|A_n|=\infty$.
Assume for contradiction that $\nu_{\mathbf L}$ has an atom $\{v\}$ (i.e. $\nu_{\mathbf L}(\{v\})>0$). Then 
$\langle\phi,\mathbf L\rangle\geq \nu_{\mathbf L}(\{v\})^2>0$, contradicting 
$\lim_{n\rightarrow\infty}\langle\phi,\mathbf A_n\rangle=0$. Hence $\nu_{\mathbf L}$ is atomless.

Otherwise, $|L|=\lim_{n\rightarrow\infty}|A_n|<\infty$ (as $\mathbf L$ is an elementary limit of $(\mathbf A_n)_{n\in\bbbn}$).
Let $N=|L|$.
Label $v_1,\dots,v_N$ the elements of $L$ and let $p_i=\nu_{\mathbf L}(\{v_i\})$. Then
\begin{align*}
\frac{1}{N}\sum_{i=1}^N p_i^2-\Bigl(\frac{1}{N}\sum_{i=1}^N p_i\Bigr)^2&=\frac{\langle\phi,\mathbf L\rangle}{N}-\frac{1}{N^2}\\
&=\frac{\lim_{n\rightarrow\infty}\langle\phi,\mathbf A_n\rangle}{N}-\frac{1}{N^2}\\
&=0
\end{align*}
Thus $p_i=1/N$ for every $i=1,\dots,N$.
\myqed\end{proof}
\begin{corollary}
Every {\rps} ${\rm FO}_2^{\rm local}$-limit of finite
structures is weakly uniform.
\end{corollary}

\begin{lemma}
\label{lem:L2W}
Let $X$ be a fragment that includes all quantifier free formulas.

Assume $\mathbf L$ is a modeling $X$-limit of a sequence $(G_n)_{n\in\bbbn}$ of graphs with $|G_n|\rightarrow\infty$.
Let $\overline{\nu}_{\mathbf L}$ be the completion of the measure $\nu_{\mathbf L}$.

Then there is at least one mod $0$ isomorphism $f:[0,1]\rightarrow (L,\overline{\nu}_{\mathbf L})$, and
for every such $f$ the graphon $W$ defined by
$$W(x,y)=\mathbf 1_{\Omega_{(x_1\sim x_2)}(\mathbf L)}(f(x),f(y))$$
(for $x,y$ in the domain of $f$, and $W(x,y)=0$ elsewhere)
is a random-free graphon {\rm L}-limit of $(G_n)_{n\in\bbbn}$.
\end{lemma}
\begin{proof}
Considering the formula $x_1=x_2$, we infer that $\nu_{\mathbf L}$ is atomless.
This measure is also atomless and turns $L$ into a standard probability space. 
According to the isomorphism theorem, all atomless standard probability spaces are mutually mod $0$ isomorphic hence
there is at least one mod $0$ isomorphism $f:[0,1]\rightarrow (L,\overline{\nu}_{\mathbf L})$ ($[0,1]$ is considered with Lebesgue measure).

Fix such a mod $0$ isomorphism $f$, defined on $[0,1]\setminus N_1$, with value on $L\setminus N_2$ (where $N_1$ and $N_2$ are nullsets). 
 For every Borel measurable function $g:L^n\rightarrow [0,1]$, define $g_f$ by
$g_f(x_1,\dots,x_n)=g(f(x_1),\dots,f(x_n))$ if $x_i\notin N_1$ for every $1\leq i\leq n$ and 
$g_f(x_1,\dots,x_n)=0$ otherwise. Then
 it holds that
\begin{align*}
\int_{[0,1]^n}g_f(x_1,\dots,x_n)\,{\rm d}x_1\dots{\rm d}x_n&=\int_{L^n}g(v_1,\dots,v_n)\,{\rm d}\overline{\nu}_{\mathbf L}(v_1)\dots{\rm d}\overline{\nu}_{\mathbf L}(v_n)\\
&=\int_{L^n}g(v_1,\dots,v_n)\,{\rm d}\nu_{\mathbf L}(v_1)\dots{\rm d}\nu_{\mathbf L}(v_n).
\end{align*}
It follows that for every finite graph $F$ with vertex set $\{1,\dots,n\}$, denoting by $\phi_F$ the formula 
$\bigwedge_{ij\in E(F)}(x_i\sim x_j)$, it holds that
\begin{align*}
t(F,W)&=\int_{[0,1]^n}\prod_{ij\in E(F)}W(x_i,x_j)\,{\rm d}x_1\dots{\rm d}x_n\\
&=\int_{L^n}\prod_{ij\in E(F)}\mathbf 1_{\Omega_{(x_1\sim x_2)}(\mathbf L)}(v_i,v_j)\,{\rm d}\nu_{\mathbf L}(v_1)\dots{\rm d}\nu_{\mathbf L}(v_n)\\
&=\int_{L^n}\mathbf 1_{\Omega_{\phi_F}(\mathbf L)}(v_1,\dots,v_n)\,{\rm d}\nu_{\mathbf L}(v_1)\dots{\rm d}\nu_{\mathbf L}(v_n)\\
&=\langle\phi_F,\mathbf L\rangle\\
&=\lim_{n\rightarrow\infty}\langle\phi_F,G_n\rangle\\
&=\lim_{n\rightarrow\infty}t(F,G_n).
\end{align*}
Hence $W$ is a graphon L-limit of $(G_n)_{n\in\bbbn}$. As $W$ is $\{0,1\}$-valued, it is (by definition) random-free.
\end{proof}

We deduce the following limitation of modelings as limit objects.
\begin{corollary}
\label{cor:rfl}
Let $X$ be a fragment that includes all quantifier free formulas.

Assume $(G_n)_{n\in\bbbn}$ is an $X$-convergent sequence of graphs with unbounded order, which is {\rm L}-convergent to some non random-free graphon $W$. 
Then $(G_n)_{n\in\bbbn}$ has no modeling $X$-limit.
\end{corollary}

Let us now give some example stressing that the nullsets of the mod $0$ isomorphism $f$ can be quite large, making $\mathbf L$ and
$W$ look quite different. We now give an example in the more general setting of directed graphs and non-symmetric graphons.
\begin{example}
Let $\vec{T}_n$ be the transitive tournament of order $n$, that is the directed graph on $\{1,\dots,n\}$ defined from the
natural linear order $<_n$ on $\{1,\dots,n\}$ by $i\rightarrow j$ if $i<j$. This sequence is obviously ${\rm FO}$-convergent.

It is not difficult to construct a modeling ${\rm FO}$-limit of $(\vec T_n)_{n\in\bbbn}$:
Let 
$$L=(\{0\}\times\bbbz^+)\cup(\,]0,1[\times\bbbz)\cup (\{1\}\times\bbbz^-),$$
 with the Borel $\sigma$-algebra $\Sigma$ generated by 
the product topology of $\bbbz$ (with discrete topology) and $\bbbr$ (with usual topology). On $L$ we define a linear order $<_L$ by
$(\alpha,i)<_L(\beta,j)$ if $\alpha<\beta$ or $(\alpha=\beta)$ and $(i<j)$.
That $(L,\Sigma)$ is a {\rss} follows from the o-minimality of $([0,1],<)$. 
The measure $\nu_{\mathbf L}$ can be defined as the product of Lebesgue measure on $[0,1]$ by any probability measure on $\bbbz$.
For instance, for every $B\in\Sigma$ we let
$\nu_{\mathbf L}(B)=\lambda(B\cap ([0,1]\times\{0\}))$, where $\lambda$ is Lebesgue measure.
It is not difficult to check that $\mathbf L$ is indeed a modeling ${\rm FO}$-limit of $(\{1,\dots,n\},<_n)\simeq \vec T_n$.

In this case, a mod $0$ isomorphism $f$ defined on $[0,1]$ with values in $L\setminus N_2$ (where $N_2$ is a nullset) can be defined by $f(x)=(x,0)$.
The null set $N_2$, although very large, is clearly a $\nu_{\mathbf L}$-nullset, and the obtained
(non symmetric) random-free graphon $W:[0,1]\times[0,1]\rightarrow [0,1]$ is simply defined by
$$
W(x,y)=\begin{cases}
1&\text{if }x<y\\
0&\text{otherwise.}
\end{cases}
$$
Note that $W$ corresponds to $[0,1]$ with its natural order $<$. This order is clearly an L-limit of $<_n$
(but not an elementary limit, as it is dense although no finite order is).
\end{example}

In the spirit of Lemma~\ref{lem:rss}, we propose the following problems:
\begin{problem}
\label{pb:root1}
Let $\mathbf L$ be a modeling ${\rm FO}$-limit of a sequence $(\mathbf A_n)_{n\in\bbbn}$ of $\sig$-structures, and let $v\in L$.
Does there exist a sequence $(v_n)_{n\in\bbbn}$ such that $v_n\in A_n$ and such that 
the rooted modeling $(\mathbf L,v)$ is a modeling ${\rm FO}$-limit of the rooted structures $(\mathbf A_n,v_n)$? 
\end{problem}

\begin{problem}
\label{pb:rootk}
Let $\mathbf L$ be a modeling ${\rm FO}$-limit of a sequence $(\mathbf A_n)_{n\in\bbbn}$ of $\sig$-structures. Does there exist
$f:L\rightarrow\prod_{i\in\bbbn}A_n$ such that for every $v_1,\dots,v_k\in L$,
the $k$-rooted modeling $(\mathbf L,v_1,\dots,v_k)$ is a 
modeling ${\rm FO}$-limit of the $k$-rooted structures $(\mathbf A_n,f(v_1)_n,\dots,f(v_k)_n)$? 
\end{problem}
\subsection{Interpretation Schemes  applied to Modelings}

Basic interpretation schemes will be an efficient tool to handle {\rps}s.
Let $\mathsf I$ be an interpretation scheme ${\mathsf I}$ of $\sig$-structures
in $\kappa$-structures. We have seen that $\mathsf I$ can be extended to a mapping from 
$\kappa$-{\rss} to $\sig$-\rss. In the case where $\mathsf I$ is a basic interpretation scheme, we further extend
$\mathsf I$ to a mapping  from $\kappa$-{\rps} to $\sig$-\rps:
For a $\kappa$-{\rps} $\mathbf A$, the $\sig$-{\rps} $\mathbf B=\mathsf I(\mathbf A)$ is the {\rps} on the image {\rss}
of $\mathbf A$ with the probability measure $\nu_{\mathbf B}=\nu_{\mathbf A}$. This is formalized as follows:

\begin{lemma}
\label{lem:measint}
Let $\mathsf I$ be a basic interpretation scheme ${\mathsf I}$ of $\sig$-structures
in $\kappa$-structures with exponent $k$. Extend the definition of\/ $\mathsf I$ to a mapping 
of $\kappa$-{\rps} to $\sig$-{\rps} by setting $\nu_{\mathsf I(\mathbf A)}=\nu_{\mathbf A}^{k}$.
Then
for every $\kappa$-{\rps} $\mathbf A$ and every $\phi\in{\rm FO}(\sig)$ the following equation holds:
$$
\langle\phi,\mathsf I(\mathbf A)\rangle=\langle\tilde{\mathsf I}(\phi),\mathbf
A\rangle.$$ 
\end{lemma}
\begin{proof}
Let $\mathbf A$ be a $\kappa$-{\rps}.  
For every $\phi\in{\rm FO}_p(\sig)$
 the following equation holds:
 $$
\Omega_{\phi}(\mathsf I(\mathbf A))=\Omega_{\tilde{\mathsf I}(\phi)}(\mathbf A)$$
thus $\langle\phi,\mathsf I(\mathbf A)\rangle
=\nu_{\mathsf I(\mathbf A)}^{p}(\Omega_{\phi}(\mathsf I(\mathbf A)))
=\nu_{\mathbf A}^{kp}(\Omega_{\tilde{\mathsf I}(\phi)}(\mathbf A))
=\langle \tilde{\mathsf I}(\phi),\mathbf A\rangle$.
\myqed
\end{proof}
\begin{remark}
If the basic interpretation scheme $\mathsf I$ is defined by quantifier free formulas only, then it is possible
to define $\widetilde{\mathsf I}$ in such a way that for every $\phi\in{\rm FO}(\sig)$ it holds that 
${\rm qrank}(\widetilde{\mathsf I}(\phi))\leq{\rm qrank}(\phi)$.
\end{remark}

The following strengthening of Proposition~\ref{prop:intlim}
in the case where we consider a basic interpretation scheme
 is a clear consequence of Lemma~\ref{lem:measint}.

\begin{mytinted}
\begin{proposition}
\label{prop:bintlim}
Let $\mathsf I$ be a basic interpretation scheme  of $\sig$-structures
in $\kappa$-structures.

If $\mathbf L$ is a {\rps} ${\rm FO}$-limit
of a sequence $(\mathbf A_n)_{n\in\bbbn}$ of 
$\kappa$-{\rps}s then $\mathsf{I}(\mathbf L)$
is a {\rps} ${\rm FO}$-limit of 
the sequence 
$(\mathsf{I}(\mathbf A_n))_{n\in\bbbn}$. 
\end{proposition}
\end{mytinted}

\begin{remark}
Let us mention that interpretations can be used to generalize graphings to
bounded degree $k$-regular hypergraphs, and even to bounded degree  relational structures:
Define  the degree of an element of a $\sig$-structure as the number of relations it belongs to.
Following the lines of Proposition~\ref{prop:bintlim} and considering a natural interpretation of $\sig$-structures with maximum degree $D$ 
in colored (multi)graphs with maximum degree $\max(D,r)$ (where $r$ is the maximum arity of a relation in $\sig$) we deduce from Corollary~\ref{cor:BSE} that classes of relational structures with bounded maximum degree have modeling FO-limits.
\end{remark}

\begin{mytinted}
\begin{lemma}
\label{lem:pullmu}
 Let $p\in\bbbn$ be a positive integer,
let $\mathbf L$ be a modeling, and
 let $\Type{p}{L}:L^p\rightarrow S(\mathcal B({\rm FO}_p(\sig)))$ be the
 function mapping $(v_1,\dots,v_p)\in L^p$ to the complete theory of $(\mathbf L,v_1,\dots,v_p)$
 (that is the set of the formulas $\varphi\in{\rm FO}_p(\sig)$ such that
 $\mathbf L\models\varphi(v_1,\dots,v_p)$).

 Then $\Type{p}{L}$ is a measurable map from $(L^p,\Sigma_{\mathbf L}^p)$ 
 to $S(\mathcal B({\rm FO}_p(\sig)))$ (with its Borel
 $\sigma$-algebra).
 
\noindent Let
 $(\mathbf A_n)_{n\in\bbbn}$ be an ${\rm FO}_p(\sig)$-convergent sequence
 of finite $\sig$-structures, and  let $\mu_p$ be the  associated limit measure 
(as in Theorem~\ref{thm:genmodx}).

Then $\mathbf L$ is an ${\rm FO}_p(\sig)$-limit modeling of $(\mathbf A_n)_{n\in\bbbn}$ 
if and only if $\mu_p$ is the pushforward of the product measure
$\nu_{\mathbf L}^p$ by the measurable map $\Type{p}{L}$, that is:
$$\Type{p}{L}_*(\nu_{\mathbf L}^p)=\mu_p.$$
\end{lemma}
\end{mytinted}
\begin{proof}
Recall that the clopen sets of $S(\mathcal B({\rm FO}_p(\sig)))$
are of the form $K(\phi)$ for $\phi\in{\rm FO}_p(\sig)$ and that they generate
the topology of $S(\mathcal B({\rm FO}_p(\sig)))$ hence also its Borel $\sigma$-algebra.

That $\Type{p}{L}$ is measurable follows from the fact 
that for every $\phi\in{\rm FO}_p$ the preimage of $K(\phi)$,
 that is 
$\Type{p}{L}^{-1}(K(\phi))=\Omega_\phi(\mathbf L)$,
 is measurable.
 
Assume that  $\mathbf L$ is an ${\rm FO}_p(\sig)$-limit modeling of $(\mathbf A_n)_{n\in\bbbn}$.
In order to prove that $\Type{p}{L}_*(\nu_{\mathbf L}^p)=\mu_p$, it is sufficient to
check it on sets $K(\phi)$:
$$
\mu_p(K(\phi))=\lim_{n\rightarrow\infty}\langle\phi,\mathbf A_n\rangle
=\langle\phi,\mathbf L\rangle=\nu_{\mathbf L}^p(\Type{p}{L}^{-1}(K(\phi))).
$$ 
Conversely, if $\Type{p}{L}_*(\nu_{\mathbf L}^p)=\mu_p$ then for every $\phi\in{\rm
FO}_p(\sig)$ the following equation holds:
$$
\langle\phi,\mathbf L\rangle=\nu_{\mathbf L}^p(\Type{p}{L}^{-1}(K(\phi)))=
\mu_p(K(\phi))=\lim_{n\rightarrow\infty}\langle\phi,\mathbf A_n\rangle,
$$
hence $\mathbf L$ is an ${\rm FO}_p(\sig)$-limit modeling of $(\mathbf A_n)_{n\in\bbbn}$.
\myqed\end{proof}

If $(X,\Sigma)$ is a Borel space with a probability measure $\nu$, 
it is standard to define the product $\sigma$-algebra $\Sigma^\omega$
on the infinite product space $X^\bbbn$,
 which is generated by cylinder sets of the form 
 $$
 R=\{f\in X^\bbbn:\ f(i_1)\in A_{i_1},\dots,f(i_k)\in A_{i_k}\}
 $$
 for some $k\in\bbbn$ and $A_{i_1},\dots,A_{i_k}\in\Sigma$.
 The measure $\nu^\omega$ of the cylinder $R$ defined above
 is then
 $$
 \nu^\omega(R)=\prod_{j=1}^k \nu(A_{i_j}).
 $$
 By Kolmogorov's Extension Theorem, this extends to a unique probability
 measure on $\Sigma^\omega$ (which we still denote by $\nu^\omega$).
 We summarize this as the following (see also Fig.~\ref{fig:push}).
 
 \begin{mytinted}
 \begin{theorem}
 \label{thm:pullback}
let $\mathbf L$ be a modeling, and
 let $\Type{\omega}{L}:L^\bbbn\rightarrow S(\mathcal B({\rm FO}(\sig)))$ be
 the function which assigns to $f\in L^\bbbn$ 
 the point of $S(\mathcal B({\rm FO}(\sig)))$ corresponding to 
 the set $\{\phi: \mathbf L\models\phi(f(1),\dots,f(p)), \text{ where }{\rm Fv}(\phi)\subseteq\{1,\dots,p\}\}$.
 
 Then $\Type{\omega}{L}$ is a measurable map.
 
\noindent Let
 $(\mathbf A_n)_{n\in\bbbn}$ be an ${\rm FO}(\sig)$-convergent sequence
 of finite $\sig$-structures, and  let $\mu$ be the  associated limit measure 
(see Theorem~\ref{thm:genmod}).

Then $\mathbf L$ is an ${\rm FO}(\sig)$-limit modeling of $(\mathbf A_n)_{n\in\bbbn}$ 
if and only if
$$\Type{\omega}{L}_*(\nu_{\mathbf L}^\omega)=\mu.$$
\myqed \end{theorem}
\end{mytinted}

Fig.~\ref{fig:push} visualizes Lemma~\ref{lem:pullmu} and Theorem~\ref{thm:pullback}.
\begin{figure}[ht]
\label{fig:push}
\begin{center}
\tiny\xy
\xymatrix@C=3mm@R=12mm{
&\myb{L}{\nu_{\mathbf L}}
\ar@/_1ex/[d]_{\Type{1}{L}} &\myb{L^2}{\nu_{\mathbf L}^2}\ar[l]\ar@/_1ex/[d]_{\Type{2}{L}}\ar@(u,ur)^{\Sym{2}}
&\ar[l]\dots
&\myb{L^p}{\nu_{\mathbf L}^p}\ar[l]\ar@/_1ex/[d]_{\Type{p}{L}}\ar@(u,ur)^{\Sym{p}}
&\myb{L^\bbbn}{\nu_{\mathbf L}^\omega}
\ar@/_1ex/[d]_{\Type{\omega}{L}}\ar[l]\ar@(u,ur)^{\Sym{\omega}}\\
\myb{S(\mathcal B({\rm FO}_0))}{\delta_T}
&\myb{S(\mathcal B({\rm FO}_1))}{\mu_1}\ar[l]
&\myb{S(\mathcal B({\rm FO}_2))}{\mu_2}\ar[l]\ar@/_1ex/[d]\ar@(u,ur)^{\Sym{2}}
&\ar[l]\dots
&\myb{S(\mathcal B({\rm FO}_p))}{\mu_p}\ar[l]\ar@/_1ex/[d]\ar@(u,ur)^{\Sym{p}}
&\myb{S(\mathcal B({\rm FO}))}{\mu}\ar[l]\ar@/_1ex/[d]\ar@(u,ur)^{\Sym{\omega}}\\
&
&\myb{S(\mathcal B({\rm FO}_1))^2}{\mu_1^2}\ar[ul]\ar@(u,ur)^{\Sym{2}}&\ar[l]\dots
&\myb{S(\mathcal B({\rm FO}_1))^p}{\mu_1^p}\ar[l]\ar@(u,ur)^{\Sym{p}}
&\myb{S(\mathcal B({\rm FO}_1))^\bbbn}{\mu_1^\omega}\ar[l]\ar@(u,ur)^{\Sym{\omega}}\\
}
\endxy
\end{center}
\caption{Pushforward of measures}\label{fig:push}
\end{figure}
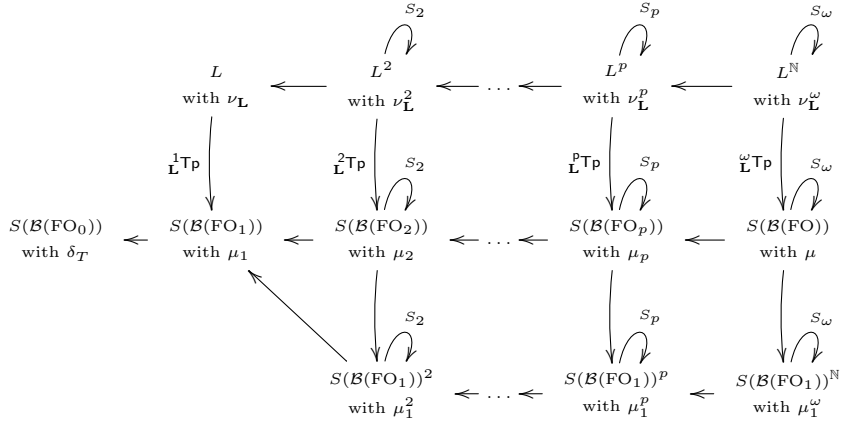
\begin{remark}
We could have considered free variables to be indexed by $\bbbz$ instead of
$\bbbn$. In such a context, natural shift operations $S$ and $T$
act respectively on the Stone space $\mathcal S$ of the Lindenbaum--Tarski
algebra of ${\rm FO}(\sig)$, and on the space $\mathbf L^\bbbz$ of the mappings from
$\bbbz$ to a $\sig$-{\rps} $\mathbf L$. If $(\mathbf A_n)_{n\in\bbbn}$ is
an ${\rm FO}$-convergent sequence with limit measure $\mu$ on $\mathcal S$, then
$(\mathcal S,\mu,S)$ is a measure-preserving dynamical system.
Also, if $\nu^\bbbz$ is the product measure on $\mathbf A$,
$(\mathbf A^\bbbz,\nu,T)$ is a Bernoulli scheme.
Then, the condition of Theorem~\ref{thm:pullback} can be restated
as follows: the modeling $\mathbf L$ is a modeling ${\rm FO}$-limit
of the sequence $(\mathbf A_n)_{n\in\bbbn}$ if and only if
 $(\mathcal S,\mu,S)$ is a factor of $(\mathbf A^\bbbz,\nu^\bbbz,T)$. 
 This setting leads to yet another interpretation of our result, which
 we hope will be treated elsewhere.
\end{remark}


\subsection{Component-Local Formulas}
The basic observation is that for $\varpi$-local formulas, we can reduce the
Stone pairing to components.

\begin{lemma}
\label{lem:stlocpair}
Let $\mathbf A$ be a $\sig$-{\rps} and
component relation $\varpi$.
Let $\psi\in{\rm FO}_p(\sig)$ be a $\varpi$-local formula of $\mathbf A$.

Assume $\mathbf A$ has countably many connected components $\{\mathbf
B_i\}_{i\in \Gamma}$.
Let $\Gamma_+$ be the
set of indexes $i$ such that $\nu_{\mathbf A}(B_i)>0$.
For $i\in \Gamma_+$ we equip $\mathbf B_i$ with the $\sigma$-algebra
$\Sigma_{\mathbf B_i}$  and the probability measure $\nu_{\mathbf B_i}$, where
$\Sigma_{\mathbf B_i}$ is restriction of $\Sigma_{\mathbf A}$ to $B_i$ and, for $X\in\Sigma_{\mathbf B_i}$, 
$\nu_{\mathbf B_i}(X)=\nu_{\mathbf A}(X)/\nu_{\mathbf A}(B_i)$. 
Then
$$
\langle\psi,\mathbf A\rangle=\sum_{i\in\Gamma}\nu_{\mathbf
A}(B_i)^p\,\langle\psi,\mathbf B_i\rangle.
$$ 
\end{lemma}
\begin{proof}
First note that each connected component of $\mathbf A$ is measurable:
let $\mathbf B_i$ be a connected component of $\mathbf A$ and let $a\in B_i$. Then
$B_i=\{x\in A: \mathbf A\models \varpi(x,a)\}$ hence $B_i$ is measurable as $\mathbf A$ is 
a \rss.
Let $Y=\{(v_1,\dots,v_p)\in A^p: \mathbf A\models\psi(v_1,\dots,v_p)\}$. 
Then $\langle\psi,\mathbf A\rangle=\nu_{\mathbf A}^{p}(Y)$.
As $\psi$ is $\varpi$-local, it also holds that 
$Y=\bigcup_{i\in \Gamma}Y_i$, where $Y_i=\{(v_1,\dots,v_p): \mathbf
B_i\models\psi(v_1,\dots,v_p)\}=Y\cap B_i^p$.
As $B_i\in\Sigma_{\mathbf A}$ and
 $Y\in\Sigma_{\mathbf A}^{p}$, it follows that $Y_i\in\Sigma_{\mathbf A}^{p}$ and (by countable additivity)
it holds that 
$$\langle\psi,\mathbf A\rangle=\nu_{\mathbf A}^{p}(Y)
=\sum_{i\in \Gamma}\nu_{\mathbf A}^{p}(Y_i)
=\sum_{i\in \Gamma_+}\nu_{\mathbf A}(\mathbf B_i)^p\nu_{\mathbf B_i}^{p}(Y_i)=
\sum_{i\in \Gamma}\nu_{\mathbf A}(\mathbf B_i)^p\,\langle\psi,\mathbf
B_i\rangle.
$$
\myqed
\end{proof}

\begin{corollary}
\label{cor:stlocpair}
Let $\mathbf A$ be a finite $\sig$-structure with
component relation $\varpi$.
Let $\psi\in{\rm FO}_p(\sig)$ be a $\varpi$-local formula of $\mathbf A$.

Let $\mathbf B_1,\dots,\mathbf B_n$ be the connected components of $\mathbf A$. 
Then
$$
\langle\psi,\mathbf A\rangle=\sum_{i=1}^n\biggl(\frac{|B_i|}{|A|}\biggr)^p\,\langle\psi, \mathbf B_i\rangle.
$$ 
\end{corollary}
We are now ready to reduce Stone pairing of local formulas
to Stone pairings with $\varpi$-local formulas on connected components.

\begin{mytinted}
\begin{theorem}
 \label{thm:sl}
Let $p\in\bbbn$ and $\phi\in{\rm FO}_p^{\rm local}(\sig)$.

Then there exist $\varpi$-local formulas $\xi_{i,j}\in{\rm FO}_{q_{i,j}}^{\rm local}$ 
($1\leq i\leq n$, $j\in I_i$) with ${\rm qrank}(\xi_{i,j})\leq {\rm qrank}(\phi)$ such that
for each $i$, $\sum_{j\in I_i} q_{i,j}=p$
and, for every {\rps} $\mathbf A$ with component relation $\varpi$
and countable set of connected components $\{\mathbf B_k\}_{k\in \Gamma}$,
 the following equation holds:
$$\langle\phi,\mathbf A\rangle=\sum_{i=1}^n \prod_{j\in I_i}\sum_{k\in\Gamma}\nu_{\mathbf A}(B_k)^{q_{i,j}}
\langle\xi_{i,j},\mathbf B_k\rangle.$$
\end{theorem}
\end{mytinted}
\begin{proof}
This is a direct consequence of Lemmas~\ref{lem:sl} and~\ref{lem:stlocpair}. 
\end{proof}

The case of sentences can be handled easily by limited counting. 
For a set $X$ and an integer $m$, define
$${\rm Big}_m(X)=\begin{cases}
1&\text{if }|X|\geq m\\
0&\text{otherwise}.
\end{cases} 
$$
\begin{mytinted}
\begin{lemma}
\label{lem:scount}
Let $\theta\in{\rm FO}_0(\sig)$. 

Then there exist formulas $\psi_1,\dots,\psi_s\in{\rm FO}_1^{\rm local}$ with quantifier rank at most
$q({\rm qrank}(\theta))$, integers $m_1,\dots,m_s\leq {\rm qrank}(\theta)$, and
a Boolean function $F$ such that for every $\sig$-structure $\mathbf A$ with component relation
$\varpi$ and connected components $\mathbf B_i$ ($i\in I$), the property $\mathbf A\models\theta$ is 
equivalent to

$$
F({\rm Big}_{m_1}(\{i,\mathbf B_i\models(\exists x)\psi_1(x)\}),\dots,
{\rm Big}_{m_s}(\{i,\mathbf B_i\models(\exists x)\psi_s(x)\}))=1.
$$
\end{lemma}
\end{mytinted}
\begin{proof}
Indeed, it follows from the Gaifman locality theorem (Theorem~\ref{thm:gaifman}) that ---
in the presence of a component relation $\varpi$ ---
every sentence $\theta$ with quantifier rank $r$ can be written as a Boolean combination of sentences $\theta_k$ 
of the form
$$\exists
y_1\dots\exists y_{m_k}\biggl(\bigwedge_{1\leq i<j\leq
m_k}\neg\varpi(y_i,y_j)\wedge\bigwedge_{1\leq i\leq m_k}\psi_k(y_i)\biggr)
$$
where $\psi_k$ is $\varpi$-local,  $m_k\leq {\rm qrank}(\theta)$, 
and ${\rm qrank}(\psi_k)\leq q({\rm qrank}(\theta))$, for some fixed function $q$.
As $\mathbf A\models\theta_k$ if and only if 
${\rm Big}_{m_k}(\{i,\mathbf B_i\models(\exists x)\psi_k(x)\})=1$, the lemma
follows.
\myqed\end{proof}

%
\subsection{Convex Combinations of Modelings}
In several contexts, it is clear when the disjoint union of converging sequences forms a converging sequence.
If two graph sequences $(G_n)_{n\in\bbbn}$ and $(H_n)_{n\in\bbbn}$ are {\rm L}-convergent or
BS-convergent, it is clear
that the sequence $(G_n\cup H_n)_{n\in\bbbn}$ is also convergent, provided that 
the limit 
$$\lim_{n\rightarrow\infty} |G_n|/(|G_n|+|H_n|)$$ exists. The same applies if we 
merge a countable set of {\rm L}-convergent (resp.\ BS-convergent) sequences
$(H_{n,i})_{n\in\bbbn}$ (where $i\in\bbbn$), with the obvious restriction  that
for each $i\in\bbbn$ all but finitely many $H_{n,i}$ are empty graphs.

We shall see that the possibility of merging  a countable set of converging
sequences to ${\rm FO}^{\rm local}$-convergence will need a further assumption, namely the following equality:
$$
\sum_{i}\lim_{n\rightarrow\infty}\frac{|G_{n,i}|}{|\bigcup_j G_{n,j}|}=1.
$$
The importance of this assumption is illustrated by the next example.
\begin{example}
Let $N_n=2^{2^n}$ (so that $N(n)$ is divisible by $2^i$ for every $1\leq i\leq 2^n$).
Consider sequences $(H_{n,i})_{n\in\bbbn}$ of edgeless black and white colored graphs
where $H_{n,i}$ is 
\begin{itemize}
  \item empty if $i>2^n$,
  \item the edgeless graph with $(2^{-i}+2^{-n})N_n$ white vertices and $2^{-i}N_n$
  black vertices if $n$ is odd,
  \item the edgeless graph with $(2^{-i}+2^{-n})N_n$ black vertices and $2^{-i}N_n$
  white vertices if $n$ is even.
\end{itemize} 
For each $i\in\bbbn$, the sequence $(H_{n,i})_{n\in\bbbn}$ is obviously {\rm L}-convergent (and even ${\rm FO}$-convergent)
as the proportion of white vertices in $H_{n,i}$ tends to $1/2$ as $n\rightarrow\infty$.
The order of $G_n=\bigcup_{i\in\bbbn}H_{n,i}$ is $3N_n$ and 
$|H_{n,i}|/|G_n|$ tends to $\frac{2}{3}\cdot 2^{-i}$ as $n$ goes to infinity. 
However the sequence $(G_n)_{n\in\bbbn}$ is not {\rm L}-convergent (hence not ${\rm FO}^{\rm local}$-convergent).
Indeed, the proportion of white vertices in $G_n$ is $2/3$ if $n$ is odd and $1/3$ is $n$ is even.
\end{example}  

Hence, we are led to the following definition.
\begin{definition}[Convex combination of {\Rps}s]
Let $\mathbf H_i$ be $\sig$-{\rps}s for $i\in I\subseteq\bbbn$ and let $(\alpha_i)_{i\in I}$ be positive real numbers
such that $\sum_{i\in I}\alpha_i=1$.

Let $\mathbf H=\coprod_{i\in I}\mathbf H_i$ be the {\rss} obtained as the disjoint union of the $\mathbf H_i$.
We endow $\mathbf H$ with the probability measure $\nu_{\mathbf H}(X)=\sum_i\alpha_i\nu_{\mathbf H_i}(X\cap H_i)$.

Then $\mathbf H$ is the {\em convex combination} of {\rps}s
$\mathbf H_i$ with {\em weights} $\alpha_i$ and we denote it by 
$\coprod_{i\in I}(\mathbf H_i,\alpha_i)$.
\end{definition}
\begin{lemma}
\label{lem:unionfom}
Let $\mathbf H_i$ be $\sig$-{\rps}s for $i\in I\subseteq\bbbn$
and let $(\alpha_i)_{i\in I}$ be positive real numbers such that $\sum_{i\in I}\alpha_i=1$.
Let $\mathbf H=\coprod_{i\in I}(\mathbf H_i,\alpha_i)$
Then
\begin{enumerate}
  \item $\mathbf H$ is a {\rps}, each $H_i$ is measurable and 
    $\nu_{\mathbf H}(H_i)=\alpha_i$ holds for every $i\in I$;
  \item if all the $\mathbf H_i$ are weakly uniform and either all the $H_i$ are infinite or all the $H_i$ are finite, $I$ is finite,
 and $\alpha_i=|H_i|/\sum_{i\in I}|H_i|$, then $\mathbf H$ is weakly uniform. 
\end{enumerate}
\end{lemma}
\begin{proof}
According to Lemma~\ref{lem:rss}, $\mathbf H$ is a {\rss}, in which each $H_i$ is measurable.
That $\nu_{\mathbf H}(H_i)=\alpha_i$ immediately follows from the definition of $\nu_{\mathbf H}$.

Assume that all the $\mathbf H_i$ are weakly uniform.
If all the $\mathbf H_i$ are finite, $I$ is finite, and $\alpha_i=|H_i|/\sum_{i\in I}|H_i|$, then
$\mathbf H$ is the {\rps} associated to the union of the $\mathbf H_i$ hence
it is weakly uniform. Otherwise all the $H_i$ are infinite, hence all the $\nu_{\mathbf H_i}$ are atomless,
$\nu_{\mathbf H}$ is atomless, and $\mathbf H$ is weakly uniform. 
\myqed\end{proof}

\begin{lemma}
\label{lem:slcomb}
Let $p\in\bbbn$ and $\phi\in{\rm FO}_p^{\rm local}(\sig)$.

Then there exist local formulas $\xi_{i,j}\in{\rm FO}_{q_{i,j}}^{\rm local}$ 
($1\leq i\leq n$, $j\in I_i$) with ${\rm qrank}(\xi_{i,j})\leq {\rm qrank}(\phi)$ such that
for each $i$, $\sum_{j\in I_i} q_{i,j}=p$
and, for every countable set of {\rps}s $\mathbf A_k$ and weights $\alpha_k$ ($k\in
\Gamma\subseteq\bbbn$ and $\sum_k\alpha_k=1$) the following equation holds, denoting $\mathbf A=\coprod_{i\in\Gamma}(\mathbf A_i,\alpha_i)$:

$$\langle\phi,\mathbf A\rangle=\sum_{i=1}^n \prod_{j\in I_i}\sum_{k\in\Gamma}\alpha_k^{q_{i,j}}
\langle\xi_{i,j},\mathbf A_k\rangle.$$
\end{lemma}
\begin{proof}
Considering, as above, the combination $\mathbf A^+=\coprod_{i\in \Gamma}(\mathbf A_i^+,\alpha_i)$, where $\mathbf A_i^+$ is
obtained by the basic interpretation scheme adding a full binary relation
$\varpi$, the result is an immediate consequence of Theorem~\ref{thm:sl}.
\myqed
\end{proof}

For $\sig$-{\rps}s $\mathbf A$ and $\mathbf B$, and $p,r\in\bbbn$ define
$$
\|\mathbf A-\mathbf B\|_{p,r}^{\rm local}=\sup\{|\langle\phi,\mathbf A\rangle-\langle\phi,\mathbf B\rangle|:\ \phi\in{\rm FO}_p^{\rm local}(\sig), {\rm qrank}(\phi)\leq r\}.
$$
The following lemma relates precisely how close  Stone pairings on two combinations of modelings are, when the modelings and weights involved in the combinations
define close Stone pairings.
\begin{lemma}
\label{lem:ploc}
Let $p,r\in\bbbn$, and let $\Gamma\subseteq\bbbn$.
For $k\in\Gamma$, let $\mathbf A_k,\mathbf B_k$ be $\sig$-{\rps}s, and let $\alpha_k,\beta_k$ be non-negative weights with
 $\sum_k\alpha_k=\sum_k\beta_k=1$.

Let $\mathbf A=\coprod_{i\in\Gamma}(\mathbf A_i,\alpha_i)$ and $\mathbf B=\coprod_{i\in\Gamma}(\mathbf B_i,\beta_i)$.
Then there exists a constant $c_{r,p}$ (which depends only on $\sig,r$, and $p$) such that it holds that 
\begin{align*}
\|\mathbf A-\mathbf B\|_{p,r}^{\rm local}
&\leq c_{r,p}\bigl(\|\alpha-\beta\|_1+\sum_{k\in\Gamma}\alpha_k\|\mathbf A_k-\mathbf B_k\|_{p,r}^{\rm local}\bigr)\\
&\leq c_{r,p}\bigl(\|\alpha-\beta\|_1+\sup_{i\in\Gamma}\|\mathbf A_i-\mathbf B_i\|_{p,r}^{\rm local}\bigr).
\end{align*}
\end{lemma}
\begin{proof}
Let $\phi\in{\rm FO}_p^{\rm local}(\sig)$ with ${\rm qrank}(\phi)\leq r$. According to Lemma~\ref{lem:slcomb}
there exist local formulas $\xi_{i,j}\in{\rm FO}_{q_{i,j}}^{\rm local}(\sig)$ 
($1\leq i\leq n$, $j\in I_i$) with ${\rm qrank}(xi_{i,j})\leq r$ such that
for each $i$, $\sum_{j\in I_i} q_{i,j}=p$
and, for every countable set of {\rps}s $\mathbf C_k$ and weights $\gamma_k$ ($k\in
\Gamma$ and $\sum_k\gamma_k=1$) the following equation holds, denoting $\mathbf C=\coprod_{i\in \Gamma}(\mathbf C_i,\gamma_i)$:
$$\langle\phi,\mathbf C\rangle=\sum_{i=1}^n \prod_{j\in I_i}\langle\xi_{i,j},\mathbf C\rangle,\quad\text{with }
\langle\xi_{i,j},\mathbf C\rangle=\sum_{k\in\Gamma}\gamma_k^{q_{i,j}}\langle\xi_{i,j},\mathbf C_k\rangle.$$
As there are only finitely many non-equivalent formulas in ${\rm FO}_p^{\rm local}(\sig)$ with quantifier rank at most $r$, there is a constant
$N_{r,p}$ such that $n\leq N_{r,p}$.

We have
\begin{align*}
|\langle\phi,\mathbf A\rangle-\langle\phi,\mathbf B\rangle|
&\leq 
\sum_{i=1}^n\biggl|\prod_{j\in I_i}\langle\xi_{i,j},\mathbf A\rangle-\prod_{j\in I_i}\langle\xi_{i,j},\mathbf B\rangle\biggr|.
\intertext{Note that if $a_i,b_i\in[0,1]$ then we get easily} 
\biggl|\prod_{i=1}^k a_i-\prod_{i=1}^k b_i\biggr|&=\biggl|(a_1-b_1)\prod_{i=2}^k a_i+b_1\bigl(\prod_{i=2}^k a_i-\prod_{i=2}^k b_i\bigr)\biggr|\\
&\leq |a_1-b_1|+ \biggl|\prod_{i=2}^k a_i-\prod_{i=2}^k b_i\biggr|\\
&\leq \sum_{i=1}^k |a_i-b_i|.
\intertext{Hence, as for every $1\leq i\leq n$ and every $j\in I$ it holds that  
$0\leq \langle\xi_{i,j},\mathbf A\rangle\leq 1$ and $0\leq\langle\xi_{i,j},\mathbf B\rangle\leq 1$, we have
}
|\langle\phi,\mathbf A\rangle-\langle\phi,\mathbf B\rangle|&\leq \sum_{i=1}^n\sum_{j\in I_i}\biggl|\langle\xi_{i,j},\mathbf A\rangle-\langle\xi_{i,j},\mathbf B\rangle\biggr|\\
&\leq \sum_{i=1}^n\sum_{j\in I_i}\biggl|\sum_{k\in\Gamma}\alpha_k^{q_{i,j}}\langle\xi_{i,j},\mathbf A_k\rangle-\sum_{k\in\Gamma}\beta_k^{q_{i,j}}
\langle\xi_{i,j},\mathbf B_k\rangle\biggr|\\
&\leq \sum_{i=1}^n\sum_{j\in I_i}\sum_{k\in\Gamma}\bigl|\alpha_k^{q_{i,j}}\langle\xi_{i,j},\mathbf A_k\rangle-\beta_k^{q_{i,j}}
\langle\xi_{i,j},\mathbf B_k\rangle\bigr|.
\intertext{Thus, as $q_{i,j}\geq 1$ and as Stone pairings $\langle\,\cdot\,,\,\cdot\,\rangle$ have value in $[0,1]$, the following inequality holds
(denoting $c_{r,p}=pN_{r,p}$):}
\|\mathbf A-\mathbf B\|_{p,r}^{\rm local}&\leq c_{r,p}\biggl(\sum_{k\in\Gamma}|\alpha_k-\beta_k|+\sum_{k\in\Gamma}\alpha_k\|\mathbf A_k-\mathbf B_k\|_{p,r}^{\rm local}\biggr).
\end{align*}
\end{proof}

\begin{lemma}
\label{lem:resloc}
Let $p,r\in\bbbn$, let $\mathbf A,\mathbf B$ be $\sig$-{\rps}, with connected components $\mathbf A_k, k\in\Gamma_{\mathbf A}$ 
and $\mathbf B_k, k\in\Gamma_{\mathbf B}$ (where $\Gamma_{\mathbf A}$ and $\Gamma_{\mathbf B}$ can be infinite non-countable).

Then the following inequality holds
$$\|\mathbf A-\mathbf B\|_{p,r}^{\rm local}<c_{r,p}\bigl(\sup_{k\in\Gamma_{\mathbf A}}\nu_{\mathbf A}(\mathbf A_k)+\sup_{k\in\Gamma_{\mathbf B}}\nu_{\mathbf B}(\mathbf B_k)+\|\mathbf A-\mathbf B\|_{1,r}^{\rm local}\bigr).$$
\end{lemma}
\begin{proof}
Let $\phi\in{\rm FO}_p^{\rm local}(\sig)$ with ${\rm qrank}(\phi)\leq r$.
The following equation holds
$$\langle\phi,\mathbf A\rangle=\sum_{i=1}^n \prod_{j\in I_i}\langle\xi_{i,j},\mathbf A\rangle.$$
It is clear that if $\zeta_{i,j}$ is component-local and $q_{i,j}>1$ then
$$\langle\xi_{i,j},\mathbf A\rangle<\sup_{k\in\Gamma_{\mathbf A}}\nu_{\mathbf A}(\mathbf A_k).$$

Let $X$ be the set of the integers $1\leq i\leq n$ such that there is $j\in I_i$ such that $q_{i,j}>1$, and let $Y$
be the complement of $X$ in $\{1,\dots,n\}$.
Then 
\begin{align*}
\bigl|\langle\phi,\mathbf A\rangle-\sum_{i\in Y}\prod_{j\in I_i}\langle\xi_{i,j},\mathbf A\rangle\bigr|&<c_{r,p} \sup_{k\in\Gamma_{\mathbf A}}\nu_{\mathbf A}(\mathbf A_k).
\intertext{Similarly, it holds that}
\bigl|\langle\phi,\mathbf B\rangle-\sum_{i\in Y}\prod_{j\in I_i}\langle\xi_{i,j},\mathbf B\rangle\bigr|&<c_{r,p} \sup_{k\in\Gamma_{\mathbf B}}\nu_{\mathbf B}(\mathbf B_k).
\intertext{Thus the statement follows from}
\biggl|\sum_{i\in Y}\prod_{j\in I_i}\langle\xi_{i,j},\mathbf A\rangle-\sum_{i\in Y}\prod_{j\in I_i}\langle\xi_{i,j},\mathbf B\rangle\biggr|
&\leq\sum_{i\in Y}\biggl|\prod_{j\in I_i}\langle\xi_{i,j},\mathbf A\rangle-\prod_{j\in I_i}\langle\xi_{i,j},\mathbf B\rangle\biggr|\\
&\leq \sum_{i\in Y}\sum_{j\in I_i}\bigl|\langle\xi_{i,j},\mathbf A\rangle-\langle\xi_{i,j},\mathbf B\rangle\bigr|\\
&\leq c_{r,p} \|\mathbf A-\mathbf B\|_{1,r}^{\rm local}.
\end{align*}
\end{proof}

\begin{mytinted}
\begin{theorem}
\label{thm:contcomb}
Let $p\in\bbbn$, let $I\subseteq \bbbn$ and, for each $i\in I$ let 
$(\mathbf A_{i,n})_{n\in\bbbn}$ be an
${\rm FO}_p^{\rm local}(\sig)$-convergent sequence of $\sig$-modelings
 and let 
$(a_{i,n})_{n\in\bbbn}$ be a convergent sequence of non-negative real numbers,
such that $\sum_{i\in I}a_{i,n}=1$ holds for every $n\in\bbbn$, and such that
$\sum_{i\in I}\lim_{n\rightarrow\infty}a_{i,n}=1$.

Then the sequence of convex combinations
$\coprod_{i\in I}(\mathbf A_{i,n},a_{i,n})$ is 
${\rm FO}_p^{\rm local}(\sig)$-convergent.
\end{theorem}
\end{mytinted}
\begin{proof}
If $I$ is finite, then the result follows from Lemma~\ref{lem:slcomb}.
Hence we can assume $I=\bbbn$.
 
Let $\phi\in{\rm FO}_p^{\rm local}$, let $q\in\bbbn$, and let $\epsilon>0$ be a
positive real.
Assume that for
each $i\in\bbbn$ the sequence $(\mathbf A_{i,n})_{n\in\bbbn}$ 
is ${\rm FO}_p^{\rm local}$-convergent
 and that $(a_{i,n})_{n\in\bbbn}$ is a convergent sequence
of non-negative real numbers, such that $\sum_{i}a_{i,n}=1$ holds for every $n\in\bbbn$.
Let $\alpha_i=\lim_{n\rightarrow\infty}a_{i,n}$, 
 let $d_i=\lim_{n\rightarrow\infty}\langle\phi,\mathbf A_{i,n}\rangle$,
and let
$C$ be such that $\sum_{i=1}^C\alpha_i>1-\epsilon/4$.
There exists $N$ such that for every $n\geq N$ and every $i\leq C$ it
holds that $|a_{n,i}-\alpha_i|<\epsilon/4C$ 
and
$|a_{i,n}^q\langle\phi,\mathbf A_{i,n}\rangle-
\alpha_{i}^q d_i|<\epsilon/2C$.
Thus
$\left|\sum_{i=1}^Ca_{i,n}^q\langle\phi,\mathbf A_{i,n}\rangle-
\sum_{i=1}^C\alpha_{i}^q d_i\right|<\epsilon/2$ and
 $\sum_{i>C+1}a_{i,n}<\epsilon/2$.
It follows that for any $n\geq N$ the following inequality holds
$$\left|\sum_{i>C+1}a_{i,n}^q\langle\phi,\mathbf A_{i,n}\rangle-
\sum_{i>C+1}\alpha_{i}^q d_i\right|
\leq \max\left(\sum_{i>C+1}a_{i,n}^q,\sum_{i>C+1}\alpha_{i}^q
d_i\right)<\epsilon/2$$
hence $|\sum_{i}a_{i,n}^q\langle\phi,\mathbf A_{i,n}\rangle-
\sum_{i}\alpha_{i}^q d_i|<\epsilon$. 

For every $\psi\in{\rm FO}_p^{\rm local}$, the expression
appearing in Lemma~\ref{lem:slcomb} for the expansion
of $\langle\phi,\coprod_{i}(\mathbf A_{i,n},a_{i,n})\rangle$
is  a finite combination of terms
of the form \linebreak $\sum_{i}a_{i,n}^{q_{i,n}}\langle\phi,\mathbf A_{i,n}\rangle$, 
where $q_{i,n}\in\bbbn$ and $\phi\in{\rm FO}_p^{\rm local}$.
It follows that  the 
value $\langle\phi,\coprod_{i}(\mathbf A_{i,n},a_{i,n})\rangle$
converges as $n$ grows to infinity. Hence 
$(\coprod_{i}(\mathbf A_{i,n},a_{i,n}))_{n\in\bbbn}$ 
is ${\rm FO}_p^{\rm local}$-convergent. 
\myqed
\end{proof}
\begin{mytinted}
\begin{corollary}
\label{cor:merge}
Let $p\geq 1$ and let $(\mathbf A_n)_{n\in\bbbn}$ be a sequence
of finite $\sig$-structures.
 
Assume $\mathbf A_n$ be the disjoint union of 
$\mathbf B_{n,i}$ ($i\in\bbbn$) where all but a finite number
of $\mathbf B_{n,i}$ are empty.
Let $a_{n,i}=|B_{n,i}|/|A_n|$.
Assume further that:
\begin{itemize}
  \item for each $i\in\bbbn$, the limit $\alpha_i=\lim_{n\rightarrow\infty}a_{n,i}$ exists,
  \item for each $i\in\bbbn$ such that $\alpha_i\neq 0$, 
  the sequence $(\mathbf B_{n,i})_{n\in\bbbn}$ is ${\rm FO}_p^{\rm local}$-convergent,
  \item the following equation holds:
  $$\sum_{i\geq 1}\alpha_{i}=1.
  $$
\end{itemize}

Then, the sequence $(\mathbf A_n)_{n\in\bbbn}$ is ${\rm FO}_p^{\rm local}$-convergent.

Moreover, if $\mathbf L_i$ is a modeling ${\rm FO}_p^{\rm local}$-limit
of $(\mathbf B_{n,i})_{n\in\bbbn}$ when $\alpha_i\neq 0$ then 
$\coprod_i (\mathbf L_i,\alpha_i)$ is a modeling ${\rm FO}_p^{\rm local}$-limit
of $(\mathbf A_{n})_{n\in\bbbn}$. 
\end{corollary}
\end{mytinted}
\begin{proof}
This follows from Theorem~\ref{thm:contcomb}, as
 $\mathbf A_n=\coprod_i (\mathbf B_{n,i},a_{n,i})$.
\myqed\end{proof}
\begin{definition}
A family of sequence $(\mathbf A_{i,n})_{n\in\bbbn}\ (i\in I)$ 
of $\sig$-structures is
{\em uniformly elementarily convergent} if, for every formula
$\phi\in{\rm FO}_1(\sig)$ there is an integer $N$ such that the following implication holds
$$
\forall i\in I,\ \forall n'\geq n\geq N,\quad 
(\mathbf A_{i,n}\models (\exists x)\phi(x))\Longrightarrow
(\mathbf A_{i,n'}\models (\exists x)\phi(x)).
$$
\end{definition}
First notice that if a family $(\mathbf A_{i,n})_{n\in\bbbn}\ (i\in I)$
of sequences is uniformly elementarily convergent, then
each sequence $(\mathbf A_{i,n})_{n\in\bbbn}$ is elementarily convergent.

\begin{lemma}
\label{lem:elunion}
Let $I\subseteq N$, and let $(\mathbf A_{i,n})_{n\in\bbbn}\ (i\in I)$ be sequences forming
a uniformly elementarily convergent family.

Then $(\bigcup_{i\in I}\mathbf A_{i,n})_{n\in\bbbn}$ is elementarily convergent.

Moreover, if $(\mathbf A_{i,n})_{n\in\bbbn}$ is elementarily convergent to $\widehat{\mathbf A}_i$
then $(\bigcup_{i\in I}\mathbf A_{i,n})_{n\in\bbbn}$ is elementarily convergent to $\bigcup_{i\in I}\widehat{\mathbf A}_i$.
\end{lemma}
\begin{proof}
Let $\sig^+$ be the signature $\sig$ augmented by a binary relational symbol
$\varpi$. Let $I_1$ be the basic interpretation scheme of $\sig^+$-structures in
$\sig$-structures defining $\varpi(x,y)$ for every $x,y$. Let $\mathbf A_{i,n}^+=I_1(\mathbf A_{i,n})$.
According to Lemma~\ref{lem:scount}, for every sentence $\theta\in{\rm FO}_0(\sig)$ there
exist formulas $\psi_1,\dots,\psi_s\in{\rm FO}_1^{\rm local}$, an integer $m$, and a Boolean
function $F$ such that the property 
$\bigcup_{i\in I}\mathbf A_{i,n}^+\models\theta$ is equivalent to
$$F({\rm Big}_{m_1}(\{i, {\mathbf A}_{i,n}\models (\exists x)\psi_1(x)\}),\dots,
{\rm Big}_{m_s}(\{i, {\mathbf A}_{i,n}\models (\exists x)\psi_s(x)\}))=1.$$
According to the definition of a uniformly elementarily convergent family
there is an integer $N$ such that, for every $1\leq j\leq s$, the value
 ${\rm Big}_{m_j}(\{i, {\mathbf A}_{i,n}\models (\exists x)\psi_j(x)\})$
 is a function of $n$, which is non-decreasing for $n\geq N$. It follows
 that this function admits a limit for every $1\leq j\leq s$ hence the
 exists an integer $N'$ such that either $\bigcup_{i\in I}\mathbf A_{i,n}^+\models\theta$
 holds for every $n\geq N'$ or it holds for no $n\geq N'$. It follows
 that $(\bigcup_{i\in I}\mathbf A_{i,n}^+)_{n\in\bbbn}$ is elementarily convergent.
 Thus (by means of the basic interpretation scheme deleting $\varpi$) 
 $(\bigcup_{i\in I}\mathbf A_{i,n})_{n\in\bbbn}$ is elementarily convergent
 
If $I$ is finite, it is easily checked that if $(\mathbf A_{i,n})_{n\in\bbbn}$ is elementarily convergent to $\widehat{\mathbf A}_i$
then $(\bigcup_{i\in I}\mathbf A_{i,n})_{n\in\bbbn}$ is elementarily convergent to $\bigcup_{i\in I}\widehat{\mathbf A}_i$.

Otherwise, we can assume $I=\bbbn$.
 Following the same lines, it is easily checked that 
 $(\bigcup_{i=1}^n\widetilde{\mathbf A}_i)_{n\in\bbbn}$ converges elementarily
 to $(\bigcup_{i\in\bbbn}\widetilde{\mathbf A}_i)_{n\in\bbbn}$.
 For $i,n\in\bbbn$, let $\mathbf B_{i,2n}=\mathbf A_{i,n}$ and
 $\mathbf B_{i,2n+1}=\widetilde{\mathbf A_i}$. As, for each $i\in\bbbn$, $\widetilde{\mathbf A}_i$
 is an elementary limit of $(\mathbf A_{i,n})_{n\in\bbbn}$ it is easily checked that
 the family of the sequences $(\mathbf B_{i,n})_{n\in\bbbn}$ is uniformly elementarily convergent.
It follows that $(\bigcup_{i\in\bbbn}\mathbf B_{i,n})_{n\in\bbbn}$ is elementarily convergent
thus the elementary limit of  $(\bigcup_{i\in I}\mathbf A_{i,n})_{n\in\bbbn}$ and
$(\bigcup_{i=1}^n\widetilde{\mathbf A}_i)_{n\in\bbbn}$ are the same, that is
$\bigcup_{i\in I}\widetilde{\mathbf A}_i$.
\myqed\end{proof}

From Corollary~\ref{cor:merge} and Lemma~\ref{lem:elunion} then follows the
next general result.
\begin{mytinted}
\begin{corollary}
\label{cor:mergefo}
Let $(\mathbf A_n)_{n\in\bbbn}$ be a sequence
of finite $\sig$-structures.
 
Assume $\mathbf A_n$ be the disjoint union of 
$\mathbf B_{n,i}$ ($i\in\bbbn$) where all but a finite number
of $\mathbf B_{n,i}$ are empty.
Let $a_{n,i}=|B_{n,i}|/|A_n|$.
Assume that:
\begin{itemize}
  \item for each $i\in\bbbn$, the limit
  $\alpha_i=\lim_{n\rightarrow\infty}a_{n,i}$ exists and
the following equation holds:
  $$\sum_{i\geq 1}\alpha_{i}=1,
  $$
    \item for each $i\in\bbbn$ such that $\alpha_i\neq 0$, 
  the sequence $(\mathbf B_{n,i})_{n\in\bbbn}$ is ${\rm FO}^{\rm
  local}$-convergent,
  \item the family $\{(\mathbf B_{n,i})_{n\in\bbbn}\ (i\in\bbbn)\}$ 
  is uniformly elementarily convergent.
\end{itemize}

Then, the sequence $(\mathbf A_n)_{n\in\bbbn}$ is ${\rm FO}$-convergent.

Moreover, if $\mathbf L_i$ is a modeling ${\rm FO}$-limit
of $(\mathbf B_{n,i})_{n\in\bbbn}$ when $\alpha_i\neq 0$ and an elementary
limit of $(\mathbf B_{n,i})_{n\in\bbbn}$ when $\alpha_i=0$ then $\coprod_i
(\mathbf L_i,\alpha_i)$ is a modeling ${\rm FO}$-limit of $(\mathbf
A_{n})_{n\in\bbbn}$.
\end{corollary}
\end{mytinted}

\subsection{Random-free graphons and Modelings}
A graphon is {\em random-free} if it is $\{0,1\}$-valued almost everywhere. Moreover, if two graphons
represent the same L-limit of finite graphs, then either they are both random-free or none of them are 
(see for instance \cite{Ja2013}). Several properties of random-free graph limits have been studied.

For example, a graph limit $\Gamma$ is random-free if and only if the random graph $G(n,\Gamma)$ of order $n$ sampled from $\Gamma$
	has entropy $o(n^2)$ \cite{Aldous1985,Ja2013} (see also \cite{Hatami2012}).
	
 A sequence of graphs $(G_n)_{n\in\bbbn}$ is {\rm L}-convergent to a random-free graphon if and only if the sequence
	$(G_n)_{n\in\bbbn}$ is convergent for the stronger metric $\delta_1$ \cite{Pikhurko2010}, where the distance 
	$\delta_1(G,H)$ of graphs $G$ and $H$ with respective vertex sets $\{x_1,\dots,x_m\}$ and $\{y_1,\dots,y_m\}$ is the minimum 
	over all non-negative $m\times n$ matrices $A=(\alpha_{i,j})$ with row sums $1/m$ and column sums $1/n$ of
	$\sum_{(i,j,g,h)\in\Delta} \alpha_{i,g}\alpha_{j,h}$, where $\Delta$ is the set of quadruples $(i,j,g,h)$ such that
	either $\{x_i,x_j\}\in E(G)$ or $\{y_g,y_h\}\in E(H)$ (but not both).

Lov\'asz and Szegedy \cite{Lovasz2010} defined a graph property (or equivalently a class of graphs) $\mathcal C$ 
to be {\em random-free} if every L-limit of graphs in $\mathcal C$ is random-free. They prove the following:
\begin{theorem}[Lov\'asz and Szegedy \cite{Lovasz2010}]
\label{thm:RFLS}
A hereditary class $\mathcal C$ is random-free if and only if there exists a bipartite graph $F$ with bipartition
$(V_1,V_2)$ such that no graph obtained from $F$ by adding edges within $V_1$ and $V_2$ is in $\mathcal C$.
\end{theorem}
From this result, one deduce for instance that the class of $m$-partite cographs is random-free (see \cite{QFTSL-arxiv} for a related study of quantifier-free limits of tree-semilattices), 
generalizing the particular cases of threshold graphs \cite{Diaconis2009} and (more general) cographs \cite{Janson2013}.

Recall that the
{\em Vapnik--Chervonenkis dimension} (or simply {\em VC-dimension}) ${\rm VC}(G)$ {\em of a graph} $G$ is the maximum
integer $k$ such that there exists in $G$ disjoint vertices $u_i$ ($1\leq i\leq k$) and $v_I$ ($\emptyset\subseteq I\subseteq \{1,\dots,k\}$) such that $u_i$ is adjacent to $v_I$ exactly if $i\in I$.
We now rephrase Lov\'asz and Szegedy Theorem~\ref{thm:RFLS} in terms of VC-dimension.
\begin{theorem}
A hereditary class $\mathcal C$ is random-free if and only if ${\rm VC}(\mathcal C)<\infty$, where
$${\rm VC}(\mathcal C)=\sup_{G\in\mathcal C} {\rm VC}(G).$$
\end{theorem}
\begin{proof}
Let $B_k$ be the bipartite graph with vertices $u_i$ ($1\leq i\leq k$) and $v_I$ ($\emptyset\subseteq I\subseteq \{1,\dots,k\}$) such that $u_i$ is adjacent to $v_I$ exactly if $i\in I$.

If ${\rm VC}(\mathcal C)<k$ then no graph obtained from $B_k$ by adding edges within the $u_i$'s and the $v_I$'s is in $\mathcal C$ hence, according to Theorem~\ref{thm:RFLS}, the class $\mathcal C$ is random-free.

Conversely, if the class $\mathcal C$ is random-free there exists, according to Theorem~\ref{thm:RFLS},
a bipartite graph $F$ with bipartition $(V_1,V_2)$ (with $|V_1|\leq |V_2|$)
such that no graph obtained from $F$ by adding edges 
within $V_1$ and $V_2$ is in $\mathcal C$. It is easily
checked that $F$ is an induced subgraph of $B_{|V_1|+\log_2 |V_2|}$ so 
${\rm VC}(\mathcal C)<\frac{|F|}{2}+\log_2|F|$.
\myqed\end{proof}

The VC-dimension of classes of graphs can also be related to the nowhere dense/somewhere dense dichotomy.
Recall that a class $\mathcal C$ is {\em somewhere dense} if there exists an integer $p$ such that
for every integer $n$ the $p$-subdivision of $K_n$ is a subgraph of a graph in $\mathcal C$, and that the class
$\mathcal C$ is {\em nowhere dense}, otherwise \cite{ND_logic,ND_characterization,Sparsity}. 
This dichotomy can also be characterized in quite a number of different ways, see \cite{Sparsity}.
Based on Laskowski \cite{Laskowski1992}, another characterization has been proved, which relates this dichotomy to VC-dimension:

\begin{theorem}[Adler and Adler \cite{Adler2013}]
\label{thm:AA}
For a monotone class of graphs $\mathcal C$, the following are equivalent:
\begin{enumerate}
	\item For every interpretation scheme ${\mathsf I}$ of graphs in graphs, the class ${\mathsf I}(\mathcal C)$ has bounded VC-dimension;
	\item For every basic interpretation scheme ${\mathsf I}$  of graphs in graphs with exponent $1$, the class ${\mathsf I}(\mathcal C)$ has bounded VC-dimension;
	\item The class $\mathcal C$ is nowhere-dense.
\end{enumerate} 
\end{theorem}

From Theorem~\ref{thm:RFLS} and~\ref{thm:AA} we deduce the following:

\begin{theorem}
\label{thm:RFND}
Let $\mathcal C$ be a monotone class of graphs. Then the following are equivalent:
\begin{enumerate}
	\item\label{enum:RF1} For every interpretation scheme ${\mathsf I}$  of graphs in graphs, the class ${\mathsf I}(\mathcal C)$ is random-free;
	\item\label{enum:RF2} For every basic interpretation scheme ${\mathsf I}$ of graphs in graphs with exponent $1$ and built using local formulas, the class ${\mathsf I}(\mathcal C)$ is random-free;
	\item\label{enum:RF3} The class $\mathcal C$ is nowhere-dense.
\end{enumerate}
\end{theorem}
\begin{proof}
Obviously, condition \eqref{enum:RF1} implies condition \eqref{enum:RF2}.
Assume that \eqref{enum:RF2} and assume for contradiction that \eqref{enum:RF3} does not hold.
Then, as $\mathcal C$ is monotone and somewhere dense, 
there is an integer $p\geq 1$ such that for every graph $n$, the $p$-subdivision ${\rm Sub}_p(K_n)$ of the complete
graph $K_n$ is in $\mathcal C$. To every finite graph $G$ we associate a graph $G'$ 
by considering an arbitrary orientation of $G$ and then building $G'$ as shown on the
figure bellow.

\begin{center}
	\includegraphics[width=\textwidth]{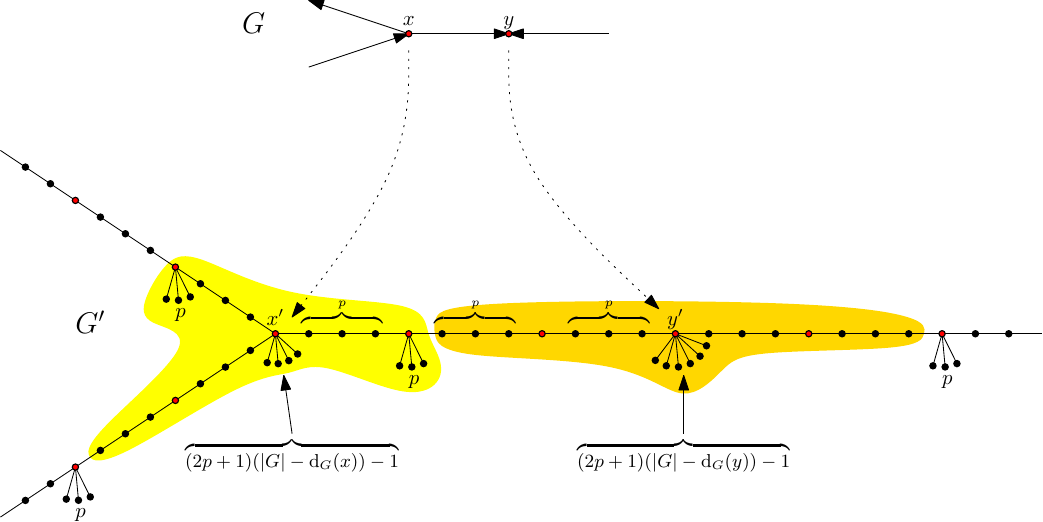}
\end{center}

Note that $G'\in\mathcal C$ as it is obviously a subgraph of the $p$-subdivision of the complete graph of order 
$(2p+1)|G|^2$. It is easily checked that there is a basic interpretation scheme ${\mathsf I}_p$  of graphs in graphs with exponent $1$ (which definitions only depends on $p$) defined using local formulas only, such
that ${\mathsf I}_p(G')=G[(2p+1)|G|]$, where $G[(2p+1)|G|]$ denotes the graph obtained from $G$ by blowing each vertex to
an independent set of size $(2p+1)|G|$.

Let $(G_i)_{i\in\bbbn}$ be a sequence of graph that is {\rm L}-convergent to a non random-free graphon $W$.
As $t(F,G)=\frac{{\rm hom}(F,G)}{|G|^{|F|}}$ is invariant by uniform blow-up of the vertices of $G$, 
for every finite graph $F$ the following equation holds:
$$t(F,{\mathsf I}_p(G_i'))=t(F,G_i[(2p+1)|G_i|])=t(F,G_i).$$
Hence $({\mathsf I}_p(G_i'))_{i\in\bbbn}$ is {\rm L}-convergent to $W$. Then the condition~\eqref{enum:RF2} contradicts the hypothesis
that $W$ is not random-free. It follows (by contradiction) that \eqref{enum:RF2} implies \eqref{enum:RF3}.

Assume condition~\eqref{enum:RF3} holds, and let $\mathsf I$ be an interpretation scheme of graphs in graphs.
 Then according to Theorem~\ref{thm:AA} the class ${\mathsf I}(\mathcal C)$ has bounded VC-dimension, hence
the hereditary closure of ${\mathsf I}(\mathcal C)$ has bounded VC-dimension thus is random-free, whence
${\mathsf I}(\mathcal C)$ is random-free.
\myqed\end{proof}

We derive the following corollary concerning existence of modeling ${\rm FO}^{\rm local}$-limits, which completes Corollary~\ref{cor:rfl} and implies Theorem~\ref{thm:modnd0}.
\begin{mytinted}
\begin{theorem}
\label{thm:modnd}
Let $\mathcal C$ be a monotone class of graphs. 

If every ${\rm FO}^{\rm local}$-convergent sequence of graphs in $\mathcal C$ has a modeling ${\rm FO}^{\rm local}$-limit then $\mathcal C$ is nowhere dense.
\end{theorem}
\end{mytinted}
\begin{proof}
Let ${\mathsf I}$ be a basic interpretation scheme of graphs in graphs built using local formulas, and let $(G_i)_{i\in\bbbn}$ be a sequence of graphs in $\mathcal C$ such that 
$|G_i|$ is unbounded, and the sequence $({\mathsf I}(G_i))_{i\in\bbbn}$ is {\rm L}-convergent. 

 By compactness, the sequence $(G_i)_{i\in\bbbn}$ has a subsequence $(G_{n_i})_{i\in\bbbn}$ that is ${\rm FO}^{\rm local}$-convergent. Hence, by hypothesis, 
$(G_{n_i})_{i\in\bbbn}$ has a modeling ${\rm FO}^{\rm local}$-limit $\mathbf L$.
According to Proposition~\ref{prop:bintlim}, 
the sequence $(\mathsf I(G_{n_i}))_{i\in\bbbn}$ has modeling ${\rm FO}^{\rm local}$-limit $\mathsf I(\mathbf L)$.
By Lemma~\ref{lem:L2W}, $\mathbf L$ defines a random-free graphon $W$ that is the L-limit of $(\mathsf I(G_{n_i}))_{i\in\bbbn}$. 
Of course, the L-limit of an {\rm L}-convergent sequence $(G_i)_{i\in\bbbn}$ with $|G_i|$ bounded is also random-free.
Hence
the class ${\mathsf I}(\mathcal C)$ is random-free. As this conclusion holds for every basic interpretation scheme 
${\mathsf I}$ built using local formulas we deduce from Theorem~\ref{thm:RFND} that $\mathcal C$ is nowhere dense.
\myqed\end{proof}


\subsection{Modelings FO-limits for Graphs of Bounded Degrees}
\label{sec:BS2}
Nice limit objects are known for sequence of bounded degree connected graphs, both for
BS-convergence (graphing) and for ${\rm FO}_0$-convergence (countable graphs).
It is natural to ask whether a nice limit object could exist for full ${\rm
FO}$-convergence. We shall now give a positive answer to this question.
First we take time to comment on the connectivity assumption.
A first impression 
is that ${\rm FO}$-convergence of disconnected graphs could be
considered component-wise.
This is far from being true in general.
The contrast between the behaviour of graphs with 
a first-order definable component relation (like graphs with bounded diameter components) and of graphs with bounded degree
is exemplified by the following example.

\begin{example}
Consider a BS-convergent sequence $(G_n)_{n\in\bbbn}$ of planar graphs with bounded degrees such that the limit distribution
has an infinite support. Note that $\lim_{n\rightarrow\infty} |G_n|=\infty$.
Then, as planar graphs with bounded degrees form a hyperfinite class of graphs there exists, for every 
graph $G_n$ and every $\epsilon>0$ a subgraph $S(G_n,\epsilon)$ of $G_n$ obtained by deleting at most $\epsilon |G_n|$ of edges, such that
the connected components of $S(G_n,\epsilon)$ have order at most $f(\epsilon)$. By considering a subsequence $G_{s(n)}$ we can assume
$\lim_{n\rightarrow\infty} |G_{s(n)}|/f(1/n)=\infty$. Then note that the sequences $(G_s(n))_{n\in\bbbn}$ and
$(S(G_{s(n)},1/n))_{n\in\bbbn}$ have the same BS-limit. By merging these sequences, we conclude that there
exists an ${\rm FO}^{\rm local}$ convergent sequence of graphs with bounded degrees $(H_n)$ such that
$H_n$ is connected if $n$ is even and such that the number of connected components of $H_n$ for $n$ odd tends
to infinity. 
\end{example}

\begin{example}
Using Fig.~\ref{fig:AB}, consider four sequences $(A_n)_{n\in\bbbn}$, $(B_n)_{n\in\bbbn}$,
$(C_n)_{n\in\bbbn}$,$(D_n)_{n\in\bbbn}$ of ${\rm FO}$-converging sequences where 
$|A_n|=|B_n|=|C_n|=|D_n|$ grows to infinity, and where these sequences have distinct limits.

Consider a sequence $(G_n)_{n\in\bbbn}$ defined as follows: for each $n$, $G_n$ has  two connected components 
 denoted by $H_{n,1}$ and $H_{n,2}$ obtained by joining $A_n, C_n$ and $B_n,D_n$ by a path of length $n$ (for $n$ odd),
and by joining $A_n, D_n$ and $B_n, C_n$ by a path of length $n$ (for $n$ even).
Then $(G_n)_{n\in\bbbn}$ is ${\rm FO}$-convergent. However, there is no
choice of a mapping $f:\bbbn\rightarrow\{1,2\}$ such that
$(H_{n,f(n)})$ is ${\rm FO}$-convergent (or even BS-convergent).
\begin{figure}[ht]
$$\begin{array}{l|r}
  \includegraphics[width=.4\textwidth]{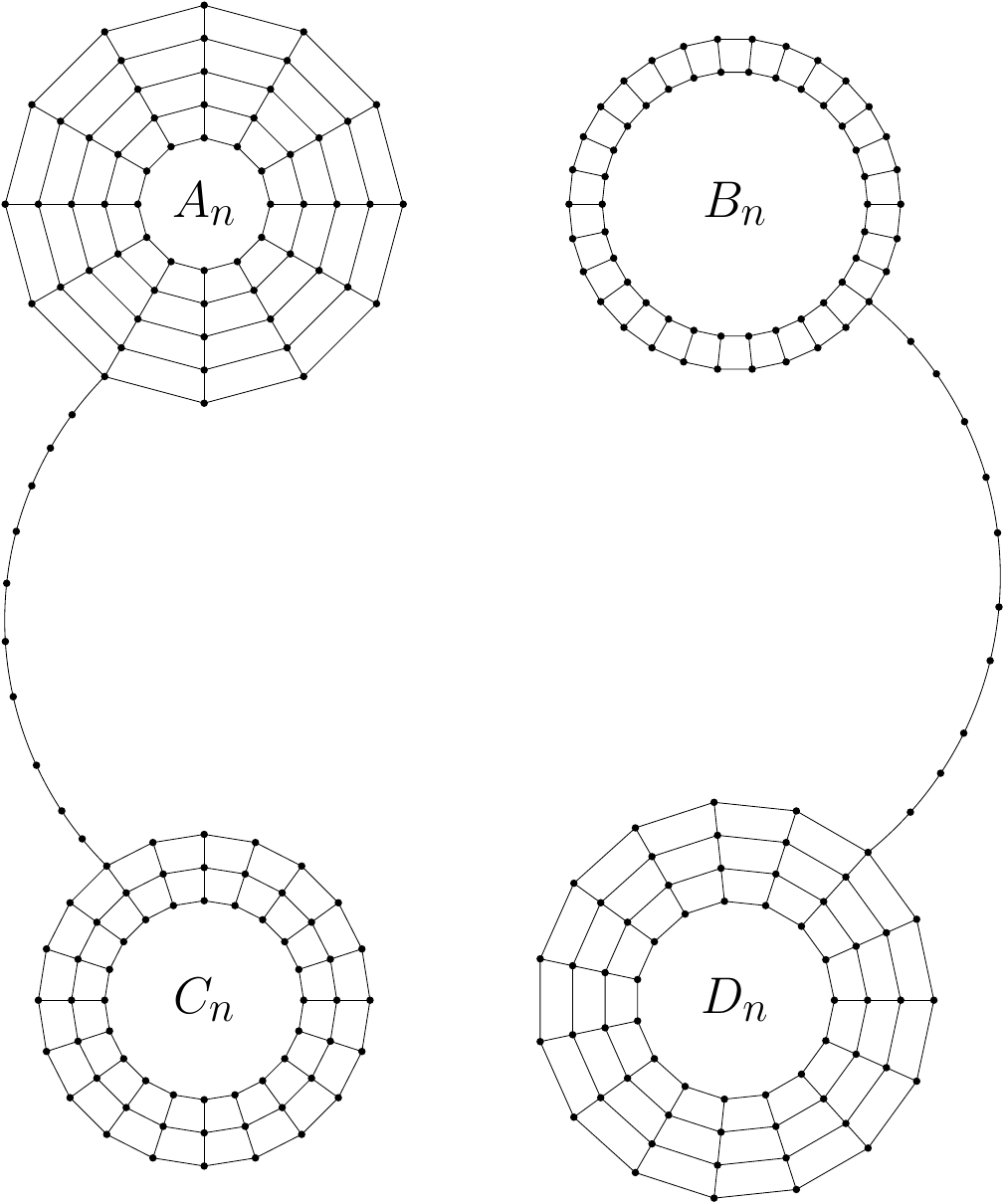}\quad&\quad
  \includegraphics[width=.4\textwidth]{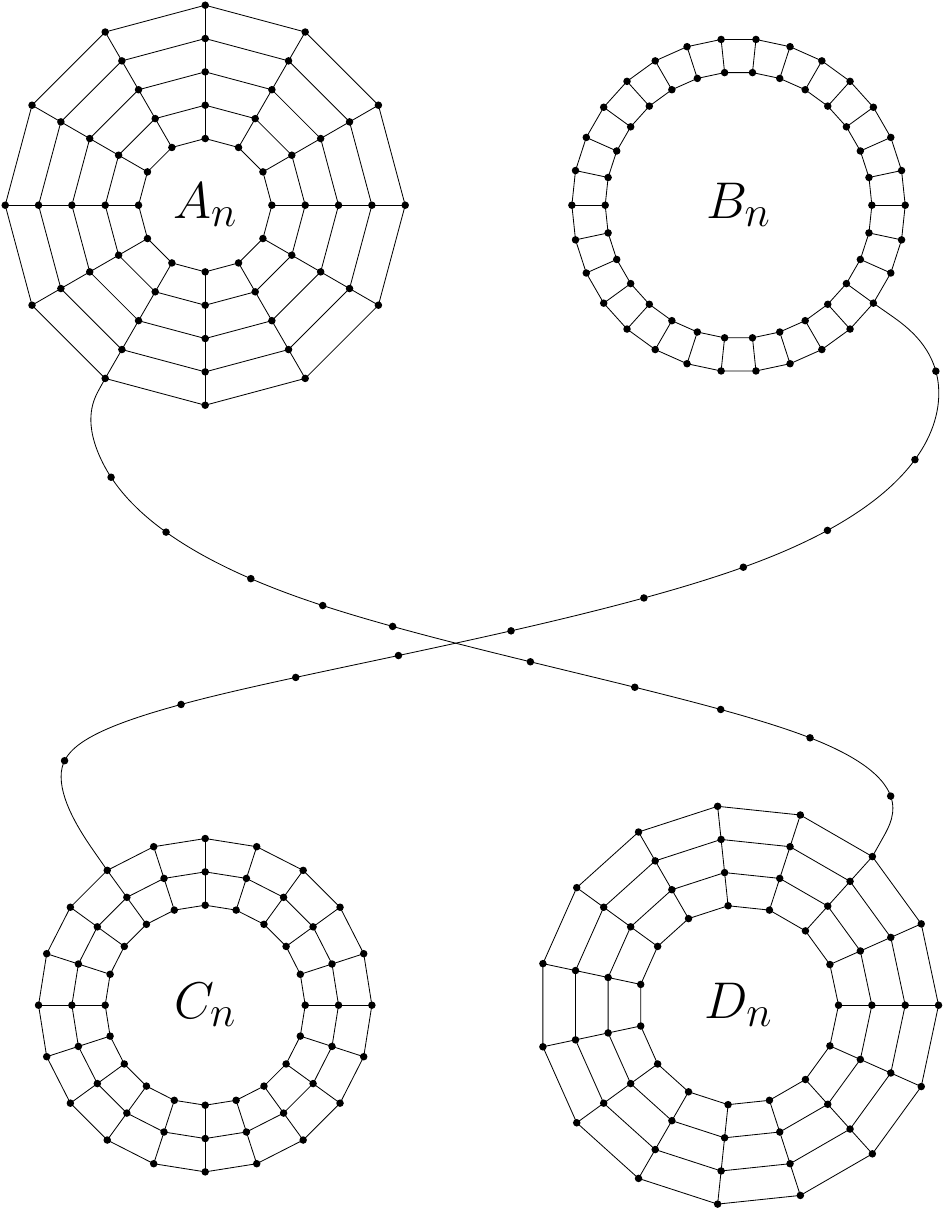}\\
  \end{array}$$
  \caption{An ${\rm FO}$-converging sequence with no component selection}
  \label{fig:AB}
\end{figure}
\end{example}

This situation is indeed related to the fact that the diameter of the graph $G_n$ in the sequence
tend to infinity as $n$ grows and that the belonging to the same connected component cannot be defined by a first-order formula.
 This situation is standard when one consider BS-limits of connected
graphs with bounded degrees: it is easily checked that, as a limit of
connected graphs, a graphing may have uncountably many connected components. 

\begin{figure}[ht]%
\includegraphics[width=.75\columnwidth]{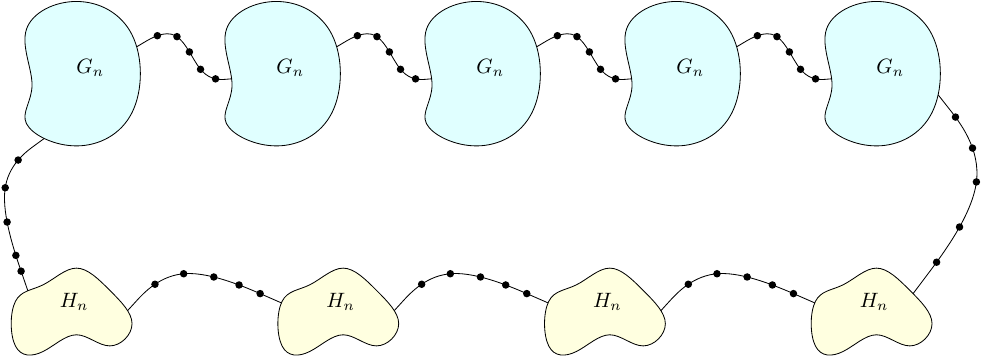}%
\caption{Construction of the graph $M_n$}%
\label{fig:convexBS}%
\end{figure}

\begin{remark}
In the spirit of the construction shown Fig.~\ref{fig:AB}, we can prove that the set of measure $\mu$ which are BS-limits of connected graphs with maximum degree $d\geq 2$ and order going to infinity is convex: Assume $(G_n)_{n\in\bbbn}$ and $(H_n)_{n\in\bbbn}$ are convergent sequences with limits $\mu_1$ and $\mu_2$, and let $0<\alpha<1$. We construct graph $M_n$ as follows: let $c_n=\min(|G_n|,|H_n|)$. We consider $\alpha |H_n|$ copies of $G_n$ and 
$(1-\alpha) |G_n|$ copies of $H_n$ linked by paths of length $\lfloor\log c_n\rfloor$ (see Fig.~\ref{fig:convexBS}). It is easily checked that the statistics of the neighborhoods of $M_n$ tend to $\alpha\mu_1+(1-\alpha)\mu_2$.
\end{remark}

Let $V$ be a standard Borel space with a measure $\mu$.
Suppose that $T_1, T_2,\dots, T_k$ are
measure preserving Borel involutions of $X$.
Then the system 
$$\mathbf G = (V, T_1, T_2, \dots, T_k, \mu)$$
 is
called a {\em measurable graphing} (or simply a {\em graphing}) \cite{Adams1990}.
A graphing $\mathbf G$ determines an equivalence relation
on the points of $V$. Simply, $x\sim_{\mathbf G} y$ if there exists
 a sequence of points $(x_1, x_2,\dots, x_m)$ of $X$
such that
\begin{itemize}
  \item $x_1 = x, x_m = y$
 \item  $x_{i+1} = T_j(x_i)$ for some $1\leq j\leq k$.
\end{itemize}
Thus there exist a natural  simple graph structure on the
equivalence classes, the {\em leafgraph}.
Here $x$ is adjacent to $y$, if $x \neq y$ and $T_j(x) = y$
 for some $1\leq j \leq k$. Now if
If $V$ is a compact metric space with a Borel measure $\mu$ and
$T_1, T_2,\dots, T_k$ are continuous
measure preserving involutions of $V$, then
$\mathbf G = (V, T_1, T_2, \dots, T_k, \mu)$ is
a {\em topological graphing}.
It is a consequence of \cite{Benjamini2001} and \cite{Gaboriau2005} that
every local weak limit of finite connected graphs with maximum
degree at most $D$ can be represented
as a measurable graphing. Elek \cite{Elek2007b} further proved the representation
can be required to be a topological graphing.

A graphing defines an edge coloration, where $\{x,y\}$ is colored
by the set of the indexes $i$ such that $y=T_i(x)$.
For an integer $r$, a graphing ${\mathbf G}=(V,T_1,\dots,T_k,\mu)$ and
 a finite rooted edge colored graph $(F,o)$ we define the set 

$$D_r(\mathbf G,(F,o))=\{x\in\mathbf G, B_r(\mathbf G,x)\simeq (F,o)\}.$$

It is easily checked that $D_r(\mathbf G,(F,o))$ is measurable.

Considering a $k$-edge colored graphing allows us to describe a vertex $x$ in a
distance-$r$ neighborhood of a given vertex $v$ by the sequence of the colors
of the edges of a path linking $v$ to $x$. Taking, among the minimal length sequences,
the one which is lexicographically minimum, it is immediate that for every 
vertex $v$ and every integer $r$ there is a injection 
$\iota_{v,r}$ from $B_r(\mathbf G,v)$
to the set of the sequences of length at most $r$ with values in $[k]$.  
Moreover, if  $B_r(\mathbf G,v)$ and $B_r(\mathbf G,v')$ are isomorphic
as edge-colored rooted graphs, then there exists a unique isomorphism
$f:B_r(\mathbf G,v)\rightarrow B_r(\mathbf G,v')$ and this isomorphism 
as the property that for every $x\in B_r(\mathbf G,v)$ it holds that 
$\iota_{v,r'}(f(x))=\iota_{v,r}(x)$.

\begin{mytinted}
\begin{lemma}
\label{lem:graphingsb}
Every graphing is a \rps.
\end{lemma}
\end{mytinted}
\begin{proof}
Let $\mathbf G=(V, T_1,\dots,T_d, \mu)$ be a graphing.
We color the edges of $G$ according to the the involutions involved.

For $r\in\bbbn$, we denote by $\mathcal F_r$ the finite set of all the 
colored rooted graphs that arise as $B_r(\mathbf G,v)$ for some
$v\in V$.
To every vertex $v\in V$ and integer $r\in\bbbn$ we associate $t_r(v)$, 
which is the isomorphism
type of the edge colored ball $B_r(\mathbf G,v)$.

According to Gaifman's locality theorem, in order to prove that 
$\mathbf G$ is a modeling, it is sufficient to prove that
for each $\phi\in {\rm FO}_{p}^{\rm local}$, the set
$$
X=\{(v_1,\dots,v_p)\in V^q:\quad\mathbf G\models\phi(v_1,\dots,v_p)\}
$$
is measurable (with respect to the product $\sigma$-algebra of $V^p$).

Let $L\in\bbbn$ be such that $\phi$ is $L$-local.
For every $\mathbf v=(v_1,\dots,v_p)\in X$ we define the  graph 
$\Gamma(\mathbf v)$ with vertex set $\{v_1,\dots,v_p\}$ such that
two vertices of $\Gamma(\mathbf v)$ are adjacent if their distance in $\mathbf G$ is 
at most $L$.  We define a partition 
$\mathcal P(\mathbf v)$ of $[p]$ as follows:
$i$ and $j$ are in a same part if $v_i$ and $v_j$ belong to a same
connected component of $\Gamma(v)$. To each part $P\in\mathcal P(\mathbf v)$,
 we associate the tuple formed by $T_P=t_{(|P|-1)L}(v_{\min P})$ and, for each
$i\in P-\{\min P\}$, a composition $F_{P,i}=T_{i_1}\circ\dots\circ T_{i_j}$ 
with $1\leq j\leq (|P|-1)L$, such that $v_i=F_{P,i}(v_{\min P})$.
We also define $F_{P,\min P}$ as the identity mapping. 
According to the locality of $\phi$, if $\mathbf v'=(v_1',\dots,v_p')\in V^p$
defines the same partition, types, and compositions, then
$\mathbf v'\in X$.  For fixed partition $\mathcal P$, types $(T_P)_{P\in\mathcal P}$,
 and compositions $(F_{P,i})_{i\in P\in\mathcal P}$, 
the corresponding subset $X'$ of $X$ is included in a (reshuffled) product $Y$
of sets of tuples of the form $(F_{P,i}(x_{\min P}))$ for 
$v_{\min P}\in W_{P}$, and is the set of all $v\in G$ such that $B_{(|P|-1)L}(\mathbf G,v)=T_P$.
Hence $W_P$ is measurable and (as each $F_{P,i}$ is measurable)
$Y$ is a measurable subset of $G^{|P|}$. Of course, this product may contain
tuples $\mathbf v$ defining another partition. A simple induction and inclusion/exclusion
argument shows that $X'$ is measurable. As $X$ is the union of a finite number of
such sets, $X$ is measurable.
\myqed\end{proof}

We now relate graphings to ${\rm FO}$-limits of bounded degree graphs.
We shall make use of the following lemma which reduces a graphing to its essential support.

\begin{lemma}[Cleaning Lemma]
Let $\mathbf G=(V, T_1,\dots,T_d, \mu)$ be a graphing.

Then there exists a subset $X\subset V$ with $0$ measure such that
$X$ is globally invariant by each of the $T_i$ and
$\mathbf G'=(V-X,T_1,\dots,T_d, \mu)$ is a graphing
such that  for every finite rooted colored graph $(F,o)$
 and integer $r$ the following equation holds:
$$
\mu(D_r(\mathbf G',(F,o)))=\mu(D_r(\mathbf G,(F,o)))
$$
(which means that $\mathbf G'$ is equivalent to $\mathbf G$) and
$$D_r(\mathbf G',(F,o))\neq\emptyset\quad\iff\quad\mu(D_r(\mathbf G',(F,o)))>0.$$
\end{lemma}
\begin{proof}
For a fixed $r$, define $\mathcal F_r$ has the set of all (isomorphism types of)
finite rooted $k$-edge colored graphs $(F,o)$ with radius at most $r$ such that 
\linebreak $\mu(D_r(\mathbf G,(F,o)))=0$.
Define
$$
X=\bigcup_{r\in\bbbn}\bigcup_{(F,o)\in \mathcal F_r}D_r(\mathbf G,(F,o)).
$$
Then $\mu(X)=0$, as it is a countable union of $0$-measure sets. 

We shall now prove that $X$ is a union of connected components of $\mathbf G$, and thus
$X$ is globally invariant by each of the $T_i$.
Namely, if $x\in X$ and $y$ is adjacent to $x$, then $y\in X$.
Indeed: if $x\in X$ then there exists an integer $r$ such that
$\mu(D(\mathbf G,B_r(\mathbf G,x)))=0$. But it is easily checked that

$$
\mu(D(\mathbf G,B_{r+1}(\mathbf G,y)))\leq d \cdot\mu(D(\mathbf G,B_r(\mathbf
G,x))).
$$
Hence $y\in X$.
It follows that for every $1\leq i\leq d$ we have $T_i(X)=X$. 
So we can define the graphing $\mathbf G'=(V-X,T_1,\dots,T_d,\mu)$.

Let $(F,o)$ be a rooted finite colored graph. Assume there exists
$x\in\mathbf G'$ such that $B_r(\mathbf G',r)\simeq (F,o)$. As $X$ is a union
of connected components, we also have $B_r(\mathbf G,r)\simeq (F,o)$ and
$x\notin X$. It follows that $\mu(D(\mathbf G,(F,o)))>0$  hence it holds that 
$\mu(D_r(\mathbf G',(F,o)))>0$.
\myqed\end{proof}

The cleaning lemma allows us a clean description of ${\rm FO}$-limits in the bounded
degree case:

\begin{mytinted}
\begin{theorem}
\label{thm:Bounded}
Let $(G_n)_{n\in\bbbn}$ be a FO-convergent sequence of finite graphs with
maximum degree $d$, with $\lim_{n\rightarrow\infty}|G_n|=\infty$. Then
there exists a graphing $\mathbf G$, which is the disjoint union
of a graphing $\mathbf G_0$ and a countable graph $\hat G$ such that
\begin{itemize}
  \item The graphing $\mathbf G$ is a {\rps} ${\rm FO}$-limit of the sequence $(G_n)_{n\in\bbbn}$.
  \item The graphing $\mathbf G_0$ is a BS-limit of the sequence $(G_n)_{n\in\bbbn}$ such that
  $$D_r(\mathbf G_0,(F,o))\neq\emptyset\quad\iff\quad\mu(D_r(\mathbf G_0,(F,o)))>0.$$
  \item The countable graph $\hat G$ is an elementary limit of the sequence $(G_n)_{n\in\bbbn}$.
\end{itemize}
\end{theorem}
\end{mytinted}
\begin{proof}
Let $\mathbf G_0$ be a BS-limit, which has been ``cleaned'' using the previous
lemma, and let $\hat G$ be an elementary limit of $G$. It is clear that 
$\mathbf G=\mathbf G_0\cup\hat G$ is also a BS-limit of the sequence, so the lemma amounts
in proving that $\mathbf G$ is elementarily equivalent to $\hat G$.

According to Hanf's theorem \cite{Hanf1965}, it  is sufficient to prove that for
all integers $r,t$ and for every rooted finite graph $(F,o)$ (with
maximum degree $d$) the following equality holds:

$$\min(t,|D_r(\mathbf G,(F,o))|)=\min(t,|D_r(\hat G,(F,o))|).
$$
Assume for contradiction that this is not the case. Then 
$|D_r(\hat G,(F,o))|<t$ and $D_r(\mathbf G_0,(F,o))$ is not empty. However, as
$\mathbf G_0$ is clean, this implies \linebreak $\mu(D_r(\mathbf G_0,(F,o)))=\alpha>0$.
It follows that for every sufficiently large $n$ it holds that 
$|D_r(G_n,(F,o))|>\alpha/2\, |G_n|>t$. Hence $|D_r(\hat G,(F,o))|>t$, contradicting
our hypothesis.

That $\mathbf G$ is a {\rps} then follows from Lemma~\ref{lem:graphingsb}.
\myqed\end{proof}

\begin{remark}
Not every graphing with maximum degree $2$ is an ${\rm FO}$-limit modeling of a
sequence of finite graphs (as it needs not be an elementary limit of finite graphs). Indeed: let $\mathbf G$ be a graphing that is an
${\rm FO}$-limit modeling of the sequence of cycles. The 
disjoint union of $\mathbf G$ and a ray is a graphing
$\mathbf G'$, which has the property that all its vertices  but one have degree
$2$, the exceptional vertex having degree $1$. As this property is not satisfied
by any finite graph, $\mathbf G'$ is not the ${\rm FO}$-limit of a sequence
of finite graphs.
\end{remark}

Let us finish this section by giving an interesting example, which shows
that the cleaning lemma sometimes applies in a non-trivial way:
\begin{example}
Consider the graph $G_n$ obtained from a de Bruijn sequence of
length $2^n$ as shown Fig~\ref{fig:debruijn}.

\begin{figure}[ht]
\begin{center}
\includegraphics[width=.5\textwidth]{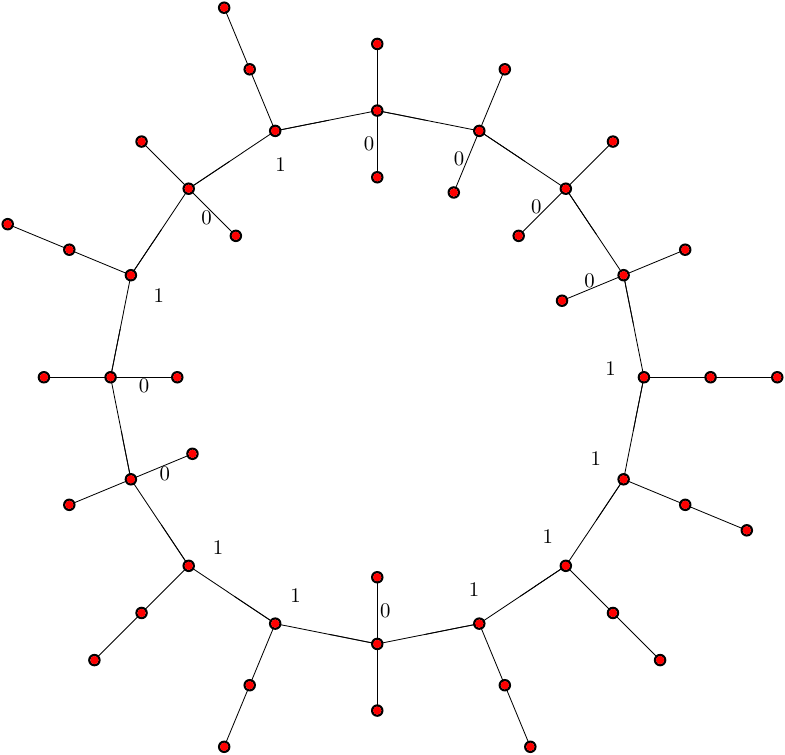}
\end{center}
\caption{The graph $G_n$ is constructed from a de Bruijn sequence of length
$2^n$.}
\label{fig:debruijn}
\end{figure}

It is easy to define a graphing $\mathbf G$, which is the limit of the sequence
$(G_n)_{n\in\bbbn}$: as vertex set, we consider the rectangle
$[0;1)\times[0;3)$. We define a measure preserving function $f$ and two
measure preserving involutions $T_1,T_2$ as follows:
\begin{align*}
f(x,y)&=\begin{cases}
(2x,y/2)&\text{if }x<1/2\text{ and }y<1\\
(2x-1,(y+1)/2)&\text{if }1/2\leq x\text{ and }y<1\\
(x,y)&\text{otherwise}
\end{cases}\\
T_1(x,y)&=\begin{cases}
(x,y+1)&\text{if }y<1\\
(x,y-1)&\text{if }1\leq y<2\\
(x,y)&\text{otherwise}
\end{cases}
\end{align*}
\begin{align*}
T_2(x,y)&=\begin{cases}
(x,y+1)&\text{if }x<1/2\text{ and }1\leq y<2\\
(x,y+2)&\text{if }1/2\leq x\text{ and }y<1\\
(x,y-1)&\text{if }x<1/2\text{ and }2\leq y\\
(x,y-2)&\text{if }1/2\leq x\text{ and }2\leq y\\
(x,y)&\text{otherwise}
\end{cases}\\
\end{align*}

Then the edges of $\mathbf G$ are the pairs $\{(x,y),(x',y')\}$ such that
$(x,y)\neq(x',y')$ and either $(x',y')=f(x,y)$, or $(x,y)=f(x',y')$, or
$(x',y')=T_1(x,y)$, or $(x',y')=T_2(x,y)$.

If one considers a random root $(x,y)$ in $\mathbf G$, then the connected
component of $(x,y)$ will almost surely be a rooted line with some decoration,
as expected from what is seen from a random root in a sufficiently large $G_n$.
However, special behaviour may happen when $x$ and $y$ are rational. Namely, it
is possible that the connected component of $(x,y)$ becomes finite. For
instance, if $x=1/(2^n-1)$ and $y=2^{n-1}x$ then the orbit of $(x,y)$ under the
action of $f$  has length $n$; thus the connected component of $(x,y)$ in
$\mathbf G$ has order $3n$. Of course, such finite connected components do
not appear in $G_n$.
Hence, in order to clean $\mathbf G$, infinitely many components have to be
removed.
\end{example}

Let us give a simple example exemplifying the distinction between BS and ${\rm
FO}$-convergence for graphs with bounded degree.
\begin{example}
Let $G_n$ denote the $n\times n$ grid. The BS-limit object is a probability
distribution concentrated on the infinite grid with a specified root. A limit
graphing can be described as the Lebesgue measure on $[0,1]^2$, where $(x,y)$ is adjacent
to $(x\pm\alpha \bmod 1,y\pm\alpha \bmod 1)$ for some irrational number $\alpha$.

This graphing, however, is not an ${\rm FO}$-limit of the sequence $(G_n)_{n\in\bbbn}$ as 
every ${\rm FO}$-limit has to contain four vertices of degree $2$. An ${\rm FO}$-limit 
graphing can be described as the above graphing restricted to $[0,1)^2$ (obtained by deleting all 
vertices with $x=1$ or $y=1$). One checks for instance that this graphing contains four 
vertices of degree $2$ (the vertices $(\alpha,\alpha)$, $(1-\alpha,\alpha)$, $(\alpha,1-\alpha)$, and $(1-\alpha,1-\alpha)$) and infinitely many vertices of degree $3$. 
\end{example}

We want to stress that our general and unifying approach to structural limits 
was not developed for its own sake and that it 
 provided a proper setting (and, yes, encouragement) for the study of classes of sparse graphs.
 So far classes of graphs with bounded degree are the only classes of sparse graphs where the structural 
 limits were constructed  efficiently.
 (Another example of limits of sparse graphs is provided by scaling limits of transitive graphs
 \cite{Benjamini2012} which proceeds in different direction and is not considered here.)

\section{Decomposing Sequences: the Comb Structure}
\label{sec:decompose}
The combinatorics of limits of equivalence relations (such as components) is
complicated. We start this analysis by considering the
combinatorics of ``large'' equivalence classes. This leads to the notion
of spectrum, which will be analyzed in this section.

\subsection{Spectrum of a First-order Equivalence Relation}
\begin{definition}[$\varpi$-spectrum]
Let $\mathbf A$ be a $\sig$-{\rps} (with measure $\nu_{\mathbf A}$), 
and let $\varpi\in{\rm FO}_{2}(\sig)$ be a formula expressing a component relation on $\mathbf A$ (see Definition~\ref{def:pirel}).
Let $\{C_i: i\in\Gamma\}$ be the set of all the $\varpi$-equivalence classes of 
$A$, and let $\Gamma_+$ be the (countable) subset of 
$\Gamma$ of the indexes $i$ such that $\nu_{\mathbf A}(C_i)>0$.

The {\em $\varpi$-spectrum} ${\rm Sp}_{\varpi}(\mathbf A)$ 
of $\mathbf A$ is the (countable) sequence of the values 
$\nu_{\mathbf A}(C_i)$ (for $i\in\Gamma_+$) ordered in non-increasing order.

\end{definition}
\begin{lemma}
\label{lem:spectr}
For $k\in\bbbn$, let $\varpi^{(k)}$ be the formula
$\bigwedge_{i=1}^{k}\varpi(x_i,x_{i+1})$. Then the following equation holds:

$$\sum_{i\in\Gamma_+}\nu_{\mathbf A}(C_i)^{k+1}=
\langle\varpi^{(k)},\mathbf A\rangle.$$
\end{lemma}
\begin{proof}
Let $k\in\bbbn$. Define
$$D_{k+1}=\{(x_1,\dots,x_{k+1})\in A^{k+1}:\
A\models\varpi^{k}(x_1,\dots,x_{k+1})\}.$$

According to Lemma~\ref{lem:cc}, each $C_i$ is measurable, thus
$\bigcup_{i\in\Gamma_+}C_i$ is measurable and so is
$R=A\setminus\bigcup_{i\in\Gamma_+}C_i$.

Considering the indicator function $\mathbf 1_{D_{k+1}\cap R^{k+1}}$
of $D_{k+1}\cap R^{k+1}$ and applying Fubini's theorem, we get
$$
\int_{A^{k+1}}\mathbf 1_{D_{k+1}\cap R^{k+1}}\,{\rm d}\nu_{\mathbf A}^{k+1}=
\idotsint \mathbf 1_{R}(x_1,\dots,x_{k+1})\ 
{\rm d}\nu_{\mathbf A}(x_1,\dots, {\rm d}\nu_{\mathbf A}(x_{k+1})=0.$$
as for every fixed $a_1,\dots,a_k$ (with $a_1\in C_\alpha$, for some
$\alpha\in\Gamma\setminus\Gamma_+$) we have 

$$0\leq \int \mathbf
1_{R}(a_1,\dots,a_k,x_{k+1})\,{\rm d}\nu_{\mathbf A}(x_{k+1})\leq \nu_{\mathbf A}(C_\alpha)=0.$$

It follows (by countable additivity) that 
$$\langle\varpi^{(k)},\mathbf A\rangle=\nu_{\mathbf A}^{k+1}(D_{k+1})=
\nu_{\mathbf A}^{k+1}(\bigcup_{i\in\Gamma_+}C_i^{k+1})
=\sum_{i\in\Gamma_+}\nu_{\mathbf A}(C_i)^{k+1}.
$$
\myqed\end{proof}
It follows from Lemma~\ref{lem:spectr} that the spectrum ${\rm
Sp}_{\varpi}(\mathbf A)$ is computable from the sequence of (non-increasing) values
$(\langle \varpi^{(k)},\mathbf A\rangle)_{k\in\bbbn}$. 

We assume that every finite sequence $\mathbf x=(x_1,\dots,x_n)$ of positive reals 
is implicitly embedded in an infinite sequence by defining $x_i=0$ for $i>n$.
Recall the usual $\ell_k$ norms:
$$\|\mathbf x\|_k=\Bigl(\sum_{i} |x_i|^k\Bigr)^{1/k}.$$
Hence above equations rewrite as
\begin{equation}
\|{\rm Sp}_{\varpi}(\mathbf A)\|_{k+1}=\langle \varpi^{(k)},
\mathbf A\rangle^{1/(k+1)}.
\end{equation}

We shall prove that the spectrum is, in a certain sense, defined by a continuous function.
We need the following technical lemma.

\begin{lemma}
\label{lem:0}
For each $n\in\bbbn$, let $\mathbf a_{n}=(a_{n,i})_{i\in\bbbn}$ be a non-increasing sequence of positive
real numbers with bounded sum (i.e $\|\mathbf a_{n}\|_1<\infty$ for every $n\in\bbbn$). 

Assume that for every integer $k\geq 1$ the limit
$s_k=\lim_{n\rightarrow\infty}\|\mathbf a_{n}\|_k$
exists.

Then $(\mathbf a_n)_{n\in\bbbn}$ converges in the space
$c_0$ of all sequences converging to zero 
(with norm $\|\,\cdot\,\|_\infty$).
\end{lemma}
\begin{proof}
We first prove that the sequences converge pointwise, that is that 
 there exists a sequence $\mathbf x=(x_i)_{i\in\bbbn}$ such that for every $i\in\bbbn$ the following equation holds:
 $$x_i=\lim_{n\rightarrow\infty} a_{n,i}.$$

For every $\epsilon>0$, if $s_k<\epsilon$ then $a_{n,1}<2\epsilon$ for
all sufficiently large values of $n$.
Thus if $s_k=0$ for some $k$, the limit
$\lim_{n\rightarrow}a_{n,i}$ exists for every $i$ and is null.
Thus, we can assume that $s_k$ is strictly positive for every $k\in\bbbn$.  

Fix $k\in\bbbn$. There exists $N\in\bbbn$ such that for
every $n\geq N$ it holds that 
$|s_k^k-\|\mathbf a_{n}\|_k^k|<s_k^k/k$.
As $(a_{n,i})_{i\in\bbbn}$ is a non-increasing sequence of positive real
numbers, for every $n\neq N$ the following inequality holds 
$$a_{n,1}^k\leq \|\mathbf a_{n}\|_k^k<s_k^k (1+1/k)$$
and
$$a_{n,1}^{k-1}\geq \|\mathbf a_{n}\|_k^k>s_k^k (1-1/k).$$
Hence
$$
\log s_k+\frac{\log (1+1/k)}{k}\geq
\log a_{n,1}\geq \bigl(1+\frac{1}{k-1}\bigr)\bigl(\log s_k+\frac{\log
(1-1/k)}{k}\bigr).$$
Thus $x_1=\lim_{n\rightarrow\infty}a_{n,1}$ exists and
$x_1=\lim_{k\rightarrow\infty}s_k$.
Inductively, we get that for each $i\in\bbbn$, the limit 
$x_i=\lim_{n\rightarrow\infty}a_{n,i}$ exists and that

$$
x_i=\lim_{k\rightarrow\infty}(s_k^k-\sum_{j<i}x_j^k)^{1/k}.
$$ 

We now prove that the converge is uniform, that is that
for every $\epsilon>0$ there exists $N$ such that for every $n\geq N$ the following inequality holds:
$$
\|\mathbf x-\mathbf a_{n}\|_\infty<\epsilon.
$$
As $\mathbf a_n\in \ell_1$ and $\|\mathbf a_n\|_1$ converges there exists $M$ such that 
$\|\mathbf a_n\|_1\leq M$ for every $n\in\bbbn$.
Let $\epsilon>0$. Let $A=\min\{i: x_i\leq\epsilon/3\}$. 
(Note that $A\leq 3M/\epsilon$.)
There exists $N$ such that for every $n\geq N$ it holds that 
$\sup_{i\leq A}|x_i-a_{n,i}|<\epsilon/3$.
Moreover, for every $i>A$ the following inequality holds:
$$
0\leq a_{n,i}\leq a_{n,A}<x_A+\epsilon/3<2\epsilon/3. 
$$
As $0\leq x_i\leq \epsilon/3$ for every $i>A$ the following inequality holds:
$$
|x_i-a_{n,i}|<\epsilon
$$
for every $i>A$ (hence for every $i$).
Thus $(\mathbf a_n)_{n\in\bbbn}$ converges in $\ell_\infty$.
As obviously each $\mathbf a_n$ has $0$ limit,
$(\mathbf a_n)_{n\in\bbbn}$ converges in $c_0$. 
\myqed\end{proof}

\begin{mytinted}
\begin{lemma}
\label{lem:contsp}
Let $\sig$ be a signature.
The mapping $\mathbf A\mapsto{\rm Sp}_\varpi(\mathbf A)$ is
a continuous mapping from the space of $\sig$-modelings with
a component relation $\varpi$ (with the topology
of ${\rm FO}^{\rm local}(\sig)$-convergence) to
the space $c_0$ of all sequences converging to zero (with $\|\,\cdot\,\|_\infty$
norm).
\end{lemma}
\end{mytinted}
\begin{proof}
Assume $\mathbf A_n$ is an ${\rm FO}^{\rm local}(\sig)$-convergent
 sequence of $\sig$-modelings.

Let $(\spect_{n,1},\dots,\spect_{n,i},\dots)$
be the $\varpi$-spectrum of $\mathbf A_n$ (extended by zero values if finite),
and let $\mathbf a_n=(a_{n,i})_{i\in\bbbn}$ be the sequence
defined by $a_{n,i}=\spect_{n,i}^2$. Then
for every integer $k\geq 1$ it holds that 
$$
\|a_n\|_k=\|{\rm Sp}_\varpi(\mathbf A_n)\|_{2k}^{2}
= \langle \varpi^{(2k-1)},
\mathbf A_n\rangle^{1/k}.
$$
Hence $s_k=\lim_{n\rightarrow\infty}\|a_n\|_k$ exists.
According to Lemma~\ref{lem:0}, $(\mathbf a_n)_{n\in\bbbn}$ converges in
$c_0$, thus so does $({\rm Sp}_\varpi(\mathbf A_n))_{n\in\bbbn}$.
\myqed\end{proof}

\begin{definition}
Let $(\mathbf A_n)_{n\in\bbbn}$ be a sequence of finite $\sig$-structures.
Let $\varpi$ be a component relation, and for simplicity assume $\varpi\in\sig$.
In the following, we assume that $\varpi$-spectra are 
extended to infinite sequences by adding zeros if necessary.

\begin{itemize}
  \item  The sequence $(\mathbf A_n)_{n\in\bbbn}$ is {\em $\varpi$-nice}
  if ${\rm Sp}_\varpi(\mathbf A_n)$ converges pointwise;
\item The {\em limit $\varpi$-spectrum} of a $\varpi$-nice sequence $(\mathbf
A_n)_{n\in\bbbn}$ is the pointwise limit of ${\rm Sp}_\varpi(\mathbf A_n)$;
\item the {\em $\varpi$-support}  is the set $I$
 of the indexes $i$ for which the limit $\varpi$-spectrum is non-zero;
\item the sequence has {\em full} $\varpi$-spectrum if, for every index $i$ not in the $\varpi$-support,
there is some $N$ such that the $i$th value of ${\rm Sp}_\varpi(\mathbf A_n)$ is zero for every
$n>N$.   
\end{itemize}
\end{definition}

As proved in Lemma~\ref{lem:contsp}, every ${\rm FO}^{\rm local}$-convergent
sequence is $\varpi$-nice. 

\begin{lemma}
\label{lem:loc1}
Let $(\mathbf A_n)$ be a $\varpi$-nice sequence of $\sig$-structures
with empty $\varpi$-support.

 Then the following conditions are equivalent:
 \begin{enumerate}
  \item the sequence $(\mathbf A_n)$ is ${\rm FO}^{\rm local}$-convergent;
  \item the sequence $(\mathbf A_n)$ is ${\rm FO}_1^{\rm local}$-convergent.
\end{enumerate}
Moreover, for every $\varpi$-local formula $\phi$ with $p>1$ free variables it
holds that 
$$
\lim_{n\rightarrow\infty}\langle\phi,\mathbf A_n\rangle=0.
$$ 
\end{lemma}
\begin{proof}
${\rm FO}^{\rm local}$-convergence obviously implies ${\rm FO}_1^{\rm
local}$-convergence.
So, assume that $(\mathbf A_n)_{n\in\bbbn}$ is ${\rm FO}_1^{\rm local}$-convergent, and
let $\phi$ be a $\varpi$-local first-order formula with $p>1$ free variables.
For $n\in\bbbn$, let $\mathbf B_{n,i}$ ($i\in\Gamma_n$) denote the connected
components of $\mathbf A_n$. As $(\mathbf A_n)$ is $\varpi$-nice and has empty
$\varpi$-support, there exists for every $\epsilon>0$ an integer $N$ such that
for $n>N$ and every $i\in\Gamma_n$ it holds that  $|B_{n,i}|<\epsilon |A_n|$.
Then, according to Corollary~\ref{cor:stlocpair}, for $n>N$

\begin{align*}
\langle\phi,\mathbf
A_n\rangle&=\sum_{i\in\Gamma_n}\biggl(\frac{|B_{n,i}|}{|A_n|}\biggr)^p\langle\phi,\mathbf B_{n,i}\rangle\\
&\leq \sum_{i\in\Gamma_n}\biggl(\frac{|B_{n,i}|}{|A_n|}\biggr)^p\\
&< \sum_{i\in\Gamma_n} \frac{|B_{n,i}|}{|A_n|}\epsilon^{p-1}=\epsilon^{p-1}.
\end{align*}

Hence $\langle\phi,\mathbf A_n\rangle$ converges (to $0$) as $n$ grows to infinity. It
follows that $(\mathbf A_n)_{n\in\bbbn}$ is ${\rm FO}^{\rm local}$-convergent,
according to Theorem~\ref{thm:sl}.
\myqed\end{proof}

\begin{lemma}
\label{lem:grp}
Let $(\mathbf A_n)_{n\in\bbbn}$ be an ${\rm FO}^{\rm local}$-convergent sequence
of finite $\sig$-struc\-tures, with component relation
$\varpi\in\sig$ and limit $\varpi$-spectrum $(\spect_i)_{i\in I}$.
For $n\in\bbbn$, let $\mathbf B_{n,i}$ be the connected components
of $\mathbf A_n$ order in non-decreasing order (with $\mathbf B_{n,i}$ empty if $i$ is greater 
than the number of connected components of $\mathbf A_n$).
Let $a\leq b$ be the first and last occurrence of 
$\spect_a=\spect_b$ in the $\varpi$-spectrum and 
let $\mathbf A_n'$ be the union of all the $\mathbf B_{n,i}$ for
$a\leq i\leq b$.

Then $(\mathbf A_n')_{n\in\bbbn}$ is ${\rm FO}$-convergent if $\spect_a>0$ and
${\rm FO}^{\rm local}$-convergent if $\spect_a=0$.

Assume moreover that $(\mathbf A_n)_{n\in\bbbn}$ has a {\rps} ${\rm FO}^{\rm local}$-limit ${\mathbf L}$.
Let $\mathbf L'$ be the union of the connected components $\mathbf L_i$ of $\mathbf L$ 
with $\nu_{\mathbf L}(L_i)=\spect_a$. Equip $\mathbf L'$ with the $\sigma$-algebra $\Sigma_{\mathbf L'}$ which
is the restriction of $\Sigma_{\mathbf L}$ to $L'$ and the probability measure $\nu_{\mathbf L'}$ defined
by $\nu_{\mathbf L'}(X)=\nu_{\mathbf L}(X)/\nu_{\mathbf L}(L')$ (for $X\in\Sigma_{\mathbf L'}$).

Then $\mathbf L'$ is a {\rps} ${\rm FO}$-limit of $(\mathbf A_n')_{n\in\bbbn}$ if $\spect_a>0$ and
a {\rps} ${\rm FO}^{\rm local}$-limit of $(\mathbf A_n')_{n\in\bbbn}$ if $\spect_a=0$.
\end{lemma}
\begin{proof}
Extend the sequence $\spect$  to the null index by defining $\spect_0=2$.
Let $r=\min(\spect_{a-1}/\spect_a, \spect_b/\spect_{b+1})$ (if
$\spect_{b+1}=0$ simply define $r=\spect_{a-1}/\spect_a$). 
Notice that $r>1$.  
Let $\phi$ be a $\varpi$-local formula with $p$ free variables.
According to Corollary~\ref{cor:stlocpair} the following equation holds:

$$
\langle\phi,\mathbf A_n\rangle=\sum_{i}\biggl(\frac{|B_{n,i}|}{|A_n|}\biggr)^p
\langle\phi,\mathbf B_{n,i}\rangle.
$$
In particular, it holds that 
$$\langle \varpi^{(p)},\mathbf A_n\rangle=\sum_i \left(\frac{|B_{n,i}|}{|A_n|}\right)^p.
$$
Let $\alpha>1/(1-r^p)$.
Define 
$$w_{n,i}=\biggl(\frac{|B_{n,i}|}{|A_n|}\biggr)^p(\alpha+\langle\phi,\mathbf B_{n,i}\rangle).$$

From the definition of $r$ it follows that for each 
 $n\in\bbbn$, $w_{n,i}>w_{n,j}$ if $i<a$ and $j\geq a$ or $i\leq b$ and $j>b$.
Let $\sigma\in \Sym{\omega}$ be a permutation of $\bbbn$, 
such that
$a_{n,i}=w_{n,\sigma(i)}$ is non-increasing. It holds that 
$$\sum_i a_{n,i}=\sum_i w_{n,i}=\alpha \langle
\varpi^{(p)},\mathbf A_n\rangle+\langle\phi,\mathbf A_n\rangle.$$ 

Hence
$$\lim_{n\rightarrow\infty}\sum_i a_{n,i}^p$$
exists.
According to Lemma~\ref{lem:0} it follows that
for every $i\in\bbbn$ the limit $\lim_{n\rightarrow\infty}a_{n,i}$ exists.
Moreover, as $\sigma$ globally preserves the set $\{a,\dots,b\}$ it follows that
the limit
$$
d=\lim_{n\rightarrow\infty}\sum_{i=a}^b\biggl(\frac{|B_{n,i}|}{|A_n|}\biggr)^p(\alpha+\langle\phi,\mathbf B_{n,i}\rangle)
$$
exists. As for every $i\in\{a,\dots,b\}$ it holds that 
$\lim_{n\rightarrow\infty}|B_{n,i}|/|A_n|=\spect_a$ and as
$\langle\phi,\mathbf A_{n}'\rangle=\sum_{i=a}^b(|B_{n,i}|/|A_n|)^p\langle\phi,\mathbf B_{n,i}\rangle$
we deduce
$$
\lim_{n\rightarrow\infty}\langle\phi,\mathbf A_{n}'\rangle=d-(b-a+1)\alpha.
$$
Hence $\lim_{n\rightarrow\infty}\langle\phi,\mathbf A_{n}'\rangle$ exists for every
$\varpi$-local formula and, according to Theorem~\ref{thm:sl}, the sequence $(\mathbf A_n')_{n\in\bbbn}$ is
${\rm FO}^{\rm local}$-convergent.

Assume $\spect_a>0$. Let $N=b-a+1$. To each sentence $\theta$ we associate the formula $\widetilde{\theta}\in{\rm FO}_N^{\rm local}$ 
that asserts that the substructure induced by the closed neighborhood of $x_1,\dots,x_N$ satisfies $\theta$
and that $x_1,\dots,x_N$ are pairwise distinct and non-adjacent. For sufficiently large $n$, the structure
$\mathbf A_n'$ has exactly $N$ connected components.
It is easily checked that if $\mathbf A_n'$ does not satisfy $\theta$ then
$\langle\widetilde{\theta},\mathbf A_n'\rangle=0$, although if  $\mathbf A_n'$ does satisfy $\theta$
then 
$$\langle\widetilde{\theta},\mathbf A_n'\rangle\geq \Bigl(\frac{\min_{a\leq i\leq b} |B_{n,i}|}{\sum_{i=a}^b |B_{n,i}|}\Bigr)^N,$$
hence $\langle\widetilde{\theta},\mathbf A_n'\rangle>(2N)^{-N}$ for all sufficiently large $n$. As
$\langle\widetilde{\theta},\mathbf A_n'\rangle$ converges for every sentence $\theta$, we deduce that
the sequence $(\mathbf A_n')_{n\in\bbbn}$ is elementarily convergent.
According to Theorem~\ref{thm:fole},  the sequence $(\mathbf A_n')_{n\in\bbbn}$ is thus ${\rm FO}$-convergent.

Now assume that $(\mathbf A_n)_{n\in\bbbn}$ has a {\rps} ${\rm FO}^{\rm local}$-limit ${\mathbf L}$.
First note that $L_i$ being an equivalence class of $\varpi$ it holds that  $L_i\in\Sigma_{\mathbf L}$,
hence $L'\in\Sigma_{\mathbf L}$ and $\nu_{\mathbf L}(L')$ is well defined. For every $\varpi$-local formula
$\phi\in{\rm FO}_p(\sig)$ it holds, according to Lemma~\ref{lem:stlocpair}, that
\begin{align*}
\langle\phi,\mathbf L'\rangle&=\sum_{i=a}^b \nu_{\mathbf L'}(L_i)^p \langle\phi,\mathbf L_i\rangle\\
&=\frac{1}{\nu_{\mathbf L}(L')^p}\sum_{i=a}^b \nu_{\mathbf L}(L_i)^p \langle\phi,\mathbf L_i\rangle.\\
\end{align*}
We deduce that
$$
\langle\phi,\mathbf L'\rangle=\lim_{n\rightarrow\infty}\langle\phi,\mathbf A_n'\rangle.
$$ 
According to Theorem~\ref{thm:sl}, it follows that the same equality holds for every 
$\phi\in{\rm FO}^{\rm local}(\sig)$ hence
$\mathbf L'$ is a {\rps} ${\rm FO}^{\rm local}$-limit of the sequence
$(\mathbf A_n')_{n\in\bbbn}$.

As above, for $\spect_a>0$, if $\mathbf L'$ is a {\rps} ${\rm FO}^{\rm local}$-limit
of $(\mathbf A_n')_{n\in\bbbn}$ then it is a {\rps} ${\rm FO}$-limit.

\myqed\end{proof}

\begin{lemma}
\label{lem:refine}
Let $(\mathbf A_n)_{n\in\bbbn}$ be an ${\rm FO}$-convergent sequence
of finite $\sig$-structures, with component relation
$\varpi$ (expressing usual notion of connected components).
Assume all the $\mathbf A_n$ have at most $k$ connected components.
Denote by $\mathbf B_{n,1},\dots,\mathbf B_{n,k}$ these components (adding
 empty $\sig$-structures if necessary).

Assume that for each $1\leq i\leq k$ it holds that 
$\lim_{n\rightarrow\infty} |B_{n,i}|/|A_n|=1/k$.

Then there exists a sequence $(\sigma_n)_{n\in\bbbn}$ of permutations of $[k]$
such that for each $1\leq i\leq k$ the sequence
$(\mathbf B_{n,\sigma_n(i)})_{n\in\bbbn}$ is ${\rm FO}$-convergent.
\end{lemma}
\begin{proof}
To a formula $\phi\in{\rm FO}_p(\sig)$ we associate the $\varpi$-local formula
$\widetilde{\phi}\in{\rm FO}_p^{\rm local}(\sig)$ asserting that all the free variables
are $\varpi$-adjacent and that their closed neighborhood (that is their connected component)
satisfies $\phi$. Then essentially the same proof as above allows to refine
$\mathbf A_n$ into sequences such that $\langle\phi,\mathbf A_{n,i}'\rangle$ is constant
on the connected components of each of the $\mathbf A_n'$.  
Considering formulas allowing to split at least one of the sequences, we repeat this process
 (at most $k-1$ times) until
each $\mathbf A_{n,i}'$ contains equivalent connected components. Then, $\mathbf A_{n,i}'$ can be
split into connected components in an arbitrary order, thus obtaining the sequences $\mathbf B_{n,i}$.   
\myqed\end{proof}

So we have proved that a ${\rm FO}$-convergent can be decomposed by isolines of
the $\varpi$-spectrum (that is by groups of connected components with same asymptotic measure). In the next sections, we shall investigate how to refine this 
further.

\subsection{Sequences with Finite Spectrum}
For every $\varpi$-nice sequence \linebreak $(\mathbf A_n)_{n\in\bbbn}$ with finite support $I$, 
we define the {\em residue} $\mathbf R_n$ of $\mathbf A_n$ as the union of the connected components
$\mathbf B_{n,i}$ of $\mathbf A_n$ such that $i\notin I$.

\begin{figure}[ht]
\begin{center}
\includegraphics[width=.7\textwidth]{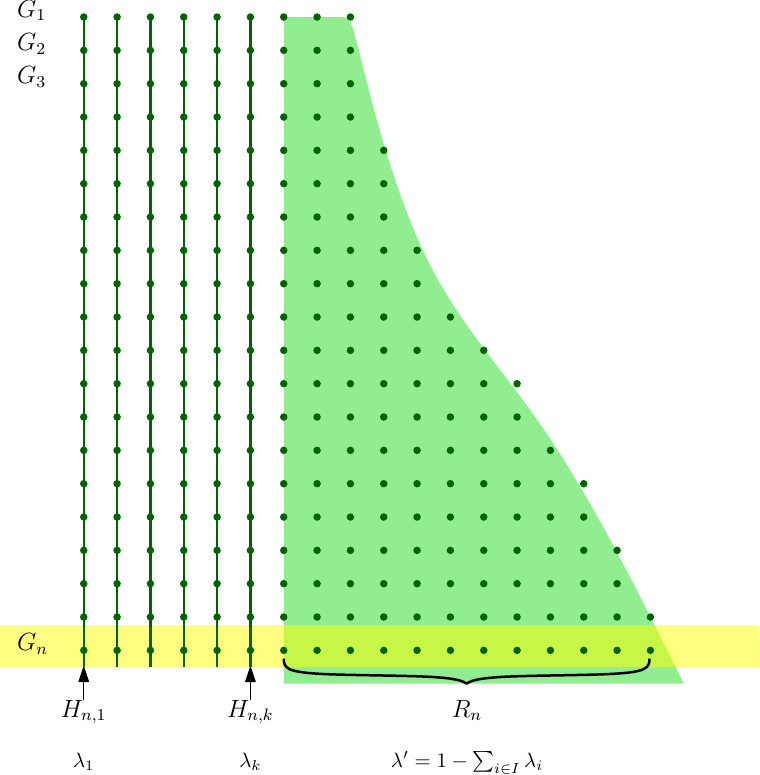}
\end{center}
\caption{Illustration of the Comb structure for sequences with finite support}
\label{fig:comb0}
\end{figure}

When one considers an ${\rm FO}^{\rm local}$-convergent sequence $(\mathbf A_n)$ with a finite support then 
the sequence of the residues forms a sequence which is either ${\rm FO}^{\rm local}$-convergent or
``negligible'' in the sense that $\lim_{n\rightarrow\infty} |R_n|/|A_n|=0$. This is formulated as follows:   

\begin{lemma}
\label{lem:residue0}
Let $(\mathbf A_n)_{n\in\bbbn}$ be a sequence of $\sig$-structures with component relation $\varpi$.
For each $n\in\bbbn$ and $i\in\bbbn$, let $\mathbf B_{n,i}$ be the $i$-th largest
connected component of $\mathbf A_n$.

Assume that $(\mathbf A_n)_{n\in\bbbn}$ is ${\rm FO}^{\rm local}$-convergent and
has finite spectrum $(\spect_i)_{i\in I}$. 
Let $\mathbf R_n$ be the residue of $\mathbf A_n$.

Then $\spect'=\lim_{n\rightarrow\infty} |R_n|/|A_n|$ exists and
either $\spect'=0$ or $(\mathbf R_n)_{n\in\bbbn}$ is ${\rm FO}^{\rm local}$-convergent.
\end{lemma}
\begin{proof}
Clearly, $\spect'=1-\sum_i\spect_i$. Assume $\spect'>0$.
First notice that for every $\epsilon>0$ there exists $N$ such that for every
$i>N$, the $\sig$-structure $\mathbf R_n$ has no connected component of size at least
$\epsilon/2\spect'|A_n|$ and $\mathbf R_n$ has order at least $\spect'/2 |A_n|$.
Hence, for every $i>N$, the $\sig$-structure $\mathbf R_n$ has no connected component of size at least
$\epsilon |R_n|$. According to Lemma~\ref{lem:loc1}, proving that
$(\mathbf R_n)_{n\in\bbbn}$ is ${\rm FO}^{\rm local}$-convergent reduces to proving that
$(\mathbf R_n)_{n\in\bbbn}$ is ${\rm FO}_1^{\rm local}$-convergent.

Let $\phi\in{\rm FO}_1^{\rm local}$. 
We group the $\sig$-structures $\mathbf B_{n,i}$ (for $i\in I$) by values of $\spect_i$
as $\mathbf A_{n,1}',\dots,\mathbf A_{n,q}'$. Denote by $c_j$ the common value of $\spect_i$
for the connected components $\mathbf B_{n,i}$ in $\mathbf A_{n,j}'$. 
According to Corollary~\ref{cor:stlocpair} it holds (as $\phi$ is clearly $\varpi$-local) that
\begin{align*}
\langle\phi,\mathbf A_n\rangle&=\sum_{i}\frac{|B_{n,i}|}{|A_n|}\langle\phi,\mathbf B_{n,i}\rangle\\
&=\sum_{i\in I}\frac{|B_{n,i}|}{|A_n|}\langle\phi,\mathbf B_{n,i}\rangle
+\sum_{i\notin I}\frac{|B_{n,i}|}{|A_n|}\langle\phi,\mathbf B_{n,i}\rangle\\
&=\sum_{j=1}^q \frac{|A_{n,j}'|}{|A_n|}\langle\phi,\mathbf A_{n,j}'\rangle
+\frac{|R_{n}|}{|A_n|}\langle\phi,\mathbf R_{n}\rangle.
\end{align*}
According to Lemma~\ref{lem:grp}, each sequence
$(\mathbf A_{n,j}')_{n\in\bbbn}$ is ${\rm FO}$-convergent. 
Hence the limit $\lim_{n\rightarrow\infty} \langle\phi,\mathbf R_{n}\rangle$ exists and we have
$$
\lim_{n\rightarrow\infty}\langle\phi,\mathbf R_{n}\rangle=\frac{1}{\spect'}\biggl(
\lim_{n\rightarrow\infty}\langle\phi,\mathbf A_n\rangle-
\sum_{j=1}^q c_j\ \lim_{n\rightarrow\infty}\langle\phi,\mathbf A_{n,j}'\rangle
\biggr).
$$
It follows that the sequence $(\mathbf R_n)_{n\in\bbbn}$ is ${\rm FO}^{\rm local}$-convergent.
\myqed\end{proof}
The following result finally determines the structure of converging sequences of (disconnected) $\sig$-structures
with finite support. This structure is called {\em comb structure}, see Fig~\ref{fig:comb0}. 
\begin{mytinted}
\begin{theorem}[Comb structure for $\sig$-structure sequences with finite spectrum]
\label{thm:comb0}
Let $(\mathbf A_n)_{n\in\bbbn}$ be an ${\rm FO}^{\rm local}$-convergent sequence of
finite $\sig$-structures with component relation $\varpi$ and finite spectrum $(\spect_i)_{i\in I}$.
Let $\mathbf R_n$ be the residue of $\mathbf A_n$.

Then there exists, for each $n\in\bbbn$, 
a permutation $f_n:I\rightarrow I$ such that the following holds
\begin{itemize}
  \item $\lim_{n\rightarrow\infty} \max_{i\notin I}|B_{n,i}|/|A_n|=0$;
  \item $\lim_{n\rightarrow\infty} |R_n|/|A_n|$ exists;
  \item for every $i\in I$, the sequence $(B_{n,f_n(i)})_{n\in\bbbn}$ is ${\rm
  FO}$-convergent and $\lim_{n\rightarrow\infty} |B_{n,f_n(i)}|/|A_n|=\spect_i$;
  \item either $\lim_{n\rightarrow\infty} |R_n|/|A_n|=0$, 
  or the sequence $(\mathbf R_n)_{n\in\bbbn}$ is ${\rm FO}^{\rm local}$-convergent.
\end{itemize}
Moreover, if $(\mathbf A_n)_{n\in\bbbn}$ is ${\rm FO}$-convergent 
then $(\mathbf R_n)_{n\in\bbbn}$ is elementary-convergent. 
\end{theorem}
\end{mytinted}
\begin{proof}
This lemma is a direct consequence of Lemmas~\ref{lem:grp},~\ref{lem:refine} and~\ref{lem:residue0}, 
except that we still have to prove ${\rm FO}$-convergence of $(\mathbf
R_n)_{n\in\bbbn}$ in the case where $(\mathbf A_n)_{n\in\bbbn}$ is ${\rm FO}$-convergent.
As $I$ is finite,   the elementary convergence 
of $(\mathbf R_n)_{n\in\bbbn}$ easily follows from the one of $(\mathbf A_n)$ and the one of the 
$(\mathbf B_{n,f_n(i)})$ for $i\in I$. 
\myqed\end{proof}

\subsection{Sequences with Infinite Spectrum}
Let $(\mathbf A_n)_{n\in\bbbn}$ be a $\varpi$-nice sequence with infinite spectrum (and support $I=\bbbn$).
In such a case, the notion of a residue becomes more tricky and will need
some technical definitions.
Before this, let us take the time to give an example illustrating the
difficulty of the determination of the residue $\mathbf R_n$ in the comb structure of 
sequences with infinite spectrum.
\begin{example}
Consider  the sequence
$(G_n)_{n\in\bbbn}$ where $G_n$ is the union of $2^n$ stars
$H_{n,1},\dots,H_{n,2^n}$, where the $i$-th star $H_{n,i}$
has order $2^{2^n}(2^{-i}+2^{-n})/2$. Then it holds that 
$$\spect_i=\lim_{n\rightarrow\infty}|H_{n,i}|/|G_n|=2^{-(i+1)}$$
 hence
$\sum_i\spect_i=1/2$ thus the residue asymptotically should contain half of the vertices 
of $G_n$!
An ${\rm FO}$-limit of this sequence is shown Fig.~\ref{fig:example1}.

\begin{figure}[h!t]
\begin{center}
\includegraphics[width=.5\textwidth]{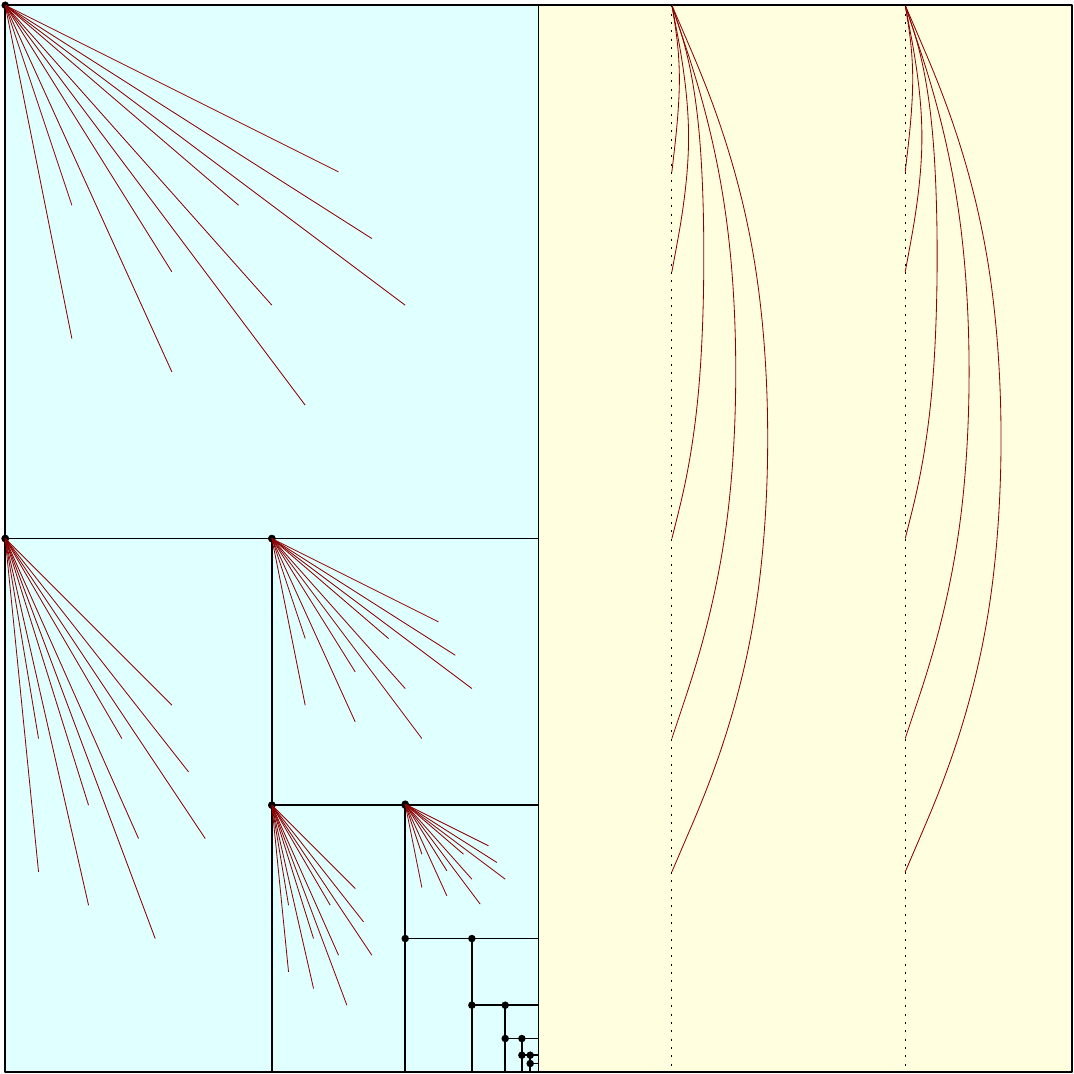}
\end{center}
\caption{An ${\rm FO}^{\rm local}$-limit. On the left side, each rectangle correspond
to a star with the upper left point as its center; on the right side, each
vertical line is a star with the upper point as its center.}
\label{fig:example1}
\end{figure}
\end{example}

This example is not isolated. In fact it is quite frequent in many of its
variants. To decompose such examples we need a convenient separation. This is
provided by the notion of clip.
\begin{definition}
\begin{itemize}
\item A {\em clip} of a $\varpi$-nice sequence $(\mathbf A_n)_{n\in\bbbn}$ with support $\bbbn$ is a
non-decreasing function $C:\bbbn\rightarrow\bbbn$ such that
$\lim_{n\rightarrow\infty}C(n)=\infty$ and
$$
\forall n'\geq n\quad \sum_{i=1}^{C(n)}
\biggl|\frac{|B_{n',i}|}{|A_{n'}|}-\spect_i\biggr|\leq
\sum_{i>C(n)}\spect_i$$
\item  The {\em residue} $\mathbf R_n$ of $\mathbf A_n$ with respect to a
 clip $C(n)$ is the disjoint union of the $\mathbf B_{n,i}$ for $i>C(n)$.
\end{itemize}
\end{definition} 
\begin{proposition}
\label{prop:clip}
Every $\varpi$-nice sequence $(\mathbf A_n)_{n\in\bbbn}$ with infinite support has a clip $C_0$, which is defined by 
$$
C_0(n)=\sup \biggl\{M,\quad (\forall n'\geq n)\ \sum_{i=1}^M
\biggl|\frac{|B_{n',i}|}{|A_{n'}|}-\spect_i\biggr|\leq
\sum_{i>M}\spect_i\biggr\}.$$
Moreover, $\lim_{n\rightarrow}C_0(n)=\infty$ and a non decreasing function $C$
is a clip of $(\mathbf A_n)_{n\in\bbbn}$ if and only if $C\leq C_0$ and $\lim_{n\rightarrow}C(n)=\infty$. 
\end{proposition}
\begin{proof}
Indeed, for each $n\in\bbbn$, the value $z_l(M)=\sup_{n'\geq n}\sum_{i=1}^M 
\biggl|\frac{|B_{n',i}|}{|A_{n'}|}-\spect_i\biggr|$ is non-decreasing function
of $C$ with $z_l(0)=0$, and $z_r(M)=\sum_{i>M}\spect_i$ is 
a decreasing function of $C$ with $z_r(0)=\sum_i\spect_i>0$ hence $C_0$ is well defined.
Moreover, for every integer $M$, let $\alpha=\sum_{i>M}\spect_i>0$. Then,
as $\lim_{n\rightarrow\infty}|B_{n',i}|/|A_{n'}|=\spect_i$ there exists
$N$ such that for every $n'\geq N$ and every $1\leq i\leq M$ it holds that  
$||B_{n',i}|/|A_{n'}|-\spect_i|\leq\alpha/M$ thus
for every $n'\geq N$ it holds that
$$\sum_{i=1}^{M}\biggl|\frac{|B_{n',i}|}{|A_{n'}|}-\spect_i\biggr|\leq\alpha=\sum_{i>M}\spect_i.$$
It follows that $C_0(N)\geq M$. Hence $\lim_{n\rightarrow\infty}C_0(n)=\infty$.
 
That a non decreasing function $C$ is a clip of $(\mathbf A_n)_{n\in\bbbn}$ if and only
if $C\leq C_0$ and $\lim_{n\rightarrow}C(n)=\infty$ follows directly from the definition.
\myqed\end{proof}
\begin{lemma}
\label{lem:lambdap}
Let $(\mathbf A_n)_{n\in\bbbn}$ be a $\varpi$-nice sequence with support $\bbbn$, and let 
$C$ be a clip of $(\mathbf A_n)_{n\in\bbbn}$.

Then the limit $\spect'=\lim_{n\rightarrow\infty}\frac{|R_n|}{|A_n|}$
exists and  $\spect'=1-\sum_i\spect_i$.
\end{lemma}
\begin{proof}
As $C$ is a clip, the following inequality holds for every $n\in\bbbn$  

$$\sum_i\spect_i-2\sum_{i>C(n)}\spect_i\leq
\sum_{i=1}^{C(n)}\frac{|B_{n,i}|}{|A_n|}\leq \sum_i\spect_i.
$$

Also, for every $\epsilon>0$ there exists $n$ such that
$|\sum_{i=1}^{C(n)}\spect_i-\sum_i\spect_i|<\epsilon$, that is:
$\sum_{i>C(n)}\spect_i<\epsilon.$
It follows that

$$
\lim_{n\rightarrow\infty}\sum_{i=1}^{C(n)}\frac{|B_{n,i}|}{|A_n|}=
\sum_i\spect_i.
$$
Hence the limit $\spect'=\lim_{n\rightarrow\infty}\frac{|R_n|}{|A_n|}$
exists and  $\spect'=1-\sum_i\spect_i$.
\myqed\end{proof}

\begin{lemma}
\label{lem:residue}
Let $(\mathbf A_n)_{n\in\bbbn}$ be a sequence of $\sig$-structures with component relation $\varpi$.
For each $n\in\bbbn$ and $i\in\bbbn$, let $\mathbf B_{n,i}$ be the $i$-th largest
connected component of $\mathbf A_n$ (if $i$ is at most equal to the number of connected components of $\mathbf A_n$,
the empty $\sig$-structure otherwise).

Assume that $(\mathbf A_n)_{n\in\bbbn}$ is ${\rm FO}$-convergent.

Let $C:\bbbn\rightarrow\bbbn$ be a clip of $(\mathbf A_n)_{n\in\bbbn}$, and let
 $\mathbf R_n$ be the residue of $\mathbf A_n$ with respect to $C$.

Let $\spect'=\lim_{n\rightarrow\infty} |R_n|/|A_n|$. Then
either $\spect'=0$ or $(\mathbf R_n)_{n\in\bbbn}$ is ${\rm FO}^{\rm local}$-convergent.
\end{lemma}
\begin{proof}
According to Lemma~\ref{lem:lambdap}, $\lim_{n\rightarrow\infty} |R_n|/|A_n|$ exists and $\spect'=1-\sum_i\spect_i$.
 Assume $\spect'>0$.
First notice that for every $\epsilon>0$ there exists $N$ such that for every
$i>N$, the $\sig$-structure $\mathbf R_n$ has no connected component of size at least
$\epsilon/2\spect'|\mathbf A_n|$ and $\mathbf R_n$ has order at least $\spect'/2 |A_n|$.
Hence, for every $i>N$, the $\sig$-structure $\mathbf R_n$ has no connected component of size at least
$\epsilon |R_n|$. According to Lemma~\ref{lem:loc1}, proving that
$(\mathbf R_n)_{n\in\bbbn}$ is ${\rm FO}^{\rm local}$-convergent reduces to proving that
$(\mathbf R_n)_{n\in\bbbn}$ is ${\rm FO}_1^{\rm local}$-convergent.

Let $\phi\in{\rm FO}_1^{\rm local}$ (thus $\phi$ is $\varpi$-local). 
Let $\epsilon>0$.
There exists $k\in\bbbn$ such that $\sum_{i\leq
k}\spect_i>1-\spect'-\epsilon/3$ and such that $\spect_{k+1}<\spect_k$.
We group the $\sig$-structures $\mathbf B_{n,i}$ (for $1\leq i\leq k$) by values of $\spect_i$
as $\mathbf A_{n,1}',\dots,\mathbf A_{n,q}'$. Denote by $c_j$ the common value of $\spect_i$
for the connected components $\mathbf B_{n,i}$ in $\mathbf A_{n,j}'$. 
According to Lemma~\ref{lem:grp}, each sequence
$(\mathbf A_{n,i}')_{n\in\bbbn}$ is ${\rm FO}$-convergent. Define
$$
\mu_i=\lim_{n\rightarrow\infty}\langle\phi,\mathbf A_{n,i}'\rangle.
$$
There exists
$N$ such that for every $n>N$ the following inequality holds
$$
\sum_{i=1}^q|\langle\phi,\mathbf A_{n,i}'\rangle-\mu_i|<\epsilon/3.
$$
According to Corollary~\ref{cor:stlocpair} it holds, 
for every $n>N$, that
\begin{align*}
\langle\phi,\mathbf A_n\rangle&=\sum_{i}\frac{|B_{n,i}|}{|A_n|}\langle\phi,\mathbf B_{n,i}\rangle\\
&=\sum_{i=1}^k\frac{|B_{n,i}|}{|A_n|}\langle\phi,\mathbf B_{n,i}\rangle
+\sum_{i=k+1}^{C(n)}\frac{|A_{n,i}|}{|A_n|}\langle\phi,\mathbf B_{n,i}\rangle
+\sum_{i>C(n)}\frac{|B_{n,i}|}{|A_n|}\langle\phi,\mathbf B_{n,i}\rangle\\
&=\sum_{i=1}^q c_i\langle\phi,\mathbf A_{n,i}'\rangle
+\sum_{i=k+1}^{C(n)}\frac{|B_{n,i}|}{|A_n|}\langle\phi,\mathbf B_{n,i}\rangle
+\frac{|R_{n}|}{|A_n|}\langle\phi,\mathbf R_{n}\rangle.
\end{align*}
Thus we have
\begin{align*}
\biggl|\spect'\langle\phi,\mathbf R_{n}\rangle-\bigl(\langle\phi,\mathbf A_n\rangle-\sum_{i=1}^q c_i\mu_i\bigr)\biggr|
&\leq \sum_{i=1}^q|\langle\phi,\mathbf A_{n,i}'\rangle-\mu_i|
+\sum_{i=k+1}^{C(n)}|B_{n,i}|/|A_n|\\
&\quad+\bigl||R_{n}|/|A_n|-\spect'\bigr|\\
&\leq\epsilon.
\end{align*}

It follows that $\lim_{n\rightarrow\infty}\langle\phi,\mathbf R_{n}\rangle$ exists.
By sorting the $C(n)$ first connected components of each $\mathbf A_n$ according to both 
$\spect_i$ and Lemma~\ref{lem:refine} we obtain the following expression 
for the limit:
$$
\lim_{n\rightarrow\infty}\langle\phi,\mathbf R_{n}\rangle=\frac{1}{\spect'}
(\lim_{n\rightarrow\infty}\langle\phi,\mathbf A_{n}\rangle-\sum_{i<\widehat{C}}
\spect_i\lim_{n\rightarrow\infty}\langle\phi,\mathbf B_{n,i}\rangle).
$$
\myqed\end{proof}
Finally, we obtain the main results of this section.
\begin{mytinted}
\begin{theorem}[Comb structure for $\sig$-structure sequences with infinite spectrum (local convergence)]
\label{thm:comb}
Let $(\mathbf A_n)_{n\in\bbbn}$ be an ${\rm FO}^{\rm local}$-convergent sequence of
finite $\sig$-structures with component relation $\varpi$,
support $\bbbn$, and spectrum $(\spect_i)_{i\in\bbbn}$.
Let $C:\bbbn\rightarrow\bbbn$ be a clip of $(\mathbf A_n)_{n\in\bbbn}$, and let
 $\mathbf R_n$ be the residue of $\mathbf A_n$ with respect to $C$.

Then there exists, for each $n\in\bbbn$, 
a permutation  $f_n:[C(n)]\rightarrow [C(n)]$ such that, extending $f_n$ to $\bbbn$ by
putting $f(i)$ to be the identity for $i>C(n)$, the following holds:
\begin{itemize}
  \item $\lim_{n\rightarrow\infty} \max_{i>C(n)}|B_{n,i}|/|A_n|=0$;
  \item $\spect'=\lim_{n\rightarrow\infty} |R_n|/|A_n|$ exists;
  \item for every $i\in\bbbn$, $(B_{n,f_n(i)})_{n\in\bbbn}$ is ${\rm
  FO}$-convergent;
  \item either $\spect'=0$ or the sequence $(\mathbf R_n)_{n\in\bbbn}$ is ${\rm FO}^{\rm local}$-convergent.
\end{itemize}
\end{theorem}
\end{mytinted}
\begin{proof}
This lemma is a direct consequence of the previous lemmas.
\myqed\end{proof}
\begin{figure}[ht]
\begin{center}
\includegraphics[width=.7\textwidth]{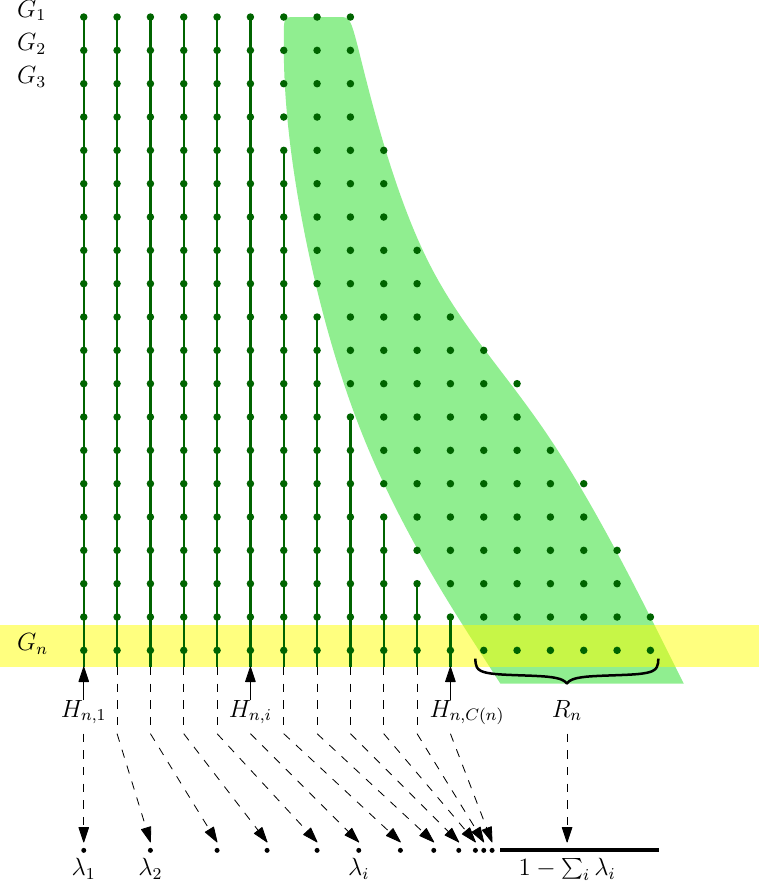}
\end{center}
\caption{Illustration of the Comb structure theorem}
\label{fig:comb}
\end{figure}

We shall now extend the Comb structure theorem to full ${\rm
FO}$-convergence.
In contrast with the case of a finite $\varpi$-spectrum, 
the elementary convergence aspects  will be non trivial and will require
a careful choice of a clip for the sequence.
\begin{lemma}
\label{lem:elemres}
Let $(\mathbf A_n)_{n\in\bbbn}$ be an ${\rm FO}^{\rm local}$-convergent sequence of
finite $\sig$-struc\-tures with component relation $\varpi$,
such that $\lim_{n\rightarrow\infty}|A_n|=\infty$.
Let $\mathbf B_{n,i}$ be the connected components of $\mathbf A_n$. 
Assume that the connected components  with same $\spect_i$ have been
reshuffled according to Lemma~\ref{lem:refine}, so
that $(\mathbf B_{n,i})_{i\in\bbbn}$ is ${\rm FO}$-convergent for each $i\in\bbbn$.

For $i\in\bbbn$, let $\widehat{\mathbf B_i}$ be 
an elementary limit of 
$(\mathbf B_{n,i})_{n\in\bbbn}$.
Then there exists a clip $C$ such that the sequence $(\mathbf R_n)_{n\in\bbbn}$ of the residues
is elementarily convergent. Moreover, if $\widehat{\mathbf R}$ is an elementary limit
of $(\mathbf R_n)_{n\in\bbbn}$, then $\bigcup_i\widehat{\mathbf B_i}\cup\widehat{\mathbf R}$ is an elementary
limit of $(\mathbf A_n)_{n\in\bbbn}$. 

Let $\mathbf B_{n,i}'$ be either $\mathbf B_{n,i}$ if $C(n)\geq i$
or the empty $\sig$-structure if $C(n)<i$.
Then the family 
consisting in the sequences $(\mathbf B_{n,i}')_{i\in\bbbn}$ ($i\in\bbbn$) and
of the sequence $(\mathbf R_n)_{n\in\bbbn}$ is uniformly elementarily convergent.
\end{lemma}
\begin{proof}
Let $\widehat{\mathbf A}$ be 
an elementary limit of 
 $(\mathbf A_n)_{n\in\bbbn}$.
 
For $\theta\in{\rm FO}_1^{\rm local}$ and $m\in\bbbn$ we 
denote by $\theta^{(m)}$ the sentence
$$
\theta^{(m)}:\quad (\exists x_1\dots\exists x_m)\ \biggl(\bigwedge_{1\leq i<j\leq m}
 \neg\varpi(x_i,x_j)\ \wedge\ \bigwedge_{i=1}^m \theta(x_i)\biggr).
$$

According to Theorem~\ref{thm:gaifman}, elementary convergence
of a sequence of $\sig$-struc\-tures with component relation $\varpi$
can be checked by considering sentences of the form $\theta^{(k)}$ for
$\theta\in{\rm FO}_1^{\rm local}$ and $k\in\bbbn$.

Note that for every $k<k'$ and every $\sig$-structure $\mathbf A$, if it holds that
$\mathbf A\models \theta^{(k')}$ then it holds that  $\mathbf A\models \theta^{(k)}$.
 Define
\begin{align*}
M(\theta)&=\sup\{k\in\bbbn:\quad \widehat{\mathbf A}\models\theta^{(k)}\}\\
\Omega(\theta)&=\{i\in\bbbn:\quad \widehat{\mathbf B_i}\models(\exists x)\theta(x)\}.
\end{align*}

Note that obviously $|\Omega(\theta)|\leq M(\theta)$.

For $r\in\bbbn$, let 
$\theta_1,\dots,\theta_{F(r)}$ be an enumeration of the local first-order formulas 
with a single free variable with quantifier rank at most $r$ (up to logical equivalence).
Define
$$
A(r)=\max(r,\max_{a\leq F(r)}\max\Omega(\theta_a)).
$$

Let  
$$C_0(n)=\sup \biggl\{K\in\bbbn:\quad (\forall n'\geq n)\ \sum_{i=1}^K
\biggl|\frac{|B_{n',i}|}{|A_{n'}|}-\spect_i\biggr|\leq
\sum_{i>K}\spect_i\biggr\}$$
be the standard (maximal) clip on $(\mathbf A_n)_{n\in\bbbn}$
(see Proposition~\ref{prop:clip}).

Let $B(r)$ be the minimum integer such that
\begin{enumerate}
  \item it holds that $C_0(B(r))\geq A(r)$ (note that
  $\lim_{n\rightarrow\infty}C_0(n)=\infty$, according to
  Proposition~\ref{prop:clip});
  \item for every $n\geq B(r), a\leq F(r)$ and every $k\leq r$ it holds that  
 $\mathbf A_n\models \theta_a^{(k)}$ if and only if $M(\theta_a)\geq k$
 (note that this holds for sufficiently large $n$ as $\widehat{\mathbf A}$ is
 an elementary limit of $(\mathbf A_n)_{n\in\bbbn}$);
 \item for every $i\leq A(r)$ and $a\leq F(r)$ the following equivalence holds:
 $$
 \mathbf B_{n,i}\models(\exists x)\theta_a(x)\quad\iff\quad \widehat{\mathbf B_i}\models(\exists x)\theta_a(x).
 $$
(note that this holds for sufficiently large $n$ as $\widehat{\mathbf B_i}$ is
 an elementary limit of $(\mathbf B_{n,i})_{n\in\bbbn}$ and as we consider
 only finitely many values of $i$);
 
\end{enumerate}

we define the non-decreasing function $C:\bbbn\rightarrow\bbbn$ by 
$$
C(n)=\max\{A(r): B(r)\leq n\}.
$$
As $\lim_{r\rightarrow\infty}A(r)=\infty$ and as $C_0(B(r))\geq A(r)$
it holds that $\lim_{r\rightarrow\infty}B(r)=\infty$.
Moreover, for every $r\in\bbbn$ it holds that $C_0(B(r))\geq A(r)$ hence
$C_0(n)\geq C(n)$. According to Proposition~\ref{prop:clip}, it follows
that the function $C$ is a clip on $(\mathbf A_n)_{n\in\bbbn}$.

Let $(\mathbf R_n)_{n\in\bbbn}$ be the residue of $(\mathbf A_n)_{n\in\bbbn}$
with respect to the clip $C$, and let $\mathbf B_{n,i}'$ be defined
as $\mathbf B_{n,i}$ if $i\leq C(n)$ and the empty $\sig$-structure otherwise. 
Then it is immediate from the definition of the clip $C$ that 
the family $\{(B_{n,i}')_{n\in\bbbn}:\ i\in\bbbn\}$ is uniformly elementarily
convergent.
Using Lemma~\ref{lem:scount}, it is also easily checked that 
 the residue $(\mathbf R_n)_{n\in\bbbn}$ of $(\mathbf A_n)_{n\in\bbbn}$
with respect to the clip $C$ is elementarily convergent and thus, that the
family $\{(B_{n,i}')_{n\in\bbbn}:\ i\in\bbbn\}\cup\{(\mathbf R_n)_{n\in\bbbn}\}$
is uniformly elementarily convergent.
 \myqed\end{proof}

The extension of the Comb structure theorem to ${\rm FO}$-convergence now follows directly. 
\begin{mytinted}
\begin{theorem}[Comb structure for $\sig$-structure sequences with infinite spectrum]
\label{thm:fullcombinf}
Let $(\mathbf A_n)_{n\in\bbbn}$ be an ${\rm FO}$-convergent sequence of
finite $\sig$-structures with component relation $\varpi$ and infinite spectrum
$(\spect_i)_{i\in\bbbn}$.

Then there exists
a clip $C:\bbbn\rightarrow\bbbn$ with residue
 $\mathbf R_n$ and, for each $n\in\bbbn$, 
a permutation  $f_n:[C(n)]\rightarrow [C(n)]$ such that, extending $f_n$ to $\bbbn$ by
putting $f(i)$ to be the identity for $i>C(n)$, and 
letting $\mathbf B_{n,i}'$ be either $\mathbf B_{n,f_n(i)}$ if $C(n)\geq i$
or the empty $\sig$-structure if $C(n)<i$, the following holds:
\begin{itemize}
  \item $\mathbf A_n=\mathbf R_n\cup\bigcup_{i\in\bbbn}\mathbf B_{n,i}'$;
  \item $\lim_{n\rightarrow\infty} \max_{i>C(n)}|B_{n,i}'|/|A_n|=0$;
  \item $\lim_{n\rightarrow\infty} |R_n|/|A_n|$ exists;
  \item for every $i\in\bbbn$, $(\mathbf B_{n,i}')_{n\in\bbbn}$ is ${\rm
  FO}$-convergent;
  \item 
  either $\lim_{n\rightarrow\infty} |R_n|/|A_n|=0$ and $(\mathbf R_n)_{n\in\bbbn}$ is elementarily convergent, 
  or the sequence $(\mathbf R_n)_{n\in\bbbn}$ is ${\rm FO}$-convergent;
   \item the family $\{(\mathbf B_{n,i}')_{n\in\bbbn}:
   i\in\bbbn\}\cup\{(\mathbf R_n)_{n\in\bbbn}\}$ is uniformly elementarily convergent.
\end{itemize}
\myqed\end{theorem}
\end{mytinted}
This ends the (admittedly very technical and complicated) analysis of the component structure of limits.
This was not developed for its own sake, but it will be all needed in the Part~3 of this paper, to construct
modeling ${\rm FO}$-limits for convergent sequences of trees with bounded height and, 
 by means of a fitting basic interpretation scheme, to graphs with bounded 
 tree-depth (defined in \cite{Taxi_tdepth}), or graphs
 with bounded SC-depth (defined in \cite{Ganian2012}).

In a broader sense, this detailed analysis was a cradle of much of the further research (see {\em Addendum} Section~\ref{sec:add}).
\chapter{Limits of Graphs with Bounded Tree-depth}%
In this part, we mainly consider 
	the signature $\sig$, which consists in a binary 
relation $\sim$ (symmetric adjacency relation),
a unary relation $R$ (property of being a root), 
and $c$ unary relations $C_i$ (the coloring).
Colored rooted trees with height at most $h$ are particular $\sig$-structures, and the
class of (finite or infinite) colored rooted trees with height at most $h$ will be denoted by
$\mathcal Y^{(h)}$. (Here we shall be only concerned with trees that are either finite, countable, or of size continuum.)

\section{${\rm FO}_1$-limits of Colored Rooted Trees with Bounded Height}
In this section, we explicitly define {\rps} ${\rm FO}_1$-limits for ${\rm FO}_1$-convergent sequences
of colored rooted trees with bounded height and characterize modelings which are ${\rm FO}_1$-limits for ${\rm FO}_1$-convergent sequences of (finite) colored rooted trees with bounded height.

\subsection{Preliminary Observations}
We take some time for some preliminary observations on the logical
structure of rooted colored trees with bounded height.
These observations will use arguments based on
Ehrenfeucht-Fra{\"\i}ss{\'e} games and strategy stealing. (For definitions of
$\equiv^n$ and Ehrenfeucht-Fra{\"\i}ss{\'e} games, see 
Section~\ref{sec:fop}.)

For a rooted colored tree $\mathbf Y\in\mathcal Y^{(h)}$ and a vertex $x\in Y$, we denote by 
$\mathbf Y(x)$ the {\em subtree of $\mathbf Y$ rooted at $x$} --- 
that is the subtree of $\mathbf Y$  with root $x$ induced by $x$ and all its descendants ---
and (for a non-root $x$) 
by $\mathbf Y\setminus\mathbf Y(x)$ the rooted
tree obtained from $\mathbf Y$ by removing all the vertices
in $\mathbf Y(x)$. 

The following two lemmas show that, like for isomorphism, equivalence between two colored rooted trees can be
reduced to equivalence of branches.

\begin{lemma}
\label{lem:joinr}
Let $\mathbf Y,\mathbf Y'\in\mathcal Y^{(h)}$, 
let $s,s'$ be sons of the roots of $\mathbf Y$ and $\mathbf Y'$, respectively.

Let $n\in\bbbn$.
If  
$\mathbf Y(s)\equiv^{n}\mathbf Y'(s')$ and 
$\mathbf Y\setminus\mathbf Y(s)\equiv^{n}\mathbf Y'\setminus\mathbf Y'(s')$,
then $\mathbf Y\equiv^{n}\mathbf Y'$.
\end{lemma}
\begin{proof}
Assume $\mathbf Y(s)\equiv^{n}\mathbf Y'(s')$ and 
$\mathbf Y\setminus\mathbf Y(s)\equiv^{n}\mathbf Y'\setminus\mathbf Y'(s')$.
In order to prove  $\mathbf Y\equiv^{n}\mathbf Y'$ we play
an $n$-step Ehrenfeucht-Fra{\"\i}ss{\'e}-game ${\rm EF}_0$ on $\mathbf Y$ and
$\mathbf Y'$ as Duplicator. Our strategy will be based on two
auxiliary $n$-step Ehrenfeucht-Fra{\"\i}ss{\'e}-games, ${\rm EF}_1$ and ${\rm
EF}_2$, on $\mathbf Y(s)$ and $\mathbf Y'(s')$ and on 
$\mathbf Y\setminus\mathbf Y(s)$ and $\mathbf Y'\setminus\mathbf Y'(s')$,
respectively, against Duplicators following a winning strategy.
 Each time
Spoiler selects a vertex in game ${\rm EF}_0$, we play the same vertex in the game ${\rm EF}_1$ or ${\rm EF}_2$ (depending on the
tree the vertex belongs to), then we mimic the selection of the Duplicator
of this game. It is easily checked that this strategy is a winning
strategy.
 \myqed\end{proof}
\begin{lemma}
\label{lem:splitr}
Let $\mathbf Y,\mathbf Y'\in\mathcal Y^{(h)}$, 
let $s,s'$ be sons of the roots of $\mathbf Y$ and $\mathbf Y'$, respectively.

Let $n\in\bbbn$.
If $\mathbf Y\equiv^{n+h}\mathbf Y'$ and 
$\mathbf Y(s)\equiv^{n}\mathbf Y'(s')$, then 
$\mathbf Y\setminus\mathbf Y(s)\equiv^{n}\mathbf Y'\setminus\mathbf Y'(s')$.
\end{lemma}
\begin{proof}
Assume $\mathbf Y\equiv^{n+h}\mathbf Y'$ and 
$\mathbf Y(s)\equiv^{n}\mathbf Y'(s')$.

We first play (as Spoiler) $s$ in $\mathbf Y$ then $s'$ in $\mathbf Y'$.
Let $t'$ and $t$ be the corresponding plays of
Duplicator.
Then the further $n$ steps of the game have to map vertices
in $\mathbf Y(s)$, $\mathbf Y(t)$, 
$\mathbf Y\setminus(\mathbf Y(s)\cup\mathbf Y(t))$ to 
$\mathbf Y'(t')$, $\mathbf Y'(s')$, 
$\mathbf Y'\setminus(\mathbf Y'(t')\cup\mathbf Y'(s'))$
(and converse), for otherwise $h-2$ steps would allow Spoiler to win the 
game. Also, by restricting our play to one of these pairs of
trees, we deduce
$\mathbf Y(s)\equiv^n \mathbf Y'(t')$,
$\mathbf Y(t)\equiv^n \mathbf Y'(s')$,
and $\mathbf Y\setminus(\mathbf Y(s)\cup\mathbf Y(t))\equiv^n
\mathbf Y\setminus(\mathbf Y'(s')\cup\mathbf Y'(t'))$.  
As $\mathbf Y'(s')\equiv^n\mathbf Y(s)$ it follows
$$\mathbf Y(t)\equiv^n \mathbf Y'(s')\equiv^n\mathbf Y(s)\equiv^n \mathbf
Y'(t').$$
Hence, according to Lemma~\ref{lem:joinr}, 
as $\mathbf Y\setminus(\mathbf Y(s)\cup\mathbf Y(t))=(\mathbf Y\setminus\mathbf
Y(s))\setminus\mathbf Y(t)$ and 
$\mathbf Y'\setminus(\mathbf Y'(s')\cup\mathbf Y'(t'))=(\mathbf
Y'\setminus\mathbf Y'(s'))\setminus\mathbf Y'(t')$,
we deduce
$\mathbf Y\setminus\mathbf Y(s)\equiv^n\mathbf Y'\setminus\mathbf Y'(s')$.
\myqed\end{proof}

Let $\sig^\bullet$ denote the signature obtained from $\sig$ by adding
a new unary relation $S$ (marking a {\em special} vertex, which is not necessarily the root).
Let $\theta_\bullet$ be the sentence
$$
(\exists x)(S(x)\wedge(\forall y\ S(y)\rightarrow (y=x))),
$$
which states that a $\sig^\bullet$ contains a unique special vertex.
We denote by $\mathcal Y^{(h)}_\bullet$ the class obtained
by marking as special a single vertex of a colored rooted tree with height at most $h$.
Let ${\rm Unmark}$ be the interpretation of $\sig$-structures in
$\sig^\bullet$-structures consisting in forgetting $S$ (so that
${\rm Unmark}$ projects $\mathcal Y^{(h)}_\bullet$ onto $\mathcal Y^{(h)}$).

\begin{lemma}
\label{lem:reducem}
Let $\mathbf Y,\mathbf Y'\in\mathcal Y^{(h)}_\bullet$ be such that $\mathbf Y$ (resp.\ $\mathbf Y'$) has
special vertex $m$ (resp.\ $m'$). Assume that both $m$ and $m'$
have height $t>1$ (in $\mathbf Y$ and $\mathbf Y'$, respectively).
Let $v$ (resp.\ $v'$) be son of the root of $\mathbf Y$ (resp.\ $\mathbf Y'$) that is an ancestor of $m$ (resp.\ $m'$).

Then for every $n\in\bbbn$, if ${\rm Unmark}(\mathbf Y)\equiv^{n+h}{\rm Unmark}(\mathbf Y')$ and
$\mathbf Y(v)\equiv^n\mathbf Y'(v')$, then $\mathbf Y\equiv^n\mathbf Y'$. 
\end{lemma}
\begin{proof}
Assume ${\rm Unmark}(\mathbf Y)\equiv^{n+h}{\rm Unmark}(\mathbf Y')$ and
$\mathbf Y(v)\equiv^n\mathbf Y'(v')$. Then it holds that 
${\rm Unmark}(\mathbf Y(v))\equiv^{n}{\rm Unmark}(\mathbf Y'(v'))$ thus,
according to Lemma~\ref{lem:splitr},
\begin{align*}
\mathbf Y\setminus \mathbf Y(v)&={\rm Unmark}(\mathbf Y)\setminus{\rm
Unmark}(\mathbf Y(v))\\
&\equiv^n {\rm Unmark}(\mathbf Y')\setminus{\rm
Unmark}(\mathbf Y'(v'))=\mathbf Y'\setminus \mathbf Y'(v').
\end{align*}
Hence, according to Lemma~\ref{lem:joinr}, it holds that 
$\mathbf Y\equiv^n\mathbf Y'$ (as the marking could be considered as a coloring).
\myqed\end{proof}

\begin{figure}[ht]
\begin{center}
\includegraphics[width=\textwidth]{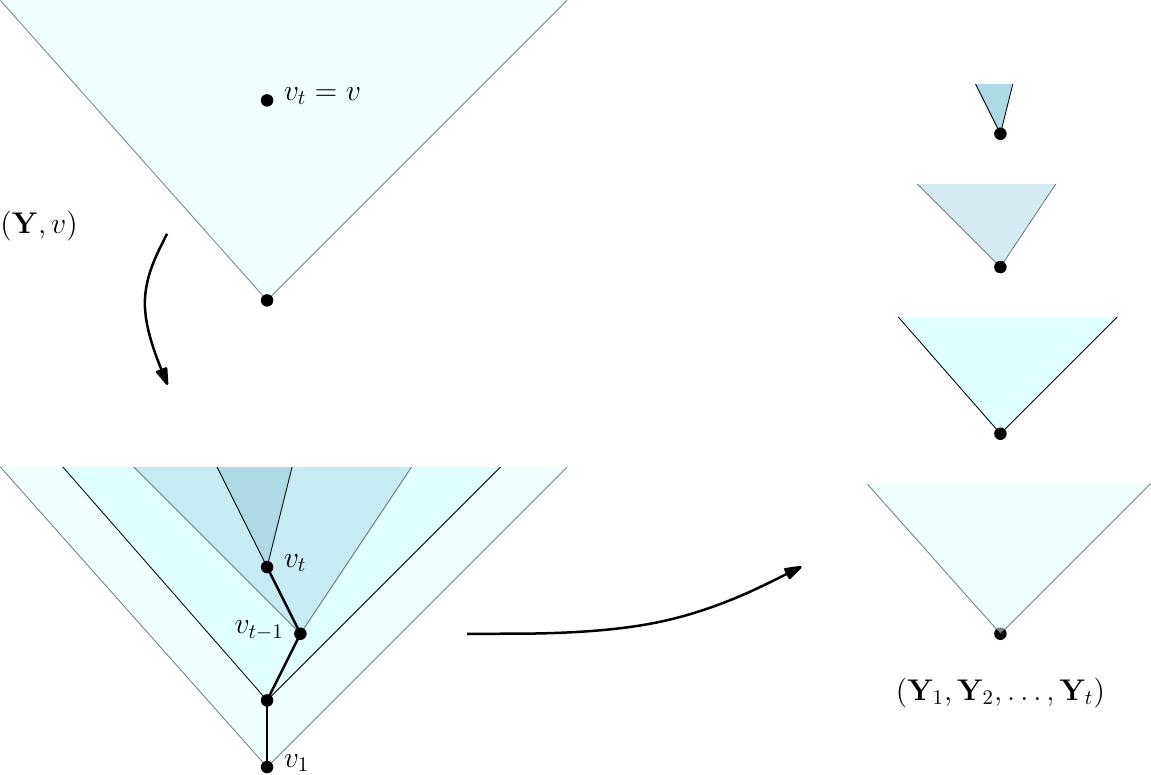}
\end{center}
\caption{Transformation of a rooted tree with a distinguished vertex $(\mathbf Y,v_t)$ into
a tuple of rooted trees $(\mathbf Y_1,\dots,\mathbf Y_t)$.}
\label{fig:1to0}
\end{figure}

The next lemma states that the properties of a colored rooted
trees with a distinguished vertex $v$ (which is not necessarily the root) can be retrieved from the properties
of the subtree rooted at $v$, the subtree rooted at the father of $v$, etc.\   
(see Fig.~\ref{fig:1to0}). 

\begin{lemma}
\label{lem:0to1}
Let $\mathbf Y,\mathbf Y'\in \mathcal Y^{(h)}$,
$v_t\in Y$ and $v_t'\in Y'$ be vertices with height $t$. 
For $1\leq i<t$, let $v_i$ (resp.\ $v_i'$) be the ancestor of $v_t$ (resp.\ of $v_t'$) at height $i$.
 
Then for every integer $n$ it holds that 
\begin{align*}
(\forall 1\leq i\leq t)\ \mathbf Y(v_i)\equiv^{n+h+1-i} \mathbf
Y'(v_i')\ &\Longrightarrow\ (\mathbf Y,v_t)\equiv^n (\mathbf Y',v_t')\\
(\mathbf Y,v_t)\equiv^{n+(t-1)h} (\mathbf
Y',v_t')\ &\Longrightarrow\ (\forall 1\leq i\leq t)\ 
\mathbf Y(v_i)\equiv^{n+(t-i)h} \mathbf Y'(v_i')
\end{align*}
\end{lemma}
\begin{proof}
We proceed by induction over $t$. If $t=1$, then the statement obviously holds.
So, assume $t>1$ and that the statement holds for $t-1$.

Let $\mathbf Y_\bullet,\mathbf Y_\bullet'\in\mathcal Y^{(h)}_\bullet$  be the marked rooted colored
trees obtained from $\mathbf Y$ and $\mathbf Y'$
by marking $v_t$ (resp.\ $v_t'$) as a  special vertex.

Assume $(\forall 1\leq i\leq t)\ \mathbf Y(v_i)\equiv^{n+h+1-i} \mathbf
Y'(v_i')$. By induction, 
$(\forall 2\leq i\leq t)\ \mathbf Y(v_i)\equiv^{n+(h-1)+1-(i-1)}
\mathbf Y'(v_i')$ implies $(\mathbf Y(v_2),v_t)\equiv^n (\mathbf
Y'(v_2'),v_t')$, that is \linebreak $\mathbf Y_\bullet(v_2)\equiv^n \mathbf Y_\bullet'(v_2')$.
As $\mathbf Y\equiv^{n+h}\mathbf Y'$, it follows from Lemma~\ref{lem:reducem}
that $\mathbf Y_\bullet\equiv^n\mathbf Y_\bullet'$, that is: $(\mathbf Y,v_t)\equiv^n
(\mathbf Y',v_t')$.

Conversely, if $(\mathbf Y,v_t)\equiv^{n+(t-1)h} (\mathbf Y',v_t')$
(i.e. $\mathbf Y_\bullet\equiv^{n+(t-1)h} \mathbf Y_\bullet'$)
an repeated application of Lemma~\ref{lem:splitr} gives
$\mathbf Y_\bullet(v_i)\equiv^{n+(t-i)h}\mathbf Y_\bullet'(v_i)'$ hence
(by forgetting the marking)
$\mathbf Y(v_i)\equiv^{n+(t-i)h}\mathbf Y'(v_i)'$.
\myqed\end{proof}

This lemma allows to encode the complete theory of a colored rooted
tree $\mathbf Y$ of height at most $h$ with special vertex $v$ as
a tuple of complete theories of colored rooted trees with height at most $h$.

In fact what follows could be described as a fine analysis of the Stone dual of the Boolean algebra of all the formulas having a model in $\mathcal Y^{(h)}$.
As the height $h$ is bounded, the classes $\mathcal Y^{(h)}$ can be axiomatized by finitely many axioms (hence by some single axiom $\eta_{\mathcal Y^{(h)}}$), it
is a basic elementary class.
For an integer $p\geq 0$, we introduce a short notation for the Stone
space associated to the Lindenbaum--Tarski algebra of
formulas on  $\mathcal Y^{(h)}$ with
$p$ free variables:
$$\mathfrak Y_p^{(h)}=S(\mathcal B({\rm FO}_p(\sig),\eta_{\mathcal Y^{(h)}})).$$

We shall now move from models to theories, specifically from $\mathcal Y_\bullet^{(h)}$ (colored rooted trees with height at most $h$ and a special vertex)
to the Stone space $\mathfrak Y_1^{(h)}$ and from $\mathcal Y^{(h)}$ (colored rooted trees with height at most $h$) to 
the Stone space $\mathfrak Y_0^{(h)}$.

In that direction, we first show how the notion of ``property of the subtree $\mathbf Y(v)$ of $\mathbf Y$ rooted at the vertex $v$''
translates into a relativization homomorphism $\varrho:\mathcal B({\rm FO}_0(\sig))\rightarrow\mathcal B({\rm FO}_1(\sig))$.

We consider the simple interpretation ${\mathsf I}_\bullet$ of $\sig$-structures
in $\sig^\bullet$-structures, which maps a $\sig^\bullet$-structure $\mathbf Y_\bullet$
to the $\sig$-structure defined as follows: let $x\simeq y$ be defined as $(x\sim y)\vee (x=y)$.
Then 
\begin{itemize}
  \item the domain of ${\mathsf I}_\bullet(\mathbf Y_\bullet)$ is defined by the formula
  $$S(x_1)\vee(\forall y_1,\dots,y_h) \Bigl(\bigl(R(y_1)\wedge \bigwedge_{i=1}^{h-1} \neg S(y_i)\wedge (y_i\simeq y_{i+1})\bigr)\rightarrow (y_h\neq x_1)\Bigr);
  $$
	\item the adjacency relation $\sim$ is defined as in $\mathbf Y_\bullet$ (i.e. by the formula $(x_1\sim x_2)$);
  \item the relation $R$ of ${\mathsf I}_\bullet(\mathbf Y_\bullet)$ is defined by the formula $S(x_1)$.
\end{itemize}
Although ${\mathsf I}_\bullet$ maps general $\sig_\bullet$-structures to $\sig$-structures, we shall be only concerned
by the specific property that ${\mathsf I}_\bullet$ maps a rooted tree $\mathbf Y_\bullet\in{\mathcal Y}^{(h)}_\bullet$ with special vertex
$v$ to the rooted tree ${\rm Unmark}(\mathbf Y_\bullet)(v)$.

In a sake for simplicity, for $\mathbf Y\in{\mathcal Y}^{(h)}$ we denote by $(\mathbf Y,v)$ (where $\mathbf Y$ is 
a $\sig$-structure) the $\sig^\bullet$-structure obtained by adding the new relation $S$ with
$v$ being the unique special vertex. 

\begin{lemma}
There is a Boolean algebra homomorphism  
$$\varrho:\mathcal B({\rm FO}_0(\sig),\eta_{\mathcal Y^{(h)}})\rightarrow\mathcal B({\rm FO}_1(\sig),\eta_{\mathcal Y^{(h)}})$$
 (called {\em relativization}), such that
for every sentence $\phi\in{\rm FO}_0(\sig)$, every $\mathbf Y\in\mathcal Y^{(h)}$, and every $v\in Y$ the following equivalence holds
$$\mathbf Y(u)\models \phi\quad\iff\quad \mathbf Y\models \varrho(\phi)(u).$$
\end{lemma}
\begin{proof}
The lemma
follows from the property
$$\mathbf Y(u)\models \phi\quad\iff\quad {\mathsf I}_\bullet(\mathbf Y,u)\models\phi \quad\iff\quad (\mathbf Y,u)\models {\tilde{\mathsf I}}_\bullet(\phi).$$
The formula $\rho(\phi)$ is obtained from the sentence ${\tilde{\mathsf I}}_\bullet(\phi)$ by replacing 
each occurrence of $S(y)$ by $y=x_1$.
\myqed\end{proof}

Using relativization and Lemma~\ref{lem:0to1}, we can translate the transformation shown on Figure~\ref{fig:1to0}
to a encoding of elements of $\mathfrak Y_1^{(h)}$ into tuples of elements $\mathfrak Y_0^{(h)}$. Intuitively, a element 
$T\in \mathfrak Y_1^{(h)}$ defines the properties of a colored rooted tree $\mathbf Y$ with special vertex $x_1$, and
the relativization $\rho$ allows us to extract from $T$ the tuple of the complete theories of the subtrees of $\mathbf Y$ rooted at 
 $x_1$, the father of $x_1$, etc.\  Moreover, the meaning of Lemma~\ref{lem:0to1} is that what we
obtain only depends on the complete theory of $(\mathbf Y,x_1)$, that is only on $T$. 

\begin{definition}
For $1\leq i\leq h$, let $\eta_i\in{\rm FO}_1(\sig)$ be the formula 
stating that the height of $x_1$ is $i$.

We define the mapping ${\rm Encode}: \mathfrak Y_1^{(h)}\rightarrow \biguplus_{k=1}^h (\mathfrak Y_0^{(h)})^k$
as follows:

For $T\in \mathfrak Y_1^{(h)}$, let $k$ be the (unique) integer such that $\eta_k\in T$. 
Then ${\rm Encode}(T)$ is the $k$-tuple $(T_0,\dots,T_{k-1})$, where
\begin{itemize}
	\item $T_{k-1}$ is the set of sentences $\theta\in{\rm FO}_0(\sig)$ such that $\rho(\theta)\in T$ (intuitively, the complete theory
	of the subtree rooted at $x_1$);
  \item $T_{k-2}$ is the set of sentences $\theta\in{\rm FO}_0(\sig)$ such that
	$$\bigl((\exists y_1)(\eta_{k-1}(y_1)\wedge y_1\sim x_1\wedge\rho(\theta)(y))\bigr)\in T$$
	(intuitively, the complete theory	of the subtree rooted at the father of $x_1$);
	\item $T_{k-1-i}$ is the set of sentences $\theta\in{\rm FO}_0(\sig)$ such that
	$$\bigl((\exists y_1\dots y_i)(\bigwedge_{j=1}^i\eta_{k-j}(y_j)\wedge\bigwedge_{j=1}^{i-1}(y_j\sim y_{j+1}\wedge y_1\sim x_1\wedge\rho(\theta)(y_i)\bigr)\in T$$
(intuitively, the complete theory	of the subtree rooted at the ancestor of $x_1$ which has height $k-i$);
	\item $T_{0}=T\cap \in{\rm FO}_0(\sig)$ (intuitively, the complete theory of the whole rooted tree).
\end{itemize}
\end{definition}

\begin{lemma}
\label{lem:encode}
${\rm Encode}$ is a homeomorphism of $\mathfrak Y_1^{(h)}$ and
${\rm Encode}(\mathfrak Y_1^{(h)})$, which is
a closed subspace of $\biguplus_{k=1}^h(\mathfrak Y_0^{(h)})^k$.
\end{lemma}
\begin{proof}
This lemma is a direct consequence of Lemma~\ref{lem:0to1}.
\end{proof}

\subsection{The Universal Relational Sample Space $\mathbb Y_h$}
The aim of this section is to construct a rooted colored forest on
a standard Borel space $\mathbb Y_h$ that is {\em ${\rm FO}_1$-universal}, in the
sense that every ${\rm FO}_1$-convergent sequence of colored rooted trees will have a
{\rps} ${\rm FO}_1$-limit obtained by assigning an adapted probability measure 
to one of the connected components of $\mathbb Y_h$.

\begin{definition} 
For theories $T,T'\in\mathfrak Y_0^{(h)}$, we define
$w(T,T')\geq k$ if and only if there exists a model $\mathbf Y$ of $T$, such that the root of $\mathbf Y$
 has $k$ (distinct) sons
$v_1,\dots,v_k$ with ${\rm Th}(\mathbf Y(v_i))=T'$. 
\end{definition}
\begin{lemma}
\label{lem:w}
For $k\in\bbbn$ and $\phi\in{\rm FO}_0$, let $\zeta_k(\phi)$ be the sentence 
$(\exists^{\geq k}y)\ \rho(\phi)(y)$.
Then $w(T,T')\geq k$ if and only if $\zeta_k(\phi)\in T$ holds for every $\phi\in T'$.
\end{lemma}
\begin{proof}
If $w(T,T')\geq k$, then $\zeta_k(\phi)\in T$ holds for every $\phi\in T'$, hence we only have to
prove the opposite direction.
Assume that $\zeta_k(\phi)\in T$ holds for every $\phi\in T'$, but that there is $\phi_0\in T'$ such that
$\zeta_{k+1}(\phi_0)\notin T$. Let $Y$ be a model of $T$, and let $v_1,\dots,v_k$ be the sons of the root of $Y$ such that
$Y(v_i)\models\phi_0$. For every $r\in\bbbn$, $r\geq{\rm qrank}(\phi_0)$, let $\psi_r$ be the conjunction of the
sentences in $T'$ with quantifier rank $r$. Obviously, $\psi_r\in T'$. Moreover, as $\zeta_k(\psi_r)\in T$, it holds that $Y\models\zeta_k(\psi_r)$. As $\psi_r\vdash \phi_0$, it follows that for every $1\leq i\leq k$ it holds that $Y(v_i)\models\psi_r$ (only possible choices). As this holds for every $r$, we infer that for every $1\leq i\leq k$, $Y(v_i)$ is a model of $T'$ hence
$w(T,T')\geq k$. Now assume that for every $k\in\bbbn$ and every $\phi\in T'$ it holds that 
$\zeta_k(\phi)\in T$. Let $Y$ be a model of $T$, let $Y'$ be a model of $T'$, and let $\tilde{Y}$ be obtained from $Y$ by adding 
(at the root of $Y$) a son $u$ with subtree $\tilde{Y}(u)$ isomorphic to $Y'$. By an easy application of an 
Ehrenfeucht-Fra{\"\i}ss\'e game, we get that $\tilde{Y}$ is elementarily equivalent to $Y$, hence a model of $T$. Thus $w(T,T')\geq k$.
\end{proof}
Let $\overline{\bbbn}$ be the one point compactification of $\bbbn$, that is 
$\overline{\bbbn}=\bbbn\cup\{\infty\}$ with open sets generated by complements of finite sets.
\begin{lemma}
The function $w:\mathfrak Y_0^{(h)}\times\mathfrak Y_0^{(h)}\rightarrow\overline{\bbbn}$ is upper semicontinuous (with respect to product topology
of Stone space $\mathfrak Y_0^{(h)}$).
\end{lemma}
\begin{proof}
For $r\in\bbbn$ define the function $w_r:\mathfrak Y_0^{(h)}\times\mathfrak Y_0^{(h)}\rightarrow\overline{\bbbn}$ by:
$$w_r(T,T')=\sup\{k\in\bbbn:\ \forall \psi\in T'\ ({\rm qrank}(\psi)\leq r)\Rightarrow\zeta_k(\psi)\in T\}.$$
It follows from Lemma~\ref{lem:w} that the following equation holds
$$w(T,T')=\inf_{r\in\bbbn}w_r(T,T').$$
Hence, in order to prove that 
the function $w$ is upper semicontinuous, it is sufficient to 
prove that the functions $w_r$ are continuous. Let $(T_0,T_0')\in\mathfrak Y_0^{(h)}\times\mathfrak Y_0^{(h)}$.
We distinguish two cases:
	
	\noindent -- Firstly, assume $w_r(T_0,T_0')=k$. 
	
	\noindent If ${\rm dist}(T',T_0')<2^{-r}$ and ${\rm dist}(T,T_0)<2^{-\max\{{\rm qrank}(\zeta_{k+1}(\psi)):\ {\rm qrank}(\psi)\leq r\}}$,
	then it holds that $w_r(T,T')=w_r(T_0,T_0')$;
	
	\noindent -- Secondly, assume $w_r(T_0,T_0')=\infty$, and let $k\in\bbbn$.
	
\noindent If ${\rm dist}(T',T_0')<2^{-r}$ and ${\rm dist}(T,T_0)<2^{-\max\{{\rm qrank}(\zeta_{k+1}(\psi)):\ {\rm qrank}(\psi)\leq r\}}$, then it holds that  $w_r(T,T')>k$.
\end{proof}

For $z=(z_1,\dots,z_a)\in\overline{\bbbn}^a$ define the subset $F_z$ of
$(\mathfrak Y_0^{(h)})^{a+1}$ by
$$F_z=\{(T_0,\dots,T_{a}): w(T_{i-1},T_{i})=z_i\}.$$

For $t\in\overline{\bbbn}$, define
$$X_t=\begin{cases}
\{1,\dots,t\},&\text{if }t\in\bbbn,\\
[0,1],&\text{if }t=\infty.
\end{cases}
$$

For $z=(z_1,\dots,z_a)\in\overline{\bbbn}^a$, define
$X_z=\prod_{i=1}^a X_{z_i}$. 
Let
$$V_h=\mathfrak Y_0^{(h)}\uplus\biguplus_{a=1}^{h-1}\biguplus_{z\in\overline{\bbbn}^a} (F_z\times X_z).$$

\begin{definition}
The {\em universal forest} $\mathbb Y_{h}$ has vertex set $V_h$.
The roots of $\mathbb Y_{h}$ are the elements in $\mathfrak Y_0^{(h)}$.
The edges of $\mathbb Y_{h}$ are the pairs of the form
$$\{((T_0,T_1,\dots,T_a),(\alpha_1,\dots,\alpha_a)),((T_0,T_1,\dots,T_{a+1}),(\alpha_1,\dots,\alpha_{a+1}))\}$$
where $T_i\in\mathfrak Y_0^{(h)}$, $\alpha_i\in[0,1]$ and $a\in\{0,\dots,h-1\}$.

Moreover, the vertex set $V_h$ inherits the topological structure of
$\biguplus_{i=1}^h(\mathfrak Y_0^{(h)})^i\times[0,1]^{i-1}$, which defines
a $\sigma$-algebra $\Sigma_h$ on $V_h$ (as the trace on $V_h$ of the Borel $\sigma$-algebra
of  $\biguplus_{i=1}^h(\mathfrak Y_0^{(h)})^i\times[0,1]^{i-1}$).
\end{definition}

\begin{remark}
Let $T_0$ be the complete complete theory of a colored rooted tree with height at most $h$.
Then, by construction, $T_0$ is the complete theory of the connected component of $\mathbb Y_h$ rooted at $T_0$.
In particular, no two connected components of $\mathbb Y_h$ are elementarily equivalent.
\end{remark}

The remaining of this section will be devoted to the proof
of Theorem~\ref{thm:yrss}, which states
that $\mathbb Y_h$ is a \rss. In order to prove this result,
we shall need a preliminary lemma, which  expresses that the property of a tuple
of vertices in a colored rooted tree with bounded height is 
completely determined by the individual properties of the vertices 
in the tuple and the heights of the lowest common ancestors of 
every pair of vertices in the tuples.

\begin{lemma}
\label{lem:pto1}
Fix rooted trees $\mathbf Y,\mathbf Y'\in\mathcal Y^{(h)}$. 
Let $u_1,\dots,u_p$ be $p$ vertices of $\mathbf Y$, let $u_1',\dots,u_p'$
be $p$ vertices of $\mathbf Y'$, and let $n\in\bbbn$.

Assume that for every $1\leq i\leq p$ it holds that  
$(\mathbf Y,u_i)\equiv^{n+h} (\mathbf Y',u_i')$ and
that for every $1\leq i,j\leq p$ the height of $u_i\wedge u_j$ in $\mathbf Y$ is the same
as the height of $u_i'\wedge u_j'$ in $\mathbf Y'$
(where $u\wedge v$ denotes the lowest common ancestor of $u$ and $v$).

Then $(\mathbf Y,u_i,\dots,u_p)\equiv^n (\mathbf Y',u_1',\dots,u_p')$.
\end{lemma}
\begin{proof}
In the proof we consider $p+1$ simultaneous
Ehrenfeucht-Fra{\"\i}ss{\'e} games (see Fig.~\ref{fig:mgame}).
\begin{figure}[htp]
\begin{center}
  \includegraphics[width=.6\textwidth]{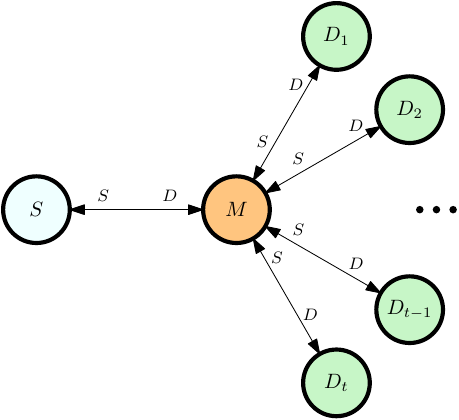}
  \caption{Schematic representation of how  a winning strategy for ${\rm EF}((Y,u_1,\dots,u_p),(Y',u_1',\dots,u_p'),n)$ is built using
  $p$ auxiliary games ${\rm EF}((Y,u_i),(Y',u_i'),n+h)$.}
  \label{fig:mgame}
\end{center}
\end{figure}

Consider an $n$-step Ehrenfeucht-Fra{\"\i}ss{\'e} ${\rm EF}((\mathbf Y,u_1,\dots,u_p),(\mathbf Y',u_1',\dots,u_p'),n)$ 
 on $(\mathbf Y,u_1,\dots,u_p)$ and 
$(\mathbf Y',u_1',\dots,u_p')$.
We build a strategy for Duplicator by considering $p$ auxiliary Ehrenfeucht-Fra{\"\i}ss{\'e} games
 ${\rm EF}((\mathbf Y,u_i),(\mathbf Y',u_i'),n+h)$
on $(\mathbf Y,u_i)$ and $(\mathbf Y',u_i')$ (for $1\leq i\leq p$) where we play the role of Spoiler
against Duplicators having a winning strategy for $(n+h)$-step games.

For every vertex $v\in Y$ (resp.\ $v'\in Y'$) let $p(v)$ (resp.\ $p'(v)$) be the maximum 
ancestor of $v$ (in the sense of the furthest from the root) such that $p(v)\leq u_i$ 
(resp.\ $p'(v)\leq u_i'$) for some $1\leq i\leq p$. 
We partition $Y$ and $Y'$ as follows: for every vertex $v\in Y$ (resp.\ $v'\in Y'$)
we put $v\in V_i$ (resp.\ $v'\in V_i'$) if $i$ is the minimum integer such that
$p(v)\leq u_i$ (resp.\ such that $p'(v)\leq u_i'$), see Fig~\ref{fig:multigame}.

\begin{figure}[htp]
\begin{center}
  \includegraphics[width=.7\textwidth]{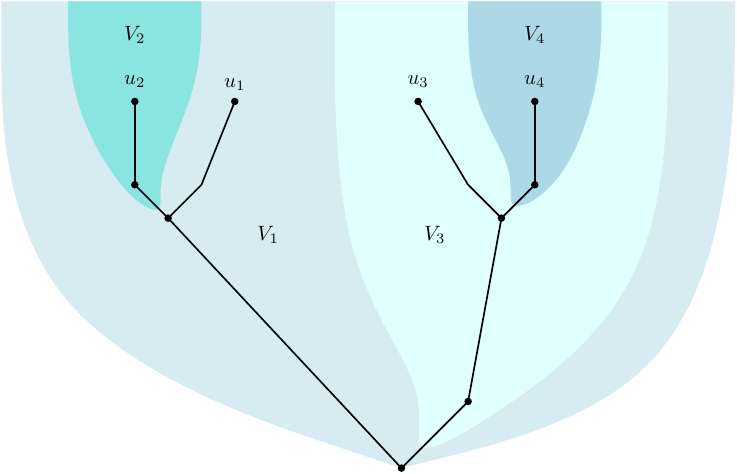}
  \caption{The partition $(V_1,V_2,V_3,V_4)$ of $Y$ induced by $(u_1,u_2,u_3,u_4)$.}
  \label{fig:multigame}
\end{center}
\end{figure}
 
Note that each $V_i$ (resp.\ $V_i'$ induces a connected subgraph
of $\mathbf Y$ (resp.\ of $\mathbf Y'$).

Assume that at round $j\leq n$, Spoiler plays a vertex $v\in (\mathbf Y,u_1,\dots,u_p)$ 
(resp.\ a vertex $v'\in (\mathbf Y',u_1',\dots,u_p')$).

If $v\in V_i$ (resp.\ $v'\in V_i'$) for some $1\leq i\leq p$ then we play $v$ (resp.\ $v'$) on
$(\mathbf Y,u_i)$ (resp.\ $(\mathbf Y',u_i')$).
We play Duplicator on $(\mathbf Y',u_1',\dots,u_p')$ (resp.\ on $(\mathbf Y,u_1,\dots,u_p)$) with the same move as our Duplicator
opponent did on $(\mathbf Y',u_i)$ (resp.\ on $(\mathbf Y,u_i)$).
If such a move is not legal (i.e. does not define a local isomorphism)
then it is easily checked that $h$ additional moves (at most)
are sufficient for at least one of the Spoilers to win one of the $p$ games, contradicting the hypothesis
 of $p$ winning strategies for Duplicators. 
 It follows that $(\mathbf Y,u_i,\dots,u_p)\equiv^n (\mathbf Y',u_1',\dots,u_p')$.
\myqed\end{proof}

\begin{mytinted}
\begin{theorem}
\label{thm:yrss}
The rooted colored forest $\mathbb Y_h$ (equipped with the 
$\sigma$-algebra $\Sigma_h$) is a \rss.
\end{theorem} 
\end{mytinted}
\begin{proof}
It suffices to prove that for every $p\in\bbbn$ and every $\varphi\in{\rm FO}_p$ the set 
$$
\Omega_\varphi(\mathbb Y_h)=\{(v_1,\dots,v_p)\in V_h^p:\ \mathbb
Y_h\models\varphi(v_1,\dots,v_p)\}
$$
is measurable.

Let $\varphi\in{\rm FO}_{p}$ and let $n={\rm qrank}(\varphi)$. 

We partition $V_h$ into equivalence classes modulo
$\equiv^{n+h}$, which we denote \linebreak $C_1,\dots,C_N$.

Let $i_1,\dots,i_p\in[N]$ and, for $1\leq j\leq p$, let
$v_{j}$ and $v_{j}'$ belong to $C_{i_j}$.

According to
Lemma~\ref{lem:pto1}, if the heights of the lowest common ancestors
of the pairs in $(v_1,\dots,v_p)$ coincide with
the heights of the lowest common ancestors
of the pairs in $(v_1',\dots,v_p')$ 
then it holds that 
$$(\mathbb Y_h,v_1,\dots,v_p)\equiv^n(\mathbb Y_h,v_1',\dots,v_p')$$
thus $(v_1,\dots,v_p)\in \Omega_\varphi(\mathbb Y_h)$ if and only if
$(v_1',\dots,v_p')\in\Omega_\varphi(\mathbb Y_h)$.

It follows from Lemma~\ref{lem:encode} (and the definition of $\mathbb V_h$ and $\Sigma_h$) that each $C_j$ is measurable.
According to Lemma~\ref{lem:encode} and the encoding of the vertices of $V_h$, the conditions
on the heights of lowest common ancestors rewrite as equalities and inequalities of coordinates.
It follows that $\Omega_\varphi(\mathbb Y_h)$ is measurable.    
\myqed\end{proof}
\subsection{Modeling ${\rm FO}_1$-limits of Colored Rooted Trees with Bounded Height}
Let $(\mathbf Y_n)_{n\in\bbbn}$ be an ${\rm FO}_1$-convergent of colored rooted trees with height at most $h$, and let
$\widetilde{\mathbf Y}$ be the connected component of $\mathbb Y_h$ that is an elementary limit of 
$(\mathbf Y_n)_{n\in\bbbn}$. According to Lemma~\ref{lem:cc}, $\widetilde{\mathbf Y}$ is a {\rss}.
We have to transfer the measure $\mu$ we obtained in
Theorem~\ref{thm:genmodxc} on $S(\mathcal B({\rm FO_1}))$ to $\widetilde{\mathbf Y}$.

\begin{definition}
\label{def:nu}
Let $\mu$ be a measure on $\mathfrak Y_1^{(h)}$. We define $\nu$ on $\mathbb Y_{h}$ as follows:
let $\widetilde{\mu}={\rm Encode}_*(\mu)$ be the pushforward of $\mu$ by ${\rm Encode}$ (see page~\pageref{def:pushforward}).
For $t\in\overline{\bbbn}$ we equip $X_t$ with uniform discrete probability
measure if $t<\infty$ and the Haar probability measure if $t=\infty$. For
$z\in\overline{\bbbn}^a$, $X_z$ is equipped with the corresponding product measure, which we denote by $\leb_z$ (not to be confused with
signature $\sig$).

We define the measure $\nu$ as follows: let $A$ be a measurable subset of 
$V_h$, let $A_0=A\cap\mathfrak Y_0^{(h)}$, and let $A_z=A\cap (F_z\times X_z)$.
Then
$$
\nu(A)=\widetilde{\mu}(A_0)+\sum_{a=1}^{h-1}\sum_{z\in\overline{\bbbn}^a} (\widetilde{\mu}\otimes\leb_z)(A_z).
$$
(Notice that the sets $A_z$ are measurable as $F_z\times X_z$ is measurable for every $z$.) 
\end{definition}

\begin{lemma}
\label{lem:push1}
The measure $\mu$ is the push-forward of $\nu$ by the projection  
$P:\mathbb Y_{h}\rightarrow\mathfrak Y_1^{(h)}$ defined by
$$
P((T_0,T_1,\alpha_1,\dots,T_a,\alpha_a))={\rm Encode}^{-1}(T_0,\dots,T_a),
$$
that is: $\mu=P_*(\nu)$.
\end{lemma}
\begin{proof}
First notice that $P$ is continuous, as {\rm Encode} is a homeomorphism (by Lemma~\ref{lem:encode}).
Let $B$ be a measurable set of $\mathfrak Y_1^{(h)}$.
Let $A=P^{-1}(B)$. Then $A\cap (F_z\times X_z)=({\rm Encode}(B)\cap F_z)\times X_z$ hence
$$(\widetilde{\mu}\otimes\leb_z)(A\cap (F_z\times X_z))=\nu({\rm Encode}(B)\cap F_z)\leb_z(X_z)=\widetilde{\mu}({\rm Encode}(B)\cap F_z).$$ 
It follows that
\begin{align*}
P_*(\nu)(B)&=\nu(A)\\
&=\widetilde{\mu}(A\cap\mathfrak Y_0^{(h)})+\sum_{a=1}^{h-1}\sum_{z\in\overline{\bbbn}^a} (\widetilde{\mu}\otimes\leb_z)(A\cap (F_z\times X_z))\\
&=\widetilde{\mu}({\rm Encode}(B)\cap\mathfrak Y_0^{(h)})+\sum_{a=1}^{h-1}\sum_{z\in\overline{\bbbn}^a} \widetilde{\mu}({\rm Encode}(B)\cap F_z)\\
&=\widetilde{\mu}\bigl({\rm Encode}(B)\cap(\mathfrak Y_0^{(h)}\uplus\biguplus_{a=1}^{h-1}\biguplus_{z\in\overline{\bbbn}^a} F_z)\bigr)\\
&=\widetilde{\mu}\circ{\rm Encode}(B)\\
&=\mu(B).
\end{align*}
(as $z$ ranges over a countable set and as all the $F_z$ are measurable). 
Hence $\mu=P_*(\nu)$.
\myqed\end{proof}

\begin{lemma}
\label{lem:FO1muY}
Let $\mu$ be a pure measure on $\mathfrak Y_1^{(h)}$ and
let $T_0$ be the complete theory of $\mu$ (see Definition~\ref{def:pure}).
Let $\nu$ be the measure defined from $\mu$ by Definition~\ref{def:nu}.
Let $\widetilde{\mathbf Y}$ be the connected component of $\mathbb Y_{h}$ containing the support of $\nu$.
Let $\nu_{\widetilde{\mathbf Y}}$ be the restriction of $\nu$ to $\widetilde{\mathbf Y}$.

Then $\widetilde{\mathbf Y}$, equipped with the probability
 measure $\nu_{\widetilde{\mathbf Y}}$ is a {\rps} such that for every $\varphi\in{\rm FO}_1$  the following equality holds
$$
\langle\varphi,\widetilde{\mathbf Y}\rangle=\mu(K(\varphi)).
$$
Let $X\subset \mathfrak Y_1^{(h)}$ be the set of all $T\in\mathfrak Y_1^{(h)}$ such that $x_1$ is not the root
(i.e. $X=\{T\in \mathfrak Y_1^{(h)}: R(x_1)\notin T\}$).
Let $f: X\rightarrow \mathfrak Y_0^{(h)}$ be the second projection of ${\rm Encode}$
(if ${\rm Encode}(T)=(T_0,\dots,T_i)$ then $f(T)=T_1$).
Let $\kappa=f_*(\mu)$ be the pushforward of $\mu$ by $f$.
Intuitively, for $T\in \mathfrak Y_0^{(h)}$, $\kappa(\{T\})$ is the global measure of all the subtrees 
with complete theory $T$ that are rooted at a son of the root.

Let $r_{\widetilde{\mathbf Y}}$ be the root of $\widetilde{\mathbf Y}$. Then it holds that 
$$\sup_{v\sim r_{\widetilde{\mathbf Y}}}
\nu_{\widetilde{\mathbf Y}}(\widetilde{\mathbf Y}(v))=\sup_{T\in X}\frac{\kappa(\{T\})}{w(T_0,T)}.
$$
\end{lemma}
\begin{proof}
 As $\mu$ is pure, the complete theory of $\mu$ is the theory $T_0$ to which every point 
 of the support of $\mu$ projects. Hence the support of $\mu$ defines a unique
 connected component $\widetilde{\mathbf Y}$ of $\mathbb Y_h$. 
That for every $\varphi\in{\rm FO}_1$  it holds that 
$$
\langle\varphi,\widetilde{\mathbf Y}\rangle=\mu(K(\varphi))
$$
is a direct consequence of Lemma~\ref{lem:push1}. 

The second equation is a direct consequence of the construction of $\nu_{\widetilde{\mathbf Y}}$.
\myqed\end{proof}
\begin{theorem}
\label{thm:FO1Y}
Let $\mathbf Y_n$ be an ${\rm FO}_1$-convergent sequence of colored rooted trees with height at most $h$,
and let $\mu$ be the limit measure of $\mu_{\mathbf Y_n}$ on $\mathfrak Y_1^{(h)}$.
Let $\nu$ be the measure defined from $\mu$ by Definition~\ref{def:nu}.
Let $\widetilde{\mathbf Y}$ be the connected component of $\mathbb Y_{h}$ containing the support of $\nu$.
Let $\nu_{\widetilde{\mathbf Y}}$ be the restriction of $\nu$ to $\widetilde{\mathbf Y}$.

Then $\widetilde{\mathbf Y}$, equipped with the probability
 measure $\nu_{\widetilde{\mathbf Y}}$, is a {\rps} ${\rm FO}_1$-limit of 
$(\mathbf Y_n)_{n\in\bbbn}$. 

Moreover, it holds that 
$$\sup_{v\sim r_{\widetilde{\mathbf Y}}}
\nu_{\widetilde{\mathbf Y}}(\widetilde{\mathbf Y}(v))
\leq\liminf_{n\rightarrow\infty}\max_{v\sim r_{\mathbf Y_n}}\frac{|\mathbf Y_n(v)|}{|\mathbf Y_n|}.
$$
\end{theorem}
\begin{proof}
 As $(\mathbf Y_n)_{n\in\bbbn}$ is elementarily convergent, the complete theory of the elementary limit
 of this sequence is the theory $T_0$ to which every point 
 of the support of $\mu$ projects, hence $\mu$ is pure. 
According to Lemma~\ref{lem:FO1muY}, $\widetilde{\mathbf Y}$ is an ${\rm FO}_1$-modeling limit of
$(\mathbf Y_n)_{n\in\bbbn}$. 

Let $\kappa$ be defined as in Lemma~\ref{lem:FO1muY}.
If $\kappa$ is atomless, then $\sup_{v\sim r_{\widetilde{\mathbf Y}}}\nu_{\widetilde{\mathbf Y}}(\widetilde{\mathbf Y}(v))=0$
hence the inequation holds.

Let $T$ be such that $\kappa(\{T\})>0$.
For every $\epsilon>0$ there exists $\theta_\epsilon\in T$ such that
$$\kappa(\{T\})\leq \kappa(\{T':\ T'\ni\theta_\epsilon\})\leq\kappa(\{T\})+\epsilon.$$
Moreover, it follows from Lemma~\ref{lem:w} that
\begin{align*}
w(T_0,T)&=\lim_{\epsilon\rightarrow 0}\sum\{w(T_0,T'):\ T'\ni\theta_\epsilon\}.
\intertext{Then, as $(\mathbf Y_n)_{n\in\bbbn}$ is elementarily convergent to a rooted tree with theory $T_0$ it holds that }
w(T_0,T)&=\lim_{\epsilon\rightarrow 0}\lim_{n\rightarrow\infty} |\{v\in Y_n:\ v\sim r_{\mathbf Y_n}\text{ and }\mathbf Y_n(v)\models\theta_\epsilon\}|\\
&=\lim_{\epsilon\rightarrow 0}\lim_{n\rightarrow\infty} |\{v\in Y_n:\ \mathbf Y_n\models (v\sim r_{\mathbf Y_n})\wedge\rho(\theta_\epsilon)(v)\}|.
\end{align*}
As $\lim_{n\rightarrow\infty} |\{v\in Y_n:\ \mathbf Y_n\models (v\sim r_{\mathbf Y_n})\wedge\rho(\theta_\epsilon)(v)\}|$ is non-increasing when $\epsilon\rightarrow 0$, and is a non-negative integer or $\infty$, there exists $\epsilon_0$ such that
for every $0<\epsilon<\epsilon_0$ it holds that 
$$
w(T_0,T)=\lim_{n\rightarrow\infty} |\{v\in Y_n:\ \mathbf Y_n\models (v\sim r_{\mathbf Y_n})\wedge\rho(\theta_\epsilon)(v)\}|.
$$

For $\epsilon>0$, let $\phi_\epsilon$ be the formula stating that the subtree rooted at a son of the root that contains $x_1$ satisfies $\theta_\epsilon$.
Then it holds that 
\begin{align*}
\kappa(\{T':\ T'\ni\theta_\epsilon\})&=\mu(K(\phi_\epsilon))\\
&=\lim_{n\rightarrow\infty}\langle\phi_\epsilon,\mathbf Y_n\rangle\\
&=\lim_{n\rightarrow\infty}\frac{\sum\bigl\{|\mathbf Y_n(v)|:\ \mathbf Y_n\models (v\sim r_{\mathbf Y_n})\wedge\rho(\theta_\epsilon)(v)\bigr\}}{|\mathbf Y_n|}\\
\end{align*}
Hence, for every $0<\epsilon<\epsilon_0$ it holds that 
\begin{align*}
\frac{\kappa(\{T\})}{w(T_0,T)}&\leq \lim_{n\rightarrow\infty}\frac{\epsilon+\sum\biggl\{\frac{|\mathbf Y_n(v)|}{|\mathbf Y_n|}:\ \mathbf Y_n\models (v\sim r_{\mathbf Y_n})\wedge\rho(\theta_\epsilon)(v)\biggr\}}{|\{v\in Y_n:\ \mathbf Y_n\models (v\sim r_{\mathbf Y_n})\wedge\rho(\theta_\epsilon)(v)\}|}\\
&\leq \epsilon+\liminf_{n\rightarrow\infty}\max\biggl\{\frac{|\mathbf Y_n(v)|}{|\mathbf Y_n|}:\ \mathbf Y_n\models (v\sim r_{\mathbf Y_n})\wedge\rho(\theta_\epsilon)(v)\biggr\}\\
&\leq \epsilon+\liminf_{n\rightarrow\infty}\max_{v\sim r_{\mathbf Y_n}}\frac{|\mathbf Y_n(v)|}{|\mathbf Y_n|}.
\intertext{Hence}
\frac{\kappa(\{T\})}{w(T_0,T)}&\leq\liminf_{n\rightarrow\infty}\max_{v\sim r_{\mathbf Y_n}}\frac{|\mathbf Y_n(v)|}{|\mathbf Y_n|}.
\end{align*}

\myqed\end{proof}

\subsection{Inverse Theorems for ${\rm FO}_1$-limits of Colored Rooted Trees with Bounded Height}
We characterize here the measures $\mu$ on $S(\mathcal B({\rm FO}_1))$, which are weak limits
of measures $\mu_{Y_n}$ for some ${\rm FO}_1$-convergent sequence $(Y_n)_{n\in\bbbn}$ 
of colored rooted trees with height at most $h$.

\begin{fact}
If $(Y_n)_{n\in\bbbn}$ is an ${\rm FO}_1$-convergent sequence
of colored rooted trees with height at most $h$, then $\mu$ is pure
and its complete theory is the limit in $S(\mathbf B({\rm FO}_0))$ of the complete theories
of the rooted trees $Y_n$.
\end{fact}

\begin{figure}[ht]
	\centering
		\includegraphics[width=.75\textwidth]{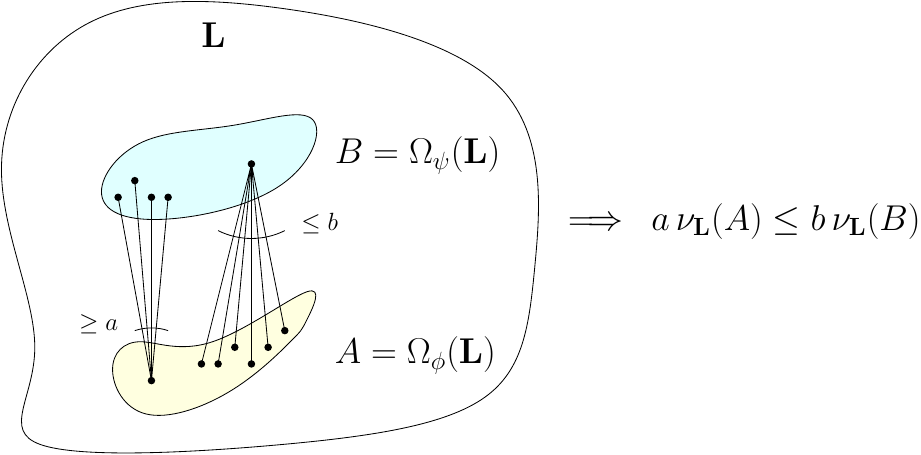}
	\caption{
A modeling $\mathbf L$ satisfies the FMTP if, for every first-order definable subsets
$A,B$ of $L$ and all integers $a,b$ with the property that every element in $A$ has at least $b$ neighbors in $B$
and every element in $B$ has at most $b$ neighbors in $A$, it holds that 
$a\,\nu_{\mathbf L}(A)\leq b\,\nu_{\mathbf L}(B)$ .
	}
	\label{fig:FMTP}
\end{figure}

We now define a Finitary Mass Transport Principle for probability measures on a the Stone dual of $\mathcal B({\rm FO}_p(\lambda))$, 
in a similar way that a Finitary Mass Transport Principle was
 introduced for modelings in Section~\ref{sec:limfin}.
\begin{mytinted}
\begin{definition}
A probability measure $\mu$ on $S(\mathcal B({\rm FO}_p(\sig)))$ ($p\geq 1$) or
$S(\mathcal B({\rm FO}(\sig)))$ satisfies the 
	{\em Finitary Mass Transport Principle} (FMTP) if
for every $\phi,\psi\in{\rm FO}_1(\sig)$ and all integers $a,b$ such that
$$\begin{cases}
\phi\entails (\exists^{\geq a}y)\,(x_1\sim y)\wedge\psi(y)\\
\psi\entails (\exists^{\leq b}y)\,(x_1\sim y)\wedge\psi(y)
\end{cases}$$
it holds that 
$$a\,\mu(K(\phi))\leq b\,\mu(K(\psi)).$$
~\\

Similarly, a modeling $\mathbf L$ satisfies the FMTP if, for every $\phi,\psi,a,b$ as above
the following holds (see Fig.~\ref{fig:FMTP}):
$$a\langle\phi,\mathbf L\rangle\leq b\langle\psi,\mathbf L\rangle.$$
\end{definition}
\end{mytinted}

\begin{fact}
Every finite structure $\mathbf A$ satisfies the FMTP and, consequently, the measures $\mu_{\mathbf A}$ associated
to $\mathbf A$ on $S(\mathcal B({\rm FO}_p))$ ($p\geq 1$) and $S(\mathcal B({\rm FO}))$ satisfy the FMTP.
\end{fact}

Let $r\in\bbbn$. We denote by ${\rm FO}_1^{(r)}$ the fragment of ${\rm FO}_1$ with formulas having quantifier-rank at most $r$.
Note that $\mathcal B({\rm FO}_1^{(r)})$ is a finite Boolean algebra, hence $S(\mathcal B({\rm FO}_1^{(r)}))$ is a finite space.

The following approximation lemma lies in the center of our inverse argument. 
\begin{lemma}
\label{lem:approxY1}
Let $\mu$ be a pure measure on $S(\mathcal B({\rm FO}_1(\sig)))$ with support in $\mathfrak Y_1^{(h)}$ that satisfies the FMTP.
Then, for every integer $r\geq 1$ there exist integer $C=C(\sig,r)$ such that for every $N\in\bbbn$ there is a colored
rooted tree $Y_N$ with the following properties:
\begin{enumerate}
	\item $N\leq |Y_N|\leq N+C$;
	\item for every $\varphi\in{\rm FO}_1$ with quantifier rank at most $r$ it holds that 
$$
\bigl|\langle\varphi,Y_N\rangle-\mu(K(\varphi))\bigr|\leq C/N.
$$
\item the trees $Y_N$ (with root $r_N$) are balanced in the following sense: for every
modeling $\mathbf L$ (with root $r_{\mathbf L}$)
such that $\langle\phi,\mathbf L\rangle=\mu(K(\phi))$ holds for every $\phi\in {\rm FO}_1$, we have
$$\max_{v\sim r_N}\frac{|Y_N(v)|}{|Y_N|}\leq \max\Bigl(\frac{1}{r+h},\sup_{v\sim r_{\mathbf L}}\nu_{\mathbf L}(\mathbf L(v))\Bigr)+C/N.$$

\end{enumerate}
\end{lemma}
\begin{proof}
Note that it is sufficient to prove that there exists $C$ such that for every $N\geq C$ the statement holds (as the initial statement
obviously holds with constant $2C$ instead of $C$).

For integers $k,r$ and a sentence $\phi\in{\rm FO}_0(\sig)$, let $\zeta_k(\phi)$ be the sentence
$$(\exists^{\geq k}y)\ \varrho(\phi)(y),$$
and let
$$c(s,k)=k+1+\max\{{\rm qrank}(\zeta_k(\phi)):\ \phi\in{\rm FO}_0(\sig)\text{ and }{\rm qrank}(\phi)\leq s\}.$$

For formulas $\phi,\psi$ we define
$$w'(\phi,\psi)=\begin{cases}
0&\text{ if }\phi\entails\nexists y\ \varrho(\psi)(y)\\
k&\text{ if }0<k<r+h,\ \phi\entails\zeta_k(\psi),\text{ and }\phi\nvdash\zeta_{k+1}(\psi)\\
r+h&\text{ otherwise.}
\end{cases}
$$
Let $T,T'\in\mathfrak Y^{(h)}_0$ be complete theories of rooted trees, let $a,b$ are integers such that $a\geq c(b,r+h)$,
let $\phi=\bigwedge (T\cap{\rm FO}_0^{(a)})$, and let $\psi=\bigwedge (T'\cap{\rm FO}_0^{(b)})$.
Then either $w'(\phi,\psi)<r+h$ or $\phi\entails\zeta_{r+h}(\psi)$.
This means that for any model $\mathbf Y$ of $T$, either $w'(\phi,\psi)<r+h$ and the root of $\mathbf Y$
has exactly $w'(\phi,\psi)$ sons $v$ such that ${\rm Th}(\mathbf Y(v))\cap {\rm FO}_0^{(b)}=T'\cap {\rm FO}_0^{(b)}$, or
$w'(\phi,\psi)=r+h$ and the root of $\mathbf Y$
has at least $r+h$ sons $v$ such that ${\rm Th}(\mathbf Y(v))\cap {\rm FO}_0^{(b)}=T'\cap {\rm FO}_0^{(b)}$.

Let $\widetilde{\mu}={\rm Encode}_*(\mu)$ (see Lemma~\ref{lem:encode}) be the pushforward of $\mu$
on  $\biguplus_{i=1}^h(\mathfrak Y_0^{(h)})^i$.
For a given integer $r$, we define integers $a_{r,0},a_{r,1},\dots,a_{r,h-1}$ by
$$a_{r,h-1}=r+h,\ a_{r,h-2}=c(a_{r,h-1},r+h),\dots,\ a_{r,0}=c(a_{r,1},r+h).$$
Let $F$ be the mapping defined on $\biguplus_{i=1}^h(\mathfrak Y_0^{(h)})^i$ by
$$F(T_0,\dots,T_i)=(T_0\cap {\rm FO}_0^{(a_{r,0})},\dots,T_i\cap {\rm FO}_0^{(a_{r,i})}).$$
We note that $\mathfrak X=F\bigl(\biguplus_{i=1}^h(\mathfrak Y_0^{(h)})^i\bigr)$ is a finite space,
and we endow $\mathfrak X$ with the discrete topology. (Note that $F$ is continuous.)
We define the probability measure
$\widetilde{\mu}^{(r)}=F_*(\widetilde{\mu})$ on $\mathfrak X$ as the pushforward of $\widetilde{\mu}$ by $F$.

We will construct disjoint sets $V_{\hat T_0,\dots,\hat T_i}$ indexed by the elements $(\hat T_0,\dots,\hat T_i)$ of $\mathfrak X$.
To construct these sets, it will be sufficient to define their cardinalities and the unary relations that apply to their elements.
We proceed inductively on the length of the index tuple. 
As $\mu$ is pure, $\mathfrak X$ contains a unique $1$-tuple $(\hat T_0)$, and
we let the set $V_{\hat T_0}$ to be a singleton. The unique element $r$ of $V_{\hat T_0}$ will be the root of the approximation
tree $Y_N$. Hence we let $R(v)$ and for every color relation $C_i$ we let $C_i(r)$ if $(\forall x) R(x)\rightarrow C_i(x)$ belongs to $\hat T_0$. 
Assume sets $V_{\hat T_0,\dots,\hat T_j}$ have been constructed for every $0\leq j\leq i$ and every $(\hat T_0,\dots,\hat T_j)\in\mathfrak X$. 
Let $(\hat T_0,\dots,\hat T_{i+1})\in\mathfrak X$.
Then of course $(\hat T_0,\dots,\hat T_{i})\in\mathfrak X$.
\begin{itemize}
	\item If $\widetilde{\mu}^{(r)}(\{(\hat T_0,\dots,\hat T_i)\})=0$ and $\widetilde{\mu}^{(r)}(\{(\hat T_0,\dots,\hat T_{i+1})\})=0$
then 	
$$|V_{\hat T_0,\dots,\hat T_{i+1}}|=w'(\bigwedge \hat T_i,\bigwedge \hat T_{i+1})\,|V_{\hat T_0,\dots,\hat T_i}|;$$
	\item If $\widetilde{\mu}^{(r)}(\{(\hat T_0,\dots,\hat T_i)\})>0$ and $\widetilde{\mu}^{(r)}(\{(\hat T_0,\dots,\hat T_{i+1})\})=0$
then 	(according to FMTP) $w'(\bigwedge \hat T_i,\bigwedge \hat T_{i+1})=0$ and we let $V_{\hat T_0,\dots,\hat T_i}=\emptyset$;
	\item If $\widetilde{\mu}^{(r)}(\{(\hat T_0,\dots,\hat T_i)\})=0$ and $\widetilde{\mu}^{(r)}(\{(\hat T_0,\dots,\hat T_{i+1})\})>0$ 
then 	(according to FMTP) $w'(\bigwedge \hat T_i,\bigwedge \hat T_{i+1})=r+h$ and we let	
	$$|V_{\hat T_0,\dots,\hat T_{i+1}}|=\max((r+h)|V_{\hat T_0,\dots,\hat T_i}|,\lfloor \widetilde{\mu}^{(r)}(\{(\hat T_0,\dots,\hat T_{i+1})\})N\rfloor).$$
	\item Otherwise $\widetilde{\mu}^{(r)}(\{(\hat T_0,\dots,\hat T_i)\})>0$ and $\widetilde{\mu}^{(r)}(\{(\hat T_0,\dots,\hat T_{i+1})\})>0$.
	In this case, according to FMTP, it holds that 
	$$w'(\bigwedge \hat T_i,\bigwedge \hat T_{i+1})=\min\left(r+h, \frac{\widetilde{\mu}^{(r)}(\{(\hat T_0,\dots,\hat T_{i+1})\})}{\widetilde{\mu}^{(r)}(\{(\hat T_0,\dots,\hat T_i)\})}\right),$$
		Then, if $w'(\bigwedge \hat T_i,\bigwedge \hat T_{i+1})<r+h$ we let
		\begin{align*}
	|V_{\hat T_0,\dots,\hat T_{i+1}}|&=w'(\bigwedge \hat T_i,\bigwedge \hat T_{i+1})|V_{\hat T_0,\dots,\hat T_i}|
	\intertext{and otherwise we let}
	|V_{\hat T_0,\dots,\hat T_{i+1}}|&=\max((r+h)|V_{\hat T_0,\dots,\hat T_i}|,\lfloor \widetilde{\mu}^{(r)}(\{(\hat T_0,\dots,\hat T_{i+1})\})N\rfloor).
	\end{align*}
	\end{itemize}
The colors of the elements of $V_{\hat T_1,\dots,\hat T_i}$ are easily defined: for $v\in V_{\hat T_1,\dots,\hat T_i}$ and color relation $C_i$ we let
$C_i(v)$ if $(\forall x) R(x)\rightarrow C_i(x)$ belongs to $\hat T_i$.

The tree $Y_N$ has vertex set $\bigcup V_{\hat T_1,\dots,\hat T_i}$.
Each set $V_{\hat T_0,\dots,\hat T_{i+1}}$ is partitioned as equally as possible into
$|V_{\hat T_0,\dots,\hat T_{i}}|$ parts, each part being adjacent to a single vertex in $V_{\hat T_0,\dots,\hat T_{i}}$.
It follows that the degree in $V_{\hat T_0,\dots,\hat T_{i+1}}$ of a vertex in $V_{\hat T_0,\dots,\hat T_{i}}$ lies 
 between $\lceil|V_{\hat T_0,\dots,\hat T_{i+1}}|/|V_{\hat T_0,\dots,\hat T_{i}}|\rceil$ and 
$\lfloor|V_{\hat T_0,\dots,\hat T_{i+1}}|/|V_{\hat T_0,\dots,\hat T_{i}}|\rfloor$, and that
(by construction and thanks to FMTP)
this coincides with $w'(\bigwedge \hat T_i,\bigwedge \hat T_{i+1})$ (when $<r+h$) or is at least 
$w'(\bigwedge \hat T_i,\bigwedge \hat T_{i+1})$ (when $=r+h$).

For $(\hat T_0,\dots,\hat T_i)\in\mathfrak X$,  it is easily checked that
$$\bigl||V_{\hat T_0,\dots,\hat T_i}|-\widetilde{\mu}(\{(\hat T_0,\dots,\hat T_i)\})N\bigr|\leq (r+h)^{i}.$$

For $\phi\in{\rm FO}_1(\sig)$, let $\mathcal F_\phi=\{F\circ{\rm Encode}(T): T\in K(\phi)\cap\mathfrak Y_1^{(h)}\}$.
Let $C=(r+h)^h|\mathfrak X|$.
Then, by summing up the above inequality, we get
$$
0\leq \biggl(\sum_{(\hat T_0,\dots,\hat T_i)\in\mathcal F_\varphi} |V_{\hat T_0,\dots,\hat T_i}|\biggr)-\mu(K(\varphi))N \leq C.
$$
In particular, if $\phi$ is the true statement, we get 
$$
N\leq |Y_N|\leq N+C.
$$
Also, we deduce that for every $(\hat T_0,\hat T_1)\in\mathfrak X$ and every $v_1,v_2\in V_{\hat T_0,\hat T_1}$ the following inequality holds
$$\bigl||Y_N(v_1)|-|Y_N(v_2)|\bigr|\leq C.$$

Let $\mathfrak Z=\{T\in\mathfrak Y_1^{(h)}: T\cap{\rm FO}_0(\sig)=T_0\}$, let $T\in\mathfrak Z$, let $(T_0,\dots,T_i)={\rm Encode}(T)$, 
and let $(\hat T_0,\dots,\hat T_i)=F(T_0,\dots,T_i)$.
We now prove that if $v\in V_{\hat T_0,\dots,\hat T_i}$ and $(T_0',\dots,T_i')={\rm Encode}({\rm Th}(Y_N,v))$ then it holds that 
$T_j\cap{\rm FO}_0^{(r+h)}=T_j'\cap{\rm FO}_0^{(r+h)}$ for every $1\leq j\leq i$ (see Fig~\ref{fig:approx}).
\begin{figure}[ht]
	\centering
		\includegraphics[width=\textwidth]{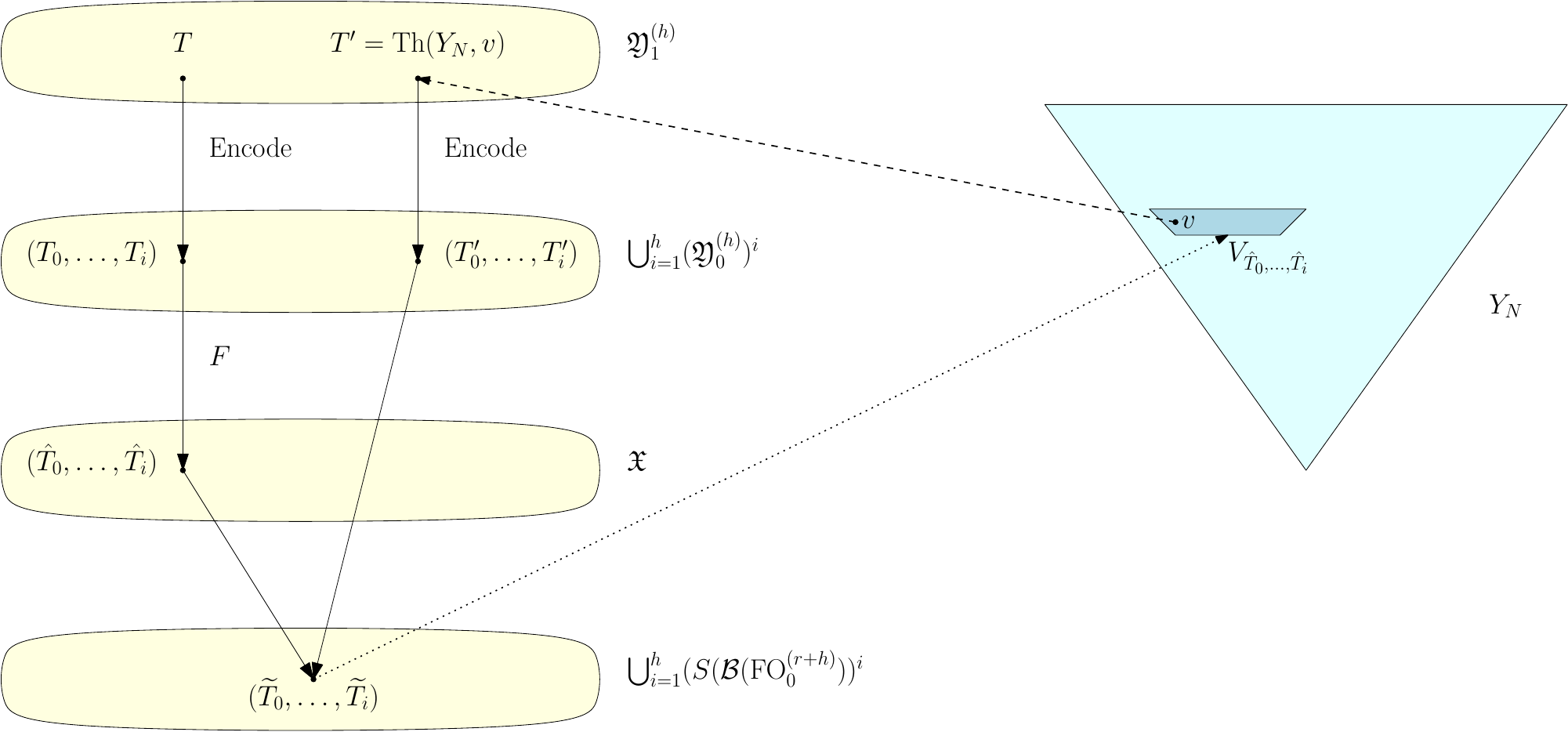}
	\label{fig:approx}
	\caption{Comparing $T$ with its approximation in $Y_N$}
\end{figure}

First note that it is sufficient to prove $T_i\cap{\rm FO}_0^{(r+h)}=T_i'\cap{\rm FO}_0^{(r+h)}$, as the other equalities follow
by considering the ancestors of $v$. If $T_i$ (hence $T_i'$) is the complete theory of a single vertex tree, then
by construction of $V_{\hat T_0,\dots,\hat T_i}$, it holds that  $T_i=T_i'$. Assume now that 
$i$ is such that for every  $(T_0,\dots,T_{i+1})\in {\rm Encode}(T)$ with  $T\in\mathfrak Z$,
it holds that $T_{i+1}\cap{\rm FO}_0^{(r+h)}=T_{i+1}'\cap{\rm FO}_0^{(r+h)}$. 
Let $A$ be a model of $T_i$ and let $B$ be a model of $T_i'$. It follows from the induction that 
the roots of $A$ and $B$ have the same number of sons (up to $r+h$) with subtrees, which are $(r+h)$-equivalent to 
a fixed rooted tree. By an easy argument based on an Ehrenfeucht-Fra{\"\i}ss\'e game, it follows that $A$ and $B$ are
$(r+h)$-equivalent hence $T_{i}\cap{\rm FO}_0^{(r+h)}=T_{i}'\cap{\rm FO}_0^{(r+h)}$. 

According to Lemma~\ref{lem:0to1}, we deduce that for $T\in\mathfrak Z$ and the corresponding vertex $v\in V_{\hat T_0,\dots,\hat T_i}$
it holds that ${\rm Th}(Y_N,v)\cap {\rm FO}_1^{(r)}=T\cap {\rm FO}_1^{(r)}$.

It follows that for every $\varphi\in{\rm FO}_1^{(r)}$ the following equation holds:
$$
\langle\varphi,Y_N\rangle=\sum_{(\hat T_0,\dots,\hat T_i)\in\mathcal F_\varphi} \frac{|V_{\hat T_0,\dots,\hat T_i}|}{|Y_N|}.
$$
Hence
$$
\bigl|\langle\varphi,Y_N\rangle-\mu(K(\varphi))\bigr|\leq C/N.
$$

Let $r_N$ be the root of $Y_N$. Define
$$\alpha_N=\max_{v\sim r_N}\frac{|Y_N(v)|}{|Y_N|}.$$
Assume $\mathbf L$ is a modeling with ${\rm FO}_1$ statistics $\mu$ and root $r_{\mathbf L}$.

Let $(\hat T_0,\hat T_1)\in\mathfrak X$ (vertices in $V_{\hat T_0,\hat T_1}$ are sons of the root of $Y_N$). 
By construction, all the subtrees rooted at a vertex in $V_{\hat T_0,\hat T_1}$ have
almost the same number of vertices (the difference being at most $C$). If $w'(\hat T_0,\hat T_1)=r+h$, it follows
that for every $v\in V_{\hat T_0,\hat T_1}$ it holds that $|Y_N(v)|\leq C+|Y_N|/(r+h)$, i.e. $|Y_N(v)|/|Y_N|\leq 1/(r+h)+C/N$.
Otherwise, $w'(\hat T_0,\hat T_1)=k<r+h$ hence if $\psi$ is the formula stating that the ancestor of $x_1$ which is a son of the root
satisfies $\bigwedge \hat T_1$, then 
$$\mu(K(\psi))=\langle\psi,\mathbf L\rangle=\sum_{v\sim r_{\mathbf L}, \mathbf L\models\psi(v)}\nu_{\mathbf L}(\mathbf L(v))\leq k\sup_{v\sim r_{\mathbf L}}\nu_{\mathbf L}(\mathbf L(v)).$$
Also 
\begin{align*}
\mu(K(\psi))+\frac{C}{N}\geq \langle\psi,Y_N\rangle&=\sum_{(\hat T_0,\dots,\hat T_i)\in\mathcal F_\psi}\frac{|V_{\hat T_0,\dots,\hat T_i}|}{|Y_N|}\\
&=\sum_{v\sim r_N, Y_n\models\psi(v)}\frac{|Y_N(v)|}{|Y_N|}
\intertext{hence $\mu(K(\psi))+\frac{C}{N}\geq\max_{v\sim r_N, Y_n\models\psi(v)}\frac{|Y_N(v)|}{|Y_N|}$ if $k=1$, and otherwise}
\mu(K(\psi))+\frac{C}{N}&\geq k\max_{v\sim r_N, Y_n\models\psi(v)}\frac{|Y_N(v)|-C}{|Y_N|}\\
&\geq k\max_{v\sim r_N, Y_n\models\psi(v)}\frac{|Y_N(v)|}{|Y_N|}-C/N.
\intertext{Hence, considering the case $k=1$ and the case $k\geq 2$ (where $2C/k\leq C$) we get}
\max_{v\sim r_N, Y_n\models\psi(v)}\frac{|Y_N(v)|}{|Y_N|}&\leq \sup_{v\sim r_{\mathbf L}}\nu_{\mathbf L}(\mathbf L(v))+C/N.
\end{align*}
And we deduce
$$\alpha_N\leq \max\Bigl(\frac{1}{r+h},\sup_{v\sim r_{\mathbf L}}\nu_{\mathbf L}(\mathbf L(v))\Bigr)+C/N.$$
\myqed\end{proof}
Summarizing, we get the following two inverse results: 
\begin{mytinted}
\begin{theorem}
\label{thm:mu1}
A measure $\mu$ on $S(\mathcal B({\rm FO}_1))$ 
is the weak limit of a sequence of measures $\mu_{Y_n}$ associated to an ${\rm FO}_1$-convergent sequence $(Y_n)_{n\in\bbbn}$
of finite colored rooted trees with height at most $h$ (i.e. of finite $Y_n\in\mathcal Y^{(h)}$) if and only if
\begin{itemize}
	\item $\mu$ is pure and its complete theory belongs to $\mathfrak Y_0^{(h)}$,
	\item $\mu$ satisfies the FMTP.
\end{itemize}
\end{theorem}
\end{mytinted}
\begin{proof}
Assume that $(Y_n)_{n\in\bbbn}$ is an ${\rm FO}_1$-convergent sequence of finite $Y_n\in\mathcal Y^{(h)}$, and that 
$\mu_{Y_N}\Rightarrow\mu$. According to Remark~\ref{rem:pure}, $\mu$ is pure. As $(Y_n)_{n\in\bbbn}$ is elementarily 
convergent, the complete theory of $\mu$ is the complete theory of the elementary limit of $(Y_n)_{n\in\bbbn}$.
Also, $\mu$ satisfies the FMTP (see Section~\ref{sec:limfin}).

Conversely, assume $\mu$ is pure, that its complete theory belongs to $\mathfrak Y_0^{(h)}$, and that
it satisfies the FMTP. According to Lemma~\ref{lem:approxY1} we can construct a sequence $(Y_n)_{n\in\bbbn}$
of finite $Y_n\in\mathcal Y^{(h)}$ (considering for instance $r=n$ and $N=10^C$ where $C$ is the constant defined 
from $r,h,c$) such that for every formula $\phi\in{\rm FO}_1(\sig)$ it holds that  
$\bigl|\langle\phi,Y_n\rangle-\mu(K(\phi))\bigr|\rightarrow 0$ as $n\rightarrow\infty$, i.e. $\mu_{Y_n}\Rightarrow\mu$.
\end{proof}
From this we deduce
\begin{mytinted}
\begin{theorem}
\label{thm:mod1}
A modeling $\mathbf L$ is the ${\rm FO}_1$-limit of an ${\rm FO}_1$-convergent sequence $(Y_n)_{n\in\bbbn}$
of finite colored rooted trees with height at most $h$ (i.e. of finite $Y_n\in\mathcal Y^{(h)}$) if and only if
\begin{itemize}
	\item $\mathbf L$ is a colored rooted tree with height at most $h$ (i.e. $\mathbf L\in\mathcal Y^{(h)}$),
	\item ${\mathbf L}$ satisfies the FMTP.
\end{itemize}
\end{theorem}
\end{mytinted}
\begin{proof}
That an ${\rm FO}_1$-convergent sequence of finite rooted colored trees $Y_n\in\mathcal Y^{(h)}$ has a modeling
${\rm FO}_1$-limit is a direct consequence of Theorem~\ref{thm:FO1Y}. 
That ${\mathbf L}$ satisfies the FMTP is immediate (as the associated measure 
$\mu=\Type{1}{L}_*(\nu_{\mathbf L})$ does).

Conversely, that a colored rooted tree modeling $\mathbf L\in\mathcal Y^{(h)}$ that satisfies the FMTP
is the ${\rm FO}_1$-limit of a sequence of finite rooted colored trees $Y_n\in\mathcal Y^{(h)}$ is a direct consequence
of Theorem~\ref{thm:mu1}.
\end{proof}

\section{${\rm FO}$-limits of Colored Rooted Trees with Bounded Height}
\label{sec:trees}
In this section we explicitly define {\rps} ${\rm FO}$-limits for ${\rm FO}$-convergent sequences
of colored rooted trees with bounded height.  

We first sketch our 
method.

We consider the signature $\sig^+$, which is
the signature $\sig$ augmented by a new unary relation $P$.
Particular $\sig^+$-structures are
colored rooted forests with a {\em principal} connected component, whose root will be
marked by relation $P$ instead of $R$ (no other vertex gets $P$). The class of colored rooted
forests with a principal connected component and height at most $h$
will be denoted by $\mathcal F^{(h)}$.  
 
We consider three basic interpretation schemes:
\begin{enumerate}
  \item  $\mathsf I_{Y\rightarrow F}$ is a basic interpretation scheme 
  of $\sig^+$-structures in $\sig$-structures defined as follows:
	for every $\sig$-structure $\mathbf A$, the domain of $\mathsf I_{Y\rightarrow F}(\mathbf A)$
	is the same as the domain of $\mathbf A$, and the following holds (for every $x,y\in A$):
\begin{align*}
\mathsf I_{Y\rightarrow F}(\mathbf A)\models x\sim y&\quad\iff\quad\mathbf A\models (x\sim y)\wedge\neg R(x)\wedge\neg R(y)\\
\mathsf I_{Y\rightarrow F}(\mathbf A)\models R(x)&\quad\iff\quad\mathbf A\models (\exists z)\ R(z)\wedge (z\sim x)\\
\mathsf I_{Y\rightarrow F}(\mathbf A)\models P(x)&\quad\iff\quad\mathbf A\models R(x)
\end{align*}
In particular, 	$\mathsf I_{Y\rightarrow F}$ maps
a colored rooted tree $Y\in\mathcal Y^{(h)}$ into a colored rooted forest
$\mathsf I_{Y\rightarrow F}(\mathbf A)(Y)\in\mathcal F^{(h-1)}$, formed by
 the subtrees rooted at the sons of the former root
and a single vertex rooted principal component
(the former root);
\item $\mathsf I_{F\rightarrow Y}$ is a basic interpretation scheme 
  of $\sig$-structures in $\sig^+$-structures defined as follows:
	for every $\sig^+$-structure $\mathbf A$, the domain of $\mathsf I_{F\rightarrow Y}(\mathbf A)$
	is the same as the domain of $\mathbf A$, and the following holds (for every $x,y\in A$):
\begin{align*}
\mathsf I_{F\rightarrow A}(\mathbf A)\models x\sim y&\quad\iff\quad\mathbf A\models (x\sim y)\vee R(x)\wedge P(y)\vee R(y)\wedge P(x)\\
\mathsf I_{F\rightarrow A}(\mathbf A)\models R(x)&\quad\iff\quad\mathbf A\models P(x)
\end{align*}
In particular, 	$\mathsf I_{F\rightarrow Y}$ maps
a colored rooted forest $F\in\mathcal F^{(h)}$ into a colored rooted tree
$\mathsf I_{F\rightarrow Y}(F)\in\mathcal Y^{(h+1)}$ by making each non-principal root a son of the principal root;
\item $\mathsf I_{R\rightarrow P}$ is a basic interpretation scheme 
  of $\sig^+$-structures in $\sig$-structures	defined as follows:
	for every $\sig^+$-structure $\mathbf A$, the domain of $\mathsf I_{R\rightarrow P}(\mathbf A)$
	is the same as the domain of $\mathbf A$, adjacencies are the same in $\mathbf A$ and $\mathsf I_{R\rightarrow P}(\mathbf A)$,
	no element of $\mathsf I_{R\rightarrow P}(\mathbf A)$ is in $R$, and for every $x\in A$ the following equivalence holds:
$$I_{R\rightarrow P}(\mathbf A)\models P(x)\quad\iff\quad \mathbf A\models R(x).$$
(Roughly speaking, the relation $R$ becomes the relation $P$.) 
In particular, 	$\mathsf I_{R\rightarrow P}$ maps	
a colored rooted tree $Y\in\mathcal Y^{(h)}$ into a 
 colored rooted forest $\mathsf I_{R\rightarrow P}(Y)\in\mathcal F^{(h)}$ having a single (principal) component. 
\end{enumerate}

We now outline our proof strategy.
Let $(\mathbf Y_n)_{n\in\bbbn}$ be an ${\rm FO}$-convergent sequence 
of finite rooted colored trees ($\mathbf Y_n\in\mathcal Y^{(h)}$) such that $\lim_{n\rightarrow\infty} |Y_n|=\infty$.

For each $n$, $\mathsf I_{Y\rightarrow F}(\mathbf Y_n)$ is a forest $\mathbf F_n$, and 
$(\mathbf Y_n)_{n\in\bbbn}$ is an ${\rm FO}$-convergent sequence.
According to the Comb Structure Theorem, there exists a countable set $(\mathbf Y_{n,i})_{n\in\bbbn}$
of ${\rm FO}$-convergent sequences  of colored rooted trees $\mathbf Y_{n,i}\in\mathcal Y^{(h)}$ 
and a ${\rm FO}$-convergent sequence $(\mathbf R_n)_{n\in\bbbn}$ of residues $\mathbf R_n\in\mathcal F^{(h)}$, which are
 special colored rooted forests (as the isolated principal root obviously belongs to
$\mathbf R_n$), so that
\begin{itemize}
	\item the sequences $(\mathbf Y_{n,i})_{n\in\bbbn}$ and the sequence $(\mathbf R_n)_{n\in\bbbn}$
	 form a uniformly convergent family of sequences;
	\item for each $n\in\bbbn$ it holds that $\mathsf I_{Y\rightarrow F}(\mathbf Y_n)=\mathbf R_n\cup\bigcup_{i\in I}\mathbf Y_{n,i}$.
\end{itemize}

If the limit spectrum of $(\mathsf I_{Y\rightarrow F}(\mathbf Y_n))_{n\in\bbbn}$
is empty (i.e. $I=\emptyset$), the sequence $(\mathbf Y_n)_{n\in\bbbn}$ of colored rooted trees
is called {\em residual}, and in this case we deduce directly that a residual sequence of colored rooted trees
admit a {\rps} ${\rm FO}$-limit from our results on ${\rm FO}_1$-convergent sequences.

Otherwise, we proceed by induction over the height bound $h$.
Denote by $(\spect_i)_{i\in I}$ the limit spectrum of 
$(\mathsf I_{Y\rightarrow F}(\mathbf Y_n))_{n\in\bbbn}$, let $\spect_0=1-\sum_{i\in I}\spect_i$,
and let
$\mathbf Y_{n,0}=\mathsf I_{R\rightarrow P}\circ \mathsf I_{F\rightarrow
Y}(\mathbf R_n)$.
As $(\mathsf I_{F\rightarrow Y}(\mathbf R_n))_{n\in\bbbn}$ is residual,
$(\mathbf Y_{n,0})_{n\in\bbbn}$ has a {\rps} FO-limit $\widetilde{\mathbf Y}_0$.
By induction, each $(\mathbf Y_{n,i})_{n\in\bbbn}$ has a {\rps} FO-limit $\widetilde{\mathbf Y}_i$.
As $\mathbf Y_n=\mathsf I_{F\rightarrow Y}(\bigcup_{i\in I\cup\{0\}}{\mathbf
Y}_{n,i})$, we deduce (using uniform elementary convergence) that $(\mathbf Y_n)_{n\in\bbbn}$
has {\rps} ${\rm FO}$-limit 
$\mathsf I_{F\rightarrow Y}(\coprod_{i\in I\cup\{0\}}(\widetilde{\mathbf
Y}_{i},\spect_i))$.

This finishes the outline of our construction. Now we provide details.

\subsection{The Modeling FO-limit of Residual Sequences}
We start by a formal definition of residual sequences of colored rooted trees.
\begin{definition}
Let $(\mathbf Y_n)_{n\in\bbbn}$ be a sequence of finite colored rooted trees,
let $N_n$ be the set of all sons of the root of $\mathbf Y_n$, and let 
$\mathbf Y_n(v)$ denote (for $v\in Y_n$) the subtree of $\mathbf Y_n$
rooted at $v$.

The sequence $(\mathbf Y_n)_{n\in\bbbn}$ is {\em residual}
if  
$$
\limsup_{n\rightarrow\infty}\max_{v\in N_n}\frac{|\mathbf Y_n(v)|}{|\mathbf Y_n|}=0.
$$
\end{definition}

We extend this definition to single infinite {\rps}s.

\begin{definition}
 A {\rps} colored rooted tree $\widetilde{\mathbf Y}$ 
 with height at most $h$   is {\em residual}
 if, denoting by $N$ the neighbor set of the root, it holds that  
$$
\sup_{v\in N}\nu_{\widetilde{\mathbf Y}}(\widetilde{\mathbf Y}(v))=0.
$$
 \end{definition}
Note that the above definition makes sense as belonging to  some
$\widetilde{\mathbf Y}(v)$ (for some $v\in N$) is first-order definable hence,
as $\widetilde{\mathbf Y}$ is a {\rss},
each  $\widetilde{\mathbf Y}(v)$ is $\Sigma_{\widetilde{\mathbf Y}}$-measurable.

We first prove that  for a {\rps} colored rooted tree to be a
{\rps} ${\rm FO}$-limit of a residual sequence $(\mathbf Y_n)_{n\in\bbbn}$ of
rooted colored trees with bounded height, it is sufficient that it is a
 {\rps} ${\rm FO}_1$-limit of the sequence.
 
\begin{lemma}
\label{lem:resid1}
Assume $(\mathbf Y_n)_{n\in\bbbn}$ is a residual ${\rm FO}_1$-convergent sequence
of finite rooted colored trees with bounded height with residual {\rps} ${\rm FO}_1$-limit
$\widetilde{\mathbf Y}$. 

Then  $(\mathbf Y_n)_{n\in\bbbn}$ is ${\rm FO}$-convergent and
has {\rps} ${\rm FO}$-limit $\widetilde{\mathbf Y}$.  
\end{lemma}
\begin{proof}
Let $h$ be a bound on the height of the rooted trees $\mathbf Y_n$.
Let $\mathbf F_n=I_{Y\rightarrow F}(\mathbf Y_n)$.
Let $\varpi$ be the formula asserting ${\rm dist}(x_1,x_2)\leq 2h$.
Then $\mathbf F_n\models\varpi(u,v)$ holds if
and only if $u$ and $v$ belong to a same connected component of $\mathbf F_n$.
According to Lemma~\ref{lem:loc1}, we get that
$(\mathbf F_n)_{n\in\bbbn}$ is ${\rm FO}^{\rm local}$-convergent.
As it is also ${\rm FO}_0$-convergent, it is  ${\rm FO}$-convergent
(according to Theorem~\ref{thm:fole}). 
As $\mathbf Y_n=I_{F\rightarrow Y}(\mathbf F_n)$, we deduce
that $(\mathbf Y_n)_{n\in\bbbn}$ is ${\rm FO}$-convergent.

That $\widetilde{\mathbf Y}$ is a {\rps} FO-limit
of $(\mathbf Y_n)_{n\in\bbbn}$ then follows from
 Theorem~\ref{thm:sl}.
 \myqed\end{proof}

\begin{lemma}
\label{lem:residueY}
Let $\mathbf Y_n$ be a residual ${\rm FO}_1$-convergent sequence of colored rooted trees with height at most $h$,
let $\mu$ be the limit measure of $\mu_{\mathbf Y_n}$ on $\mathfrak T_1^{(h)}$,
and let $\widetilde{\mathbf Y}$ be the connected component of $\mathbb Y_{h}$ containing the support of $\nu$.
Then $\widetilde{\mathbf Y}$, equipped with the probability
 measure $\nu_{\widetilde{\mathbf Y}}=\nu$, is a {\rps} ${\rm FO}$-limit of 
$(\mathbf Y_n)_{n\in\bbbn}$. 
\end{lemma}
\begin{proof}
That $\widetilde{\mathbf Y}$ is a residual ${\rm FO}_1$-modeling limit of
$(\mathbf Y_n)_{n\in\bbbn}$ is a consequence of Theorem~\ref{thm:FO1Y}. 
That it is then an ${\rm FO}$-modeling limit of
$(\mathbf Y_n)_{n\in\bbbn}$ follows from Lemma~\ref{lem:resid1}
\myqed\end{proof}

\subsection{The Modeling ${\rm FO}$-Limit of a Sequence of Rooted Trees}

For an intuition of how the structure of a {\rps} FO-limit of a sequence 
of colored rooted trees with height at most $h$ could look like, consider 
a {\rps} rooted colored tree $\mathbf Y$. 
Obviously, the $\mathbf Y$ contains
two kind of vertices: the {\em heavy} vertices $v$ such that the subtree $\mathbf Y(v)$ of $\mathbf Y$
 rooted at $v$ has
positive $\nu_{\mathbf Y}$-measure  and the {\em light} vertices 
for which $\mathbf Y(v)$ has zero $\nu_{\mathbf Y}$-measure.
 It is then immediate that heavy vertices of $\mathbf Y$ induce a countable rooted subtree with same root as 
 $\mathbf Y$.

This suggest the following definitions.

\begin{definition}
A {\em rooted skeleton} is a countable rooted tree $\mathbf S$  together with a
{\em mass function} $m:S\rightarrow(0,1]$ such that
$m(r)=1$ ($r$ is the root of $\mathbf S$) and 
for every non-leaf vertex $v\in S$ it holds that 
$$m(v)\geq\sum_{u\text{ son of }v}m(u).$$
\end{definition}

\begin{figure}[ht]
	\centering
		\includegraphics[width=0.70\textwidth]{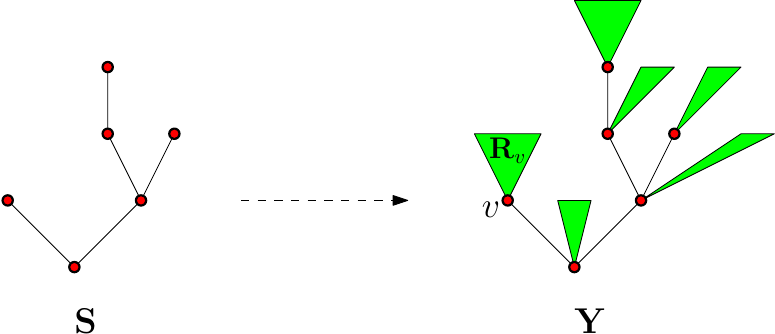}
	\caption{Grafting of trees on a skeleton}
	\label{fig:grafting}
\end{figure}

\begin{definition}
Let $(\mathbf S,m)$ be a rooted skeleton, let 
$S_0$ be the subset of $S$ with vertices $v$ such that
$m(v)=\sum_{u\text{ son of }v}m(u)$, let 
$(\mathbf R_v)_{v\in S\setminus S_0}$ be a countable sequence of 
non-empty residual $\sig$-{\rps}
indexed by $S\setminus S_0$, and let $(\mathbf R_v)_{v\in S_0}$ be a countable 
sequence of
non empty countable colored rooted trees indexed by $S_0$. 
The {\em grafting} of $(\mathbf R_v)_{v\in S\setminus S_0}$ and 
$(\mathbf R_v)_{v\in S_0}$
on $(\mathbf S,m)$ is the {\rps} $\mathbf Y$ defined as follows:
As a graph, $\mathbf Y$ is obtained by taking the disjoint union of 
$\mathbf S$ with the colored rooted trees $\mathbf R_v$ and
then identifying $v\in S$ with the root of $\mathbf R_v$ (see Fig.~\ref{fig:grafting}). 
The sigma algebra $\Sigma_{\mathbf Y}$ is defined as

$$\Sigma_{\mathbf Y}=\Bigl\{\bigcup_{v\in S\setminus S_0}M_v\cup\bigcup_{v\in S_0}M_v':\ M_v\in\Sigma_{\mathbf R_v}, M_v'\subseteq R_v\Bigr\}$$
and the measure $\nu_{\mathbf Y}(M)$ of $M\in\Sigma$ is defined by
$$\nu_{\mathbf Y}(M)=\sum_{v\in S\setminus S_0}\Bigl(m(v)-\hspace{-2mm}\sum_{u\text{ son of }v}\hspace{-2mm}m(u)\Bigr)\nu_{\mathbf R_v}(M_v),$$
where $M=\bigcup_{v\in S\setminus S_0} M_v\cup\bigcup_{v\in S}M_v'$ with $M_v\in\Sigma_{\mathbf R_v}$ and $M_v'\subseteq R_v$.
\end{definition}

\begin{lemma}
Let $\mathbf Y$ be obtained by grafting a countable sequence of non-empty 
{\rps} colored rooted trees $\mathbf R_v$ 
on a rooted skeleton $(\mathbf S,m)$. Then $\mathbf Y$ is 
a {\rps}.
\end{lemma}
\begin{proof}
We prove the statement by induction over the height of the rooted skeleton.
The statement obviously holds if $\mathbf S$ is a single vertex rooted tree
(that is if ${\rm height}(\mathbf S)=1$). Assume that the statement holds
for rooted skeletons with height at most $h$, and let $(\mathbf S,m)$ be
a rooted skeleton with height $h+1$.

Let $s_0$ be the root of $\mathbf S$ and let $\{s_i: i\in I\subseteq\bbbn\}$
be the set of the sons of $s_0$ in $\mathbf S$. For $i\in I$, 
$\mathbf Y_i=\mathbf Y(s_i)$ be the subtree of $\mathbf Y$ rooted at
$s_i$, let $\spect_i=\sum_{x\in Y_i}m(x)$, and let
$m_i$ be the mass function on $\mathbf S_i$ defined by
$m_i(v)=m(v)/\spect_i$. Also, let $\spect_0=1-\sum_{i\in I}\spect_i$.

For each $i\in I\cup\{0\}$, if $\spect_i=0$ (in which case $\mathbf R_{s_i}$
is only assumed to be a {\rss}) we turn $\mathbf R_{s_i}$  into a
{\rps} by defining a probability measure on $\mathbf R_{s_i}$ concentrated 
on $s_i$.

For $i\in I$, let $\mathbf Y_i$ be obtained by grafting the $\mathbf R_v$
on $(\mathbf S_i,m_i)$ (for $v\in S_i$), and let $\mathbf Y_0$ be the
$\sig^+$-modeling consisting in a rooted colored forest with single
(principal) component $\mathbf R_{s_0}$ (that is: $\mathbf
Y_0=\mathsf I_{R\rightarrow P}(\mathbf R_{s_0}))$. 
According to Lemma~\ref{lem:measint}, $\mathbf Y_0$ is a {\rps},
and by induction hypothesis each $\mathbf Y_i$ ($i\in I$) is a {\rps}.
 According to Lemma~\ref{lem:unionfom}, it follows that 
 $\mathbf F=\amalg_{i\in I\cup\{0\}} (\mathbf Y_i,\spect_i)$ is
 a {\rps}. Hence, according to Lemma~\ref{lem:measint},
 $\mathbf Y=\mathsf I_{F\rightarrow Y}(\mathbf F)$ is a {\rps}.
\myqed\end{proof}

Our main theorem is the following.
\begin{mytinted}
\begin{theorem}
\label{thm:limy}
Let $(\mathbf Y_n)_{n\in\bbbn}$ be an ${\rm FO}$-convergent sequence of 
finite colored rooted trees with height at most $h$.

Then there exists a skeleton $(\mathbf S,m)$ and a family
$(\mathbf R_v)_{v\in S}$ --- where 
$\mathbf R_v$ is (isomorphic to) a connected component of $\mathbb Y_{h}$, 
$\Sigma_{\mathbf R_v}$ is the induced $\sigma$-algebra
on $R_v$  --- with the property
that the grafting $\mathbf Y$ of the $\mathbf R_v$ on $(\mathbf S,m)$ 
is a {\rps} ${\rm FO}$-limit of the sequence
$(\mathbf Y_n)_{n\in\bbbn}$.
 \end{theorem}
\end{mytinted}
\begin{proof}
First notice that the statement obviously holds if
$\lim_{n\rightarrow\infty}|Y_n|<\infty$ as
then the sequence is eventually constant to a finite colored rooted 
tree $\mathbf Y$: we can let $\mathbf S$ be $\mathbf Y$ (without the colors),
$m$ be the uniform weight ($m(v)=1/|Y|$), and $\mathbf R_v$ be single vertex rooted tree whose
root's color is the color of $v$ in $\mathbf Y$.
So, we can assume that $\lim_{n\rightarrow\infty}|Y_n|=\infty$.

We prove the statement by induction over the height bound $h$.
For $h=1$, each $\mathbf Y_n$ is a single vertex colored rooted tree, and
the statement obviously holds.

Assume that the statements holds for $h=h_0-1\geq 1$ and let 
finite colored rooted trees with height at most $h_0$.
Let $\mathbf F_n=\mathsf I_{Y\rightarrow F}(\mathbf Y_n)$. Then
$(\mathbf F_n)_{n\in\bbbn}$ is FO-convergent (according to 
Lemma~\ref{lem:measint}). According to the Comb Structure Theorem,
there exists countably many convergent sequences $(\mathbf Y_{n,i})_{n\in\bbbn}$
of colored rooted trees (for $i\in I$) and an ${\rm FO}$-convergent sequence
$(\mathbf R_n)_{n\in\bbbn}$ of special rooted forests forming 
a uniformly convergent family of sequences, such that
$\mathsf I_{Y\rightarrow F}(\mathbf Y_n)=\mathbf R_n\cup\bigcup_{i\in I}\mathbf
Y_{n,i}$.

If the limit spectrum of $(\mathsf I_{Y\rightarrow F}(\mathbf Y_n))_{n\in\bbbn}$
is empty (i.e. $I=\emptyset$), the sequence $(\mathbf Y_n)_{n\in\bbbn}$ of colored rooted trees
is {\em residual}, and the result follows from Lemma~\ref{lem:residueY}.

Otherwise, let  $(\spect_i)_{i\in I}$ the limit spectrum of 
$(\mathsf I_{Y\rightarrow F}(\mathbf Y_n))_{n\in\bbbn}$, let
$\spect_0=1-\sum_{i\in I}\spect_i$, and let
$\mathbf Y_{n,0}=\mathsf I_{R\rightarrow P}\circ \mathsf I_{F\rightarrow
Y}(\mathbf R_n)$.
If $\spect=0$ then there is a connected component $\widetilde{\mathbf Y}_0$ of
$\mathbb Y_h$ that is an elementary limit of $(\mathbf Y_{n,0})_{n\in\bbbn}$;
otherwise, as
$(\mathsf I_{F\rightarrow Y}(\mathbf R_n))_{n\in\bbbn}$ is residual, $(\mathbf Y_{n,0})_{n\in\bbbn}$ has, according to
Lemma~\ref{lem:residueY}, a {\rps} FO-limit
$\widetilde{\mathbf Y}_0$.
By induction, each $(\mathbf Y_{n,i})_{n\in\bbbn}$ has a {\rps} FO-limit $\widetilde{\mathbf Y}_i$.
As $\mathbf Y_n=\mathsf I_{F\rightarrow Y}(\bigcup_{i\in I\cup\{0\}}{\mathbf
Y}_{n,i})$, we deduce, by Corollary~\ref{cor:merge}, Lemma~\ref{lem:elunion},
Theorem~\ref{thm:fole}, and Lemma~\ref{lem:measint},
that $(\mathbf Y_n)_{n\in\bbbn}$ has {\rps} ${\rm FO}$-limit 
$\mathsf I_{F\rightarrow Y}(\coprod_{i\in I\cup\{0\}}(\widetilde{\mathbf
Y}_{i},\spect_i))$.
 \myqed\end{proof}

So, in the case of colored rooted trees with bounded height, we have constructed
an explicit {\rss} that allows one to pullback the limit measure $\mu$ defined
on the Stone space $S({\mathcal B}({\rm FO}))$. 
\subsection{Inverse theorem for ${\rm FO}$-limits of Colored Rooted Trees with Bounded Height}
Recall that for $\sig$-{\rps}s $\mathbf A$ and $\mathbf B$, and $p,r\in\bbbn$ we defined
$$
\|\mathbf A-\mathbf B\|_{p,r}^{\rm local}=\sup\{|\langle\phi,\mathbf A\rangle-\langle\phi,\mathbf B\rangle|:\ \phi\in{\rm FO}_p^{\rm local}(\sig), {\rm qrank}(\phi)\leq r\}.
$$
\begin{lemma}
\label{lem:resmodapprox}
Let $\mathbf L\in\mathcal Y^{(h)}$ (with root $r_{\mathbf L}$) be a colored rooted tree modeling that satisfies the FMTP, let $p,r\in\bbbn$, and let $\epsilon>0$.
Then there exist $C_0=C_0(\sig,p,r,\epsilon), N_0=N_0(\sig,p,r,\epsilon)$ such that for every $N\geq N_0$ there exists a finite colored rooted tree $Y\in\mathcal Y^{(h)}$ such that it holds that $N\leq |Y|\leq N+C_0$, $Y\equiv^r\mathbf L$, and
$$\|Y-\mathbf L\|_{p,r}^{\rm local}<\max(\epsilon,2\sup_{v\sim r_{\mathbf L}}\nu_{\mathbf L}(\mathbf L(v))).
$$
\end{lemma}
\begin{proof}
Without loss of generality we can assume $\epsilon\geq 2\sup_{v\sim r_{\mathbf L}}\nu_{\mathbf L}(\mathbf L(v))$. 
Let $r'=\max(r,4c_{r,p}/\epsilon)$, where $c_{r,p}$ is the constant introduced 
in Lemma~\ref{lem:ploc}.. 

According to Lemma~\ref{lem:approxY1}, there is  $C_0=C(\sig,r')$ (hence $C_0$ depends on 
$\sig,p,r$, and $\epsilon$) 
such that for every $N\in\bbbn$ there exists $Y\in\mathcal Y^{(h)}$  with the following properties:
\begin{enumerate}
	\item $N\leq |Y|\leq N+C_0$;
	\item for every $\varphi\in{\rm FO}_1$ with quantifier rank at most $r'$ the following inequality holds
$$
\bigl|\langle\varphi,Y\rangle-\langle\varphi,\mathbf L\rangle\bigr|\leq C_0/N.
$$
(In particular $Y\equiv^{r'}\mathbf L$ as $N>C_0$.)
\item  we have
$$\max_{v\sim r_N}\frac{|Y(v)|}{|Y|}\leq \max\Bigl(\frac{1}{r'+h},\sup_{v\sim r_{\mathbf L}}\nu_{\mathbf L}(\mathbf L(v))\Bigr)+C_0/N.$$
\end{enumerate}
Let $N_0=4c_{r,p}C(\sig,r')/\epsilon$ and assume $N\geq N_0$.

Let $F={\mathsf I}_{Y\rightarrow F}(Y)$ and $\mathbf A={\mathsf I}_{Y\rightarrow F}(\mathbf L)$.
Let $F_i, i\in\Gamma_F$ and $\mathbf A_i, i\in\Gamma_{\mathbf A}$ be the connected components of $F$ and $\mathbf A$.
Then
$$\max_{i\in\Gamma_F}\frac{|F_i|}{|F|}\leq\max\Bigl(\frac{1}{r'+h},\sup_{i\in\Gamma_{\mathbf L}}\nu_{\mathbf A}(\mathbf A_i)\Bigr)+C_0/N<\frac{\epsilon}{2c_{r,p}}.
$$
As ${\mathsf I}_{Y\rightarrow F}$ is a quantifier free interpretation, for every $\varphi\in{\rm FO}_1$ with quantifier rank at most $r$ the following inequality holds
$$
\bigl|\langle\varphi,F\rangle-\langle\varphi,\mathbf A\rangle\bigr|\leq C_0/N\leq \frac{\epsilon}{4c_{r,p}}.
$$
In particular we have $\|F-\mathbf A\|_{1,r}^{\rm local}\leq \frac{\epsilon}{4c_{r,p}}$.
According to Lemma~\ref{lem:resloc}, it holds that
\begin{align*}
\|F-\mathbf A\|_{p,r}^{\rm local}&<c_{r,p}\biggl(\max_{i\in\Gamma_F}\frac{|F_i|}{|F|}+\sup_{i\in\Gamma_{\mathbf L}}\nu_{\mathbf A}(\mathbf A_i)+\|F-\mathbf A\|_{1,r}^{\rm local}\biggr)<\epsilon.
\end{align*}
and it follows that $\|Y-\mathbf L\|_{p,r}^{\rm local}<\epsilon$.
\end{proof}

\begin{lemma}
\label{lem:modapprox}
Let $\mathbf L\in\mathcal Y^{(h)}$ be an infinite colored rooted tree modeling that satisfies the FMTP, let $p,r\in\bbbn$ and let $\epsilon>0$. 
Then there exist constants $C_h,N_h$ (depending on $\sig, p,r,\epsilon$) such that for every $N\geq N_h$ there is
a finite colored rooted tree $Y_\epsilon\in\mathcal Y^{(h)}$ such that $N\leq |Y_\epsilon|\leq N+C_h$, $Y_\epsilon\equiv^r\mathbf L$, and
$\|\mathbf L-Y_\epsilon\|_{p,r}^{\rm local}<\epsilon$.
\end{lemma}
\begin{proof}
Let $\alpha=\epsilon^2/(2(3c_{r,p})^h)$, where $c_{r,p}$ is the constant which appears in Lemma~\ref{lem:ploc}.
A vertex $v\in L$ is $\alpha$-heavy if either $v$ is the root $r_{\mathbf L}$ of $\mathbf L$, or the father $u$ of $v$ 
in $\mathbf L$ is $\alpha$-heavy and $\nu_{\mathbf L}(\mathbf L(v))>\alpha \nu_{\mathbf L}(\mathbf L(u))$.
The $\alpha$-heavy vertices of $\mathbf L$ form a finite subtree $S$ rooted at $r_{\mathbf L}$ 
(each node $v$ of $S$ has at most $1/\alpha$ sons).

We prove by induction on the height $t$ of $S$ that --- assuming $\alpha\leq\epsilon/2$ ---
there exist constants $C_{t-1},N_{t-1}$ (depending on $\sig, p,r,\epsilon$) such that for every $N\geq N_{t-1}$ there is
a finite colored rooted tree $Y_\epsilon\in\mathcal Y^{(h)}$ such that $N\leq |Y_\epsilon|\leq N+C_{t-1}$, $Y_\epsilon\equiv^r\mathbf L$, and
$\|\mathbf L-Y_\epsilon\|_{p,r}^{\rm local}<\epsilon$.

If $t=1$ (i.e. $r_{\mathbf L}$ is the only $\alpha$-heavy vertex) then 
$$\sup_{v\sim r_{\mathbf L}}\nu_{\mathbf L}(\mathbf L(v))<\alpha.$$
Hence, according to Lemma~\ref{lem:resmodapprox}, there exists $N_0,C_0$ (depending on $\sig,r,p$, and $\epsilon$) such that for every $N\geq N_0$
there is a finite colored rooted tree $Y\in\mathcal Y^{(h)}$ such that $N\leq |Y|\leq N+C_0$, $Y\equiv^r\mathbf L$, and
$$\|Y-\mathbf L\|_{p,r}^{\rm local}<\max(\epsilon,2\sup_{v\sim r_{\mathbf L}}\nu_{\mathbf L}(\mathbf L(v)))=\epsilon.
$$

Now assume that the statement we want to prove by induction holds when $S$ has height at most $t\geq 1$, and let $\mathbf L$ be such that the associated subtree 
$S$ of $\alpha$-heavy vertices has height $t+1$.
Let $v_1,\dots,v_k$ (where $k$ is at most $1/\alpha$) be the $\alpha$-heavy sons of the root $r_{\mathbf L}$ of $\mathbf L$, let
$\mathbf L_i$ be the {\rss} defined by $\mathbf L_i=\mathbf L(v_i)$ for $1\leq i\leq k$, and let $\mathbf L_0$ be the colored rooted tree {\rss} obtained by removing
all the subtrees $\mathbf L_i$ from $\mathbf L$. 
Each $L_i$ is measurable in $\mathbf L$. Let $a_i=\nu_{\mathbf L}(\mathbf L_i)$, and let 
\begin{align*}
\epsilon'&=\frac{\epsilon}{3c_{r,p}}\\
C_t(\sig,p,r,\epsilon)&=\max\biggl(\frac{C_{t-1}(\sig,p,r,\epsilon')}{\alpha}, C_0(\sig,p,r,\epsilon'/3c_{r,p})\biggr)\\
N_t(\sig,p,r,\epsilon)&=\max\biggl(\frac{N_{t-1}(\sig,p,r,\epsilon')}{\epsilon'},\frac{C_{t-1}(\sig,p,r,\epsilon')}{\alpha\epsilon'},N_0(\sig,p,r,\epsilon'/3c_{r,p})\biggr).
\end{align*}
(Note that we do not change $\alpha$.)

Assume $a_0\geq\epsilon'$.
Let $\mathbf {\hat L_i}$ be the modeling with {\rss} ${\mathbf L_i}$ and probability measure $\nu_{\mathbf {\hat L_i}}$ which is
$a_i^{-1}\nu_{\mathbf L}\vert L_i$, where $\nu_{\mathbf L}\vert L_i$ stands for the restriction of $\nu_{\mathbf L}$ to $L_i$.
Let $S_i$ be the rooted subtree of $\alpha$-heavy vertices of $\mathbf {\hat L_i}$. Clearly, if $1\leq i\leq k$, then 
$S_i=S(v_i)$ (as we did not change $\alpha$) thus $S_i$ has height at most $t$. 
Let $\mathbf F\in\mathcal F^{(h)}$ be the forest defined from
$\mathbf F=\coprod_{i=0}^k (\mathbf {\hat L_i},a_i)$ by making the component $\mathbf {\hat L_0}$ special. It is clear that
$\mathbf L={\mathsf I}_{F\rightarrow Y}(\mathbf F)$.
For every $N\geq N_t(\sig,p,r,\epsilon)\geq N_{t-1}(\sig,p,r,\epsilon')/\epsilon'$ there exist,
by induction,  $Y_1,\dots,Y_k$ such that $a_i N\leq Y_i\leq a_i N+C_{t-1}(\sig,p,r,\epsilon')$, $Y_i\equiv^r\mathbf {\hat L_i}$, and
$\|Y_i-\mathbf {\hat L_i}\|_{p,r}^{\rm local}<\epsilon'$.
As the induction step is carried on at most $h$ times, it will always hold that $\alpha\leq\epsilon'^2/2$ hence
$$\sup_{v\sim r_{\mathbf {\hat L_0}}}\nu_{\mathbf {\hat L_0}}(\mathbf {\hat L_0}(v)))
\leq\frac{1}{\epsilon'}\sup_{v\sim r_{\mathbf {L_0}}}\nu_{\mathbf {L}}(\mathbf {L_0}(v)))\leq \frac{\alpha}{\epsilon'}\leq \epsilon'/2.$$
Also, according to Lemma~\ref{lem:resmodapprox}, for every $N\geq N_{t-1}(\sig,p,r,\epsilon')\geq N_0(\sig,p,r,\epsilon')$
there is a finite colored rooted tree $Y_0\in\mathcal Y^{(h)}$ such that $N\leq |Y_0|\leq N+C_0(\sig,p,r,\epsilon')\leq N+C_{t-1}(\sig,p,r,\epsilon')$, 
$Y'\equiv^r\mathbf {\hat L_0}$, and
$$\|Y_0-\mathbf {\hat L_0}\|_{p,r}^{\rm local}<\max(\epsilon',2\sup_{v\sim r_{\mathbf {\hat L_0}}}\nu_{\mathbf {\hat L_0}}(\mathbf {\hat L_0}(v)))=\epsilon'.
$$

Then 
$$\frac{a_i}{N+C_{t-1}(\sig,p,r,\epsilon')/\alpha}\leq \frac{|Y_i|}{\sum_{i=0}^k|Y_i|}\leq \frac{a_i+C_{t-1}(\sig,p,r,\epsilon')}{N}.$$
Thus
$$\biggl|a_i-\frac{|Y_i|}{\sum_{i=0}^k|Y_i|}\biggr|<\frac{C_{t-1}(\sig,p,r,\epsilon')}{\alpha N}\leq\epsilon'.$$
Let $G$ be the disjoint union of the $Y_i$. Hence the following inequality holds, according to Lemma~\ref{lem:ploc}
$$\|\mathbf F-G\|_{p,r}^{\rm local}\leq 2c_{r,p}\epsilon'<\epsilon.$$
Moreover, $N\leq |G|\leq N+C_{t-1}(\sig,p,r,\epsilon')/\alpha\leq N+C_{t}(\sig,p,r,\epsilon)$.

If $a_0<\epsilon'$ we consider $Y_1,\dots,Y_k$ as above, but $Y_0$ is chosen with the only conditions that
 $|Y_0|\leq C_{0}(\sig,p,r,\epsilon')\leq C_{t-1}(\sig,p,r,\epsilon')$ and $Y_0\equiv^r \mathbf L_0$. (Actually, $Y_0$ can be chosen so that
$|Y_0|$ is bounded by a function of $\sig,p$, and $r$ only.)
Let $G$ be the disjoint union of the $Y_i$. Let $\mathbf {\hat L_0}$ be the modeling with {\rss} $\mathbf L_0$ and 
 probability measure $\nu_{\mathbf {\hat L_0}}=a_0^{-1}\nu_{\mathbf L}\vert L_0$ if $a_0>0$, and any probability measure
if $a_0=0$ (for instance the discrete probability measure concentrated on $r_{\mathbf {\hat L_0}}$).
Let $\mathbf F\in\mathcal F^{(h)}$ be the forest defined from
$\mathbf F=\coprod_{i=0}^k (\mathbf {\hat L_i},a_i)$ by making the component $\mathbf {\hat L_0}$ special. It is clear that
$\mathbf L={\mathsf I}_{F\rightarrow Y}(\mathbf F)$.
Then, according to Lemma~\ref{lem:ploc}
\begin{align*}
\|\mathbf F-G\|_{p,r}^{\rm local}&\leq c_{r,p}\bigl(\epsilon'+\sum_{i=1}^k a_i\|\mathbf {\hat L_i}-Y_i\|_{p,r}^{\rm local}+a_0\bigr)\\
&<c_{r,p}\bigl(2\epsilon'+\sup_{1\leq i\leq k}\|\mathbf {\hat L_i}-Y_i\|_{p,r}^{\rm local}\bigr)\\
&\leq 3c_{r,p}\epsilon'=\epsilon.
\end{align*}
and, as above, $N\leq |G|\leq N+C_{t}(\sig,p,r,\epsilon)$.
Now, let $Y_\epsilon={\mathsf I}_{F\rightarrow Y}(G)$.
As ${\mathsf I}_{F\rightarrow Y}$ is basic and quantifier free, and as ${\mathsf I}_{F\rightarrow}(\mathbf Y)=\mathbf L$ it holds that 
$\|\mathbf L-Y_\epsilon\|_{p,r}^{\rm local}<\epsilon$ and $N\leq |Y|\leq N+C_{t}(\sig,p,r,\epsilon)$. 
\end{proof}

\begin{mytinted}
\begin{theorem}
\label{thm:limy2}
A modeling $\mathbf L$ is the ${\rm FO}$-limit of an ${\rm FO}$-convergent sequence $(Y_n)_{n\in\bbbn}$
of finite colored rooted trees with height at most $h$ if and only if
\begin{itemize}
	\item $L$ is a colored rooted tree with height at most $h$,
	\item $L$ satisfies the FMTP.
\end{itemize}
\end{theorem}
\end{mytinted}

\section{Limits of Graphs with Bounded Tree-depth}
\label{sec:td}
Let $Y$ be a rooted forest.
The vertex $x$ is an {\em ancestor} of $y$ in $Y$ if $x$ belongs to
the path linking $y$ and the root of the tree of $Y$ to which $y$ belongs to.
The {\em closure} $\clos(Y)$ of a rooted forest $Y$ is the graph with
vertex set $V(Y)$ and edge set $\{\{x,y\}: x\text{ is an ancestor of
  }y\text{ in }Y, x \neq y\}$. The {\em height} of a rooted forest is the maximum number
  of vertices in a path having a root as an extremity.
The {\em tree-depth} ${\rm td}(G)$ of a graph $G$ is the minimum height
of a rooted forest $Y$ such that $G\subseteq\clos(Y)$. This notion is defined in \cite{Taxi_tdepth}
and studied in detail in \cite{Sparsity}. In particular, graphs with bounded tree-depth serve
as building blocks for {\em low tree-depth decompositions}, see \cite{POMNI,POMNII,POMNIII}.
It is easily checked that for each integer $t$ the property ${\rm td}(G)\leq t$ is first-order definable.
It follows that for each integer $t$ there exists a first-order formula $\xi$ with a single free variable
such that for every graph $G$ and every vertex $v\in G$ the following equivalence holds:
$$G\models \xi(v)\quad\iff\quad {\rm td}(G)\leq t\text{ and }{\rm td}(G-v)<{\rm td(G)}.$$

Let $t\in\bbbn$. 
We define the basic interpretation scheme $\mathsf I_t$, which interprets 
 the class of
connected graphs with tree-depth at most $t$ in the class of $2^{t-1}$-colored rooted trees: given a
$2^{t-1}$-colored rooted tree $\mathbf Y$ (where colors are coded by $t-1$ unary
relations $C_1,\dots,C_{t-1}$), the vertices $u,v\in Y$ are adjacent in $\mathsf
I_t(\mathbf Y)$ if the there is an integer $i$ in $1,\dots,t-1$ such that $\mathbf Y\models C_i(v)$  
 and $u$ is the ancestor of $v$ at height $i$ or $\mathbf Y\models C_i(u)$ and
 $v$ is the ancestor of $u$ at height $i$.

Continuing this with all above results we arrive to the closing statement of this paper.

\begin{mytinted}
\begin{theorem}
\label{thm:limtd}
Let $(\mathbf G_n)_{n\in\bbbn}$ be an ${\rm FO}$-convergent sequence of 
finite colored graphs with tree-depth at most $h$.
Then there exists a colored rooted tree modeling $\mathbf L\in\mathcal Y^{(h)}$ satisfying the FMTP,
 such that the modeling $\mathbf G=\mathsf I_h(\mathbf L)$ has tree-depth at most $h$ and
is a {\rps} ${\rm FO}$-limit of the sequence
$(\mathbf G_n)_{n\in\bbbn}$.

Conversely, if there is colored rooted tree modeling $\mathbf L\in\mathcal Y^{(h)}$ satisfying the FMTP
and if $\mathbf G=\mathsf I_h(\mathbf L)$, then there is an ${\rm FO}$-convergent sequence 
$(\mathbf G_n)_{n\in\bbbn}$ of 
finite colored graphs with tree-depth at most $h$, such that $\mathbf G$ is a modeling ${\rm FO}$-limit
of $(\mathbf G_n)_{n\in\bbbn}$.
 \end{theorem}
\end{mytinted}
\begin{proof}
For each $\mathbf G_n$, there is a colored rooted tree $\mathbf Y_n\in\mathcal Y^{(h)}$
such that $\mathbf G_n=\mathsf I_h(\mathbf Y_n)$.
By compactness, the sequence $(\mathbf Y_n)_{n\in\bbbn}$ has a converging
subsequence $(\mathbf Y_{i_n})_{n\in\bbbn}$, which admits a {\rps}  ${\rm
FO}$-limit $\mathbf Y$ (according to Theorem~\ref{thm:limy2}), and it follows from
Lemma~\ref{lem:measint} that $\mathsf I_h(\mathbf Y)$ is a
{\rps} ${\rm FO}$-limit (with tree-depth at most $h$)
 of the sequence
$(\mathbf G_{i_n})_{n\in\bbbn}$, hence a
{\rps} ${\rm FO}$-limit of the sequence
$(\mathbf G_n)_{n\in\bbbn}$.

Conversely, if there is colored rooted tree modeling $\mathbf L\in\mathcal Y^{(h)}$ satisfying the FMTP
and if $\mathbf G=\mathsf I_h(\mathbf L)$ then, according to  Theorem~\ref{thm:limy2} there is an FO-convergent sequence
$(\mathbf Y_n)_{n\in\bbbn}$ of finite colored rooted trees with FO-limit $\mathbf L$. It follows from Lemma~\ref{lem:measint} that $\mathsf I_h(\mathbf Y)$ is a
{\rps} ${\rm FO}$-limit  of the sequence 
$(\mathbf G_{n})_{n\in\bbbn}$, where $G_n=\mathsf I_h(Y_n)$ is a finite graph with tree-depth at most $h$.
\myqed\end{proof}

\chapter{Concluding Remarks}
 \label{sec:cncl}
\section{Selected Problems}
We hope that the theory developed here will encourage
further researches. Here we list and summarize a sample of related problems (some of which we discussed in Section~\ref{sec:intro}). 
 
 \subsection{Modeling Limits for Nowhere Dense Classes} 
The first problem concern existence of {\rps} ${\rm FO}$-limits.
Recall that a class $\mathcal C$ is {\em nowhere dense} \cite{ND_logic,ECM2009,Nevsetvril2010a,ND_characterization}
if, for every integer $d$ there is an integer $N$ such that the $d$-subdivision of $K_N$ is not a subgraph
of a graph in $\mathcal C$.
We have proven, see Theorem~\ref{thm:modnd}, that if a monotone class $\mathcal C$ is such that
every ${\rm FO}$-convergent sequence of graphs in $\mathcal C$ has a modeling ${\rm FO}$-limit, then
$\mathcal C$ is nowhere dense. It is thus natural to ask whether the converse statement holds.

\begin{problem}
\label{pb:nd}
Let $\mathcal C$ be a nowhere dense class of graphs.
Is it true that every ${\rm FO}$-convergent sequence $(G_n)_{n\in\bbbn}$ of finite graphs
in $\mathcal C$ admit a {\rps} ${\rm FO}$-limit?
\end{problem}

 \subsection{Inverse Problems} 

The Aldous--Lyons conjecture
\cite{Aldous2006} states that
every unimodular distribution on rooted countable graphs with bounded degree is
the limit of a bounded degree graph sequence.
 One of the reformulations of this conjecture is
  that every graphing is an ${\rm
 FO}^{\rm local}$ limit of a sequence of finite graphs.
 The importance of this conjecture appears, for instance,
in the fact that it
would imply that all groups are sofic, which
would prove a number of famous conjectures which are proved for sofic groups but
still open for all groups. 

If Aldous--Lyons conjecture holds, then it follows that every graphing is local-equivalent to a graphing with the finite model property. Indeed, if $(G_n)_{n\in\bbbn}$ is BS-convergent to a graphing $\mathbf G$, then (by compactness) some subsequence of $(G_n)_{n\in\bbbn}$ is FO-convergent and (by Corollary~\ref{cor:BSE}) has a graphing FO-limit 
$\mathbf G'$, which is local-equivalent to $\mathbf G$. Hence the following problem can be seen as a natural first step towards the resolution of Aldous--Lyon conjecture:
\begin{problem}
	Is every graphing local-equivalent to a graphing with the finite model property?
\end{problem}

If the previous problem would have a positive answer then
the next problem would be a possible strengthening of Aldous--Lyons conjecture.

\begin{problem}
Is every graphing $\mathbf G$ with the finite model property 
an ${\rm FO}$-limit of a sequence of finite graphs?
\end{problem}

Although the existence of a modeling ${\rm FO}$-limit for ${\rm FO}$-convergent sequences of graphs with bounded tree-depth
follows easily from our study of ${\rm FO}$-convergent sequence
of rooted colored trees, the inverse theorem is more difficult.
Indeed, if we would like to extend the inverse theorem for rooted colored trees 
to bounded tree-depth modelings (thus removing the condition $\mathbf G=\mathsf I_h(\mathbf L)$ in Theorem~\ref{thm:limtd}), we naturally have to address the following question:
\begin{problem}
Is it true that there is a function $f:\bbbn\rightarrow\bbbn$ 
such that for every graph modeling $\mathbf L$ with tree-depth
at most $t$ there exists a rooted colored tree modeling $\mathbf Y$ with height at most $f(t)$
such that $\mathbf L=\mathsf I_{f(t)}(\mathbf Y)$, where the
$\mathsf I_h$   (for $h\in\bbbn$) are the basic interpretation schemes introduced in Section~\ref{sec:td}? 
\end{problem}

\subsection{Classes with Bounded SC-depth}
We can generalize  our main construction of limits to other tree-like
classes. For example, in a similar way that we obtained a {\rps} ${\rm
FO}$-limit for ${\rm FO}$-convergent sequences of graphs with bounded tree-depth, it is possible to get 
a {\rps} ${\rm FO}$-limit for ${\rm FO}$-convergent sequences
of graphs with bounded SC-depth, where SC-depth is defined as follows~\cite{Ganian2012}:

Let $G$ be a graph and let $X\subseteq V(G)$.
We denote by $\overline{G}^X$ the graph $G'$ with vertex set $V(G)$ where
$x\neq y$ are adjacent in $G'$ if 
(i) either $\{x,y\}\in E(G)$ and $\{x,y\}\not\subseteq X$, or
(ii) $\{x,y\}\not\in E(G)$ and $\{x,y\}\subseteq X$.
In other words, $\overline{G}^X$ is the graph obtained from $G$ 
by complementing the edges on~$X$.

\def\SC#1{{\mathcal SC}(#1)}
\begin{definition}[SC-depth]
\label{def:SC-depth}
We define inductively the class $\SC n$ as follows:
\begin{itemize}
  \item We let $\SC0=\{K_1\}$;
  \item if $G_1,\dots,G_p\in\SC n$ and $H= G_1\dot\cup\dots\dot\cup G_p$
  denotes the disjoint union of the $G_i$,
  then for every subset $X$ of vertices of $H$ we
  have $\overline{H\,}^X\in\SC{n+1}$.
\end{itemize}
The {\em SC-depth} of $G$ is the minimum integer $n$ such that $G\in\SC n$.
\end{definition}
Note that classes with bounded SC-depth 
can be seen as the first step towards moving from the study of monotone classes to the study of hereditary classes (that is classes closed under induced subgraphs). For instance, classes of graphs with bounded tree-depth are exactly those monotone classes of graphs where first-order logic and monadic second-order logic have the same expressive power, while classes with bounded SC-depth are exactly those hereditary classes where first-order logic and monadic second-order logic have the same expressive power \cite{6280445}.
\subsection{Classes with Bounded Expansion}
A graph $H$ is a {\em shallow topological minor} of a graph $G$ at depth $t$ if
some $\leq 2t$-subdivision of $H$ is a subgraph of $G$. For a class $\mathcal C$
of graphs we denote by $\mathcal C\shtm t$ the class of all shallow topological 
minors at depth $t$ of graphs in $\mathcal C$. The class $\mathcal C$ has {\em
bounded expansion} if, for each $t\geq 0$, the average degrees of the graphs in
the class $\mathcal C\shtm t$ is bounded, that is (denoting by $\overline{\rm
d}(G)$ the average degree of a graph $G$):
$$(\forall t\geq0)\quad \sup_{G\in\mathcal C\shtm t} \overline{\rm
d}(G)<\infty.$$
The notion of classes with bounded expansion were introduced by the authors 
in~\cite{ICGT05,Taxi_stoc06,POMNI}, and their properties further
studied in
\cite{POMNII,POMNIII,Dvo2007,Dvov2007,Taxi_blrn,Nevsetvril2010,ECM2009,Nevsetvril2010a,Sparsity,BEEx}
and in the monograph~\cite{Taxi_hom}. Particularly, classes with bounded
expansion include classes excluding a topological minor, like classes
with bounded maximum degree, planar graphs, proper minor closed classes, etc.\ 

Classes with bounded expansion have the characteristic property that they admit
special decompositions --- the so-called {\em low tree-depth decompositions} ---
related to tree-depth:
\begin{theorem}[\cite{Taxi_stoc06,POMNI}]
\label{thm:ltd}
 Let $\mathcal C$ be a class of graphs. Then $\mathcal C$ has bounded expansion 
if and only if for every integer $p\in\bbbn$ there exists $N(p)\in\bbbn$ such
that the vertex set of every graph $G\in\mathcal C$ can be partitioned into
at most $N(p)$ parts in such a way that the subgraph of $G$ induced by any
$i\leq p$ parts has tree-depth at most $i$.
\end{theorem}

By an inductive argument, following \cite{Grohe2011}, we can prove that for every integer $p,r$ and every 
class $\mathcal C$ of $\sig$-structure with bounded expansion, there is a
signature $\sig^+\supseteq\sig$, such that every $\sig$-structure
$\mathbf A\in\mathcal C$ can be lifted into a $\sig^+$-structure $\mathbf
A^+$ with same Gaifman graph, in such a way that for every first-order formula 
$\phi\in{\rm FO}_p(\sig)$ with quantifier rank at most $r$ there is
an existential formula $\widetilde\phi\in{\rm FO}_p(\sig^+)$ such that
for every $v_1,\dots,v_p\in A$ the following equivalence holds:
$$\mathbf A\models\phi(v_1,\dots,v_p)\quad\iff\quad\mathbf
A^+\models\widetilde\phi(v_1,\dots,v_p).$$
Moreover, by considering a slightly stronger notion of lift if necessary, we
can assume that $\widetilde\phi$ is a local formula.
We deduce that there is an integer $q=q(\mathcal C,p,r)$ such that
checking  $\phi(v_1,\dots,v_p)$ can be done by considering satisfaction 
of  $\widetilde\psi(v_1,\dots,v_p)$ in subgraphs induced by $q$ color classes
of a bounded coloration. Using a low-tree depth decomposition
(and putting the corresponding colors in the signature $\sig^+$), 
we get that there exists finitely many induced substructures $\mathbf A^+_I$
($I\in \binom{[N]}{q}$) with tree-depth at most $q$ and
the property that  for every first-order formula 
$\phi\in{\rm FO}_p(\sig)$ with quantifier rank at most $r$ there is
an existential formula $\widetilde\phi\in{\rm FO}_p(\sig^+)$ such that
for every $v_1,\dots,v_p\in A$ with set of colors $I_0\subseteq I$ the following equivalence holds:
$$\mathbf A\models\phi(v_1,\dots,v_p)\ \iff\ 
\exists I\in \binom{[N]}{q-p}:\ \mathbf
A_{I\cup I_0}^+\models\widetilde\phi(v_1,\dots,v_p).$$ Moreover, the Stone
pairing $\langle\phi,\mathbf A\rangle$ can be computed by inclusion/exclusion from
Stone pairings $\langle\phi,\mathbf A^+_I\rangle$ for $I\in\binom{[N]}{\leq q}$.

Thus, if we consider an ${\rm FO}$-convergent sequence $(\mathbf
A_n)_{n\in\bbbn}$,
the tuple of limits of the $\sig^+$-structures $(\mathbf A_n)^+_I$
behaves as a kind of approximation of the limit of the $\sig$-structures
$\mathbf A_n$. 
\section{Addendum}
\label{sec:add}
\subsection{Modeling Limits for Nowhere Dense Classes}
Since the submission of this paper a great progress has been made on Problem~\ref{pb:nd}. Based on the results of this paper, classes of graphs for which there exist explicit modeling FO-limits (satisfying the Finitary Mass Transport Principle) now include the class of forests \cite{modeling} and, more generally, classes of graphs with bounded pathwidth \cite{gajarsky2016first}.

For the general case, it has first been proved \cite{SurveyND}
that for every FO-convergent sequence of graphs $(G_n)_{n\in\bbbn}$ there exists a modeling $\mathbf L$ such that for every first-order formula $\phi$ the following equation holds
$$
\langle\phi,\mathbf L\rangle=0\quad\iff\quad\lim_{n\rightarrow\infty}\langle\phi,G_n\rangle=0.
$$
This result has been extended to prove that every FO-convergent sequence of nowhere dense graphs has a modeling limit \cite{modeling_arxiv}. (Note that this result heavily relies on this paper.)
But 
it is still open whether this modeling limit could be required to satisfy the Finitary Mass Transport Principle.
\subsection{Asymptotic Connectivity}
Some further applications include
the study of the connectivity structure of ${\rm FO}^{\rm local}$-convergent sequences we started in Section~\ref{sec:decompose} has been further refined in \cite{modeling} to study modeling limits of forests (with unbounded height). In \cite{Loclim} we deal with the important notion of clustering of a convergent sequence, and show that connectivity properties --- although not first-order definable --- can be established in ${\rm FO}^{\rm local}$-convergent sequences by means of Fourier analysis.
\subsection{Inverse Problems}
The study of existence of modeling limits for simple algebraic structures has led us to prove that FO-convergent sequences of mappings admit a modeling limit \cite{MapLim}, and we have been able to prove  inverse theorems in this case \cite{MapApprox}: every
atomless modeling mapping  that satisfies the Finitary Mass Transport Principle is the ${\rm FO}^{\rm local}$-limit of an ${\rm FO}^{\rm local}$-convergent sequence of finite mappings, and if moreover its complete theory has the Finite Model Property then if it is the ${\rm FO}$-limit of an ${\rm FO}$-convergent sequence of finite mappings.
\subsection{Rooting of Modelings}
Problems~\ref{pb:root1} and~\ref{pb:rootk} have been solved negatively in \cite{Christofides2016}, where it is nevertheless proved that Problem~\ref{pb:root1} holds for almost all rootings of the modeling limit.
\subsection{Others}
The analytic framework of $X$-convergence  has also been presented in \cite{Nevsetvril2014}.
One of the main reasons for interest in our notion of convergence is that it allows to consider structures with arbitrary (countable) signature, and interpretations of these. For instance, it led to the study of limits of mappings \cite{MapLim} (mentioned above), limits of matroids \cite{matroid_limit}, and quantifier-free convergence of tree-semilattices \cite{QFTSL-arxiv}.
\section*{Acknowledgements}
The authors would like to thank Pierre Charbit for suggesting 
to use a proof by
induction for Lemma~\ref{lem:0to1} and Cameron Freer for his help
in proving that the FO-limit of a sequence of random graphs cannot be a {\rps}.

The authors would like to express their gratitude to the anonymous referee for his careful reading of this manuscript.

\backmatter
\providecommand{\bysame}{\leavevmode\hbox to3em{\hrulefill}\thinspace}
\providecommand{\MR}{\relax\ifhmode\unskip\space\fi MR }
\providecommand{\MRhref}[2]{%
  \href{http://www.ams.org/mathscinet-getitem?mr=#1}{#2}
}
\providecommand{\href}[2]{#2}


\end{document}